
%
\input amstex
\documentstyle{amsppt}%
\NoBlackBoxes

\font\ninesy=cmsy9
\font\ninei=cmmi9
\font\twelvbf=cmbx12
\font\ninerm=cmr9
\font\twelvrm=cmr12
\font\tent=cmtt10
\font\eightmsbm=msbm8
%

\def\myl{{\vcenter{\hrule width 0.13 true in height 0.005 true in}}}
\def\myline{$\myl$}
\def\myline{\myl}
\def\mysqcuplus{{\rlap{$\sqcup$}{\raise1.5pt\hbox{\sevenbf +}}}}

\def\sqcupP{\ \!{{\rlap{$_{_{\!}}\sqcup
$}{\vbox{\moveleft2.8pt\hbox{{\ }\raise1.0pt\hbox{\fivebf P}}}}} }
}

\def\Psqcup#1{\ \!{{\rlap{$_{_{\!}}\sqcup
$}{\vbox{\moveleft2.1pt\hbox{{\ }\raise1.0pt\hbox{\fivebf
\hbox{#1}}}}}} } }
\def\myPsqcup#1{\Psqcup#1}
\def\mysubsneqq{{
{{\subset} \raise1pt\hbox{$\!\!\!\!\!\!_{_{_{ {_{\ {\not\
\!\!\!{\not}{\hbox{\sevenbf=}}}}}}}}$}}}}
\def\mysubsetneqq{$\mysubsneqq$}
\def\mmysubsetneqq{\mysubsneqq}

\def\nabl#1  
{\font\mysyfont=cmsy#1
{
{\!_{\!}}^{_{_{\hbox{\mysyfont \char"72}}}}\!}}%
\def\mynabla#1{$\nabl#1 $} 

\def\barnabl#1{\!\bar{\ \!\nabl#1 }}
\def\mybarnabla#1{\barnabl#1 }

\def\thickl{{\vcenter{\hrule width 0.1 true in height 0.015 true in }}}

\def\mthickline{\thickl}
\def\DEFdoubleH{{
\rlap{\raise3.55pt\hbox{$_{_{-\!\!\!{\tenbf -}}}\!$}}{\rlap{ }
{{^{\!}\hbox{\tenbf H}}}}}}
\def\doubleH{$\DEFdoubleH$}
\def\mdoubleH{\DEFdoubleH}

\def\eDEFdoubleH{{
\rlap{\raise2.95pt\hbox{$_{_{-}}\!$}}{\rlap{ }
{{^{\!}\hbox{\eightbf H}}}}}}

\def\medoubleH{\eDEFdoubleH}

\def\DEFminidoubleH{{
\rlap{\raise2.6pt\hbox{$_{_{\ \!{\!}{\hbox{{\sevenrm
-}$_{\!}${\sevenrm -}}\!}}}$}}{\rlap{ }{^{\!}\hbox{\eightbf H}}}}}

\def\mminidoubleH{\DEFminidoubleH}

%
\def\FFrame#1#2#3{                   %
     \vbox{\hrule height#2           
          \hbox{\vrule width#2       %
                 \hskip#1            
                 \vbox{\vskip#1{}    %
                       #3            
                       \vskip#1}     %
                 \hskip#1
                 \vrule width#2}
           \hrule height#2}}

 \hoffset=0.77 in

\topmatter
\title Augmental Homology and the K\"{u}nneth
Formula for 
Joins\endtitle
\author G\"oran Fors\endauthor
\address $\!\!{\!}$Department of Mathematics, University of Stockholm, S-$\!$106 91 Stockholm,
Sweden \endaddress

\email goranf\@matematik.su.se\endemail


\abstract
%
\def\nidoubleH{{
\rlap{\raise2.7pt\hbox {$_{_{\
_{\!}\!_{\!}{\hbox{-$\!$-}\!}}}$}}{\rlap{ }{{{\!}\hbox{\eightbf
H}}}}}}

The ``simplicial complexes" and ``join" ($\ast$) today used within
combinatorics aren't the classical concepts, cf. \cite{25} %
p.\ 108-9, but, except for {\eightsy;}, complexes having {\eightsy
f;g} as a subcomplex resp.$\
${{\eightrm{\char"06}}\lower1.0pt\hbox{\fiverm1}{\eightsy
{\char"03}}{\eightrm{\char"06}}\lower1.0pt\hbox{\fiverm2}:=}$
${{\eightsy{\char"66}}$\ _{\!}\!${\eighti{\char"1B}}$
_{\!}$\lower1.5pt\hbox{\fiverm1}{\sixsy{\char"5B}}$
${\eighti{\char"1B}}$_{\!}
$\lower1.5pt\hbox{\fiverm2}{\sixsy {\char"6A}}{\eighti
{\char"1B}}\lower1.5pt\hbox{\fivei i}{\sixsy {\char"32}}{\eightrm
{\char"06}}\lower1.5pt\hbox{\fivei i}$_{^{\!}}${\eightsy
{\char"67}}},
implying a tacit change of unit element w.r.t. the join operation,
from  {\eightsy {\char"3B}} to {\eightsy {\char"66}}{\eightsy
{\char"3B}}{\eightsy {\char"67}}.
Extending the classical realization functor to this category of
simplicial complexes we end up with a ``restricted" category of
topological spaces,
``containing" the classical and where the classical (co)homology
theory, as well as the ad-hoc invented reduced versions,
automatically becomes obsolete, in favor of a unifying and more
algebraically efficient\nobreak\ theory.
This very modest category modification greatly improves the
interaction between algebra and topology.
E.g. it makes it possible to calculate the homology groups of a
topological pair-join, expressed in the relative factor groups,
leading up to a truly simple boundary formula for joins of
manifolds:
%
{{\eightrm Bd({\eighti
X}}$^{\!}$\lower1.0pt\hbox{\fiverm1}{\eightsy
{\char"03}}{\eightrm{\eighti
X}}$^{\!}$\lower1.0pt\hbox{\fiverm2}{\eightrm)}=}%
{{\eightrm ((Bd{\eighti
X}}$^{\!}$\lower1.0pt\hbox{\fiverm1}{\eightrm)%
{\eightsy {\char"03}}%
{\eightrm{\eighti X}}$^{\!}$\lower1.0pt\hbox{\fiverm2}}{\eightrm)}%
}{\sevensy{\char"5B}}%
{\eightrm({\eighti X}}$^{\!}${\lower1.0pt\hbox{\fiverm1}%
{\eightsy {\char"03}}%
{\eightrm (Bd{\eighti X}}$^{\!}$\lower1.0pt\hbox{\fiverm2}}{\eightrm))}, %
the "product"-counterpart of which is true also classically.
It's also easily seen that no finite simplicial $n$-manifold has
an $(n-2)$-dimensional boundary, cf. Cor.\ 1 p. 26.
%
%
%
\endabstract

\endtopmatter

\document

\head Contents\endhead
{\sevenrm

{\parindent=0.1cm
{{\eightrm 1.}\ \ \ {\eightrm I}{\sevenrm NTRODUCTION}\dotfill
{\eightrm p.\ \ 2}}
\vskip-2.0pt
{{\eightrm 2.}\ \ \ {\eightrm A}{\sevenrm UGMENTAL HOMOLOGY
THEORY}\dotfill {\eightrm p.\ \ 3}}

\vbox{\moveright0.24cm\vbox{\hskip 0.3cm \FFrame{1.00pt}{0.00pt}{
\baselineskip=8.0pt
\hsize=4.67 true in
%
%
\noindent$\!_{\!}$
\dots \ %
 {\sevenrm 2.1}.
{\sevenrm Notations and Definitions of Underlying Categories}\dotfill%
{\sevenrm p.\ \ 3}$\!\!$%
\break
\dots\ {\sevenrm 2.2}. %
{\sevenrm Simplicial Augmental Homology and
Realizations}.\dotfill \ {\sevenrm p.\ \ 4}$\!\!$%
\break
\dots\ {{\sevenrm 2.3}. {\sevenrm Singular Augmental Homology
Theory}\ \dotfill {\sevenrm
p.\ \ 6}}$\!\!$}%
}}
\vskip-2.0pt
{{\eightrm 3.}\ \ \ {\eightrm 3.}1.\  {\eightrm A}$
${\sevenrm UGMENTAL HOMOLOGY MODULES FOR PRODUCTS AND JOINS}
\dotfill {\eightrm p.\ \ 7}}
\vskip-2.0pt
\vbox{\moveright0.24cm\vbox{\hskip 0.3cm \FFrame{3.0pt}{0.00pt}{
\baselineskip=8.0pt%
\hsize=4.68 true in%
\noindent$\!\!_{\!}$\dots \ %
{{\sevenrm 3.}1. Definitons of the Product and Join Operations
\dotfill {\sevenrm p.\ \ \ $_{\!\!}$7}}$_{\!}$%
\break
\noindent$\phantom{}\!\!_{\!}$
\dots \ %
{{\sevenrm 3.}2. Augmental Homology for Products and Joins
\dotfill {\sevenrm p.\ \ \ $_{\!\!}$8}}$_{\!}$%
\break
\noindent$\phantom{}\!\!_{\!}$
\dots \ %
{{\sevenrm 3.}3. Local Augmental Homology Groups for Products and Joins
\dotfill{\sevenrm p.\ $_{\!}$11}}$_{\!}$%
\break
\noindent$\phantom{}\!\!_{\!}$\dots \ %
{{\sevenrm 3.}4. Singular Homology Manifolds under Products and
Joins
\dotfill   p.\ 12}$_{\!}$}
}}
\vskip-3.0pt
{${\eightrm \hbox{\eightrm Specialization I: \ %
General Topological Properties for Realizations of Simplicial
Complexes
}} $}\dotfill {\eightrm p.\nobreak\ 14}$_{\!}$%
\vskip-1.0pt

\vbox{\moveright0.24cm\vbox{\hskip 0.3cm \FFrame{3.0pt}{0.00pt}{%
\baselineskip=8.0pt%
\hsize=4.68 true in%
\noindent{$\phantom{}\!\!_{^{\!}}$}\dots\ \ $\!${\sevenrm
{\sevenrm I.}1.\ \ Realizations and Local Homology Groups Related
to Simplicial Products and Joins \dotfill p.\ 14}$\!{\!}$
\break%
%
%
\noindent \phantom{i}{$\!\!\!_{^{\!\!}}$}\dots\ {\sevenrm
{\sevenrm I}.2.\ $_{\!}$ Simplicial Connectedness Properties under
Products and Joins\dotfill p.\ 15}$_{\!}$}
}}
\vskip-3.0pt
{\eightrm Specialization II: \ Concepts Related to Combinatorics
and Commutative Algebra\dotfill {\eightrm p. 17}}$_{\!}$%
%

\vbox{\moveright0.24cm\vbox{\hskip 0.3cm \FFrame{1.00pt}{0.00pt}{%
\baselineskip=8.0pt%
\hsize=4.68 true in%
%
\noindent$\phantom{}\!_{\!}$\dots\ {\sevenrm {\sevenrm II.}1.\ \
Definition of
Stanley-Reisner rings\dotfill p.\ 17}$\!\!$%
\break
\dots\ {\sevenrm {\sevenrm II.}2.\ \ Buchsbaum, Cohen-Macaulay and
2-Cohen-Macaulay Complexes\dotfill  p.\ 18}$\!\!$%
\break%
\dots\ {\sevenrm {\sevenrm II.}3.\ \ Segre Products\dotfill   p.\ 21}$\!\!$%
\break
\dots\ {\sevenrm {\sevenrm II.}4.\ \ Gorenstein Complexes\dotfill
p.\ 22}$\!\!$}%
}}
\vskip-2.0pt

{\eightrm Specialization III: Simplicial Manifolds \dotfill
{\eightrm p.\ 23}}$_{\!}$%
%

\vbox{\moveright0.24cm\vbox{\hskip 0.3cm  \FFrame{1.00pt}{0.00pt}{
\baselineskip=8.0pt%
\hsize=4.68 true in%
%
\noindent $\phantom{}\!_{\!}$\dots
{\sevenrm {\sevenrm III.}1.\ \
Definitions\dotfill    p.\ 23}$\!\!$%
\break
\dots{\sevenrm {\sevenrm III.}2.\ \ Auxiliaries\dotfill   p.\ 24}$\!\!$%
\break
\dots{\sevenrm {\sevenrm III.}3.\ \ Products and Joins of
Simplicial Manifolds\dotfill    p.\ 27}$\!\!$%
\break%
\dots{\sevenrm {\sevenrm III.}4.\ \ Simplicial Homology Manifolds
and Their Boundaries\dotfill    p.\ 28}$\!\!$}}}

{\eightrm Appendix: \ \ \ Simplicial Calculus and Simplicial Sets
\dotfill{\eightrm p.\ 30}$\!$%

{\sevenrm REFERENCES\dotfill   {\eightrm p.\ 32}}}$\!$%
}}


\head 1. Introduction\endhead
The importance of the {\it join} {\it operation} within algebraic
topology, becomes apparent for instance through Milnor's
construction of the universal principal fiber bundle in \cite{18} %
where he also formulate the non-relative K\"unneth formula for
joins as;
$\
\widetilde{\!\!\hbox{\tenit H}}_{_{\!q+1\!}}
        (\!X\!_{_1}\!\ast Y\!\!_{_1}\!)
\ \! {{{_{\hbox{\fivebf Z}}}}\atop{\raise3.5pt\hbox{$\cong$}}}\ \!
{\!\rlap{$_{_{_{i\!+\!j=q}}}$}} {\ \raise0.2pt\hbox{$\oplus$}} \
\hbox{\tenbf(}\ \!\widetilde{\!\!\hbox{\tenit H}}_{i}
           (X\!_{_1}\!)\!
\otimes\!\!_{_{\hbox{\fivebf Z}}}\ \!\widetilde{\!\!\hbox{\tenit
H}}_{\!j} (Y\!\!_{_1}\!)\hbox{\tenbf)}\
{\rlap{$\!\!_{_{_{i\!+\!j=q-\!1}}}$}{\ \ \!
\raise0.2pt\hbox{$\oplus$}}}\ \ {\hbox{\tenrm
Tor}}_1^{\hbox{\fivebf Z}}\hbox{\tenbf(}\
\!\widetilde{\!\!\hbox{\tenit
H}}_{i}(X\!_{_1}\!),\widetilde{\tenrm H}_{\!j}
      (Y\!\!_{_1}\!)\hbox{\tenbf)}
$
i.e. the
$``\!X\!{_{^{_{2}}}}\!\!\!=\!\!Y{_{^{_{\!\!2}}}}\!\!\!=\!\emptyset"\!$-case
in our Th.\ 4 p.\ 10.
Milnor's results, apparently, inspired\ \  G.W. Whitehead to
introduce the {\sl Augmental Total Chain Complex}\ \
${\raise2pt\hbox{$\widetilde{}$}}\!\!{\hbox{\tenbf S}}(\circ)$ and
{\sl Augmental Homology},
$\!\phantom{\ |}
{\raise2pt\hbox{$\widetilde{}$}}\!_{\!}{\hbox{{\hbox{$\!\!$\it
H}}}} _{^{_{\star\!}}}(\circ),$ in \cite{29}. 
%
This was an attempt
of one of the most prominent topologists of our time, to, within
the classical frame, extend Milnor's formula to \hbox{topological
pair joins.}

G.W. Whitehead gave the empty space, $\!\emptyset,$ the status of
a (-1)-dimensional standard simplex but, in his pair space theory
he never took into account that $\emptyset$ then would get the
identity map, Id$_{_{^{\emptyset^{\!}}}}$, as a generator for its
(-1)-dimensional singular augmental chain\nobreak\ group, which,
correctly interpreted, actually makes his pair space theory
identical to the ordinary relative homology functor.
$\emptyset$ plays a definite role\break in the Eilenberg-Steenrod
formalism, cf. \cite{7} %
p. 3-4, while the ``convention"
$\!\hbox{\tenrm H}_{{{\hbox{\seveni i}}}}(\cdot):=
%
\hbox{\tenrm H}_{{{\hbox{\seveni
i}}}}(\cdot,\emptyset;\hbox{\tenbf Z})$,
cp. \cite{7} %
pp.\ 3 + 273,
is more than a mere ``convension" in that it connects the single
and the pair space theory and both concepts \hbox{should be
handled with care.}

In contemporary combinatorics there is a (-1)-dimensional simplex
$\emptyset$ containing no vertices. A moment of reflexion on the
realization functor, $|\cdot|$, reveals the need for a new
topological join unit, $\{\wp\}=|\{\emptyset\}|\neq|\emptyset| =
\emptyset$. We're obviously dealing with a non-classical
situation, which one mustn't try to \hbox{squeeze into a classical
framework.}

The new simplex, $\emptyset$, is contained in every non-empty
simplicial complex and induces on its own, a new (-1)-dimensional
simplicial sphere $\{\emptyset\}.$ The only homology-apparatus
invented by the combinatorialists to handle their ingenious
category modification were a jargon like -``We'll use reduced
homology with
$\widetilde{\hbox{H}}\!\!\!_{_{^{-1}}}\!(\{\emptyset\})={\hbox{\tenbf
Z}}".$

\indent Classical Simplicial and Singular Homology are accompanied
by {\it Reduced} {\it Homology} {\it Functors} in a mishmash that
severely cripples Algebraic Topology, e.g. it leaves the reduced
functor without influence on the boundary definition w.r.t.
manifolds.

To construct a (unifying) homology theory one starts with a
uniquely defined category of {\it admissible sets} (of pairs) and
three naive concepts; {\it homotopy}, {\it excision} and {\it
point}. Now, it's a matter of making these concepts comply with
the formalism in \cite{7} %
p.\ 114-118.
It's basic knowledge that no two categories are equipped with the
same homology theory, since the source category is a part of the
definition as is the domain in a function definition.
So, our $\mdoubleH$-functor(s) doesn't induce just another
homology theory, it's the \underbar{first} one constructed for the
\hbox{algebraically modernized categories.}

Generalizing the boundary concept from triangulable homological
manifolds to that of any simplicial quasi-$n$-manifold $\Sigma$,
see pp.\ 11 + 23, the boundary is made up by all simplices having
a ``link" with zero reduced homology in top dimension i.e.
$\hbox{\tenrm Bd}_{_{^{\!\ }}}\!\Sigma:=$
$\{\sigma\in \Sigma |\
\widetilde{\hbox{H}}\!\!\!\!\!\!_{_{^{n-\#\sigma}}}\!(\hbox{\tenrm
Lk}_{_{^{\!\Sigma}}}\!\sigma)=0\},$
where
${\hbox{\tenrm Lk}}_{_{^{\Sigma}}}\!\sigma:= \{\tau\in \Sigma|
[\sigma\cap \tau =\emptyset]\land [\sigma\cup \tau \in
\Sigma_o]\}.$
Since $\hbox{\tenrm Lk}_{_{^{\!\Sigma}}}\!\emptyset=\Sigma$,
$\hbox{\tenrm Bd}_{_{^{\!\ }}}\!{\hbox{\tensy P}}\!_{_{^{{\tenbf \
}}}}^{^{_{\hbox{\fivebf \ 2}}}}\!\!\!=\{\emptyset\}\ne \emptyset,$
for the real projective plane ${\hbox{\tensy P}}\!\!_{_{^{{\tenbf
\ }}}}^{^{_{\tenbf \ 2}}}\!.$
Also any point $\bullet:=\{\{{\hbox{\tenrm v}}\},\emptyset\}$ has
the join-unit $\{\emptyset\}$ as boundary, i.e. the boundary of
any $0$-ball is the $-1$-sphere. $\bullet$ is the only finite
\underbar{orientable} manifold having $\{\emptyset\}$ as its
boundary.

The splitting of homology into a reduced and a relative part,
really jams up algebraic topology.
Indeed, classically, it's difficult even to see that the boundary
of ``the cone of a M\"obius band" is the real projective plane,
cf. Prop.\ 1 \hbox{p.\ 11 + Ex.\ 2 p.\nobreak\ 28.}
Moreover, classically it's a hard-motivated truth that
$\hbox{\tenrm Bd}\!\!\lower1.3pt\hbox{$_{_{\hbox{\fivebf
Z}_{{\!\hbox{\fivebf p}}}}}$}\!\!_{\!\!} (\hbox{\tensy
P}\raise4pt\hbox{{\fiverm 2q}}_{\!}\ast\!\hbox{\tensy
P}\raise4pt\hbox{{\fiverm 2r}}_{\!})
\!= \!
\hbox{\tensy P}\raise4pt\hbox{{\fiverm 2q}}_{\!}
\cup
\hbox{\tensy P}\raise4pt\hbox{{\fiverm 2r}}{\!} $
with
$\hbox{\tensy P}\raise4pt\hbox{{\fiverm n}}$ the projective
$n$-space and
$\!{{\hbox{\bf Z}_{_{\!\hbox{\fivebf p}}}}}\!\!$ the prime-number
field modulo $\hbox{\tenbf p}\neq2$, cf. Ex.\nobreak\ 3\nobreak\
p.\nobreak\ 28.


\head  2. Augmental Homology Theory
\endhead
\subhead {\tenrm 2.1} Notations and Definition of Underlying
Categories
\endsubhead
\indent\indent
The typical morphisms in the classical category {\bf $\hbox{\tensy
K}$} of simplicial complexes with vertices in {\bf W} are the
simplicial maps as defined in \cite{25} p.\ 109, implying in
particular that;
$\hbox{\tenbf Mor}\raise1pt\hbox{$_{_{\!\!_{\hbox{\fivesy
K}}}}$}\!\!(\emptyset,_{\!}{\Sigma})\!=\!
\{\emptyset\}\!=\!\{0_{\emptyset^{\!},{\Sigma}}\},$
$\hbox{\tenbf Mor}_{_{\!\hbox{\fivesy
K}}}\!({\Sigma},\emptyset)=\emptyset\ \hbox{\tenrm if}\
\Sigma\neq\emptyset\ \ \hbox{\tenrm and}\ \hbox{\tenbf
Mor}_{_{\!\hbox{\fivesy K}}}\!(\emptyset,{\emptyset})=
\{\emptyset\}=\{0_{\emptyset,{\emptyset}}\}=\{\hbox{\tenrm
id}_{\emptyset}\}$ where $0_{_{\!{\Sigma,\Sigma^{\prime}}}}\!
=\!\emptyset = \hbox{\tenrm the\ empty\ function}$ $\hbox{\tenrm
from}\ \Sigma\ \hbox{\tenrm to}\ \Sigma^{\prime}.$ So; $
0_{_{\!{\Sigma,\Sigma^{\prime}}}}\!\in\hbox{\tenbf
Mor}_{_{\!\hbox{\fivesy K}}} (\Sigma,\Sigma^\prime)
\Longleftrightarrow \Sigma=\emptyset. $\ \ If in a category
\hbox{{$\varphi_i\in \hbox{\tenbf Mor} ({ R_i,S_i}),\ i=1,2,$} we
put;}

\medskip
\noindent
$\varphi_1 \sqcup \varphi_2: R_1\sqcup R_2 \longrightarrow
S_1\sqcup S_2 :r\mapsto$
$\cases
{\varphi_1(r)} {if \hbox{$r\in R_1$}} \cr
{\varphi_2(r)}
{if \hbox{$r\in R_2$}} \cr
\endcases
$
where $ \sqcup :=$ ``disjoint union".

\definition{Definition}
\hbox{\tenrm (of the objects in  $\hbox{\tensy K}\!_{_{^{o}}}\!$)}
An $($abstract$)$ simplicial complex $\Sigma$ on a vertex set
${V_{\Sigma_{}}}$ is a collection $($empty or non-empty$)$ of
finite $($\hbox{\tenbf empty} or non-empty$)$ subsets $\sigma$  of
${V_{\Sigma_{}}}$ satisfying;
(\hbox{\tenbf a}) If $\hbox{\tenrm v}\in {V_{\Sigma_{}}}$, then
$\{\hbox{\tenrm v}\}\in \Sigma$.
(\hbox{\tenbf b}) If $\sigma\in \Sigma$ and $\tau\subset \sigma$
then $\tau\in \Sigma$.\enddefinition

So, $\!\{\emptyset\}$ is allowed as an object in \hbox{\tenbf
$\hbox{\tensy K}\!_{_{^{o}}}$}.
We will write
``concept\vbox{\moveleft1pt\hbox{${_{^{_{o}}}}\!\!$}}"
or \ ``concept\vbox{\moveleft1pt\hbox{${_{^{_{\wp}}}}\!\!$}}" when
we want to stress that a concept is to be related to our modified
categories.

If $ \vert\sigma\vert\!=\#\sigma\!:=
\hbox{\tenrm card}(_{\!}\sigma{_{^{\!\!_{\ }}}}\!_{\!})\!=\!q\
\!\!${\eightrm+}$1$ then $\dim\sigma{_{^{\!\!_{\ }}}}\!\!:=\!q$
and $\sigma{_{^{\!\!_{\ }}}}$ is said to be a
$q$-face${_{^{\!_o}}}\!$ or a
$q$-simplex${_{^{\!_o}}}$ of $\Sigma{_{^{\!_o}}}$ and
$\dim\Sigma{_{^{\!_o}}}\!\!:=$sup$\{\dim (\sigma{_{^{\!\!_{\
}}}}\!)\vert \sigma{_{^{\!\!_{\ }}}}\!\!\in\!
\Sigma{_{^{\!_o}}}\}$.
Writing\ $\emptyset{_{^{\!_o}}}\!$ when using $\emptyset$ as a
simplex, we get dim$(\emptyset)\!=\!-\infty$ and
$\dim(\{\emptyset{_{^{\!_o}}}\!\})\!=\!\dim(\emptyset{_{^{\!_o}}}\!)\!=\!-1.$

\medskip
\noindent \hbox{\tenbf Note} that 
each object in the category
\hbox{\tenbf $\hbox{\tensy K}\!_{_{^{o}}}\!$} of simplicial
complexes$_{_{^{\!o}}}\!$ except for $\emptyset$, includes
$\{\emptyset\!_{_{^{o}}}\!\}$ as a\nobreak\ subcomplex. A typical
object$_{_{^{\!o}}}\!$ in \hbox{\tenbf $\hbox{\tensy
K}\!_{_{^{o}}}$} is \hbox{\tenbf ${\Sigma} \sqcup
\{\emptyset\!_{_{^{o}}}\!\}$} or $\emptyset$ where
${\Sigma}\in\hbox{\tensy K}$ and {$\psi$} is a morphism in
\hbox{\tenbf $\hbox{\tensy K}\!_{_{^{o}}}$} if;

\smallskip
$ (\hbox{\tenrm a})\ \ \psi=\varphi \sqcup \hbox{\tenrm
id}_{\{\emptyset\}}\ \hbox{\tenrm for\ some}\ \varphi\in
\hbox{\tenbf Mor}_{_{\!\hbox{\fivesy K}}}(\Sigma,\Sigma^\prime)\ \
\ \hbox{\tenrm or} $

\smallskip
$ (\hbox{\tenrm b})\ \ \psi= 0_{\emptyset,\Sigma_o}.\
 \hbox{\tenrm (In\ particular},\
\hbox{\tenbf Mor}_{_{\!\hbox{\fivesy
K}\!_{_{^{o}}}}}(\Sigma_o,\{\emptyset\})=\emptyset $ \hbox{\tenrm
if} and only if $ \Sigma_o\neq\{\emptyset\},\emptyset.) $

\medskip
{A functor \hbox{\tenbf E}:\ \hbox{\tenbf $\hbox{\tensy
K}_o$}$\longrightarrow$\hbox{\tenbf $\hbox{\tensy K}$}:}

\smallskip
\noindent Set \hbox{\tenbf E$(\Sigma_o)=
 \Sigma_o \setminus \{\emptyset\}\in Obj(\hbox{\tensy K})$}
and given a morphism  ${\rlap{$\psi$}{\ \ :}}\ {\Sigma_o}
\rightarrow {\Sigma_o^\prime}$ we put;

\smallskip
\item{}
\hbox{\tenbf E$(\psi)=\varphi $} if $\psi$ fulfills \hbox{\tenrm
(a)} above and
\item{}
\hbox{\tenbf E$(\psi)= 0_{\emptyset,\hbox{\tenbf E}(\Sigma_o)}$}
if $\psi$ fulfills \hbox{\tenrm (b)} above.

\medskip
{A functor \hbox{\tenbf E$_o$}: $\hbox{\tenbf \hbox{\tensy
K}}\longrightarrow$\hbox{\tenbf $\hbox{\tensy K}_o$}:}

\smallskip
\noindent Set\ \hbox{\tenbf E$_o(\Sigma){\rlap{=} {\ \ \Sigma
\sqcup \{\emptyset\}}}{\rlap{$\in$}{\ \ Obj(\hbox{\tensy K}_o)}}$}
and given simplicial ${\rlap{$\varphi$}{\ \ :}}{\Sigma}
\rightarrow {\Sigma^\prime},\ \hbox{\tenrm putting}\ \psi:=\varphi
\sqcup \hbox{\tenrm id}_{\{\emptyset\}}$, gives;

\smallskip
\item {---} \hbox{\tenbf EE$_o=$id$_{\hbox{\fivesy K}}$}
\smallskip
\item {---} \hbox{\tenbf imE$_o= Obj(\hbox{\tensy K}_o)\setminus \{\emptyset\}$}
\smallskip
\item { ---}   $\hbox{\tenbf E}_o\hbox{\tenbf E}=\hbox{\tenrm id}_{\hbox{\fivesy K}}$ except for
                $\hbox{\tenbf E}_o\hbox{\tenbf E}(\emptyset)=\{\emptyset \}$

\medskip
Similarly, let \hbox{\tenbf $\hbox{\tensy C}$} be the category of
topological spaces and continuous maps. $\!$Consider\nobreak\
the\nobreak\ category \hbox{\tenbf  $\hbox{\tensy
D}\!_{_{^{\!\wp}}}\! $} with objects: $\emptyset$ together with {$
{ X} \!_{_{^{\wp}}}\!\! :=\! { X}\ \!\hbox{\tenbf +}\
\!\{\wp\}$},$ \ \hbox{\tenrm for\ all}\ { X}\in \hbox{\tenbf
\hbox{\tenrm Obj}(\hbox{\tensy C})},$ i.e. the set {$ {
X}\!_{_{^{\wp}}}\!\!:=\! { X} \sqcup \{\wp\}$} equipped with the
weak topology, $\tau\!\!_{_{X}\!_{_{^{\!\wp}}}\!\!}$, with respect
to $_{^{\!}}{ X}_{^{\!}}$ and $\{\wp\}$,
cf. \cite{6} %
Def.\ 8.4\nobreak\ p.\nobreak\ 132.

\medskip
\noindent $f\!\!_{_{^{\wp}}}\!\! \in \hbox{\tenbf
Mor}\!_{_{^{\hbox{\fivesy D}\!\!_{_{^{\wp}}}\!\! }}}({
X}\!_{_{^{\wp}}}\! ,{ Y}\!\!\!_{_{^{\wp}}}\! )$ if; $\
f\!\!_{_{^{\wp}}}\!\! =
\cases  \hbox{\hbox{\tenrm a)}}\ \ f\hbox{\tenbf +} \hbox{\tenrm
id}\!_{_{\{\!\wp\!\}}}\ (:=f\sqcup\hbox{\tenrm
id}\!_{_{\{\!\wp\!\}}})\ \hbox{\tenrm with}\ f\in \hbox{\tenbf
Mor}\!_{_{\hbox{\fivesy C}}}\!(X,Y)\ \hbox{\tenrm and}\cr
f\ \hbox{\tenrm is}\ \underline{{\hbox{\tenrm on}\ \!X}}\
\hbox{\tenrm to}\ { Y},\ \hbox{\tenrm i.e}\ \hbox{\tenrm the\
domain\ of}\ f\hbox{\tenrm \ is\ the}\ {\tenrm whole\ of}\ \!X\!\cr
\hbox{\tenrm and}\
{ X}\!_{_{^{\wp}}}\!\! ={ X}\hbox{\tenbf +} \{\wp\},\ {
Y}\!\!\!_{_{^{\wp}}}\!\! ={ Y} \hbox{\tenbf +} \{\wp\}\
\hbox{\tenrm \ \hbox{\tenbf or}}
\cr
\hbox{\hbox{\tenrm b)}}\ \ 0_{{{{\emptyset,Y\!\!_{_{^{\!\wp}}}\!\!
}}}} \ \ (=\emptyset = \hbox{\tenrm the\ empty\ function\ from}\
\emptyset\ \hbox{\tenrm to}\ { Y}\!\!_{_{^{\!\wp}}}\! ).\cr
\endcases
$

\medskip
\noindent $\!\!\hbox{\tenrm There\ are\ functors}\indent
\hbox{\tensy F}\!\!_{_{^{\!\wp}}}\!:\!\Big\{
\matrix
\hbox{\teni \hbox{\tenbf \hbox{\tensy C}}} \longrightarrow
\hbox{\tensy D}\!_{_{^{\!\wp}}}\hfill\cr
{ X}\mapsto\!{ X}\hbox{\tenbf +} \{\wp\}
 \cr
\endmatrix,
\  \hbox{\tensy F}\!:\!\Bigg\{
\matrix
\hbox{\tensy D}\!_{_{^{\!\wp}}}\longrightarrow \hbox{\teni
\hbox{\tenbf \hbox{\tensy C}}}\hfill\cr
\!{X}\hbox{\tenbf +} \{\wp\} \mapsto{X}\cr
\emptyset\mapsto\emptyset \hfill \cr
\endmatrix
$
resembling $\hbox{\tenbf E}\!_{_{{o}}}\!$ resp. $\hbox{\tenbf E}
$.

\smallskip
\noindent \hbox{\tenbf Note.} The ``$\hbox{\tensy
F}\!\!_{_{^{\wp}}}$-lift topologies"$\!$,
$$\tau_{_{\!\!{X_{_{^{\!{o}}}}}}}\!\!\!:=\! \tau_{_{\!\!{X}}}\cup
\{\!{X}_{_{^{\!{o\!}}}}\!\}\!= \!\{\hbox{\tensy O}_{_{^{\!{o\!}}}}
| \hbox{\tensy O}_{_{^{\!{o\!}}}}\!=
\!\!{X}_{_{^{\!{o\!}}}}\setminus(\hbox{\tensy N}\sqcup\{\wp\})\
\!;\ \hbox{\tensy N}\ \hbox{\tenrm closed\ in}\ {X}
\}\cup\nobreak\{\!{X_{_{^{\!{o\!}}}}}\!\}\ \ $$
\hbox{\tenrm and}
$$ \tau_{_{\!\!{X_{_{^{\!{o}}}}}}}\!\!\!:=\! \hbox{\tenbf
\hbox{\tensy F}}\!\!_{_{^{\wp}}}\!(\tau_{_{\!\!{X}}})_{\!}
\cup_{\!} \{\emptyset_{}\} \!=\! \{ \hbox{\tenbf \hbox{\tensy
O}}\!_o\!\!=\!\hbox{\tenbf \hbox{\tensy
O}}_{\!}\sqcup_{\!}\{\wp\}\ \!|\ \! \hbox{\tenbf \hbox{\tensy
O}}\!\in\!\tau_{_{\!\!{X}}}
\}_{\!}\cup_{\!} \{\emptyset_{}\} $$
would also \hbox{\tenrm give} $_{\!}\hbox{\tensy
D}\!\!_{_{^{\wp}}}\!$ due to the domain restriction in
$\hbox{\tenbf a},$
making $_{\!}\hbox{\tensy D}\!\!_{_{^{\wp}}}\!$ a link between the
two constructions of partial maps treated
in \cite{2} %
pp.\ 184-6.
No {extra} morphisms$ _{_{^{\!\wp}}}\! $\nobreak\ has
been allowed into $\hbox{\tenbf \hbox{\tensy
D}}\!_{_{^{\!\wp}}}\!\!$ $(\hbox{\tensy K}_{\!o ^{\!}} )$ in the
sense that the morphisms$_{_{^{\!\wp}}}\!$ are all targets under
$\hbox{\tensy F}\!\!_{_{^{\!\wp}}}$
$
\!\!\
(\hbox{\tenbf E}_{_{^{\!o}}}\! )$
except $0_{{{{\emptyset,Y\!\!_{_{^{\!\wp}}}\!\!   }}}} $ defined
through item b, re-establishing $\emptyset${\spaceskip2.0pt\ as
the unique initial object.}

The underlying principle for our definitions is that a concept in
$\hbox{\tensy C}\ (\hbox{\tensy K})$ is carried over to
$\hbox{\tensy D}_{\!\wp}\ (\hbox{\tensy K}_o)$ by $\hbox{\tensy
F}\!\!_{\wp}\ (\hbox{\tenbf E}_o)$ with addition of definitions of
the concept$_{\!o}$ for cases that isn't a proper image under $
\hbox{\tensy F}\!\!_{\wp}\ (\hbox{\tenbf E}_o)$. The definitions
of the product/join operations ``$\times\!_o"\!$,
$``\ast\!\!_{{^o}}"\!$, $``{\rlap{{\lower2.5pt\hbox
{\vbox{\moveright0.1pt\hbox{$^{^{\land}}$}}}}}{{
\ast}}}\!_{{^o}}"$ in page\ 7 and ``$\setminus\!_{_{^o}}{_{\!}}"$
below, certainly follows this principle.

\definition{Definition}
{${^{\hbox{$\emptyset$}}}\ \!{{\!\! /\!{^{^{_{o}}}}}}\!_{{\!\!
X}\!_{_{\!\wp}}} \!\!_{\!}:=\emptyset $ else; $ { ^{
X\!_{_{\!\wp\!1}}}}{{\!\! /\!{^{^{_{o}}}}}}\!_{{\!\!
X}\!_{_{\!\wp\!2}}}\!\!_{\!}:= {\hbox{\tensy F}}\!_\wp\bigl(
^{{\hbox{\tensy F}}( X_{\wp1})}\!\!{\rlap{/}\!{\ }}
 _{{\hbox{\tensy F}}(X_{\wp2})}
\bigr) $ if $ X\!\!_{_{\wp2}}\!\!\neq \!\emptyset$} in
${\hbox{\tensy D}}\!_{\!o\wp}$ for $
X\!_{\wp2}\!\subset\!X\!_{\!\wp1}.$

\smallskip
\noindent
where ``$/$" is the classical ``quotient" except that $
^{{\hbox{\tensy F}}( X_{\wp1})}\!{\rlap{/}{\ }}
 _{\emptyset}
\!:=\! {\hbox{\tenbf \hbox{\tensy F}}( X_{\wp1})},$ cp. \cite{2} %
p.\ 102.
\enddefinition

\definition{Definition}
$X_{\!\wp1}+\!\!_{_{^{o}}}X_{\!\wp2}:=
\cases
\!X_{\!\wp1} & \!if\
X_{\!\wp2}=\emptyset\
(or\ vice\ versa)\cr
\!{\hbox{\tensy F}\!_{\wp}}({{\hbox{\tensy F}}(X_{\!\wp1})+
{\hbox{\tensy F}}(X_{\!\wp2})}) & \!if\
X_{\!\wp1}\neq \emptyset \neq X_{\!\wp2}
\cr
\endcases
$

\smallskip
\noindent
where ``$+$" is the classical ``topological sum"$\!$,
defined in \cite{6} %
p.\ $\!$127 as the ``free union"$\!$.
\enddefinition

Proposition 1 p.\ 11 is our key motivation for introducing a
topological (-1)-object, which then imposed the following
definition of a ``setminus", ``$\setminus{_{_{^{\!o\!}}}}$", in
${\hbox{\tensy D}\!_{\wp}}$.

\definition{Definition}
${ X}_{\!\wp}\setminus\!{_{_{^{\!o\!}}}}{X}_{\!\wp}^{^\prime}:=
\cases
\emptyset & \hbox{\tenrm if}\ { X}_{\!\wp}\!=\emptyset,\ {
X}_{\!\wp}\raise2pt\hbox{$\mmysubsetneqq$} { X}_{\!\wp}^\prime\
\hbox{\tenrm or}\ { X}_{\!\wp}^\prime=\{{\wp}\} \cr
{\hbox{\tensy F}\!_{\wp}}\bigl({\hbox{\tensy F}}({X_{\!\wp}})\ \!
{\raise1.5pt\hbox{\eightmsbm \char"72}}\ \!
{\hbox{\tensy F}}({X}_{\!\wp}^\prime)\bigr) & \hbox{\tenrm else}.\
_{\!}
\hbox{\tenrm(``{\raise1.5pt\hbox{\eightmsbm \char"72}}"$\!$:=
classical ``setminus"$\!.)\!$}%
\cr%
\endcases
$
\enddefinition

\noindent\hbox{\tenbf Notations:} We have used w.r.t.:=with
respect to, and $\tau\!\!_{_{X}}$:=the topology of $X$.$\!$ We'll
also use; $\hbox{\tenbf PID}\!:=\!$ Principal Ideal Domain,
l.h.s.(r.h.s.) $\!:=\!$ left (right) hand side,\break
iff $\!:=\!$ if and only if, cp.$:=\!$ compare$_{\!}$ (cf.$:=\!$
cp.!), \hbox{\tenbf LHS} $\!:=\!$ Long $\mdoubleH$omology Sequence
and {\tenbf M-$\!$Vs} $\!:=\!$ Mayer-$\!$Vietoris sequence.
Let, here in Ch. 2, $\Delta =
\{_{\!}{\Delta}_{\!}\!^{^{_{_{\!}p}}}\!,_{\!}\partial\}$ be the
classical singular chain complex and let $\!$``$\simeq$"$\!$
denote ``homeomorphism" or ``chain\nobreak\ isomorphism".
$\bullet$ ($\bullet\bullet$) denotes the one (two) point space
$\{\bullet,\wp\}$ ($\{\bullet\bullet,\wp\}$).

\smallskip
If $X\!_{_{^{\wp}}}\!\! \neq\emptyset,\ \{\wp\}$ then
($X\!_{_{^{\wp}}}\!,\tau\!\!_{_{X\!_{_{^{\!\wp}}}\!\!}}$) is a
non-connected  space, and it therefore seems adequate to define ${
X}\!_{_{^{\wp}}}\!\!\neq \emptyset,\ \{\wp\}$ to have a certain
point set topological ``property$\!_{_{^{\wp}}}\!\!$" if
$\hbox{\tenbf \hbox{\tensy F}}({ X}\!_{_{^{\wp}}}\!\!)$ has the
``property" in question, e.g. ${ X} \!_{_{^{\wp}}}\! $ is
connected$\!_{_{^{\wp}}}\!$ \underbar{iff} {X} is connected.

\subhead {\rm 2.2} Simplicial Augmental Homology Theory and
realizations\endsubhead

\smallskip%
\indent%
$\mdoubleH$ denotes the simplicial as well as the singular
augmental (co)homology functor$\!_o.\!$

\smallskip
Choose oriented  $q$-simplices to generate
$C^{o}\!\!_{_{\!\!^q}}(\Sigma_o;\hbox{\tenbf G})$, where the
coefficient module $\hbox{\tenbf G}$ is a unital
($\leftrightarrow\!1\!_{^{_{\hbox{\fivebf A}}}}\!\circ\!g=g$)
module over any commutative ring \hbox{\tenbf A} with unit.
Now;

\smallskip
\item{-}
$\hbox{\tenbf C}^{o}(\emptyset;\hbox{\tenbf G}$) is identically 0
in all dimensions, where 0 is the additive unit-element.

\item{-}
{$\hbox{\tenbf C}^{o}(\{\emptyset_o\};\!\hbox{\tenbf G}$) is
identically 0 in all dimensions except for
$ C^{o}\!\!\!\!_{_{-1}}\!(\{\emptyset_o\};\!\hbox{\tenbf G})\ \!
\cong\ \!\hbox{\tenbf G}.$}

\noindent
\item{-}
$\hbox{\tenbf C}^{o}(\Sigma_o;\hbox{\tenbf G})\equiv
\widetilde{\hbox{\tenbf C}}(\hbox{\tenbf E}(\Sigma_o);\hbox{\tenbf
G})\equiv$ the classical
``$\{\emptyset\!_{_{^{o\!}}}\!\}$"-augmented chain.


\medskip By just hanging on to the
``$\{\emptyset\!_{_{^{o\!}}}\!\}$-augmented chains", also when
defining relative chains$_o$, we get the {\it Relative} {\it
Simplicial} {\it Augmental} {\it Homology} {\it Functor} {\it for}
$\hbox{\tensy K}_{_{^{o}}}\!$-{\teni pairs}, denoted
$\mdoubleH_{\ast}$ and fulfilling;

$$\mdoubleH_{i}(\Sigma_{o1},\Sigma_{o2};\hbox{\tenbf G})=
\cases
\hbox{\tenrm H}_i(\hbox{\tenbf E}(\Sigma_{o1}),\hbox{\tenbf
E}(\Sigma_{o2});\hbox{\tenbf G}) & if\ \hbox{$\Sigma_{o2}\neq
\emptyset$}\cr
\widetilde{\hbox{\tenrm H}}_i(\hbox{\tenbf
E}(\Sigma_{o1}),\hbox{\tenbf G}) & if\ \hbox{$\Sigma_{o1}\neq\ \!
\{\emptyset_o\}, \emptyset$,}\ and\
 \hbox{$\Sigma_{o2}=\ \emptyset$}\cr
$%
$\cases
\cong \hbox{\tenbf G}  & if\ \hbox{$i=-1$} \cr
=0 & if\ \hbox{$i\neq -1$}\cr
\endcases
& when\ \hbox{$\Sigma_{o1}\!=\!\{\emptyset_o\}$}\ \ and\ \
\hbox{$\!\!\Sigma_{o2}=\ \emptyset$}\cr
0 & \ for\ all\ $i$\ when\ \hbox{$\Sigma_{o1}=\Sigma_{o2}=\
\emptyset$}.\cr\endcases $$

\medskip
$\emptyset\!\ne\!\{\emptyset\!_{_{^{o\!}}}\!\}$ and both lacks
final sub-objects, which under any useful definition of the {\it
realization of a simplicial  complex} implies that
$|\emptyset|\ne|\{\emptyset\!_{_{^{o\!}}}\!\}|$ demanding the
addition of a non-final object
$\{\wp\}=|\{\emptyset_{_{^{\!\!o}}}\!\}|$ into the classical
category of topological spaces as {\it join-unit} and
(-1)-dimensional standard simplex.
This approach conforms Homology $\!$Theory $\!$and considerably
$\!$simplifies the study of manifolds, cf. p.\ 13.

We will use Spanier's definition of the ``function space
realization" ${ |\Sigma_{_{^{o}}} |}$ as given
in \cite{25} %
p. 110, unaltered, except for the `` $_{_{^{o}}}$"'s and the
underlined $\underline{\underline{addition}}$ where
$\wp:=\alpha_0$ and $\alpha_0(\hbox{\tenrm v})\equiv 0\ \forall\
\hbox{\tenrm v}\in V\!\!\!_{_{^{\Sigma}}}$:

\medskip
\indent{\hbox{\tenbf -}\lower1.0pt\hbox{\twelvbf ``}} We now
define a covariant functor from the category of simplicial
complexes$_{_{^{\!o}}}$ and simplicial maps$_{_{^{\!o}}}$ to the
category of topological spaces$_{_{^{\!o}}}$ and continuous
maps$_{_{^{\!o}}}$. Given a nonempty simplicial
complex$_{_{^{_{\!}o}}}$ ${\Sigma_{_{^{\!o}}} }$, let ${
|\Sigma_{_{^{\!o}}}|}$ be the set of all functions $\alpha$ from
the set of vertices of ${\Sigma_{_{^{\!o}}}}$ to $\hbox{\tenbf
I}:=[0,1]\ \hbox{\tenrm such\ that;}$

\item{(a)} For any $\alpha$, $\{\hbox{\tenrm v}\in V_{\Sigma_o } |
         \alpha(\hbox{\tenrm v})\neq 0\}$
          is a simplex$_o$ of ${\Sigma_o }$ (in particular,
          $\alpha(\hbox{\tenrm v})\neq 0$ for only a
          finite set of vertices).
\item{(b)}
For any $\alpha\ \underline{\underline{\neq\alpha_0}}$,
$\sum_{_{_{_{\!\!\!\!\!\!\!\!\!{\hbox{\fiverm v}\in
V\!\!_{_{\Sigma_o}}}}}}}\!\!\!\!\!\! \alpha(\hbox{\tenrm v})=1.$\
\

\smallskip
\nobreak \noindent If  ${\Sigma_o }=\emptyset$ , we define
         ${ |\Sigma_o |}=\emptyset$.
\indent\indent\lower1.5pt\hbox{\twelvbf "}

\medskip
The {\it barycentric} {\it coordinates}$_o$ $\alpha$, defines a
metric
$$d(\alpha,\beta)\!=\!\! {\rlap {$ \sqrt {\sum
{\lower4.5pt\hbox{$^{  [\hbox{\teni {\char"0B}}(\hbox{\tenrm
v})-\hbox{\teni {\char"0C}}(\hbox{\tenrm v})]^{^{2}} }$}} } $} {
{\phantom{\ \ }_{_{_{_{_{_{{\hbox{\fiverm v}\in{\hbox{\fivei
V}}\!_{_{^{\Sigma_o}}} }}}}}}}} } }$$
on $|\Sigma_o|$ inducing the topological space $|\Sigma_o|_d$ with
{\it the} {\it metric} {\it topology}.
We'll equip $|\Sigma_o|$ with another topology and for this
purpose we define the {\it closed} {\it simplex$_o$} $|\sigma_o|$
of $\sigma_o\in{\Sigma_o }$ i.e.
$|\sigma_o|:=\{\alpha\in|\Sigma_o|\ |\ { [\alpha(\hbox{\tenrm
v})\neq0] \Longrightarrow [\hbox{\tenrm v}\in \sigma_o] \}.}$

\definition{Definition}
For ${\Sigma_{_{^{\!o}}} }\ne\emptyset$,
$|\Sigma_{_{^{\!o}}}|$ is topologized through
${|\Sigma_{_{^{\!o\!}}}|}\!:=|\hbox{\tenbf
E}(\Sigma_{_{^{\!o\!}}})|+\{\alpha_0\}$, which is equivalent to
give ${|\Sigma_{_{^{\!o}}}|}$ the weak topology w.r.t. the
$|\sigma\!_{_{^{\!o\!}}}|$'s, naturally imbedded in $\hbox{\tenbf
R}\!^{^{_{n}}}\!\hbox{\tenbf +}\{\wp\}$ and we define
$\Sigma_{_{^{\!o}}}$ to be connected if $|\Sigma_{_{^{\!o\!}}}|$
is, i.e. if $\hbox{\tensy
F}(\vert\Sigma_{_{^{\!o\!}}}\vert)\simeq\vert\hbox{\tenbf E}
(\Sigma_{_{^{\!o\!}}})\vert\ \hbox{\tenrm is}.$
\enddefinition

\proclaim{Proposition}
$|{_{\!}}\Sigma_{{{\!o\!}}}|$ is always homotopy
equivalent to
$ |{_{\!}}\Sigma_{{{\!o\!}}}|_{_{^{\!d}}}$
\hbox{\tenrm(\cite{9} %
pp.\ 115, 226.)} and \indent
$|{_{\!}}\Sigma_{{{\!o\!}}}|$ is
homeomorphic to
$ |{_{\!}}\Sigma_{{{\!o\!}}}|_{_{^{\!d}}}$ \underbar{iff}
$\Sigma_o$ is
locally finite {\rm(\cite{25} %
p.\ 119 Th.\ 8.)}.
\qed
\endproclaim

\goodbreak%

\subhead {\rm 2.3} Singular Augmental Homology Theory\endsubhead

\medskip
$|\sigma\!_{_{^{\!\wp\!}}}|$ imbedded in $\hbox{\tenbf
R}\!^{^{_{n}}}\!\hbox{\tenbf +}\{\wp\}$ generates a satisfying set
of\nobreak\ ``standard simplices$ \!_{_{^{\wp}}}\! $" and
``singular simplices$\!_{_{^{\wp}}}\!"\!.$ This implies in
particular that the ``$p$-standard simplices$\!_{_{^{\wp}}}\!\!$",
denoted $\Delta\!^{{\wp}p}$, are defined by
$\Delta\!^{{\wp}p}\!:=\Delta\!^p\hbox{\tenbf +}\ \!\{\wp\}$ where
$\Delta\!^p$ denotes the usual $p$-dimensional standard simplex
and $\hbox{\tenbf +}$ is the topological sum, i.e.
$\Delta\!^{{\wp}p}:=\Delta\!^p\sqcup\{\wp\}$ with the weak
topology w.r.t. $\Delta\!^p$ and $\{\wp\}$.
Now, and most important:
$\Delta\!^{^{_{\wp(\!-1\!)}}}\!\!:=\{\wp\}$.

Let $T^p$ denote an arbitrary classical singular $p$-simplex
($p\geq 0$). The ``singular $p$-simplex$\!_{_{^{\wp}}}\! $",
denoted $\sigma\!\lower3pt\hbox{$^{^{{\ \!\!\wp}\ \!\!\!p}}$}$,
now stands for a function of the following kind:
$$%
\sigma\!\lower3pt\hbox{$^{^{{\ \!\!\wp}\
\!\!\!p}}$}: \Delta\!^{{\wp}p}=\Delta\!^p\hbox{\tenbf
+}\{\wp\}\longrightarrow X\hbox{\tenbf +} \{\wp\} \ \ \hbox{\tenrm
where}\ \ \sigma\!\lower3pt\hbox{$^{^{{\ \!\!\wp}\
\!\!\!p}}$}(\wp)=\wp\ \ \hbox{\tenrm and}
$$%
\smallskip
\centerline{$ \sigma\!\lower3pt\hbox{$^{^{{\ \!\!\wp}\ \!\!\!p}}$}
\!\!\!_{_{|\!\Delta\!^p}} =T^p \ \hbox{\tenrm for\ some\
ordinary}\ p\ \!\hbox{-}\ \!\hbox{\tenrm dimensional\ singular\
simplex}\ T^p\ \hbox{\tenrm for\ all}\ p\geq 0.$ }

$$\hbox{\tenrm In\ particular;} \ \ \ \ \ \ \ \ \ \
\sigma\!\lower3pt\hbox{$^{^{{\ \!\!\wp}(\!-1\!)}}$}
:\{\wp\}\longrightarrow X\!_{_{^{\wp}}}\!= X\hbox{\tenbf +}
\{\wp\} :\wp \mapsto \wp.$$

\smallskip
The boundary function $\partial_{_{{^{^{\!}}}^\wp}}\!$ is defined
by $
\partial\!_{_{{^{^{\!}}}^{\wp_{_{\!}}p}}}\!
(\sigma\!\lower1pt\hbox{$^{^{{_{_{\!}}}_{\wp_{_{\!}}p}}}$}\!)$
$:=\hbox{\tensy F}\!\!_{_{^{\!\wp}}}\!(\partial\!_p(T^p)) \ \
\hbox{\tenrm if}\ p>0$ \noindent where $\partial\!_p$ is the
ordinary  singular boundary function, and $\partial
_{\wp0}(\sigma\!\lower3pt\hbox{$^{^{{\ \!\!\wp}\ \!\!\!0}}$})
\!\equiv\! \sigma\!\lower3pt\hbox{$^{^{{\ \!\!\wp}(\!-1\!)}}$}$
for every singular 0-simplex$\!_{_{^{\wp}}}\!$\
$\sigma\!\lower3pt\hbox{$^{^{{\ \!\!\wp}\ \!\!\!0}}$}\!\!.$
Let ${\Delta{\raise1.5pt\hbox{$\!\!^{\wp}$}}}\!=
\!\{{\Delta}_{\!}\!^{^{_{\wp_{\!}p}}}\!\!,\partial_{\wp}\}$ denote
the singular augmental chain complex$_{_{\!}\wp}.\ \hbox{\tenrm
Observation;}\ |\Sigma \!_{_{^{o}}}\!|
\neq\emptyset\Longrightarrow\! |\Sigma
\!_{_{^{o}}}\!|=\hbox{\tenbf \hbox{\tensy
F}}\!\!_{\alpha_0}(|\hbox{\tenbf E} (\Sigma
\!_{_{^{o}}}\!)|)\nobreak\in\hbox{\tenbf \hbox{\tensy
D}}\!_{\alpha_0}.$

\medskip
By the strong analogy to classical homology, we omit the proof of
the $\hbox{\tenrm next\ lemma.}\!$

\proclaim{Lemma} \hbox{\tenrm (Analogously for co$\
\!_{\!}\mdoubleH$omology.)}
$$\mdoubleH_{i}(X_{\wp1},X_{\wp2};\hbox{\tenbf G})\!=\!
\cases
\hbox{\tenrm H}_i(\hbox{\tenbf \hbox{\tensy
F}}(X_{\wp1}),\hbox{\tenbf \hbox{\tensy F}}(X_{\wp2});\hbox{\tenbf
G}) & if\ \hbox{$X_{\wp2}\neq \emptyset$}\cr
\widetilde{\hbox{\tenrm H}}_i(\hbox{\tenbf \hbox{\tensy
F}}(X_{\wp1});\hbox{\tenbf G})
& if\ X_{\wp1}\neq \{\wp\},\ \emptyset\ and\
\hbox{$X_{\wp2}=\emptyset$}\cr
$
$\!\!\cases
\cong\hbox{\tenbf G} & if\ \hbox{$i=-1$} \cr =0 & if\ \hbox{$i\neq
-1$}\cr
\endcases
$
$
&
when\ X_{\wp1}\!= \{\wp\} \ and\ \hbox{$X_{\wp2}=\emptyset$}\cr
0 &  for\ all\ i\ when\ X_{\wp1}\!=X_{\wp2}\!=\emptyset.\!\!\!
\qed \cr
\endcases
$$
\endproclaim

\remark{Note} \hbox{\tenrm i.}
${\Delta}_{\!}\!^{^{_{_{\!}i}}}\!(X_{1},X_{2};\hbox{\tenbf G})
\!\cong\!
{\Delta}_{\!}\!^{^{_{\wp_{\!}i\!}}}\!(\hbox{\tensy
F}\!\!_{\wp}(X_{1}),\ \!\!\hbox{\tensy
F}\!\!_{\wp}(X_{2});\hbox{\tenbf G})$ always. $\!$So,
\hbox{$\mdoubleH\!\!_{_{^{\lower1pt\hbox{-}1}\!}}\!(\hbox{\tensy
F}\!\!_{\wp}\!(X_{1\!}),\ \!\!\hbox{\tensy
F}\!\!_{\wp}\!(X_{2});\hbox{\tenbf G}) \!\equiv\!0.\!$}

\smallskip
\item{ii.}
 $
{\Delta{\raise1.5pt\hbox{$\!\!^{\wp}$}}}(X_{\wp1},X_{\wp2}\!)
\simeq {\Delta}(\hbox{\tensy F}(X_{\!\wp1\!}), \hbox{\tensy
F}(X_{\!\wp2}))$ except \underbar{iff}
$X_{\!\wp1\!}\neq_{\!}X_{\!\wp2\!} = \emptyset$ when the only
non-isomorphisms occur for \indent
$
{\Delta}{\!\!}^{^{_{\wp_{\!}\hbox{\fivebf($_{\!}$-1$_{\!}$)}}}}\!(
X_{\!\wp1\!},\emptyset)\cong \hbox{\tenbf Z}
\ncong
\hbox{\tenbf0} \cong
{\Delta}_{\!}\!^{^{_{_{\!}\hbox{\fivebf($_{\!}$-1$_{\!}$)}}}}\!
(\hbox{\tenbf \hbox{\tensy F}}(X_{\!\wp1\!}),\emptyset) $\indent
when
$ {\mdoubleH}_{_{^{_{\!}0}}}\!(X_{\wp1},\emptyset)
\oplus\hbox{\tenbf Z}
\cong
\hbox{\tenrm H}_{_{^{_{\!}0}}}\!(\hbox{\tensy
F}(X_{\wp1}),\emptyset) $ if $X_{\!\wp1\!}\neq\!\{\emptyset\}$.

\smallskip
\item{iii.}
$\hbox{\tenbf
C}^{o}\!(\Sigma_{_{^{\!o1}}},\Sigma_{_{^{\!o2}}};\hbox{\tenbf G})
\approx {
\Delta{\raise1.5pt\hbox{$\!\!^{\wp}$}}}\!(|\Sigma_{_{^{\!o1}}}|,|\Sigma_{_{^{\!o2}}}|;\hbox{\tenbf
G})$ connects the simplicial and singular functor$_{_{^{\!o}}}\!.$
 (``$\approx$" stands for ``chain equivalence".)

\smallskip
\normalbaselines
\item{iv.}
$
\mdoubleH_{0}(X_{_{\!}\wp}\raise1pt\hbox{\eightrm+}_{_{^{\!\!\!o}}}
Y_{\!\wp},\{\wp\};\hbox{\tenbf G}) =
\mdoubleH_{0}(X_{_{\!}\wp},\{\wp\};\hbox{\tenbf G})\oplus
\mdoubleH_{0}(Y_{\!\wp},\{\wp\};\hbox{\tenbf G}) $ but $
\mdoubleH_{0}(X_{_{\!}\wp}\raise1pt\hbox{\eightrm+}_{_{^{\!\!\!o}}}
Y_{\!\wp},\emptyset;\hbox{\tenbf G})
=\break
=\ \mdoubleH_{0}(X_{_{\!}\wp},\emptyset;\hbox{\tenbf G})\oplus\
\mdoubleH_{0}(Y_{\!\wp},\{\wp\};\hbox{\tenbf G})\ =\
\mdoubleH_{0}(X_{_{\!}\wp},\{\wp\};\hbox{\tenbf G})\oplus\
\mdoubleH_{0}(Y_{\!\wp},\emptyset;\hbox{\tenbf G}).$
\endremark

\definition{Definition}
The $p$:th {\it Singular} {\it Augmental Homology Group of}
$X_{_{\!^{\wp}}}\!\!$ w.r.t.\ $\!\hbox{\tenbf G}\!= \break
=\mdoubleH_{_{^{p\!\!}}}(\!X_{_{\!^{\wp}}}\!;\hbox{\tenbf
G})\!:=\mdoubleH_{_{^{p\!\!}}}(\!X_{_{\!^{\wp}}}\!,\emptyset;\hbox{\tenbf
G}).
$ The {\it Coefficient} {\it Group}$_{o}\!$
$:={\mdoubleH}_{_{\!\!^{-1}}}\!(_{\!}\{\wp\},\emptyset;\hbox{\tenbf
G}).$

Using \hbox{\tensy F}$\!_{\!\!\wp}\ (\hbox{\tenbf E}_{\!o\!}),$ we
``lift" the concepts of \hbox{\teni homotopy}, \hbox{\teni
excision} and \hbox{\teni point} in $\hbox{\tensy C}\
(\hbox{\tensy K})$ into $\hbox{\tenbf \hbox{\tensy
D}}\!_{\wp}$-concepts $(\hbox{\tensy K}_{\!o\!}$-concepts$)$
\hbox{\teni homotopy}$_{_{\!}{\wp}}$, \hbox{\teni
excision}$_{_{\!}{\wp}}$ and \hbox{\teni point}$_{_{\!}{\wp}}$
$(=:\ \!\!\bullet)$, respectively.

\enddefinition

\smallskip
So; $f_{o},g_{o}\in\hbox{\tenbf  \hbox{\tensy D}}\!_{\wp}$ are
\hbox{\teni homotopic}$_{\wp}$ $\hbox{\tenrm if\ and\ only\ if}$
$f_{o}=g_{o}= 0_{{{{\emptyset,Y\!\!_{_{^{o}}}}}}}$ or there are
homotopic maps $f_{1},g_{1}\in\hbox{\tenbf  \hbox{\tensy C}}$ such
that $f_{o}=f_{1}+\hbox{\tenrm Id}_{\tenrm {\rlap{$
$}{\{\wp}}{\rlap{{\rlap{$\ $}{\ }}}{\} }}},
g_{o}=g_{1}+\hbox{\tenrm Id}_{\tenrm {\rlap{$
$}{\{\wp}}{\rlap{{\rlap{$\ $}{\ }}}{\} }}}.$

\smallskip
An inclusion$_o$ $(i_o,i_{o_{A_o}}): (X_o\setminus_{_{\!o}}
\!U_o,A_o\setminus_{_{\!o}\!}U_o) \longrightarrow (X_o,A_o),\
U_o\neq\{\wp\}$, is an \hbox{\teni excision}$_\wp$ if and only if
there is an excision
$(i,i_{_{A}}):(X%
{\raise1.5pt\hbox{\eightmsbm \char"72}}\ U,A\ \!
{\raise1.5pt\hbox{\eightmsbm \char"72}}\ U)\longrightarrow (X,A)$
such that $i_{o}=i_{ }+ \hbox{\tenrm Id}_{\tenrm {\rlap{$\
$}{\{\wp}}{\rlap{{\rlap{$\ $}{\ }}}{\} }}}$ and
$i_{o_{A_o}}=i_{_{A}}+ \hbox{\tenrm Id}_{\tenrm {\rlap{$\
$}{\{\wp}}{\rlap{{\rlap{$\ $}{\ }}}{\} }}}.$

\smallskip
$\{\hbox{\tenbf P},\wp\}\in\hbox{\tenbf \hbox{\tensy
D}}\!\!_{_{^{\wp}}}\!$ is a \hbox{\teni point}$_{_{^{\!\wp}}}\!$
\underbar{iff} $\{\hbox{\tenbf P}\}+\{\wp\}=\hbox{\tensy
F}\!_{\wp}(\{\hbox{\tenbf P}\})$ and $\{\hbox{\tenbf
P}\}\!\in\hbox{\tensy C}$ is a \hbox{\teni point}.\break
So,
$\{\wp\}\ \hbox{is\ \hbox{\tenbf not}\ a\ point}_{_{^{\!\wp}}}\!.$

\remark{Conclusion} $(\mdoubleH,{ {\partial}}\!_{\wp\!})$,
abbreviated $\mdoubleH$, is a {\it homology theory on the
$h$-category of pairs from} $\hbox{\tenbf \hbox{\tensy
D}}\!_{\wp}$ $(\hbox{\tensy K}_{\!o\!})$,
c.f. \cite{7} %
p.\ 117, i.e. $\mdoubleH$ fulfills the $h$-category analogues,
given in {\cite{7} %
\S\S8-9\ pp.\ 114-118}, of the seven Eilenberg-Steenrod axioms
from  {\cite{7} %
\S3\ pp.\ 10-13}.
The necessary verifications are either equivalent to the classical
or completely trivial. E.g. the {\teni dimension} {\teni axiom} is
fulfilled since $\{\wp\}\ \hbox{is\ not\ a\
point}_{_{^{\!\wp}}}\!.$
\endremark

\medskip
Since the exactness of the relative Mayer-Vietoris sequence of a
proper triad, follows from the axioms, cf.
\cite{7} %
p.\ 43 and, paying proper attention to Note iv, we'll use it
without further motivation.

$ { \widetilde{\hbox{\tenrm H}}(\!X) \!=\! \mdoubleH(\hbox{\tensy
F}\!_{\!\!\wp}(X),\emptyset) } $ explains all the ad-hoc reasoning
surrounding the $\widetilde{\hbox{\tenrm H}}$-functor.

%


\head 3. Augmental Homology Modules for Products and Joins
\endhead
\subhead {3.1} Definitions of the Product and Join Operations
\endsubhead
Let
\lower1.0pt\hbox{$^{_{_{\hbox{\sevensy \char"72} }}}$}
be one of the classical topological product/join operations
``$\times"\!$\hbox{\tenbf /}$\!\!$\hbox{\tenbf /}
$``\ast"\!$\hbox{\tenbf /}$\!\!$\hbox{\tenbf /}
$``{\rlap{{\lower2.5pt\hbox
{\vbox{\moveright0.1pt\hbox{$^{^{\land}}$}}}}}{{\ast}}}"$, defined
in
{\cite{25} %
p. 4} \hbox{\tenbf /}$\!\!$\hbox{\tenbf /} in
{\cite{6} %
p.\ 128 \ Ex.\ 3 including \cite{6} %
p.\ 135 Problem\
6:1}, \cite{21} %
p.\ 373 and \cite{29} %
p.\ 128 \hbox{\tenbf /}$\!\!$\hbox{\tenbf /} in
{\cite{2} %
pp. $\!$159-160 and \cite{18} %
}
respectively. Recall that; $X\!\ast\emptyset=\!X=\emptyset\ast X$
classically.

\definition{Definition}
$X_{\wp1}\lower1.0pt\hbox{$^{_{_{\hbox{\sevensy \char"72}
}}}$}\!\!_{{^o}} X_{\wp2}:=
\cases
\emptyset & if X_{\wp1}=\emptyset\ \hbox{\tenrm or}\
X_{\wp2}=\emptyset \cr
{\hbox{\tensy F}_{\wp}}({\hbox{\tensy
F}(X_{\wp1})\lower1.0pt\hbox{$^{_{_{\hbox{\sevensy \char"72}
}}}$}\hbox{\tensy F}(X_{\wp2})}) & \hbox{\tenrm if} X_{\wp1}\neq
\emptyset \neq X_{\wp2}.\cr
\endcases
$
\enddefinition

\smallskip
From now on we'll delete the \raise0.6pt\hbox{$\wp$/\eightrm o}\
$\!_{_{\!}}$-indices. So, e.g. $\!$``$X$ connected" now means
``$\hbox{\tensy F}(X)\ \hbox{\tenrm connected}"\!\!.$

\definition{Equivalent Join Definition}
Put $\emptyset{\myPsqcup1}\ \!X\!=\!X\!{\myPsqcup1}\
\emptyset\!:=\!\emptyset$.\
\hbox{\tenrm If} $X\!_{\!}\ne\!\emptyset;$
$\{\wp\}\!{\myPsqcup1}\ \!X\!=\!\{\wp\}$, $\!X_{\!}{\myPsqcup1}\
\{\wp\}\!_{\!}=\!X\!$.
For $X,\!Y\!_{\!}\ne\!_{\!}\emptyset, \wp$
let $X\!{_{_{^{}}}}\!{\myPsqcup1}\ Y\!{_{_{^{}}}}$ denote the set
$X\!{_{_{^{}}}}\times Y\!{_{_{^{}}}}\times(0,1]$ \hbox{\teni
pasted} to the set $X\!{_{_{^{}}}}$ by
$ \varphi\!{_{_{^{1}}}}: X\!{_{_{^{}}}}\times
Y\!{_{_{^{}}}}\!\times\{1\} \longrightarrow X\!{_{_{^{}}}};
(x{_{_{^{}}}}\!,y\!{_{_{^{}}}},1) \mapsto x{_{_{^{}}}}, $ i.e. the
quotient set of $\big(X\!{_{_{^{}}}}\times
Y\!{_{_{^{}}}}\times(0,1]\big) \sqcup X\!{_{_{^{}}}},$ under the
equivalence relation $(x{_{_{^{}}}}\!, y{_{_{^{}}}},1)\sim
x{_{_{^{}}}}\!$ and let
$p{_{_{^{_{\!}1}}}}\!:\big(X\!{_{_{^{}}}}\times Y\!{_{_{^{}}}}\!
\times(0,1]\big)\sqcup X\!{_{_{^{}}}}\! \longrightarrow
X\!{_{_{^{}}}}{\myPsqcup1}\ Y\!{_{_{^{}}}}\! $ be the quotient
function. For $X,Y=\emptyset\ \hbox{\tenrm or}\ \{\wp\}$ let $
X\!{_{_{^{}}}}\!{\myPsqcup0}\ \!Y\!{_{_{^{}}}}:=
Y\!{_{_{^{}}}}\!{\myPsqcup1}\ \!X\!{_{_{^{}}}} $ and else
the set $X\!{_{_{^{}}}}\times Y\!{_{_{^{}}}}\times[0,1)$
\hbox{\teni pasted} to the set $Y\!{_{_{^{}}}}$ by the function $
\varphi\!{_{_{^{2}}}}: X\!{_{_{^{}}}}\times
Y\!{_{_{^{}}}}\!\times\{0\} \longrightarrow Y\!{_{_{^{}}}};
(x{_{_{^{}}}}\!,y\!{_{_{^{}}}},0) \mapsto y{_{_{^{}}}}, $ and let
$ p{_{_{^{_{\!}2}}}}\!:\big(X\!{_{_{^{}}}}\!\times Y\!{_{_{^{}}}}
\times[0,1)\big)\sqcup Y\!{_{_{^{}}}}\!\! \longrightarrow
X\!{_{_{^{}}}} {\myPsqcup0}\ Y\!{_{_{^{}}}} $ be the
quotient\nobreak\ function. Put $
X\!{_{_{^{}}}}\!\circ\!Y\!{_{_{^{}}}} :=\big(X\!{_{_{^{}}}}
{\myPsqcup1}\ Y\!{_{_{^{}}}}\big) \cup
\big(X{_{_{^{}}}}\myPsqcup0\ Y\!{_{_{^{}}}}\big). $

$ (x{_{_{^{}}}}\!,y\!{_{_{^{}}}},t)\in X\!{_{_{^{}}}}\!\times\!
Y\!{_{_{^{}}}}\!\times[0,1] $
specifies the point $ (x{_{_{^{}}}}\!,y\!{_{_{^{}}}},t)\in
X\!{_{_{^{}}}}\!{\myPsqcup1}\ \!Y\!{_{_{^{}}}} \cap
X\!{_{_{^{}}}}\!{\myPsqcup0}\ \!Y\!{_{_{^{}}}}\!, $ one-to-one,
$\hbox{\tensl if}\ 0<t<1$ \hbox{\tensl and}
the equivalence class containing $x{_{_{^{}}}}$ if $t=1$
$(y\!{_{_{^{ }}}}$ if $t=0)$, which we denote $(x{_{_{^{}}}}\!,1)\
((y{_{_{^{}}}},0)).$
This allows ``coordinate functions" $\xi\!:\!
X\!{_{_{^{}}}}{\circ}Y\!{_{_{^{}}}}\!\!\rightarrow\![0,1]$,
$\eta{_{_{^{\!1}}}}\!\!:\!X\!{_{_{^{}}}}\!{\myPsqcup1}\
\!Y\!{_{_{^{}}}} \!\!\rightarrow\!X\!{_{_{^{}}}}, $
$\eta_{_{^{2}}}\!\!:\! X\!{_{_{^{}}}}\!{\myPsqcup0}\
\!Y\!{_{_{^{}}}} \!\!\rightarrow\!\! Y\!{_{_{^{}}}}\ \hbox{\tenrm
extendable\ to}\ \! X\!{_{_{^{}}}}{\circ}Y\!{_{_{^{}}}}
\hbox{\tenrm through} $ $
\eta{_{_{^{_{_{\!}1}}}}}\!(y{_{_{^{}}}}\!,0):\equiv
x\!{_{_{^{0}}}}\! \in  X\!{_{_{^{}}}} $ resp. $
\eta{_{_{^{_{\!}2}}}}\!(x{_{_{^{}}}}\!,1):\equiv y\!{_{_{^{0}}}}\!
\in  Y\!{_{_{^{}}}} $ and a projection $ p\!:X\!{_{_{^{}}}}\sqcup
\big(X\!{_{_{^{}}}}\times Y\!{_{_{^{}}}}\!\times[0,1]\big)
{\sqcup}\ \!Y\!{_{_{^{}}}}\! \!\rightarrow\!
X\!{_{_{^{}}}}{\circ}Y\!{_{_{^{}}}}\!. $

Let $X\!{_{_{^{}}}} {\rlap{{\lower2.5pt\hbox
{\vbox{\moveright0.17pt\hbox{$^{^{ \hbox{\fivesy {\char"5E}}
}}$}}}}}{{\ast}}} Y\!{_{_{^{}}}}\! $ denote $
X\!{_{_{^{}}}}{\circ}Y\!{_{_{^{}}}}\! $ equipped with the smallest
topology making $\xi, \eta{_{_{^{1}}}},\eta{_{_{^{2}}}} $
continuous and $X\!{_{_{^{}}}}{{\ast}}Y\!{_{_{^{}}}}\!,
X\!{_{_{^{}}}}{\circ}Y\!{_{_{^{}}}} $ with the quotient topology
w.r.t. $p$, i.e. the largest topology
making $p$ continuous
$(\Rightarrow \tau\!\! _{_{^{X
\!{{\hat\ast}}\!
Y\!{_{_{^{}}}}\! }}} \! \raise1.5pt\hbox{ {\sixsy {\char"1A}}
}\tau\!\! _{_{ X\!{_{_{^{}}}}\!{{\ast}}\!Y\!{_{_{^{}}}}\! }}).$
\enddefinition

\definition{Pair-definitions}
$({X}_{_{^{\!1}}},{X}_{_{^{\!2}}})\lower1.0pt\hbox{$^{_{_{\hbox{\sevensy
\char"72} }}}$}\!\!_{_{\ominus}}\ \!
({Y}\!_{_{^{\!1}}},{Y}\!_{_{^{\!2}}}) \!:= \!
({X}_{_{^{\!1}}}\lower1.0pt\hbox{$^{_{_{\hbox{\sevensy \char"72}
}}}$}\!\!_o \ \!{Y}\!_{_{^{\!1}}},
({X}_{_{^{\!1}}}\lower1.0pt\hbox{$^{_{_{\hbox{\sevensy \char"72}
}}}$}\!\!_o \ \!{Y}\!_{_{^{\!2}}}) \ominus
({X}_{_{^{\!2}}}\lower1.0pt\hbox{$^{_{_{\hbox{\sevensy \char"72}
}}}$}\!\!_o \ \! {Y}\!_{_{^{\!1}}})),$ where
${\ominus}$\ stand for ``$\cup$" or ``$\cap$"
and
if either $X_{_{^{\!2}}}\!$ or $Y\!_{_{^{\!2}}}\ \!$ is not closed
\hbox{\tenrm(}open\hbox{\tenrm)}
\hbox{\tenrm cf. \cite{6} %
p.\ 43 Cor.\ 1},
$ ({X}_{_{^{\!1}}}\!\!\ast\!_{_{^{\!\ }}}\!\! {Y}\!_{_{^{\!1}}},
({X}_{_{^{\!1}}}\!\!\ast\!_{_{^{\!\ }}}\!\!
{Y}\!_{_{^{\!2}}})\ominus ({X}_{_{^{\!2}}}\!\!\ast\!_{_{^{\!\
}}}\!\! {Y}\!_{_{^{\!1}}})) $ has to be interpreted as
$ ({X}_{_{^{\!1}}}\!\!\ast\!_{_{^{\!\ }}}\! {Y}\!_{_{^{\!1}}},
({X}_{_{^{\!1}}}\!
{\rlap{\lower3.97pt\hbox{$^{^{_{{_{\!}}\bigcirc}}}$}}
{\vbox{\moveright0.1pt\hbox{$\ast$}}}}\ \!
{Y}\!_{_{^{\!2}}})\ominus ({X}_{_{^{\!2}}}\!
{\rlap{\lower3.97pt\hbox{$^{^{_{{_{\!}}\bigcirc}}}$}}
{\vbox{\moveright0.2pt\hbox{$\ast$}}}}\ \! {Y}\!_{_{^{\!1}}})) $
i.e. $({X}_{_{^{\!1}}}\!\circ{Y}\!_{_{^{\!2}}})\ominus
({X}_{_{^{\!2}}}\!\circ{Y}\!_{_{^{\!1}}}) $ with the subspace
topology in the 2:nd component. Analogously for simplicial
complexes with ``$\times\!$" $($``$\ast_{\!}$"$)$ from
\cite{7} %
p.\ 67 Def. 8.8
(\ \!\cite{25} %
 p.\ 109 Ex.\ 7.)
\enddefinition

\remark{Note} \hbox{\tenbf i.}
$(X\!{_{_{^{}}}}\!{\myPsqcup1}\
Y\!{_{_{^{}}}})^{_{^{_{\!}t\geq0_{{\!}}._{{\!}}5}}}\!$ in
${X}\!_{_{^{\ \! }}}\!{\ast}${\teni Y}$\!_{_{^{\ }}}\!$ is
homeomorphic to the mapping cylinder w.r.t. the coordinate map
${q}_{_{^{\!1}}}\!\!:\!\!{X}_{_{^{\!\
}}}\!\!\!\times\!{Y}_{_{^{\!\ }}}\! \!\!\rightarrow\! X. $ \ \ \
\hbox{\tenbf ii.}
${X}\!_{_{^{2}}}\!{\rlap{{\lower2.5pt\hbox
{\vbox{\moveright0.1pt\hbox{$^{^{\land}}$}}
}}}{{\ast}}}\!{_{_{\!}}}$ {\teni Y}$\!_{_{^{2}}}$
is a subspace of {\teni X}$_{_{\!1}} \!{\rlap{{\lower2.5pt\hbox
{\vbox{\moveright0.1pt\hbox{$^{^{\land}}$}}}}}{{\ast}}}\!{_{_{\!}}}
$ {\teni Y}$\!_{_{1}}$ by
\cite{2} %
5.7.3 p.\ 163.$_{_{\!}}$
${X}\!_{_{^{2}}}\!{\ast}${\teni Y}$\!_{_{^{2}}}\!$
is a subspace of$\!$ {\teni X}$_{_{\!1}}\!{\ast}${\teni
Y}$\!_{_{1}}$ if {\teni X}$_{_{\!2}}$,{\teni Y}$\!_{_{2}}\!$ are
closed (open). \ \ \
{\bf iii.} $({X}_{_{^{\!1}}},\{\wp\})\times
({Y}\!_{_{^{\!1}}},{Y}\!_{_{^{\!2}}})=
({X}_{_{^{\!1}}},\emptyset)\times
({Y}\!_{_{^{\!1}}},{Y}\!_{_{^{\!1}}})\ \hbox{\tenrm if}\
{Y}\!_{_{^{\!2}}}\neq\emptyset$
and %
$({X}_{_{^{\!1}}}\!\circ{Y}\!_{_{^{\!2}}})\ \!{\cap}\ \!
({X}_{_{^{\!2}}}\!\circ{Y}\!_{_{^{\!1}}})\!=\!
{X}_{_{^{\!2}}}\!\circ{Y}\!_{_{^{\!2}}}$.
\ \ \ {\bf iv.}
${\rlap{{\lower2.5pt\hbox
{\vbox{\moveright0.17pt\hbox{$^{^{\land}}$}}}}}{{\ast}}}$ and
$\ast$ are both commutative but, while
${\rlap{{\lower2.5pt\hbox
{\vbox{\moveright0.17pt\hbox{$^{^{\land}}$}}}}}{{\ast}}}$
is associative by \cite{2} %
p.\ 161, $\ast$ isn't in general, cf. p\nobreak\ 14.
\ \ \ {\bf v.}
``$\times\!\!_{_{\cap}}\!$" is (still, cf. {\cite{5} %
p.\ 15},) the categorical \hbox{product on pairs from
$\hbox{\tensy D}\!\!_{_{^{\wp}}}\!$.}
 \ {\bf vi.}
$\lower1.0pt\hbox{$^{_{_{\hbox{\sevensy \char"72}
}}}$}\!\!:=\!\lower1.0pt\hbox{$^{_{_{\hbox{\sevensy \char"72}
}}}$}\!\!\!_{_{{{\cup}}}}.$
\endremark

\subhead  3.2 Augmental Homology for Products and Joins
\endsubhead
$\!$Through Lemma $\!${\eightbf+}$\!$ Note {\bf ii} p.\ 6 we
convert the classical K\"unneth formula cf. $\!$\cite{25} $\!$p.\
$\!$235, mimicking what Milnor did, partially$_{\!}$ ($\equiv\
\!$line\ 1),at the end of his proof of %
\cite{18} %
{p.\ 431 Lemma\ 2.1.}
The ability of a full and clear understanding of Milnor's proof
could be regarded as sufficient prerequisites for our next six
pages.
The new object $\{\wp\}$ gives the classical K\"unneth formula
($\equiv$4:th line) additional strength but much of the classical
beauty is lost - a loss which is regained in the join version i.e.
in Theorem\ 4 p.\ 10.

\proclaim{Theorem 1}
For $\{X\!_{_{1}}\!\!\times\! Y\!\!_{_{2\!}},
X\!_{_{2\!\!}}\times\!Y\!\!_{_{1}}\}$ excisive,
$\hbox{\tenbf q}\geq0$, {\bf R} a {\bf PID}, and assuming
$\hbox{\tenrm Tor}_1^{\hbox{\fivebf R}}(\hbox{\tenbf
G},\hbox{\tenbf G}^\prime)\!=\!0\ \hbox{\teni then};$

\smallskip%
\noindent {$\mdoubleH_{q}((X_{1},X_{2})\times (Y_{1},
Y_{2});\hbox{\tenbf G}\ \!{{{\otimes}} _{\!}{{_{_{\hbox{\fivebf
R}}}}\hbox{\tenbf G}^\prime}})
\raise2pt\hbox{$
{{{_{\hbox{\fivebf R}}}}\atop{\lower1pt\hbox{$\cong$}}}
$}
$

\smallskip
\noindent$\!
\raise2pt\hbox{$
{{{_{\hbox{\fivebf R}}}}\atop{\lower1pt\hbox{$\cong$}}}
$}
{\cases
[\mdoubleH_i(X_{1};\hbox{\tenbf G})
\!\otimes\!\!{_{_{_{\!\hbox{\fivebf R}}}}}\!
{
{\mdoubleH_j(Y_{1};\hbox{\tenbf
G}^\prime)]\lower2pt\hbox{$_{q}$}\!\oplus}
}
(\mdoubleH_{q}(X_{1};\hbox{\tenbf G})
\!\otimes\!\!{_{_{_{\!\hbox{\fivebf R}}}}}\!
\!\hbox{\tenbf G}^\prime)
\!\oplus\!
(\hbox{\tenbf G}
\!\otimes\!\!{_{_{_{\!\hbox{\fivebf R}}}}}\!
\mdoubleH_{q}(Y_{1};{\hbox{\tenbf G}}^\prime))
{\rlap{$\oplus$}{\hbox{\tenbf \ \ T}_1}}\ {if}\ \hbox{\tenbf C}_1\
\cr
[\mdoubleH_i(X_{1};\hbox{\tenbf G})\otimes_{_{\hbox{\fivebf R}}}
\mdoubleH_j(Y_{1},Y_{2};\hbox{\tenbf G}^\prime)]_q \oplus
(\hbox{\tenbf G}\!\otimes_{_{\hbox{\fivebf R}}}\!
\mdoubleH_{q}(Y_{1},Y_{2};\hbox{\tenbf
G}^\prime))\oplus\hbox{\tenbf T}_2 \indent\indent\ \ \ {if}\
\hbox{\tenbf C}_2\
\cr
[\mdoubleH_i(X_{1},X_{2};\hbox{\tenbf G})\otimes_{_{\hbox{\fivebf
R}}} \mdoubleH_j(Y_{1};\hbox{\tenbf G}^\prime)]_q \oplus
(\mdoubleH_{q}(X_{1},X_{2};\hbox{\tenbf G})
\otimes_{_{\hbox{\fivebf R}}}\hbox{\tenbf
G}^\prime)\oplus\hbox{\tenbf T}_3 \indent\ \ \ \ \!{if}\
\hbox{\tenbf
C}_3\ %
\cr
[\mdoubleH_i(X_{1},X_{2};\hbox{\tenbf G})\otimes_{_{\hbox{\fivebf
R}}}\mdoubleH_j(Y_{1},Y_{2};\hbox{\tenbf
G}^\prime)]_q\oplus\hbox{\tenbf T}_4\ \
\hskip4.12cm{if}\ \hbox{\tenbf C}_4\
\endcases
}
\!\!\!\!\!\!\!\!\!
\hfill(_{\!}\hbox{\tenbf 1}_{\!})\!\!
$}
\endproclaim

\smallskip
The torsion terms, i.e. the $\hbox{\tenbf T}\!_{_{^{}}}$-terms,
splits as those ahead of them, resp., e.g.
$$\hbox{\tenbf T}\!_{_{^{1}}}=[\hbox{\tenrm
Tor}{\rlap{$_{_1}$}{^{\hbox{\fivebf R}}}}\bigl(
\mdoubleH_i(X_{1};\hbox{\tenbf G}), \mdoubleH_j(Y_{1};
\hbox{\tenbf G}^\prime)\big)]{\rlap{$_{_{q-1}}$} {\ \ \ \
\!\raise1pt\hbox{$\oplus$}}}\hskip3.5cm$$
\vskip-1.0cm
$$\hskip3.5cm\oplus\ \!\hbox{\tenrm Tor}{\rlap{$_{_1}$}{^{\hbox{\fivebf R}}}}
\bigl(\mdoubleH_{q-1}(X_{1};\hbox{\tenbf G}{\rlap{)}} \ ,
{\hbox{\tenbf G}^\prime}{\rlap{$\bigr)$}{\ \ \oplus}} \ \!
\hbox{\tenrm Tor}{\rlap{$_{_1}$}{^{\hbox{\fivebf
R}}}}\bigl({\rlap{{\rlap{ }\hbox{\tenbf G}}} } \ \ \
,\mdoubleH_{q-1}(Y_{1};\hbox{\tenbf G}^\prime)\bigr),
$$
\vskip-0.5cm
and
\vskip-0.3cm
$$
\hbox{\tenbf T}\!_{_{^{4}}}=[\hbox{\tenrm
Tor}{\rlap{$_{_1}$}{^{\hbox{\fivebf R}}}}\hbox{\tenbf(}
\mdoubleH_i(X_{1},X_{2};\hbox{\tenbf G}),\
\!\mdoubleH_j(Y_{1},Y_{2}; \hbox{\tenbf
G}^\prime)\hbox{\tenbf)}]{{\lower3pt\hbox{{\fivei
q}{\fiverm-1}}}}.\hskip2.5cm
$$

\noindent $\hbox{\tenbf C}_1:=\!``X\!_{1}\!\times\!Y\!\!_{1}\neq
\emptyset,\{\wp\}\ {and}\ X_{2}\!=\!\emptyset\!=\!Y_{2}",$
$\hbox{\tenbf C}_2:=\!``X\!_{1}\!\times\!Y\!\!_{1}\neq
\emptyset,\{\wp\}\ {and}\ X_{2}\!=\!\emptyset\!\neq\!Y_{2}",$
$\hbox{\tenbf C}_3:=\!``X\!_{1}\!\times\!Y\!\!_{1}\neq
\emptyset,\{\wp\}\ { and}\ X_{2}\!\neq\!\emptyset\!=\!Y_{2}",$
$\hbox{\tenbf C}_4:=\!``X_{1}\!\!\times\!\!Y\!_{1} \!=\!
\emptyset,\!\{\wp\}\ \hbox{\tenbf or} \ X_{2}\neq \emptyset\neq
Y_{2}",$
and
[\hbox{\bf...}]$_q$ is still, i.e. as in \cite{25} %
p.\ 235 Th.\ 10, to be interpreted as$\nobreak\ \raise1.5pt\hbox{
$\bigoplus\hbox{\tenbf...}\!\!\!\!\!\!\!\!\!\!\!\!\!\!\!\!\!\!\!\!\!\!\!\!
\lower3pt\hbox{$_{_{\scriptstyle i+j=q\ \!\!\&\ \!\!\scriptstyle
i,j\geq0}}$} .$} \qed$

\proclaim{Lemma}
Let $ f\!\!:\!(X,A)\!\rightarrow \!(Y,B) $ be a relative
homeomorphism, i.e., $f\!:X\!\rightarrow\!Y$ is continuous  and
$f\!:\!X\setminus\!{_{_{^{\!o\!}}}} A\!\rightarrow
Y\setminus\!{_{_{^{\!o\!}}}} B$
is a homeomorphism.
If $F\!\!:\!N\!\times\hbox{\tenbf I}\!\rightarrow\! N\! $ is a
$(\!$strong$)$ $(\!$neighborhood$)$ deformation retraction of $N$
onto $A$ and {\teni B} and \hbox{{\teni f}$(_{\!}${\teni N}}$\
_{\!\!})$ are closed in
$N^{^{_{\prime}}}\!\!\!:=\!f(_{\!}N\setminus\!{_{_{^{\!o\!}}}}
A)\cup B,$ then $B$ is a $(\!$strong$)$ $($neighborhood$)$
deformation retract of $N^{^{_{\prime}}}$ through$;$
$$F^{^{_{\prime}}}\!\! : N^{^{_{\prime}}}\!\!\times\hbox{\tenbf
I}\!\rightarrow N^{^{_{\prime}}} \!;\!
\cases
{\!}\hbox{\tenrm({\teni y},{\teni t})\
$\!$\lower0pt\hbox{${{\mapsto}}$}\ $\!${\teni y}} & {\!\!\!}if\
{\teni y} \!\in\! {\teni B}, {\teni t} \!\in\! {\tenbf I} \cr
{\!}\hbox{\tenrm({\teni y},{\teni t})\
$\!$\lower0pt\hbox{${{\mapsto}}$}} {\hbox{\tenrm{\teni
f}$\circ${\teni F} ({\teni f} $\!^{^{_{\hbox{\fiverm-}1}}}\!
$\hbox{\ninerm(}{\ninei y}\hbox{\ninerm)}, {\teni t}) }}\!\! &
{\!\!\!} if\ {\teni y} \!\in\! {\teni f}({\teni
N})\setminus\!{_{_{^{\!o\!}}}}${\teni B} ${\!}={\!} {\teni
f}({\teni N}\setminus\!{_{_{^{\!o\!}}}}{\teni A}),\ {\teni t}
\!\in\! \hbox{\tenbf I}. \cr
\endcases
$$
\endproclaim

\demo{Proof}
$_{\!\!} F^{^{_{\prime}}}\!\!$ is continuous as being so
when\nobreak\ $\hbox{\tenrm restricted}$ to
$f(N)_{\!}\times_{\!}\hbox{\tenbf I}$ resp. ${\teni
B}_{\!}\times_{\!}\hbox{\tenbf I},$ cf. \cite{2} %
p.$\ 34;\ 2.5.12.\!\!\quad \qed$
\enddemo

\proclaim{Theorem 2}
\hbox{\rm(Analogously for ${\rlap{{\lower2.5pt\hbox
{\vbox{\moveright0.1pt\hbox{$^{^{\land}}$}}}}}{{\ast}}}_{\!}. $
\ \ ($\Rightarrow
\hbox{\tenrm{\char"01}}\!_{\!}{{\raise4.5pt\hbox{\fivei{\char"7D}}}}(
\cdot\ {\!}{\rlap{{\lower2.5pt\hbox
{\vbox{\moveright0.1pt\hbox{$^{^{\land}}$}}}}}{{\ast}}}\ {\!}
\cdot )$
is chain equivalent to
$\hbox{\tenrm{\char"01}}\!_{\!}{{\raise4.5pt\hbox{\fivei{\char"7D}}}}(
\cdot_{\!}\ast_{\!}\cdot).$
)\ )}

\smallskip
If $
(X_{1},X_{2})\!\neq\!(\{\wp\},\emptyset)\!\neq\!(Y_{1},Y_{2})\ {
and\ \hbox{\tenbf G}\ an\ \hbox{\tenbf A}\hbox{-}module};$

\smallskip
\noindent
$\mdoubleH_{q}((X_{1},X_{2})\!\times\!(Y_{1},Y_{2});\hbox{\tenbf
G})
{\lower4.0pt\hbox{$^{_{\hbox{{\fivebf A}}
\atop{{\raise0.0pt\hbox{{\eightbf =}}}
\!\!\!\!\!_{\!}{\raise2.8pt\hbox{{\eightsy {\char"27}}}} }}}$} } $

\medskip
$
{\lower4.0pt\hbox{$^{_{\hbox{{\fivebf A}}
\atop{{\raise0.0pt\hbox{{\eightbf =}}}
\!\!\!\!\!_{\!}{\raise2.8pt\hbox{{\eightsy {\char"27}}}} }}}$} }\
\mdoubleH\!\!\!{_{_{q+1}}}\!((X\!_{_1}\!,\!X\!_{_2}\!)\ast
(\!Y\!\!_{_1}\!,\!Y\!\!_{_2}\!);\hbox{\tenbf G})
\oplus\
\mdoubleH\!_{_{q}}\!((X\!_{_1}\!,\!X\!_{_2}\!)\!\ast\!
(\!Y\!\!_{_1}\!,\!Y\!\!_{_2}\!)^{^{_{\!t\geq0.5}}}\hbox{\bf+}\
(X\!_{_1}\!,\!X\!_{_2}\!) \ast
(Y\!\!_{_1}\!,\!Y\!\!_{_2}\!)^{^{_{\!t\leq0.5}}}\!;\hbox{\tenbf
G})
$ {\bf=}

$$
{\lower4.0pt\hbox{$^{_{\hbox{{\fivebf A}}
\atop{{\raise0.0pt\hbox{{\eightbf =}}}
\!\!\!\!\!_{\!}{\raise2.8pt\hbox{{\eightsy {\char"27}}}} }}}$} }
{\cases
\mdoubleH_{_{q+1}}
 (\!X\!_{_1}\!\ast Y\!\!_{_1};
\hbox{\tenbf G}) \oplus \mdoubleH_{q}(X\!_{_1};\hbox{\tenbf G})
\oplus \mdoubleH_{q}(Y\!\!_{_1};\hbox{\tenbf G})\!\!\!\!\! & {
if}\ \ \
\hbox{\tenbf C}_1
\cr
\mdoubleH_{_{q+1}}((\!X\!_{_1},\emptyset)\ast
(Y\!\!_{_1},Y\!\!_{_2});\hbox{\tenbf G})\oplus
\mdoubleH_{q}(Y\!_{_1},Y\!\!_{_2};\hbox{\tenbf G}) & {if}\ \ \
\hbox{\tenbf C}_2
\cr
\mdoubleH_{_{q+1}}((\!X\!_{_1},X\!_{_2})\ast
(Y\!\!_{_1},\emptyset);\hbox{\tenbf G})\oplus
\mdoubleH_{q}(X\!_{_1},X\!_{_2};\hbox{\tenbf G}) & {if}\ \ \
\hbox{\tenbf C}_3
\cr
\mdoubleH_{_{q+1}}( (\!X\!_{_1},X\!_{_2})\ast
(Y\!\!_{_1},Y\!\!_{_2})  ;\hbox{\tenbf G}) & {if}\ \ \
\hbox{\tenbf C}_4
\indent\indent\
\cr
\endcases
} \!
\eqno{(\hbox{\tenbf 2})}
$$

\medskip
\noindent $\hbox{\tenbf C}_1:=\!``X\!_{1}\!\times\!Y\!\!_{1}\neq
\emptyset,\{\wp\}\ {and}\ X_{2}\!=\!\emptyset\!=\!Y_{2}",$
$\hbox{\tenbf C}_2:=\!``X\!_{1}\!\times\!Y\!\!_{1}\neq
\emptyset,\{\wp\}\ {and}\ X_{2}\!=\!\emptyset\!\neq\!Y_{2}",$
$\hbox{\tenbf C}_3:=\!``X\!_{1}\!\times\!Y\!\!_{1}\neq
\emptyset,\{\wp\}\ { and}\ X_{2}\!\neq\!\emptyset\!=\!Y_{2}",$
$\hbox{\tenbf C}_4:=\!``X_{1}\!\!\times\!\!Y\!_{1} \!=\!
\emptyset,\!\{\wp\}\ \hbox{\tenbf or} \ X_{2}\neq \emptyset\neq
Y_{2}".$
\endproclaim

\demo{Proof} Split ${X}_{_{^{\!\ }}}\!\!\! \ast\!{Y}_{_{^{\!\
}}}\!$ at $t\!=\!0.5$ then; $X$ ($Y$) is a strong deformation
retract of\break $ ({X}_{_{^{\!\ }}}\!\!\ast\!{Y}_{_{^{\!\
}}}\!)^{^{_{\!t\geq0.5}}}\!\! \ \ (({X}_{_{^{\!\
}}}\!\!\ast\!{Y}_{_{^{\!\ }}}\! )^{^{_{\!t\leq0.5}}}),$
the mapping cylinder w.r.t. product projection.  The relative {\bf
M-$\!$V$\!$s} w.r.t. $\!{\!}$the excisive couple of pairs
$
\{(\!X\!_{_1}\!,\!X\!_{_2}\!)_{^{\!}}\ast_{^{\!}}
(\!Y\!\!_{_1}\!,\!Y\!\!_{_2}\!)^{^{_{\!t\geq0.5}}}\!,$
$(\!X\!_{_1}\!,\!X\!_{_2}\!)_{^{\!}}\ast_{^{\!}}
(\!Y\!\!_{_1}\!,\!Y\!\!_{_2}\!)^{^{_{\!t\leq0_{{\!}}._{{\!}}5}}}
\}$
splits since the inclusion of their topological sum into
$(\!X\!_{_1}\!,\!X\!_{_2}\!)_{\!}\ast_{\!}
(\!Y\!\!_{_1}\!,\!Y\!\!_{_2}\!)$
is pair null-homotopic, cf. \cite{21} %
p.\ 141 Ex.\ 6c, and \cite{15} p.\ 32 Prop.\ 1.6.8.
Since the 1:st (2:nd) pair is acyclic if
$Y\!\!_{_{\!2}}\ \!(X\!_{_2\!})\!\neq\!\emptyset\ \hbox{\tenrm we
get Theorem\ 2.} $
Equivalently for ${\rlap{{\lower2.5pt\hbox
{\vbox{\moveright0.17pt\hbox{$^{^{\land}}$}}}}}{{\ast}}}$ by the
Lemma.\qed
\enddemo

Milnor finished his proof of \cite{18} Lemma\ 2.1 p.\ 431 by
simply comparing the r.h.s. of the $\hbox{\tenbf C}_1$-case in
Eq.\ 1 with that of Eq.\ 2. Since we are aiming at the stronger
result of ``natural chain equivalence" in Theorem\ 3 this isn't
enough and so, we'll need the following three auxiliary results to
prove our next two theorems.
We hereby avoid explicit use of ``proof by acyclic models".
\ \ (``$\approx$" stands for ``chain equivalence".)

\proclaim{5.7.4} \hbox{\rm  (\cite{2} p.\ 164.)}
\hbox{\spaceskip2.2pt\rm(Here $\hbox{\tenbf E}^0$ is a symbol for
a point, i.e. a $0$-disc also denoted $\bullet$.)}\ \

\nointerlineskip
$$
There\ is\ a\ homeomorphism:\indent
\nu:X
{\rlap{{\lower2.5pt\hbox {\vbox{\moveright0.1pt\hbox{$^{^{\land
}}$}}}}}{{\ast}}}\ \!
Y
{\rlap{{\lower2.5pt\hbox {\vbox{\moveright0.1pt\hbox{$^{^{\land
}}$}}}}}{{\ast}}}\ \!
\hbox{\tenbf E}^0
\longrightarrow
(X
{\rlap{{\lower2.5pt\hbox {\vbox{\moveright0.1pt\hbox{$^{^{\land
}}$}}}}}{{\ast}}}\ \!
\hbox{\tenbf E}_{_1}^0)\!\times\!
(Y{\rlap{{\lower2.5pt\hbox {\vbox{\moveright0.1pt\hbox{$^{^{\land
}}$}}}}}{{\ast}}}\ \!\hbox{\tenbf E}_{_2}^0)
$$

\centerline{\ \  \hfill
which restricts to a homeomorphism{\rm :} \indent
$
X {\rlap{{\lower2.5pt\hbox {\vbox{\moveright0.1pt\hbox{$^{^{\land
}}$}}}}}{{\ast}}}\ \!Y\!\!
\rightarrow
(X
{\rlap{{\lower2.5pt\hbox {\vbox{\moveright0.1pt\hbox{$^{^{\land
}}$}}}}}{{\ast}}}\ \!
\hbox{\tenbf E}^0_{_1})\times Y\cup
X\times(Y{\rlap{{\lower2.5pt\hbox
{\vbox{\moveright0.1pt\hbox{$^{^{\land }}$}}}}}{{\ast}}}\
\!\hbox{\tenbf E}^{0}_{_2}).
$
\qed}
\endproclaim

\proclaim{Corollary 5.7.9} \hbox{\rm  (\cite{15} p.\ 210.)}
If $\phi$:\ {\bf C} $\approx$ {\bf E} with inverse $\psi$ and
$\phi^\prime:$\ {\bf C}$^\prime\approx$ {\bf C}$^\prime$ with
inverse $\psi^\prime$, then

\hskip 3 cm $\phi\otimes\phi^\prime:\hbox{\tenbf
C}\otimes\hbox{\tenbf C}^\prime\approx \hbox{\tenbf
E}\otimes\hbox{\tenbf E}^\prime$ with inverse
$\psi\otimes\psi^\prime$.
\qed
\endproclaim

\proclaim{Theorem 46.2}
\hbox{\rm (\cite{21} p.\ 279.)}
For free chain complexes\ $\!\ \hbox{\tensy C}, \hbox{\tensy
D}_{\!}$ vanishing below a certain dimension\nobreak\ and\nobreak\
if\nobreak\ a chain map $\lambda\!:\hbox{\tensy
C}\!\rightarrow\!\hbox{\tensy D}$ induces homology isomorphisms in
all dimensions, then $\lambda$ is a chain equivalence.\
$\qed$
\endproclaim

\proclaim{Theorem 3}
$\!$\hbox{\spaceskip2pt\rm(The relative Eilenberg-Zilber theorem
for joins.)}
For an excisive couple $\{X_{_{\!\ \!}}\!{\rlap{{\lower2.5pt\hbox
{\vbox{\moveright0.1pt\hbox{$^{^{\land }}$}}}}}{{\ast}}}
Y\!\!_{_{2\!}}, X\!_{_{2}}\! {\rlap{{\lower2.5pt\hbox
{\vbox{\moveright0.1pt\hbox{$^{^{\land }}$}}}}}{{\ast}}}
Y_{_{\!\!\ \!}}\}$ from the category of ordered pairs $((X\!_{_{\
\!\!}}\!,X\!_{_{2}}\!),\ \! (Y\!\!_{_{\ \!\!}},Y\!\!_{_{2}}\!))\!$
of\ topological\ pairs$\!_{_{^{_{\wp}}}}{_{\!\!}},$
$ \hbox{\tenbf s}\hbox{\tenbf (} {
\Delta{\raise1.5pt\hbox{$\!\!^{\wp}$}}} \!(\!X\!_{_{\
\!}}\!,X\!_{_2})\nobreak\otimes\nobreak{\Delta
{\raise1.5pt\hbox{$\!\!^{\wp}$}}} \!(Y_{_{\!\!\
\!}}\!,Y\!\!_{_2})\hbox{\tenbf )}$ is naturally chain equivalent
to ${ \Delta{\raise1.5pt\hbox{$\!\!^{\wp}$}}}\! \hbox{\tenbf
(}(\!X\!_{_{\ \!}}\!,X\!_{_2}) \ \!{\rlap{{\lower2.5pt\hbox
{\vbox{\moveright0.1pt\hbox{$^{^{\land }}$}}}}}{{\ast}}}\ \!
(Y\!\!_{_{\ \!}}\!,Y\!\!_{_2})\hbox{\tenbf )}.$
{\rm (``{\bf s}" stands for  suspension i.e. the
suspended\nobreak\ chain equals the original except that dimension
$i$ in the original is dimension $i\raise1pt\hbox{\fivebf+}1$ in
the suspended chain.)}
\endproclaim

\demo{Proof}
The second isomorphism is the key and
is induced by the pair homeomorphism in \cite{2} %
5.7.4\ p.\ 164.
For the 2:nd last isomorphism we use \cite{15} %
p.\ 210 Corollary\ 5.7.9 and that {\bf LHS}-homomorphisms are
``chain map"-induced.
Note that the second component in the third module is an excisive
union.

$$\lower0pt\hbox{${ \mdoubleH\!_{_{\!q\!}}(X
\ \!\!{\rlap{{\lower2.5pt\hbox
{\vbox{\moveright0.1pt\hbox{$^{^{\land }}$}}}}}{{\ast}}}\ \!
Y_{^{_{^{\!}}}})\ \
{ \lower4.0pt\hbox{$^{_{\hbox{{\fivebf Z}}
\atop{{\raise0.0pt\hbox{{\eightbf =}}}
\!\!\!\!\!_{\!}{\raise2.8pt\hbox{{\eightsy {\char"27}}}} }}}$} }
 \ \ \mdoubleH\!\!\!_{_{\!q\hbox{\fivebf+}1\!}}\!(X\!_{}
\ \!{\rlap{{\lower2.5pt\hbox
{\vbox{\moveright0.1pt\hbox{$^{^{\land }}$}}}}}{{\ast}}}\ \!
Y\!_{}
\ \!{\rlap{{\lower2.5pt\hbox
{\vbox{\moveright0.1pt\hbox{$^{^{\land }}$}}}}}{{\ast}}}\ \!
\{\!\hbox{\tenbf v}_{_{\!}}\!,\wp\!\}, X\!_{}
\ \!{\rlap{{\lower2.5pt\hbox
{\vbox{\moveright0.1pt\hbox{$^{^{\land }}$}}}}}{{\ast}}}\ \!
Y\!_{})
}$}\ \
{\lower4.0pt\hbox{$^{_{\hbox{{\fivebf Z}}
\atop{{\raise0.0pt\hbox{{\eightbf =}}}
\!\!\!\!\!_{\!}{\raise2.8pt\hbox{{\eightsy {\char"27}}}} }}}$} }
$$

$$
{ \lower4.0pt\hbox{$^{_{\hbox{{\fivebf Z}}
\atop{{\raise0.0pt\hbox{{\eightbf =}}}
\!\!\!\!\!_{\!}{\raise2.8pt\hbox{{\eightsy {\char"27}}}} }}}$} }
\lower0pt\hbox{${ \ \! \mdoubleH\!\!\!_{_{\!q\hbox{\fivebf+}1\!}}
        ((\!X\!_{}
\ \!\!{\rlap{{\lower2.5pt\hbox
{\vbox{\moveright0.1pt\hbox{$^{^{\land }}$}}}}}{{\ast}}}\ \!\!
\{\!\hbox{\tenbf u}_{_{\!}},_{\!}\wp_{\!}\}\!) \!\times\!(Y\!_{} \
\!{\rlap{{\lower2.5pt\hbox {\vbox{\moveright0.1pt\hbox{$^{^{\land
}}$}}}}}{{\ast}}}\ \!\!\{\!\hbox{\tenbf v}_{_{\!}}\!,\wp\}\!),
\hbox{\tenbf(}(\!X\!_{} \ \!{\rlap{{\lower2.5pt\hbox
{\vbox{\moveright0.1pt\hbox{$^{^{\land }}$}}}}}{{\ast}}}\ \!
\{\!\hbox{\tenbf u}_{_{\!}},_{\!}\wp_{\!}\}) \!\times\!
Y\!_{}\hbox{\tenbf)} \!\cup\! \hbox{\tenbf(}\!X\!_{} \!\times\!
(Y\!_{} \ \!{\rlap{{\lower2.5pt\hbox
{\vbox{\moveright0.1pt\hbox{$^{^{\land }}$}}}}}{{\ast}}}\ \!
\{\hbox{\tenbf v}_{_{\!}}\!,\wp\})\hbox{\tenbf)} ) { =\ \! } }$}
$$

$$
{ =\ } \mdoubleH\!\!\!_{_{\!q+1\!}}
        ((\!X\!_{}
\ \!_{\!}{\rlap{{\lower2.5pt\hbox
{\vbox{\moveright0.1pt\hbox{$^{^{\land }}$}}}}}{{\ast}}}\ \!_{\!}
\{\hbox{\tenbf u}_{_{\!}},\wp\}, X) \!\times\! (Y\!_{} \
\!_{\!}{\rlap{{\lower2.5pt\hbox
{\vbox{\moveright0.1pt\hbox{$^{^{\land }}$}}}}}{{\ast}}}\ \!_{\!}
\{\hbox{\tenbf v}_{_{\!}}\!,\wp\}, Y) )\
{ \lower4.0pt\hbox{$^{_{\hbox{{\fivebf Z}}
\atop{{\raise0.0pt\hbox{{\eightbf =}}}
\!\!\!\!\!_{\!}{\raise2.8pt\hbox{{\eightsy {\char"27}}}} }}}$} }
$$

$$
{\lower4.0pt\hbox{$^{_{\hbox{{\fivebf Z}}
\atop{{\raise0.0pt\hbox{{\eightbf =}}}
\!\!\!\!\!_{\!}{\raise2.8pt\hbox{{\eightsy {\char"27}}}} }}}$} }
\ \bigg[{{{{\hbox{\eightrm %
Motivation: The$\ $underlying$\ $chains$\ $on$\ $the$\ $l.h.s.$\
$and$\ $r.h.s.$\ $are,
}}\atop {\hbox{\eightrm $\ $by$\ $Note$\ ${\eightbf ii}$\ $p.$\
$6,$\ $isomorphic$\ _{\!}$to$\ $their$\ $classical$\
$counterparts$\ $on}}
}}\atop {\hbox{\eightrm
$\ _{\!}$\ \ which$\ $we$\ $use$\ $the$\ $classical$\
$Eilenberg-Zilber$\ $Theorem.}\hfill}}\bigg]\
{\lower4.0pt\hbox{$^{_{\hbox{{\fivebf Z}}
\atop{{\raise0.0pt\hbox{{\eightbf =}}}
\!\!\!\!\!_{\!}{\raise2.8pt\hbox{{\eightsy {\char"27}}}} }}}$} }
$$

$$
{\lower4.0pt\hbox{$^{_{\hbox{{\fivebf Z}}
\atop{{\raise0.0pt\hbox{{\eightbf =}}}
\!\!\!\!\!_{\!}{\raise2.8pt\hbox{{\eightsy {\char"27}}}} }}}$} }
\ \mdoubleH\!\!\!_{_{\!q+1\!}}
        ({\Delta{\raise1.5pt\hbox{$\!\!^{\wp}$}}}(\!X\!_{}
\ \!{\rlap{{\lower2.5pt\hbox
{\vbox{\moveright0.1pt\hbox{$^{^{\land }}$}}}}}{{\ast}}}\ \!
\{\hbox{\tenbf u}_{_{\!}},\wp\}, X)
\!\otimes\!\!\!_{_{_{\hbox{\fivebf Z}}}}\! {
\Delta{\raise1.5pt\hbox{$\!\!^{\wp}$}}}(Y\!_{} \
\!{\rlap{{\lower2.5pt\hbox {\vbox{\moveright0.1pt\hbox{$^{^{\land
}}$}}}}}{{\ast}}}\ \! \{\hbox{\tenbf v}_{_{\!}}\!,\wp\}, Y))\
{\lower4.0pt\hbox{$^{_{\hbox{{\fivebf Z}}
\atop{{\raise0.0pt\hbox{{\eightbf =}}}
\!\!\!\!\!_{\!}{\raise2.8pt\hbox{{\eightsy {\char"27}}}} }}}$} }
$$

$$
\lower0pt\hbox{${
{\lower4.0pt\hbox{$^{_{\hbox{{\fivebf Z}}
\atop{{\raise0.0pt\hbox{{\eightbf =}}}
\!\!\!\!\!_{\!}{\raise2.8pt\hbox{{\eightsy {\char"27}}}} }}}$} }
\ \hbox{\tenrm H}\!\!_{_{\!q+1\!\!}}
        (\hbox{\tenbf s}{\Delta{\raise1.5pt\hbox{$\!\!^{\wp}$}}}\!
\hbox{\tenbf (} \!X\!_{}\hbox{\tenbf )}
\!\otimes\!\!\!_{_{_{\hbox{\fivebf Z}}}}\! \hbox{\tenbf s}{
\Delta{\raise1.5pt\hbox{$\!\!^{\wp}$}}}\! \hbox{\tenbf
(}Y\!_{}\hbox{\tenbf )})
{\lower4.0pt\hbox{$^{_{\hbox{{\fivebf Z}}
\atop{{\raise0.0pt\hbox{{\eightbf =}}}
\!\!\!\!\!_{\!}{\raise2.8pt\hbox{{\eightsy {\char"27}}}} }}}$} }
\ \hbox{\tenrm H}_{_{\!q\!\!}}
        (\hbox{\tenbf s}[{\Delta{\raise1.5pt\hbox{$\!\!^{\wp}$}}}\!
\hbox{\tenbf (} \!X\!_{}\hbox{\tenbf )}
\!\otimes\!\!\!_{_{_{\hbox{\fivebf Z}}}}\! {
\Delta{\raise1.5pt\hbox{$\!\!^{\wp}$}}}\! \hbox{\tenbf
(}Y\!_{}\hbox{\tenbf )}]) }$}.
$$

\smallskip
Now the non-relative Eilenberg-Zilber Theorem for\nobreak\ joins
follows from \cite{21} p.\ 279 Th.\ 46.2 above.
\hfill$\triangleright$

\medskip
Substituting, in the $\times$-original proof$\ \hbox{\tenrm \cite{25} %
p.}\ 234,$
``$\ \!{\rlap{{\lower2.5pt\hbox
{\vbox{\moveright0.1pt\hbox{$^{^{\land }}$}}}}}{{\ast}}}\ \!$",
``${\Delta{\raise1.5pt\hbox{$\!\!^{\wp}$}}}\!$/$\!\hbox{\tenbf s}
{\Delta{\raise1.5pt\hbox{$\!\!^{\wp}$}}}$", ``Theorem\ 3, 1:st
part" for ``$\times$", ``$\Delta$", ``Theorem 6" resp. will do
since;
$$\hbox{\tenbf s}\big({
\Delta{\raise1.5pt\hbox{$\!\!^{\wp}$}}}(X\!_{_1})\otimes
{\Delta{\raise1.5pt\hbox{$\!\!^{\wp}$}}}(Y\!_{_1})\big)/
\hbox{\twelvbf{\char"28}}\hbox{\tenbf s}\big({
\Delta{\raise1.5pt\hbox{$\!\!^{\wp}$}}}(X\!_{_1})\otimes
{\Delta{\raise1.5pt\hbox{$\!\!^{\wp}$}}}(Y\!_{_2})\big)+
\hbox{\tenbf s}\big({
\Delta{\raise1.5pt\hbox{$\!\!^{\wp}$}}}(X\!_{_2})\otimes
{\Delta{\raise1.5pt\hbox{$\!\!^{\wp}$}}}(Y\!_{_1})\big)
\hbox{\twelvbf{\char"29}}= $$

$$
=\! \hbox{\tenbf s}\big\{{
\Delta{\raise1.5pt\hbox{$\!\!^{\wp}$}}}\!(\!X\!_{_1}\!)\otimes
{\Delta{\raise1.5pt\hbox{$\!\!^{\wp}$}}}\!(Y\!\!_{_1}\!)/
\hbox{\twelvbf{\char"28}}\big({
\Delta{\raise1.5pt\hbox{$\!\!^{\wp}$}}}\!(\!X\!_{_1}\!)\otimes
{\Delta{\raise1.5pt\hbox{$\!\!^{\wp}$}}}\!(Y\!\!_{_2}\!)\big)+
\big({\Delta{\raise1.5pt\hbox{$\!\!^{\wp}$}}}\!(\!X\!_{_2}\!)\otimes
{\Delta{\raise1.5pt\hbox{$\!\!^{\wp}$}}}\!(Y\!\!_{_1}\!)\big)
\hbox{\twelvbf{\char"29}}\!\big\}
=
$$

\bigskip
\hfill
$
=
\hbox{\tenbf s}\big\{\hbox{\twelvbf{\char"28}}{
\Delta{\raise1.5pt\hbox{$\!\!^{\wp}$}}}\!(\!X\!_{_1}\!)/ {
\Delta{\raise1.5pt\hbox{$\!\!^{\wp}$}}}\!(\!X\!_{_2}\!)
\hbox{\twelvbf{\char"29}}\otimes \hbox{\twelvbf{\char"28}}{
\Delta{\raise1.5pt\hbox{$\!\!^{\wp}$}}}\!(Y\!\!_{_1}\!)/ {
\Delta{\raise1.5pt\hbox{$\!\!^{\wp}$}}}\!(Y\!\!_{_2}\!)\hbox{\twelvbf{\char"29}}\!\big\}.
\hfill
$
\hfill
\qed
\enddemo

\medskip%
\cite{25} %
Cor.\ 4 p.\ 231 now gives Theorem\ 4 since;
$$
\mdoubleH\lower0pt\hbox{$\!{_{_{\star}}}$}\!(\circ) \cong
\hbox{\tenbf s}
\mdoubleH\lower0.4pt\hbox{$\!\!\!{_{_{\star+1}}}$}\!(\circ) \cong
\mdoubleH\lower0.4pt\hbox{$\!\!\!{_{_{\star+1}}}$}\!(\hbox{\tenbf
s}(\circ))
\ \ \ \hbox{\tenrm and}\ \ \
{ \Delta{\raise1.5pt\hbox{$\!\!^{\wp}$}}}\! \hbox{\tenbf
(}(\!X\!_{_{\ \!}}\!,X\!{_{\!_2}}_{\!})\ast (Y\!\!_{_{\
\!}}\!,Y\!\!{_{\!_2}}_{\!})\hbox{\tenbf )}
\approx
{\Delta{\raise1.5pt\hbox{$\!\!^{\wp}$}}}\! \hbox{\tenbf
(}(\!X\!_{_{\ \!}}\!, X\!{_{_2}}_{\!}) \
\!{\rlap{{\lower2.5pt\hbox{\vbox{\moveright0.1pt\hbox{$^{^{\land
}}$}}}}}{{\ast}}}\ \! (Y\!\!_{_{\
\!}}\!,Y\!\!{_{_2}}_{\!})\hbox{\tenbf )}
$$

\noindent
by Theorem\ 2 and \cite{21} %
p.\ 279 Th.\ 46.2.
$\{X\!_{_{1}}\!\!_{\!}\ast\!
{\!}Y\!\!_{_{2\!}}, X\!_{_{2}}\!\!_{\!}\ast _{\!}\!Y\!\!_{_{1}}\}$
is excisive \underbar{iff}
$\{X\!_{_{1}}\! \!_{\!}\times_{\!} \!Y\!\!_{_{2\!}}, X\!_{_{2}}
\!\!_{\!}\times\! _{\!}Y\!\!_{_{1}}\!\}\!$
is, which is seen through a \hbox{\tenbf M-$\!$Vs}-stuffed 9-Lemma
and Theorem\ 2 (line\nobreak\ four).)

\proclaim{Theorem 4} \hbox{\tenrm ( The K\"{u}nneth Formula for
Toplogical Joins; cp. \cite{25} %
p.\ 235.)}

\smallskip
If
$\{X\!_{_{1}}
\ \!\!{\rlap{{\lower2.5pt\hbox
{\vbox{\moveright0.1pt\hbox{$^{^{\land }}$}}}}}{{\ast}}}\ \!
Y\!_{_{\!2}}\!, X\!_{_{2}}
\ \!\!{\rlap{{\lower2.5pt\hbox
{\vbox{\moveright0.1pt\hbox{$^{^{\land }}$}}}}}{{\ast}}}\ \!
Y\!\!_{_{1}}\}$
is an excisive couple in
$X\!_{_{1}}\!{\rlap{{\lower2.5pt\hbox
{\vbox{\moveright0.1pt\hbox{$^{^{\land }}$}}}}}{{\ast}}}
Y\!\!_{_{1}}\!,$
$\hbox{\tenbf R}\ a\ \hbox{\tenbf PID},\ \hbox{\tenbf G},
\!\hbox{\tenbf G}^\prime$\
\hbox{\tenbf R}-modules\nobreak\ and

\noindent
$
\hbox{\tenrm Tor}\!_1^{_{^{\hbox{\fivebf R}}}}\!(\hbox{\tenbf
G},\hbox{\tenbf G}^{\prime})=0
$,
then the functorial sequences below are
$($non-naturally$)$\nobreak\ split\nobreak\ exact;

\medskip%
\centerline{
$0\longrightarrow {\rlap{$_{_{_{i+j=q}}}$} {\
\raise2pt\hbox{$\bigoplus$}}}\ \ [\mdoubleH_{_{i}}
           (X\!_{_1},X\!_{_2};\hbox{\tenbf G})
\otimes\!\!_{_{\hbox{\fivebf R}}}
        \mdoubleH_{_{j}}
(Y\!\!_{_{1}},Y\!\!_{_{2}};\hbox{\tenbf G}^{\prime})]$ $
\longrightarrow$
}

\hfill${{\hbox{\tenrm(\hbox{\tenbf 3})}}}\!\!$

\medskip%
\centerline{ $\hbox to 0.4cm{\rightarrowfill}
\  \mdoubleH_{_{q+1}}
        ((\!X\!_{_1},X\!_{_2})\
{\rlap{{\lower2.5pt\hbox {\vbox{\moveright0.1pt\hbox{$^{^{\land
}}$}}}}}{{\ast}}}\ (Y\!_{_1},Y\!\!_{_{2}}); \hbox{\tenbf
G}\otimes\!\!_{_{\hbox{\fivebf R}}}\hbox{\tenbf G}^{\prime})
\hbox to 0.4cm{\rightarrowfill}
{{\rlap{$\!\!\!\!_{_{_{i+j=q-1}}}$}}{\raise2pt\hbox{$\bigoplus$}}}\!\!
\ \ \hbox{\tenrm Tor}_1^{\hbox{\fivebf R}}\bigl(\mdoubleH_{_{i}}
           (X\!_{_1},X\!_{_{\!2}};\hbox{\tenbf G}),
 \mdoubleH_{_{j}}
      (Y\!\!_{_1},Y\!\!_{_{2}};
\hbox{\tenbf G}^{\prime})\bigr)
\hbox to 0.4cm{\rightarrowfill}\ 0$
}
\endproclaim

Analogously with $``\ast"$ substituted for $ _{\!}``
{\rlap{{\lower2.5pt\hbox
{\vbox{\moveright0.1pt\hbox{$^{^{\land}}$}}}}}{{\ast}}}"_{\!}
$ and \cite{25} %
p.\ 247 Th.\ 11 gives the co$\mdoubleH$omology-analog.\qed

\medskip
Putting
$(\!X\!_{_1}\!,X\!_{_2}\!)\!=\!(\!\{_{\!}\wp_{\!}\},\emptyset)$
in\ Theorem\ 4, our next theorem immediately follows.

\proclaim{Theorem 5}
\hbox{\tenrm[The Universal Coefficient Theorem for (co)\doubleH
omology.]}

$$ \mdoubleH_{_{i\!}} (Y\!\!_{_1},Y\!\!_{_2};\hbox{\tenbf G})\
{\lower4.0pt\hbox{$^{_{\hbox{{\fivebf R}}
\atop{{\raise0.0pt\hbox{{\eightbf =}}}
\!\!\!\!\!_{\!}{\raise2.8pt\hbox{{\eightsy {\char"27}}}} }}}$} }
\ \big[\hbox{\tenrm \cite{25} %
p.\ 214}\big]\
{\lower4.0pt\hbox{$^{_{\hbox{{\fivebf R}}
\atop{{\raise0.0pt\hbox{{\eightbf =}}}
\!\!\!\!\!_{\!}{\raise2.8pt\hbox{{\eightsy {\char"27}}}} }}}$} }
\ \mdoubleH_{_{i\!}} (Y\!\!_{_1},Y\!\!_{_2};\hbox{\tenbf
R}\otimes\!_{\!_{\hbox{\fivebf R}}}\!\hbox{\tenbf G})\
{\lower4.0pt\hbox{$^{_{\hbox{{\fivebf R}}
\atop{{\raise0.0pt\hbox{{\eightbf =}}}
\!\!\!\!\!_{\!}{\raise2.8pt\hbox{{\eightsy {\char"27}}}} }}}$} }
$$

$$
{\lower4.0pt\hbox{$^{_{\hbox{{\fivebf R}}
\atop{{\raise0.0pt\hbox{{\eightbf =}}}
\!\!\!\!\!_{\!}{\raise2.8pt\hbox{{\eightsy {\char"27}}}} }}}$} }
\
\hbox{\twelvbf{\char"28}}\mdoubleH_{_{i\!}}(Y\!\!_{_1},Y\!\!_{_2};\hbox{\tenbf
R}) \otimes\!_{\!_{\hbox{\fivebf R}}}\!\hbox{\tenbf
G}\hbox{\twelvbf{\char"29}}\oplus \hbox{\tenrm
Tor}^{\!\hbox{\fivebf R}}_1\! \big(\mdoubleH_{_{{i\!-\!1}}}\!
(Y\!\!_{_1},Y\!\!_{_2};\hbox{\tenbf R}), \hbox{\tenbf G}\big),
$$

\noindent for {any }{\bf R}-{\bf PID} {module}\ {\bf G}.

\medskip
If\ all\ $\mdoubleH_{{\ast}}\!
        (Y\!_{_1},Y\!_{_2};\hbox{\tenbf R})$
are of finite type or {\bf G} is finitely generated, then;

\bigskip
\noindent { $ \mdoubleH^{{i\!}}
(Y\!\!_{_1},Y\!\!_{_2};\hbox{\tenbf G})
{\lower4.0pt\hbox{$^{_{\hbox{{\fivebf R}}
\atop{{\raise0.0pt\hbox{{\eightbf =}}}
\!\!\!\!\!_{\!}{\raise2.8pt\hbox{{\eightsy {\char"27}}}} }}}$} }
\mdoubleH^{{i\!}} (Y\!\!_{_1},Y\!\!_{_2};\hbox{\tenbf
R}\otimes\!_{\!_{\hbox{\fivebf R}}}\!\hbox{\tenbf G})
{\lower4.0pt\hbox{$^{_{\hbox{{\fivebf R}}
\atop{{\raise0.0pt\hbox{{\eightbf =}}}
\!\!\!\!\!_{\!}{\raise2.8pt\hbox{{\eightsy {\char"27}}}} }}}$} }
\hbox{\twelvbf{\char"28}}\mdoubleH^{{i\!}}(Y\!\!_{_1},Y\!\!_{_2};\hbox{\tenbf
R}) \otimes\!_{\!_{\hbox{\fivebf R}}}\!\hbox{\tenbf
G}\hbox{\twelvbf{\char"29}}\oplus \hbox{\tenrm
Tor}^{\!\hbox{\fiverm R}}_1\!
\big(\mdoubleH^{{{i{\raise1pt\hbox{\fivebf{+}}1}}}}\!
(Y\!\!_{_1},Y\!\!_{_2};\hbox{\tenbf R}), \hbox{\tenbf G}\big).
\!\!\!\!$\qed}
\endproclaim

%


\subhead 3.3. Local Augmental Homology Groups for Products and
Joins
\endsubhead

\bigskip%
Proposition 1 is our key motivation for introducing a topological
(-1)-object, which then imposed the definition p.\ 4, of a
``setminus", ``$\setminus\!{_{_{^{o\!}}}}$", in ${\hbox{\tensy
D}\!_{\wp}}$, revealing the true implication of boundary
definitions w.r.t. manifolds as \hbox{given in pp.\ 13\ +\ 23.}
Somewhat specialized, it's found in
\cite{10} %
p.\ 162 %
and partially also in \cite{22} %
p.\ 116 Lemma 3.3.
``${ X}\setminus x$" usually stands for ``${
X}\setminus\{x\}"_{\!}$ and we'll write $x$ for $\{x,\wp\}$ as a
notational convention.
Recall the definition of $\alpha_{_{^{\!0}}}\!$, p.\ 5 and that
\hbox{$\dim\hbox{\tenrm
Lk}\!\lower1.2pt\hbox{$_{_{\!\Sigma_{_{\!{\
}}}}}$}\!\!\!\sigma\!=\!{\dim\!{\Sigma
\!_{\!}-\!
\hbox{\eightbf\#} \sigma\!}}.$}

\proclaim{Proposition 1}
Let $\hbox{\tenbf G}$ be any module over a commutative ring {\bf
A}$_{\!}$ with unit.\ With $\alpha\ \!\!\in\ \!\!\hbox{\tenrm
Int}_{\!}\sigma$ and $\ {\!}\alpha_{\!}=_{\!}
\alpha_{_{^{\!0}}}\!$ {\underbar{\hbox{\tenrm iff}}}\
$\sigma\!=\!\emptyset\!_{_{^{o}}}{\!}$
the following module isomorphisms are all induced by chain
equivalences,
{\rm cf. \cite{21} %
p.\ 279 Th.\ 46.2 quoted here in p.\ 9}.
$$ \mdoubleH
\vbox{\moveleft0.0cm\hbox{\lower0.0pt\hbox{$_{{{_{\!}{i
\lower1.0pt\hbox{-} \hbox{\fivebf\#}_{\!}\sigma}}}}$}}}
\vbox{\moveleft0.05cm\hbox{$(\hbox{\tenrm
Lk}_{_{\!\Sigma}}\sigma;\hbox{\tenbf G})$}}\ \!
\raise1.0pt\hbox{\lower4.0pt\hbox{$^{_{\hbox{{\fivebf A}}
\atop{{\raise0.0pt\hbox{{\eightbf =}}}
\!\!\!\!\!{\!}{\raise2.8pt\hbox{{\eightsy {\char"27}}}}}}}$}}\ \!
\mdoubleH_{i}(\Sigma,\hbox{\tenrm
cost}_{_{\!\Sigma}}\sigma;\hbox{\tenbf G})\
\raise1.0pt\hbox{\lower4.0pt\hbox{$^{_{\hbox{{\fivebf A}}
\atop{{\raise0.0pt\hbox{{\eightbf =}}}
\!\!\!\!\!{\!}{\raise2.8pt\hbox{{\eightsy {\char"27}}}}}}}$}}\ \!
\mdoubleH_{i}(\vert\Sigma\vert,\vert\hbox{\tenrm
cost}_{_{\!\Sigma}}\sigma\vert;\hbox{\tenbf G})\ \!
\raise1.0pt\hbox{\lower4.0pt\hbox{$^{_{\hbox{{\fivebf A}}
\atop{{\raise0.0pt\hbox{{\eightbf =}}}
\!\!\!\!\!{\!}{\raise2.8pt\hbox{{\eightsy {\char"27}}}}}}}$}}\ \!
\mdoubleH_{i}(|\Sigma|,|\Sigma|\setminus\!{_{_{^{\!o\!}}}}\alpha;\hbox{\tenbf
G}),
$$
$$ \mdoubleH
^{_{\!}i\raise0.5pt\hbox{\fiverm-}\hbox{\fivebf\#}_{\!}\sigma}
\!_{\!}(\hbox{\tenrm Lk}_{_{\!\Sigma}}\sigma;\hbox{\tenbf G})\
\raise1.0pt\hbox{\lower4.0pt\hbox{$^{_{\hbox{{\fivebf A}}
\atop{{\raise0.0pt\hbox{{\eightbf =}}}
\!\!\!\!\!{\!}{\raise2.8pt\hbox{{\eightsy {\char"27}}}}}}}$}}\ \!
\mdoubleH^{i}(\Sigma,\hbox{\tenrm
cost}_{_{\!\Sigma}}\sigma;\hbox{\tenbf G}) \
\raise1.0pt\hbox{\lower4.0pt\hbox{$^{_{\hbox{{\fivebf A}}
\atop{{\raise0.0pt\hbox{{\eightbf =}}}
\!\!\!\!\!{\!}{\raise2.8pt\hbox{{\eightsy {\char"27}}}}}}}$}}\ \!
\mdoubleH^{i}(\vert\Sigma\vert,\vert\hbox{\tenrm
cost}_{_{\!\Sigma}}\sigma\vert;\hbox{\tenbf G})
\raise1.0pt\hbox{\lower4.0pt\hbox{$^{_{\hbox{{\fivebf A}}
\atop{{\raise0.0pt\hbox{{\eightbf =}}}
\!\!\!\!\!{\!}{\raise2.8pt\hbox{{\eightsy {\char"27}}}}}}}$}}\ \!
\mdoubleH^{i}(|\Sigma|,|\Sigma|\setminus\!{_{_{^{\!o\!}}}}\alpha;\hbox{\tenbf
G}).
$$
\endproclaim

\demo{Proof} (Cf. $\!$definitions p.\ 30-1.) The
$``\setminus_o\!\!"\!$-definition p.\ 4 and
\cite{21} %
Th.\ $\!$46.2 p.\ 279 $\!$\raise1pt\hbox{\eightbf+} $\!$pp.\
$_{\!}$ 194-199 Lemma 35.1-35.2 $\!$\raise1pt\hbox{\eightrm+}$\!$
Lemma\ 63.1 p.\ 374 gives the two ending isomorphisms since
$ \vert\hbox{\tenrm cost}_{\Sigma}\sigma\vert $
is a deformation retract of
$ \vert\Sigma\vert\!\setminus\!_{_{o}}\alpha, $
$\hbox{\tenrm while}$ already on the chain level;
\noindent$C^{o}\!\!\!_{_{\star}}(\Sigma,\hbox{\tenrm
cost}_{\Sigma}\sigma) \ \!=\!\
C^{o}\!\!\!_{_{\star}}(\overline{\hbox{\tenrm {st}}}
_{_{\!\Sigma}}(\sigma),\ \! {
\dot{\sigma}}_{_{\!}}\ast_{_{\!}}\hbox{\tenrm
Lk}_{_{\!\Sigma}}\sigma)
\ \!=\ \! C^{o}\!\!\!_{_{\star}}(\overline{{{\sigma}}}
{_{_{\!}}}\ast{_{_{\!}}}\hbox{\tenrm Lk}_{_{\!\Sigma}}\!\sigma,\
\! {\dot{\sigma}}{_{_{\!}}}\ast{_{_{\!}}}\hbox{\tenrm
Lk}_{_{\!\Sigma}}\!\sigma) \ \!\simeq\
C^{o}\!\!\!\!\!\!\!\!\! _{_{\star\lower1pt\hbox{-}\#\sigma}}
(\hbox{\tenrm Lk}_{_{\!\Sigma}}\sigma). $
\qed
\enddemo

\proclaim{Lemma}
If $x\ (y)$ is a closed point in $X\ (Y)$, then
$\!\{\!X\!_{_{}}\!\!\times\! (Y\!\!_{_{\
}}\!\!\setminus\!\!{_{_{^{o\!}}}} y), (X\!_{_{\
}}\!\!\setminus\!\!{_{_{^{o\!}}}} x)\!\times\!Y\!\!_{_{\ }}\!\},$
and
$\!\{\!X\!_{_{}}\!\ast (Y\!\!_{_{\
}}\!\!\setminus\!\!{_{_{^{o\!}}}}  y), (X\!_{_{\
}}\!\!\setminus\!\!{_{_{^{o\!}}}}  x)\ast Y\!\!_{_{\ }}\!\}\
$%
are both excisive pairs.
\endproclaim

\demo{Proof} \cite{25} %
p.\ 188 Th.\ 3, since $X\!_{_{\ }}\!\!\times\! (Y\!\!_{_{\
}}\!\setminus\!\!{_{_{^{o\!}}}}  y\!_{_{0}}\!)$
 $\big(X\!_{_{\ }}\!\!\ast\!
(Y\!\!_{_{\ }}\!\setminus\!\!{_{_{^{o\!}}}}  y\!_{_{0}}\!)\big)$
is open in $(X\!_{_{\ }}\!\!\times\! (Y\!\!_{_{\
}}\!\setminus\!\!{_{_{^{o\!}}}} y\!_{_{0}}\!)) \cup ((X\!\!_{_{\
}}\!\setminus\!\!{_{_{^{o\!}}}}  x\!_{_{0}}\!)\times\! Y\!_{_{\
}}\!\!) $
$\big(\!(X\!_{_{\ }}\!\!\ast (Y\!\!_{_{\
}}\!\setminus\!{_{_{^{o\!}}}} y\!_{_{0}}\!)) \cup ((X\!\!_{_{\
}}\!\setminus\!{_{_{^{o\!}}}}  x\!_{_{0}}\!)\ast Y\!_{_{\
}}\!\!)\!\big) $, which proves the excisivity.
\qed
\enddemo

\proclaim{Theorem 6}
For $x\!_{\!}\in \!\!X\ (y\!\in \!Y) \ {closed}$ {and}
$(t\!{_{_{^{1}}}}\!,\widetilde{x_{{\!}}\ast_{{\!}}y},
t\!{_{_{^{2}}}}\! )\!:=\!\{\!(x,y,t)\ \!{|}\
\!0_{{\!}}<_{{\!}}t\!{_{_{^{1}}}}\!\!\le \!t\!\le\!
t\!{_{_{^{2}}}}\!\!<\!1\}${\rm;}

\medskip
\noindent \hbox{\tenbf i.}\ \ \ $ \mdoubleH\!\!\!_{_{\!q+1\!}}({X}
\ \!{\rlap{{\lower2.5pt\hbox
{\vbox{\moveright0.1pt\hbox{$^{^{\land }}$}}}}}{{\ast}}}\ \!
{Y}\!,{X} \ \!{\rlap{{\lower2.5pt\hbox
{\vbox{\moveright0.1pt\hbox{$^{^{\land }}$}}}}}{{\ast}}}\ \! {Y}
\setminus\!{_{_{^{o\!}}}}
(x\!{_{_{^{}}}},y\!{_{_{^{}}}},t);\hbox{\tenbf G})
\ \raise1.0pt\hbox{\lower4.0pt\hbox{$^{_{\hbox{{\fivebf A}}
\atop{{\raise0.0pt\hbox{{\eightbf =}}}
\!\!\!\!\!{\!}{\raise2.8pt\hbox{{\eightsy {\char"27}}}} }}}$} }\
\mdoubleH\!\!\!_{_{\!q+1\!}}({X} \ \!{\rlap{{\lower2.5pt\hbox
{\vbox{\moveright0.1pt\hbox{$^{^{\land }}$}}}}}{{\ast}}}\ \!
{Y}\!,{X} \ \!{\rlap{{\lower2.5pt\hbox
{\vbox{\moveright0.1pt\hbox{$^{^{\land }}$}}}}}{{\ast}}}\ \!
{Y}\setminus\!{_{_{^{o}}}}\ \!(t\!{_{_{^{1}}}},\widetilde{x \ast
y},t\!{_{_{^{2}}}});\hbox{\tenbf G}) \
\ \raise1.0pt\hbox{\lower4.0pt\hbox{$^{_{\hbox{{\fivebf A}}
\atop{{\raise0.0pt\hbox{{\eightbf =}}}
\!\!\!\!\!{\!}{\raise2.8pt\hbox{{\eightsy {\char"27}}}} }}}$} }\
$
$$
\raise1.0pt\hbox{\lower4.0pt\hbox{$^{_{\hbox{{\fivebf A}}
\atop{{\raise0.0pt\hbox{{\eightbf =}}}
\!\!\!\!\!{\!}{\raise2.8pt\hbox{{\eightsy {\char"27}}}} }}}$} }\
\mdoubleH_{_{\!q\!}}({X}\!\times\!{Y}\!,{X}\!\times\!{Y}
\setminus\!{_{_{^{o}}}}\! (x,y);\hbox{\tenbf G})
\ \raise1.0pt\hbox{\lower4.0pt\hbox{$^{_{\hbox{{\fivebf A}}
\atop{{\raise0.0pt\hbox{{\eightbf =}}}
\!\!\!\!\!{\!}{\raise2.8pt\hbox{{\eightsy {\char"27}}}} }}}$} }\
\Big[ {{\lower3pt\hbox{\eightrm A\ simple\
}}\atop{\raise3pt\hbox{\eightrm calculation}}} \Big]
\ \raise1.0pt\hbox{\lower4.0pt\hbox{$^{_{\hbox{{\fivebf A}}
\atop{{\raise0.0pt\hbox{{\eightbf =}}}
\!\!\!\!\!{\!}{\raise2.8pt\hbox{{\eightsy {\char"27}}}} }}}$} }\
{ \mdoubleH_{_{\!q\!}}(
        (X\!{_{_{^{}}}}\!,
X\!{_{_{^{}}}}\!\setminus\!\!{_{_{^{o\!}}}} x\!{_{_{^{}}}})
\!\times\! (Y\!{_{_{^{}}}}\!,Y\!{_{_{^{}}}}
\!\setminus\!\!{_{_{^{o\!}}}} y\!{_{_{^{}}}});\hbox{\tenbf G}) }\
\ \raise1.0pt\hbox{\lower4.0pt\hbox{$^{_{\hbox{{\fivebf A}}
\atop{{\raise0.0pt\hbox{{\eightbf =}}}
\!\!\!\!\!{\!}{\raise2.8pt\hbox{{\eightsy {\char"27}}}} }}}$} }
$$
$$
\raise1.0pt\hbox{\lower4.0pt\hbox{$^{_{\hbox{{\fivebf A}}
\atop{{\raise0.0pt\hbox{{\eightbf =}}}
\!\!\!\!\!{\!}{\raise2.8pt\hbox{{\eightsy {\char"27}}}} }}}$} }\
\Big[ {{\lower3pt\hbox{\eightrm Th.\ \!2\ p.\ 8
}}\atop{\raise3pt\hbox{\eightrm line\ four}}} \Big]
\ \raise1.0pt\hbox{\lower4.0pt\hbox{$^{_{\hbox{{\fivebf A}}
\atop{{\raise0.0pt\hbox{{\eightbf =}}}
\!\!\!\!\!{\!}{\raise2.8pt\hbox{{\eightsy {\char"27}}}} }}}$} }\
\mdoubleH\!\!\!_{_{\!q+1\!}}(
        (X\!{_{_{^{}}}}\!,
X\!{_{_{^{}}}}\!\setminus\!\!{_{_{^{o\!}}}}   x\!{_{_{^{}}}})
 \
\!{\rlap{{\lower2.5pt\hbox {\vbox{\moveright0.1pt\hbox{$^{^{\land
}}$}}}}}{{\ast}}}\ \!
(Y\!{_{_{^{}}}}\!,Y\!{_{_{^{}}}}\!\setminus\!\!{_{_{^{o\!}}}}
y\!{_{_{^{}}}});\hbox{\tenbf G}).
$$
\noindent {\bf ii.}\ \
$ \hskip1.5cm \mdoubleH\!\!\!_{_{\!q+1\!}}
        (\!X\!_{}
\ \!{\rlap{{\lower2.5pt\hbox
{\vbox{\moveright0.1pt\hbox{$^{^{\land }}$}}}}}{{\ast}}}\ \!
Y\!_{}, X\!_{} \ \!{\rlap{{\lower2.5pt\hbox
{\vbox{\moveright0.1pt\hbox{$^{^{\land }}$}}}}}{{\ast}}}\ \!
Y\!_{} \setminus\!{_{_{^{o}}}}(y\!{_{_{^{}}}},{0});\hbox{\tenbf
G})\
\ \raise1.0pt\hbox{\lower4.0pt\hbox{$^{_{\hbox{{\fivebf A}}
\atop{{\raise0.0pt\hbox{{\eightbf =}}}
\!\!\!\!\!{\!}{\raise2.8pt\hbox{{\eightsy {\char"27}}}} }}}$} }\
\mdoubleH\!\!\!_{_{\!q+1\!}} ((X\!{_{_{^{}}}}\!,\emptyset)\
\!{\rlap{{\lower2.5pt\hbox {\vbox{\moveright0.1pt\hbox{$^{^{\land
}}$}}}}}{{\ast}}}\ \!
(Y\!{_{_{^{}}}}\!,Y\!{_{_{^{}}}}\!\setminus\!\!{_{_{^{o}}}}
y\!{_{_{^{}}}}\!\ ) ;\hbox{\tenbf G})
$

\medskip
and equivalently for the $(x\!{_{_{^{}}}},{1})$-points.

\smallskip
$_{\!}$All isomorphisms are induced by chain equivalences, {\rm
cf. \cite{21} %
p.\ 279 Th.\ 46.2 quoted here in p.\ 9}.
Analogously for $``\ast"$ substituted\ for
$``{\rlap{{\lower2.5pt\hbox
{\vbox{\moveright0.1pt\hbox{$^{^{\land}}$}}}}}{{ \ast}}}"_{\!} $
and for co$\mdoubleH$omology.
\endproclaim

\demo{Proof} {\bf i.}
\indent $\cases
& %
\!\!\!\!\!\! A:=X\!{_{_{^{}}}}{\myPsqcup1}\ Y\!{_{_{^{}}}}\
\!\!\setminus\!{_{_{^{\!o}}}}\ \!\!
\{(x\!{_{_{^{0}}}},y\!{_{_{^{0}}}},t)\ \!\vert\
\!t\!{_{_{^{1}}}}\le t<1 \}
\cr
&
\!\!\!\!\!\!B:=X\!{_{_{^{}}}}{\myPsqcup0}\ Y\!{_{_{^{}}}}
\!\setminus\!{_{_{^{\!o}}}}\ \!\!
\{(x\!{_{_{^{0}}}},y\!{_{_{^{0}}}},t)\ \!\vert\ \! 0<t\le
t\!{_{_{^{2}}}}\}
\cr
\endcases
\indent
\Longrightarrow $

\medskip
$
\hskip4.0cm
\Longrightarrow\!
\indent
\cases
&%
\!\!\!\!\!\! A\cup B =\!X\!_{} \ \!{\rlap{{\lower2.5pt\hbox
{\vbox{\moveright0.1pt\hbox{$^{^{\land }}$}}}}}{{\ast}}}\ \!
Y\!_{}\setminus\!{_{_{^{\!o}}}}\ \!(t\!{_{_{^{1}}}},\widetilde{x
\ast y}, t\!{_{_{^{2}}}} )
\cr & \!\!\!\!\!\! A\cap
B=\!X\!{_{_{^{}}}}\!\times\!Y\!{_{_{^{}}}}\!\times\!(0,1)
\!\setminus\!{_{_{^{\!o}}}}\ \!\! \{x\!{_{_{^{0}}}}\}\!\times\!
\{y\!{_{_{^{0}}}}\!\}\!\times\!(0,1)\hbox{\tenbf,}
\cr
\endcases
$

\medskip
\noindent with\ \
$ x\!{_{_{^{0}}}}\times y\!{_{_{^{0}}}}\times(0,1)\!:=
\!\{x\!{_{_{^{0}}}}\!\}\times \{y\!{_{_{^{0}}}}\!\}\times \{t\
\!\vert\ \! t\in(0,1) \}\cup\{\wp\}
\ \hbox{\tenrm and}\
(x\!{_{_{^{0}}}},y\!{_{_{^{0}}}},t\!{_{_{^{0}}}})\!:=\!
\{(x\!{_{_{^{0}}}},y\!{_{_{^{0}}}},t\!{_{_{^{0}}}}),\wp\}.
$

\medskip
Now, using the null-homotopy in the relative {\bf M-$_{\!}$Vs}
w.r.t. $ \{ ( X\!{_{_{^{}}}}\!{\myPsqcup1}\ Y{_{_{^{}}}}\!\!,A),
(X\!{_{_{^{}}}}\!{\myPsqcup0}\ Y\!{_{_{^{}}}}\!,B)_{\!}\}$ and the
resulting splitting of it and the involved pair deformation
retractions as in the proof of Th.\ $\!$2, we get;

\noindent $
{\mdoubleH\!\!\!_{_{\!q+1\!}}
        (\!X\!_{}
\ \!{\rlap{{\lower2.5pt\hbox
{\vbox{\moveright0.1pt\hbox{$^{^{\land }}$}}}}}{{\ast}}}\ \!
Y\!_{}, X\!_{}
\ \!{\rlap{{\lower2.5pt\hbox
{\vbox{\moveright0.1pt\hbox{$^{^{\land }}$}}}}}{{\ast}}}\ \!
Y\!_{} \ \!\setminus\!{_{_{^{\!o}}}}\
\!\!(t\!{_{_{^{1}}}},\widetilde{x\!{_{_{^{0}}}}\! \ast
y\!{_{_{^{0}}}}\!}, t\!{_{_{^{2}}}} ) } ) \ \!
{\lower4.0pt\hbox{$^{_{\hbox{{\fivebf Z}}
\atop{{\raise0.0pt\hbox{{\eightbf =}}}
\!\!\!\!\!_{\!}{\raise2.8pt\hbox{{\eightsy {\char"27}}}} }}}$} }
\ \mdoubleH_{_{\!q\!}}
        (X\!{_{_{^{}}}}\!\times\!Y\!{_{_{^{}}}}\!\times\!(0,1),
X\!{_{_{^{}}}}\!\times\!Y\!{_{_{^{}}}}\!\times\!(0,1) \
\!\setminus{_{_{^{\!o}}}}\ \! \{x\!{_{_{^{0}}}}\!\}\!\times\!
\{y\!{_{_{^{0}}}}\!\}\!\times\!(0,1)\ ) \ \!
{\lower4.0pt\hbox{$^{_{\hbox{{\fivebf Z}}
\atop{{\raise0.0pt\hbox{{\eightbf =}}}
\!\!\!\!\!_{\!}{\raise2.8pt\hbox{{\eightsy {\char"27}}}} }}}$} }
$
$$
{\lower4.0pt\hbox{$^{_{\hbox{{\fivebf Z}}
\atop{{\raise0.0pt\hbox{{\eightbf =}}}
\!\!\!\!\!_{\!}{\raise2.8pt\hbox{{\eightsy {\char"27}}}} }}}$} }
\ \Big[ {{
\lower3.5pt\hbox{\eightrm Motivation: The underlying pair on the
r.h.s.}\atop{
{\raise3.5pt\hbox{\sevenrm is a pair deformation retract of that
on the l.h.s.$\!$}}}}}
\Big]\
{\lower4.0pt\hbox{$^{_{\hbox{{\fivebf Z}}
\atop{{\raise0.0pt\hbox{{\eightbf =}}}
\!\!\!\!\!_{\!}{\raise2.8pt\hbox{{\eightsy {\char"27}}}} }}}$} }
$$
$$
{\lower4.0pt\hbox{$^{_{\hbox{{\fivebf Z}}
\atop{{\raise0.0pt\hbox{{\eightbf =}}}
\!\!\!\!\!_{\!}{\raise2.8pt\hbox{{\eightsy {\char"27}}}} }}}$} }
\ \mdoubleH_{_{\!q\!}}
        (X\!{_{_{^{}}}}\!\times\!Y\!{_{_{^{}}}}\!\times\!
\{t\!{_{_{^{0}}}}\!,\wp\},
X\!{_{_{^{}}}}\!\times\!Y\!{_{_{^{}}}}\!\times\!\{t\!{_{_{^{0}}}}\!,\wp\}
\ \!\setminus\!{_{_{^{o}}}}\ \!\! (x\!{_{_{^{0}}}}\!,\!
y\!{_{_{^{0}}}}\!,\!t\!{_{_{^{0}}}})\ )\
{\lower4.0pt\hbox{$^{_{\hbox{{\fivebf Z}}
\atop{{\raise0.0pt\hbox{{\eightbf =}}}
\!\!\!\!\!_{\!}{\raise2.8pt\hbox{{\eightsy {\char"27}}}} }}}$} }
\ {\mdoubleH_{_{\!q\!}}
        (X\!{_{_{^{}}}}\!\times\!Y\!{_{_{^{}}}}\!,
X\!{_{_{^{}}}}\!\times\!Y\!{_{_{^{}}}}\!\setminus\!{_{_{^{\!o}}}}\!
(x\!{_{_{^{0}}}}\!,\! y\!{_{_{^{0}}}}\!)\ {\!}\!) }=
$$
$$
=
\mdoubleH_{_{\!q\!}}
        (X\!{_{_{^{}}}}\!\times\!Y\!{_{_{^{}}}}\!,
(X\!{_{_{^{}}}}\!\times\!(Y\!{_{_{^{}}}}\!\setminus\!\!{_{_{^{o}}}}
y\!{_{_{^{0}}}}\!)) \cup
((X\!{_{_{^{}}}}\!\setminus\!\!{_{_{^{o}}}}x\!{_{_{^{0}}}}\!)
\!\times\!Y\!{_{_{^{}}}})\! \ ) \ \!=
{ \mdoubleH_{_{\!q\!}}(
        (X\!{_{_{^{}}}}\!,
X\!{_{_{^{}}}}\!\setminus\!\!{_{_{^{o}}}}x\!{_{_{^{0}}}}\!)
\!\times\!
(\!Y\!{_{_{^{}}}}\!,Y\!{_{_{^{}}}}\!\setminus\!\!{_{_{^{o}}}}
y\!{_{_{^{0}}}}\!) \ ). } \eqno{\hbox{$\triangleright$}}
$$

$$
\leqno{\hbox{\tenbf ii.}}
\indent
\cases &
\!\!\!\!\!\!\! A:=X\!{_{_{^{}}}}{\myPsqcup1}\
Y\!{_{_{^{}}}}\hbox{\tenbf,}
\cr
&
\!\!\!\!\!\!\!B:=X\!{_{_{^{}}}}{\myPsqcup0}\ Y\!{_{_{^{}}}}
\!\setminus\!{_{_{^{\!o}}}}\ \!\!
X\!{_{_{^{}}}}\!\times\!\{y\!{_{_{^{0}}}}\!\}\!\times\![0,1)
\cr
\endcases
 \ \ \Longrightarrow  \ \
\cases
&
\!\!\!\!\!\!\! A\cup B =\!X\!_{} \ \!{\rlap{{\lower2.5pt\hbox
{\vbox{\moveright0.1pt\hbox{$^{^{\land }}$}}}}}{{\ast}}}\ \!
Y\!_{}\ \!\setminus\!{_{_{^{\!o}}}}\ \!(y\!{_{_{^{0}}}},0)
\cr
&
\!\!\!\!\!\!\! A\cap B= X\!{_{_{^{}}}}\!\times\!
(Y\!{_{_{^{}}}}\!\setminus\!{_{_{^{\!o}}}}\ \!\!
y\!{_{_{^{0}}}})\!\times\!(0,1)
\cr
\endcases
\ \ \hbox{\tenrm where}
$$
$ (x\!{_{_{^{0}}}}, y\!{_{_{^{0}}}},t)\in
X\!{_{_{^{}}}}\!\times\!\{y\!{_{_{^{0}}}}\!\}\!\times\![0,1) $
is independent of $x\!{_{_{^{0}}}}$ and $
(x\!{_{_{^{0}}}},y\!{_{_{^{0}}}},t\!{_{_{^{0}}}}\!)\!:=
\{(x\!{_{_{^{0}}}},y\!{_{_{^{0}}}},t\!{_{_{^{0}}}}\!), \wp\}. $

\smallskip
Now use Th.\ 1 p.\ 8 line\ 2 and that the r.h.s. is a pair
deformation retract of the l.h.s.;
$
(X\!{_{_{^{}}}}\!\times\!Y\!{_{_{^{}}}}\!\times\!(0,1) ,
X\!{_{_{^{}}}}\!\times\! (Y\!{_{_{^{}}}}\
\!\setminus{_{_{^{\!o}}}}\ \! y\!{_{_{^{0}}}}\!)\!\times\!(0,1)) \
\lower4pt\hbox{$\widetilde{\phantom{{..}}}$} \
(X\!{_{_{^{}}}}\!\times\!Y\!{_{_{^{}}}} ,
X\!{_{_{^{}}}}_{\!}\times_{\!} (Y\!{_{_{^{}}}}\
\!\setminus{_{_{^{\!o}}}}\ \! y\!{_{_{^{0}}}})) =
(X\!{_{_{^{}}}}\!,\emptyset) \!\times\!
(Y\!{_{_{^{}}}}\!,Y\!{_{_{^{}}}}\!\ \!\setminus\!{_{_{^{o}}}}\
\!y\!{_{_{^{0}}}}\!).
{\qed}
$
\enddemo

\subhead
3.4. Singular Homology Manifolds under Products and Joins
\endsubhead

\normalbaselines

\definition{Definition}
$\emptyset$ is a {weak} {homology$_{_{\!\hbox{\fivebf G}}}\!$
\raise0.8pt\hbox{\sixsy {\char"00}}$\infty$-manifold}.
Else, a
\hbox{\tenbf
T}${\!}\lower1.5pt\hbox{\fivebf{\char"31}}_{\!}$-\hbox{\teni
space} ($\Leftrightarrow\!$ all points are closed)
$X\!_{\!}\in\!{\hbox{\tensy D}}_{\!\!\wp}$ is a  weak {\teni
homology}$_{_{\!\hbox{\fivebf G}}}\!$ $n$-manifold \hbox{\tenrm
(}\hbox{\teni n}-whm$_{_{^{\!\hbox{\fivebf G}}}\!}$\hbox{\tenrm)}
if for some \hbox{{\tenbf A}-module\ {\tensy R}};
\smallskip%
$\!\!\!\!\!$
{
\FFrame{0.0pt}{0.0pt}{\hsize=0.9 \hsize
{%
{{
\noindent
$
\mdoubleH_{_{^{i}}}(X,X\setminus\!_{_{^o}}{_{\!}}x;\hbox{\tenbf
G})
\ \!=\ \!
{0\ \ if\ \hbox{$i\!\ne\! n$} \  \hbox{\sl for\ $\!$all}\
\hbox{$\wp\!\ne\! x\!\in\! X$} , \hskip2.5cm \hbox{\tenbf (4.i)}
}\hfill\break
$
${
\mdoubleH_{_{^{n}}}(X,X\setminus\!_{_{^o}}{_{\!}}x;\hbox{\tenbf
G})\
\!\cong
\hbox{\tenbf G}\oplus \hbox{\tensy R}\ f_{^{\!}}or \ some\ \hbox{$
\wp\!\ne\! x\!\in\! X$}\ \hbox{\sl if}\ X\ne\{\wp\}.\
\hskip0.5cm \hbox{\bf (4.ii$^\prime$)}
}\!$  
}}}
}%
\ {\raise0.35cm\hbox{({\tenbf 4$^\prime$})}}
\indent An {\it n}-whm$_{_{^{\!\hbox{\fivebf G}}}\!}$ $X$ is {\it
joinable}
({\it n}-jwhm$_{_{^{\!\hbox{\fivebf G}}}\!}$)
if {\bf (4.i)} holds also for $x=\wp$.
}

\indent An {\it n}-jwhm$_{_{^{\!\hbox{\fivebf G}}}\!}$ $X$
is a {\it weak} {\it homology} {\it n}-{\it
sphere}$_{_{^{\!\hbox{\fivebf G} }}\!}$ $\!(${\it n}-whsp$)$
if
$\mdoubleH\!\!\!\!_{_{^{n-1}}}\!\!(X\!\setminus\!_{_{\!^o}}{_{\!}}x;\hbox{\tenbf
G})\!=\!0\ \forall\ _{\!}x\!_{\!}\in\!\! X.$
\enddefinition

\definition{Definitions {of technical nature}}
$X$ is acyclic$\!_{{_{\hbox{\fivebf G}}}}\!\!$ if
$\mdoubleH_{i}(X,\emptyset;\hbox{\tenbf G})=0$ for all $i \in
\hbox{\tenbf Z}$.
So, $\{\wp\}\ (=|\{\emptyset_o{_{\!}}\}|)$ isn't
acyclic$\!_{{_{\hbox{\fivebf G}}}}$.
$X$ \hbox{\tenrm is} {\teni weakly} {\teni
direct}$_{_{\!^{\hbox{\fivebf G}}}}\!\!$ if
$\mdoubleH_{_{^{i}}}\!(X;\hbox{\tenbf G})
\!\cong\!
\hbox{\tenbf G}\oplus \hbox{\tensy P}$
{for some \hbox{\teni i}} {and some {\bf A}-module} \hbox{\tensy
P}.
$X$ \hbox{\tenrm is} {\teni locally} {\teni weakly} {\teni
direct}$_{_{\!^{\hbox{\fivebf G}}}}$ if
$\mdoubleH_{_{^{i}}}(X,X\setminus\!_{_{^o}}{_{\!}}x;\hbox{\tenbf
G})
\cong
\hbox{\tenbf G}\oplus \hbox{\tensy Q}$
{for some \hbox{\teni i}}, {some {\bf A}-module}\ $\hbox{\tensy
Q}$\
and some $\wp\!\ne\! x\!\in\! X$.
An $n$-whm$\!_{_{^{\!\hbox{\fivebf G}}}\!\!}$ $X$ is {\teni
ordinary}$\!_{_{\hbox{\fivebf G}}}\!$ if \
${\mdoubleH_{_{^{i}}}\!(X\!\setminus\!_{_{\!^o}}x;\hbox{\tenbf G})
\!=\!0,}\!$
\noindent ${
\ \forall\ \!i\ge n\,}$ {and}\ ${\forall\ \!x\!\in \!X
}
$
\enddefinition

\proclaim{Corollary}{\rm(to Th.\ 6)}.
For locally {weakly direct}$_{_{\!^{\hbox{\fivebf G}}}}$
$(\Rightarrow \!X_{i}\neq \emptyset, \{\wp\})$
\hbox{\tenbf T}$\!_{_{1}}\!\!$-spaces $X_{_{\!1}},X_{_{\!2}}.$

\medskip
\noindent {\bf i.} \ \
$
X_{_{\!1}}\!\!\ast\!X_{_{\!2}}\
(n_{_{{\!1}}}\!+n_{_{{\!2}}}\!+1)\hbox{\tenrm
-whm}_{_{\!{\hbox{\fivebf G}}}}\!\!
$
$ \Longleftrightarrow
X_{_{{\!1}}},X_{_{{\!2}}}\! $ {both}
n$_{_{{\!i}}}\!$-jwhm$_{_{\!\hbox{\fivebf G}}}\!\!
$
$
\Longleftrightarrow X_{_{\!1}}\!\ast\!X_{_{\!2}}\!$
jwhm$_{_{\!{\hbox{\fivebf G}}}}.
$

\smallskip
\noindent {\bf ii.} \ %
$
X\!_{_{1}}\times X\!_{_{2}} $
$(n\!_{_{1}}+n\!_{_{2}})\hbox{\tenrm -whm}_{_{\hbox{\fivebf
\!G}}}\!$
$\Longleftrightarrow
X\!_{_{1}},X\!_{_{2}} $ both \hbox{\tenrm whm}$_{_{\!\hbox{\fivebf
G}}}.
$

\smallskip
\noindent
{\bf iii.} %
If $n_{_{^1}}\!\! +n_{_{^2}}\! > n _{_{^i}}\ i=1,2, $
then;
 $[X\!_{_{1}}\!\times X\!_{_{2}}
(n\!_{_{1}}\!+n\!_{_{2}}\!)\hbox{\tenrm -jwhm}_{_{\hbox{\fivebf
\!G}}}]\!\Longleftrightarrow\ \!\! [X\!_{_{1}}\!,\
\!X\!_{_{2}}\!$\ both $n_{_{^i}}\!${\rm
-jwhm}$\!_{_{^{\hbox{\fivebf G}}}}\ \hbox{\teni and}\ \hbox{\teni
acyclic}\!_{{_{\hbox{\fivebf
G}}}}], %
$
{\tenrm (since by Eq.\ {\tenbf 2};}
{$[\mdoubleH_{_i}(X_{_{1}}\!\!\times\!X_{_{2}};\hbox{\tenbf G})=0
\ \hbox{\tenrm for}\ i\not=
n_{1}+n_{2}
]
\!\Longleftrightarrow\break
\Longleftrightarrow [X_{_{1}}, X_{_{2}}\ \!\hbox{\tenrm both}\
\hbox{\teni acyclic}\!_{_{^{\hbox{\fivebf G}}}}]
\!\Longleftrightarrow\![X_{_{1}}\!\times\! X_{_{2}}\ \hbox{\tenrm
acyclic}\!_{_{^{\hbox{\fivebf G}}}}].$
{\rm So; $X_{_{1}}\!\times\! X_{_{2}}$ is never a {\rm
whsp}$_{_{\hbox{\fivebf \!G}}\!})$.}\
}

\smallskip%
\noindent {\bf iv.}
If moreover $X_{_{1}}, X_{_{2}}\!$\ are {weakly
direct}$_{_{\!^{\hbox{\fivebf G}}}}$\ then$;\ X\!_{_{1}},
X_{_{\!2}}$\ are {both} {\rm whsp}$_{_{\hbox{\fivebf \!G}}\!}$
\underbar{iff} $X\!_{_{1}}\!\!\ast\! X_{_{\!2}}\!$ is.
\endproclaim

\demo{Proof}
\noindent%
Augmental \doubleH omology, like classical, isn't sensitive to
base ring changes. So, ignore {\bf A} and instead use the integers
{\bf Z};
({\bf i-iii.}) Use Th.\ 1, 4-6  and the {\it weak
directness}$_{_{\!^{\hbox{\fivebf G}}}}$ to transpose non-zeros
from one side to the other, using Th.\ 6.{\bf ii} only for joins,
i.e., in particular, with $\epsilon = 0\ \hbox{\tenrm or}\ 1$
depending on wether
$\lower1.0pt\hbox{$^{_{_{\hbox{\sevensy \char"72} }}}\!$}$ =
$\times\ \hbox{\tenrm or}\ \ast$\ \hbox{resp., use;}

\smallskip
\noindent
$\!\!\!\!\!$
\FFrame{0.2pt}{0.0pt}{\hsize=0.95 \hsize
$ %
\hbox{\tenbf(}\mdoubleH
\vbox{\moveleft0.0cm\hbox{\lower0.0pt\hbox{$_{{{_{\!}{p
\lower0.0pt\hbox{\fivebf+} \!\hbox{\fivebf\ }_{\!}\epsilon}}}}$}}}
(X\!_{_{^{1}}}\!\lower1.0pt\hbox{$^{_{_{\hbox{\sevensy \char"72}
}}}$}\!X\!_{_{^{2}}},
X\!_{_{^{1}}}\!\lower1.0pt\hbox{$^{_{_{\hbox{\sevensy \char"72}
}}}$}\!X\!_{_{^{2}}} \!\setminus\!{_{\!}}_{_{^o}}\!\!\!
\lower4pt\hbox{\FFrame{0.2pt}{0.0pt}{
\hsize=0.085\hsize \noindent
\rlap{$\widetilde{\phantom{rrrrrr}}$}{$
{(x\!_{_{^{1}}},x\!_{_{^{2}}}\!)}$}}} \!\!;\hbox{\tenbf G})
{\lower4.0pt\hbox{$^{_{\hbox{{\fivebf Z}}
\atop{{\raise0.0pt\hbox{{\eightbf =}}}
\!\!\!\!\!_{\!}{\raise2.8pt\hbox{{\eightsy {\char"27}}}} }}}$} }
\hbox{\tenbf)}\!\! \ \ \ \ \!\mdoubleH _{_{^{p\!\!}}}
((X\!_{_{^{1}}},X\!_{_{^{1}}}\!\setminus\!{_{\!}}_{_{^o}}x\!_{_{^{1}}})\times
(X\!_{_{^{2}}},X\!_{_{^{2}}}\!\setminus\!{_{\!}}_{_{^o}}x\!_{_{^{2}}});\hbox{\tenbf
G})
{\lower4.0pt\hbox{$^{_{\hbox{{\fivebf Z}}
\atop{{\raise0.0pt\hbox{{\eightbf =}}}
\!\!\!\!\!_{\!}{\raise2.8pt\hbox{{\eightsy {\char"27}}}} }}}$} }
$\hfill\break
$\indent\indent
{\lower4.0pt\hbox{$^{_{\hbox{{\fivebf Z}}
\atop{{\raise0.0pt\hbox{{\eightbf =}}}
\!\!\!\!\!_{\!}{\raise2.8pt\hbox{{\eightsy {\char"27}}}} }}}$} }
{\hbox{\sixrm Lemma\ p.\ \!11}\brack{\hbox{\sixrm +\ \!Eq.\ \!1\
p.\ \!8}}}\ \!
{\lower4.0pt\hbox{$^{_{\hbox{{\fivebf Z}}
\atop{{\raise0.0pt\hbox{{\eightbf =}}}
\!\!\!\!\!_{\!}{\raise2.8pt\hbox{{\eightsy {\char"27}}}} }}}$} }
$
{$
{\rlap{$\!_{_{_{i+j=p\!\!}}}$}{\raise2pt\hbox{\ $\bigoplus$}}}\ \
[\mdoubleH_{i}
({X\!_{_{^{1}}},X\!_{_{^{1}}}\!\setminus\!{_{\!}}_{_{^o}}x\!_{_{^{1}}}};\hbox{\tenbf
Z})
\!\otimes\!_{_{\hbox{\fivebf Z}}}\!
        \mdoubleH_{j}
({X\!_{_{^{2}}},X\!_{_{^{2}}}\!\setminus\!{_{\!}}_{_{^o}}x\!_{_{^{2}}}\!};\hbox{\tenbf
G})] \ \ \oplus$}
\hfill\break
\indent\hskip3.0cm %
{$ \oplus\ \ {\rlap{$\!\!\!\!_{_{_{i+j=p-\!1\!\!}}} $}
{\raise2pt\hbox{$\bigoplus$}}}\ \
\hbox{\tenrm Tor}_{^{_1}}^{_{^{\hbox{\fivebf Z}}}}\!\!\bigl(
\mdoubleH_{i}
({X\!_{_{^{1}}},X\!_{_{^{1}}}\!\setminus\!{_{\!}}_{_{^o}}x\!_{_{^{1}}}}\!;\hbox{\tenbf
Z}),
 \mdoubleH_{j}
({X\!_{_{^{2}}},X\!_{_{^{2}}}\!\setminus\!{_{\!}}_{_{^o}}x\!_{_{^{2}}}}\!;
\hbox{\tenbf G})\big).
$}
}
\hfill$\triangleright$

\noindent {\bf iv.} Use, by the Five Lemma, the chain equivalence
of the second component in the first and the last item of Th.\
6.{\bf i} and the {\bf M-$_{\!}$Vs} w.r.t.
${\!}\{(X\!_{_{\ }}\!\!\ast (Y\!\!_{_{\
}}\!\setminus\!{_{_{^{o\!}}}} y\!_{_{0}}\!)),((X\!\!_{_{\
}}\!\setminus\!{_{_{^{o\!}}}} x\!_{_{0}}\!)\ast Y\!_{_{\
}}\!\!)\}$.
\qed%
\enddemo

\definition{Definition}
$\emptyset$ is a {\teni homology}$_{_{\!\hbox{\fivebf G}}}\!$
$-\infty$-{\teni manifold} and
$X\!=\!{\bullet}{\bullet}$ is a \hbox{\teni
homology}$_{_{\!\hbox{\fivebf G}}}\!$ $0$-\hbox{\teni
manifold}.$\!$
Else, a connected, locally compact Hausdorff space
$X\!_{\!}\in\!{\hbox{\tensy D}}_{\!\!\wp}$ is a ({\teni singular})
homology$\!_{_{\!\hbox{\fivebf G}}}\!$
$n$-\hbox{\teni manifold} \hbox{\tenrm
(}$n$-hm$_{_{^{\!\hbox{\fivebf G} }}\!}$\hbox{\tenbf)} if;

\medskip
$\!\!\!\!\!\!\!$
{
\FFrame{0.0pt}{0.0pt}{\hsize=0.92 \hsize
{%
{{
\noindent
$\mdoubleH_{_{^{i}}}(X,X\setminus{_{\!}}_{_{^o}}x;\hbox{\tenbf
G})%
\ \!=\ \!$
0\ \
if\ $i\!\ne\! n$ for all
$\wp\!\ne\! x\!\in\! X,
\hskip3.36cm\hbox{({\tenbf 4.i})} %
\hfill\break
$
${
\mdoubleH_{_{^{n}}}(X,X\setminus\!{_{\!}}_{_{^o}}x;\hbox{\tenbf
G})
\!\cong\!
0\ {or}\ \hbox{\tenbf G}\ \hbox{\tenrm for\ all}\ \wp\!\ne\!
x\!\in\! X \hbox{\tenrm and $=\hbox{\tenbf G}$ for some}\
x\!\in\!X.
\hskip0.3cm\hbox{({\tenbf 4.ii})}
}$
}}}}}
\ {\raise0.35cm\hbox{({\tenbf 4})}}

\smallskip\smallskip
{The} {\teni boundary}:\ \ $\!\hbox{\tenrm Bd}_{_{\!\hbox{\fivebf
G}}}\!X\!\!:=\! \{x\!\in\!X\big\vert
\mdoubleH_{n}(X,X\setminus\!{_{\!}}_{_{^o}}x;\hbox{\tenbf G})=0\}.
$
{If} {\rm Bd}$_{_{\!\hbox{\fivebf G}}}X\!\!\ne\!\emptyset$\
(\hbox{\tenrm Bd}$_{_{\!\hbox{\fivebf G}}}X\!\!=\!\emptyset),$\
$X$ is called a \hbox{\teni homology}$_{_{\!\hbox{\fivebf G}}}\!$
$n$-{\teni manifold} {\teni with} $($without$)$ {\teni boundary}.

A compact $n$- manifold {\tensy S} is {\it
orientable}$_{_{{^{^{\!\!\hbox{\fivebf G}}}}}}\!$ {if}
$\mdoubleH\!_{_{{n}}}\!(\hbox{\tensy S},\hbox{\tenrm
Bd}\hbox{\tensy S};\hbox{\tenbf G})\ \widetilde{\hbox{\rm=}}\
\hbox{\tenbf G}.$ An $n$- manifold is \hbox{\teni
orientable}$_{_{{^{^{\!\!\hbox{\fivebf G}}}}}}\!\!$ if all its
compact $n$- submanifolds are orientable \ \hbox{\twelvrm-} else
{\teni non}-{\teni orientable}$_{_{{^{^{\!\hbox{\fivebf
G}}}}}}\!._{\!}$

Orientability is left undefined for $\emptyset$.
An $n$-hm$\!_{_{^{\!\hbox{\fivebf G} }}\!\!}$ $X$ is {\teni
joinable} if $(4)$ holds also for $x\!=\wp.$
An {\it n}-hm$_{_{^{\!\hbox{\fivebf G} }}\!}$
$X{\!}\!\not={\!}\emptyset$
is a {\it homology$_{_{\!\hbox{\fivebf G}}}{\!}{\!}$ n-sphere}
$(n$-hsp$_{_{^{\!\hbox{\fivebf G} }}\!})$ if
for all $x{\!}\in{\!} X,$
$\mdoubleH_{i}(X,X{\!}\setminus\!_{_{^o}}{{\!}}x;
\hbox{\tenbf G}){\!}={\!}\hbox{\tenbf G}$ if $i=n$ and $0$ else.
\enddefinition

\remark{Note 1}
Triangulable manifolds
$ \ne\emptyset$ are ordinary, by Note\ 1 p.\ 25, and locally {\rm
weakly direct}$_{_{\!^{\hbox{\fivebf G}}}}\!$ since\
$\mdoubleH \!\! \raise-0.55pt\hbox{${\!\!\!{_{_{_{\dim\!{\Sigma_{}
}}}}}\!}$}
(_{\!}|\Sigma|,|\Sigma|\setminus\nobreak\!_{_{\!^o}}\!\alpha;\hbox{\tenbf
G})\cong$
$\hbox{\tenbf G}$ for any $\alpha\!\in\! \hbox{\tenrm Int}\sigma$
if $\sigma$ is a maxidimensional simplex i.e. if
$\#\sigma-1=:\dim\sigma=\dim\Sigma,$ since now {\rm
Lk}$_{_{\!\Sigma}}\sigma=\{\emptyset\!_{_{^{o}}}\!\}$.
{Prop. 1. p.\ 11 and Lemma p.\ 5 now gives the claim.}
If $\alpha\!\in\!\hbox{\tenrm Int}\sigma$,
Proposition 1 p.\ 11 also implies that;
\noindent $\
\mdoubleH\!\!\!\!\!\!\!\raise0.6pt\hbox{$_{_{_{\dim\!{\Sigma}}}}$}\!
\!({\Sigma},\hbox{\tenrm cost}_{_{\!{\Sigma}}}\!{\sigma
};\hbox{\tenbf Z})
\cong\
\mdoubleH\!\!\!\!\!\!\!\raise0.6pt\hbox{$_{_{_{\dim\!{\Sigma}}}}$}\!
\!(\vert\Sigma\vert,\vert\Sigma\vert
\setminus_{\!o}\alpha;\hbox{\tenbf Z})
$
$ \cong\
\mdoubleH\!\!\!\!\!\!\!
\raise0.5pt\hbox{${\!\!\!{_{_{_{\dim\!{\Sigma
\lower1pt\hbox{-}\hbox{\fivebf\#}\sigma\!}}}}}\!}$}\!\!\!\!\!\!\!
(\ \!\hbox{\tenrm Lk}_{_{\!\Sigma_{_{\!{\
}}}}}\!\!\sigma;\hbox{\tenbf Z})$ is a direct sum of {\bf
Z}-terms.
\endremark

\smallskip
$\!\!\!$When $\lower1.0pt\hbox{$^{_{_{\hbox{\sevensy \char"72}
}}}$}\!$ in Th.\ 7, all through, is interpreted as $\times$, the
word ``manifold(s)" (on the r.h.s.) temporarily excludes
$\emptyset,\{\emptyset_o\}$ and $\bullet\bullet$, and we assume
$\epsilon:=0$.
When $\lower1.0pt\hbox{$^{_{_{\hbox{\sevensy \char"72} }}}$}\!$,
all through, is interpreted as $\ast$, put $\epsilon:=1$, and let
the word ``manifold(s)" on the right hand side be limited to ``any
\underbar{compact joinable} homology$_{_{\!\!^{\hbox{\fivebf
G\!}}}}$ $n_{_{\!i}}$-manifold".

\proclaim{Theorem 7}
For locally {\teni weakly} {\teni direct}$_{_{\!^{\!\hbox{\fivebf
G}}}}$ \hbox{\tenbf T}$\!_{_{1}}\!\!$-spaces $X\!_{_{1}},
X_{_{\!2}}$ and any \hbox{\tenbf A}-module {\bf G}$:$

\smallskip
\noindent
{\bf 7.1.}
\ $X_{_{\!1}}\!\lower1.0pt\hbox{$^{_{_{\hbox{\sevensy \char"72}
}}}$}\! X_{_{\!2}}\!$ is a \hbox{\teni
homology}$_{_{\!\hbox{\fivebf G}}}\!$
$(n_{_{\!1}}\!+n_{_{\!2}}\!+\!\epsilon)$%
$\hbox{-manifold}
\Longleftrightarrow X_{i}$
is a $n_{_{\!i\!}}$-hm$_{_{\!\hbox{\fivebf G}}},\ i=1,2.$

\smallskip
\noindent {\bf 7.2.} \
$\hbox{\tenrm Bd}_{_{\!\hbox{\fivebf G}}}{\!}_{\!} (\bullet\times
X)= \bullet\times (\hbox{\tenrm Bd}_{_{\!\hbox{\fivebf
G}}}{\!}_{\!} X)$. {\sl Else;} {\rm Bd}$_{_{\!\hbox{\fivebf
G}}}{\!}_{\!} (X_{_{\!1}}\!\lower1.0pt\hbox{$^{_{_{\hbox{\sevensy
\char"72} }}}$}\! X_{_{\!2}}\!)= ((\hbox{\tenrm
Bd}_{_{\!\hbox{\fivebf G}}}{\!}_{\!}
X_{_{\!1}}\!)\lower1.0pt\hbox{$^{_{_{\hbox{\sevensy \char"72}
}}}$}\! X_{_{\!2}}\!)\cup
(X_{_{\!1}}\!\lower1.0pt\hbox{$^{_{_{\hbox{\sevensy \char"72}
}}}$}\! (\hbox{\tenrm Bd}_{_{\!\hbox{\fivebf G}}}{\!}_{\!}
X_{_{\!2}}\!)).$

\smallskip
\noindent {\bf 7.3.} \
$X_{_{\!1}}\!\lower1.0pt\hbox{$^{_{_{\hbox{\sevensy \char"72}
}}}$} \!X_{_{\!2}}\!\ {is\ orientable}{\rlap{$_{_{\!\hbox{\fivebf
G}}}$} {\ \ \Longleftrightarrow}}\ X_{_{\!1}},X_{_{\!2}}$ {are
both orientable}$_{_{\!\hbox{\fivebf G}}}.$
\endproclaim

\demo{Proof}
Th.\ 7 is trivially true for
$\!X_{_{\!i}}\!_{\!}\times_{\!}\bullet$ and
$X_{_{\!i}}\!_{\!}\ast_{\!}\{\wp\}.$
Else, exactly as for the above Corollary, adding for {\rm 7.1.}\
that for Hausdorff-like spaces $($:= all compact subsets are
locally compact$)$, in particular for Hausdorff spaces,
$X_{_{\!1}}\!\ast X_{_{\!2}}\!$ is locally compact (Hausdorff)
\underbar{iff}
$X_{_{\!1}}\!,X_{_{\!2}}\!$ both are compact (Hausdorff), cf. \cite{4} %
p.\ 224.
\qed
\enddemo

\remark{Note 2}
\noindent$\!
(X_{_{^{_{\!}1}}}\!\lower1.0pt\hbox{$^{_{_{\hbox{\sevensy
\char"72} }}}$}\!X_{_{^{_{\!}2\!}}}, \hbox{\tenrm
Bd}\!\!{{{\lower3.5pt\hbox{\fivebf
G}}}}{\!}_{\!}(X_{_{^{_{\!}1}}}\!\lower1.0pt\hbox{$^{_{_{\hbox{\sevensy
\char"72} }}}$}\!X_{_{^{_{\!}2\!}}}))
=\![{7.2}]\!=
(X_{_{^{_{\!}1}}}\!\lower1.0pt\hbox{$^{_{_{\hbox{\sevensy
\char"72} }}}$}\!X_{_{^{_{\!}2\!}}},
X_{_{^{_{\!}1}}}\!\lower1.0pt\hbox{$^{_{_{\hbox{\sevensy \char"72}
}}}$}\!\hbox{\tenrm Bd}\!\!{{{\lower3.5pt\hbox{\fivebf
G}}}}{\!}_{\!}X_{_{^{_{\!}2\!}}}\cup \hbox{\tenrm
Bd}\!\!{{{\lower3.5pt\hbox{\fivebf
G}}}}{\!}_{\!}X_{_{^{_{\!}1}}}\!\lower1.0pt\hbox{$^{_{_{\hbox{\sevensy
\char"72} }}}$}\! X_{_{^{_{\!}2\!}}})
=
[{\hbox{\tenrm Def.} \ {\!}\hbox{\tenrm p.\ \!7}}] _{\!}=_{\!}$

$=
(X_{_{^{_{\!}1\!}}},\hbox{\tenrm
Bd}\!\!{{{\lower3.5pt\hbox{\fivebf
G}}}}{\!}_{\!}X_{_{^{_{\!}1\!}}})
\lower1.0pt\hbox{$^{_{_{\hbox{\sevensy \char"72} }}}$}\!
(X_{_{^{_{\!}2\!}}},\hbox{\tenrm
Bd}\!\!{{{\lower3.5pt\hbox{\fivebf
G}}}}{\!}_{\!}X_{_{^{_{\!}2\!}}}).
$

Equivalently;
$
\hbox{\tenrm Hip}\!\!{{{\lower3.5pt\hbox{\fivebf
G}}}}{\!}_{\!}(X_{_{^{_{\!}1}}}\!\lower1.0pt\hbox{$^{_{_{\hbox{\sevensy
\char"72} }}}$}\!X_{_{^{_{\!}2\!}}})
=
X_{_{^{_{\!}1}}}\!\lower1.0pt\hbox{$^{_{_{\hbox{\sevensy \char"72}
}}}$}\!\hbox{\tenrm Hip}\!\!{{{\lower3.5pt\hbox{\fivebf
G}}}}{\!}_{\!}X_{_{^{_{\!}2\!}}}\cup \hbox{\tenrm
Hip}\!\!{{{\lower3.5pt\hbox{\fivebf
G}}}}{\!}_{\!}X_{_{^{_{\!}1}}}\!\lower1.0pt\hbox{$^{_{_{\hbox{\sevensy
\char"72} }}}$}\! X_{_{^{_{\!}2\!}}}$,
where for any space $X$:
\ \ $\!\hbox{\tenrm Hip}\!\!{{{\lower3.5pt\hbox{\fivebf G}}}}\!X
=$ the set of {\teni Homologically}$_{_{\!\hbox{\fivebf G}}}\!$
{\teni instabile} {\teni points} $:= \{x\!\in\!X\big\vert
\mdoubleH_{i}(X,X\setminus\!{_{\!}}_{_{^o}}x;\hbox{\tenbf G})=0\
\forall\ i\ \!\in \hbox{\tenbf Z}\}.$
\endremark

\example{Example}
$\bullet$, $\bullet\bullet$ and $\bullet_{\!}\bullet_{\!}\bullet$
are all 0-whm$_{_{\!\hbox{\fivebf G}}}$ but
$\bullet_{\!}\bullet_{\!}\bullet$ isn't a 0-hm$_{_{\!\hbox{\fivebf
G}}}$.
As follows from above,
$\raise1pt\hbox{\eightsy{\char"0D}}
\simeq
(\bullet\bullet){\ast}\ \!(\bullet\bullet) $
i.e. the join of two 0-hsp$_{_{\!\hbox{\fivebf G}}}$ is a
1-hsp$_{_{\!\hbox{\fivebf G}}}$, cf. Ex.\ 4 p.\ 28.
$\raise1pt\hbox{\tensy{\char"09}}
\simeq
(\bullet\bullet){\ast}\ \!(\bullet{_{\!}}\bullet{_{\!}}\bullet) $
and
$
{{ {{ \bullet{_{\!}}{_{\!}}{_{\!}}
\raise1.4pt\hbox{$^{{\nearrow^{\!_{\!}\!_{\!}\lower1.45pt\hbox{$\bullet$}}}}$}
\!\!\!\!\!{{\myline}}\raise1.7pt\hbox{{$_{{\!\!_{\!}\hbox{$\bullet$}}}$}}
\!\!\!\!\!{_{\!}}\!{_{\!}}{_{\!}}
\lower2.1pt\hbox{${{_\searrow}}$}\!_{_{_{_{\!_{\!}\!\hbox{$\bullet$}}}}}
}} }}
\simeq
$
$
\bullet\ \!{\ast}\ \!(\bullet{_{\!}}\bullet{_{\!}}\bullet)
$
are both 1-jwhm$_{_{\!\hbox{\fivebf G}}}$ but nighter is a
1-hm$_{_{\!\hbox{\fivebf G}}}$,
which actually contradicts the statement in \cite{27} %
p.\ 122 Corollary\ 2.12(ii).

The following is a trivial example of how to ``{\it eliminate}"
neighborhood retract-subspaces.
%
{\doubleH$_{_{1}}\!($\raise1pt\hbox{\eightsy{\char"09}}$;\hbox{\tenbf
G})
= $
\doubleH$_{_{1}}\!( {{| {{ \bullet{_{\!}}{_{\!}}{_{\!}}
\raise1.4pt\hbox{$^{{\nearrow^{\!_{\!}\!_{\!}\lower1.45pt\hbox{$\bullet$}}}}$}
\!\!\!\!\!{{\myline}}\raise1.7pt\hbox{{$_{{\!\!_{\!}\hbox{$\bullet$}}}$}}
\!\!\!\!\!{_{\!}}\!{_{\!}}{_{\!}}
\lower2.1pt\hbox{${{_\searrow}}$}\!_{_{_{_{\!_{\!}\!\hbox{$\bullet$}}}}}
}} |}} \lower2pt\hbox{/} _{^{^{|\raise2pt\hbox{$_{\ \!
\!{^{^{\hbox{\tenbf .}}}\ \!\!\!\!\!{\raise0.7pt\hbox{\tenbf .}}
\!\!\!_{_{\hbox{\tenbf .}}}}}$}|}}} ;\hbox{\tenbf G})$
$=$
$ [\hbox{{\tenrm \cite{21} %
Ex.\ 2\ p.\ 230}}]$
$=$
\doubleH$_{_{1}}\!( {{ \bullet{_{\!}}{_{\!}}{_{\!}}
\raise1.4pt\hbox{$^{{\nearrow^{\!_{\!}\!_{\!}\lower1.45pt\hbox{$\bullet$}}}}$}
\!\!\!\!\!{{\myline}}\raise1.7pt\hbox{{$_{{\!\!_{\!}\hbox{$\bullet$}}}$}}
\!\!\!\!\!{_{\!}}\!{_{\!}}{_{\!}}
\lower2.1pt\hbox{${{_\searrow}}$}\!_{_{_{_{\!_{\!}\!\hbox{$\bullet$}}}}}
}}\ \!,\bullet{_{\!}}\bullet{_{\!}}\bullet ;\hbox{\tenbf G})=$

\item{=$\!$}\doubleH$_{_{1}}\!(\bullet\ \!{\ast}\
\!(\bullet{_{\!}}\bullet{_{\!}}\bullet),
\bullet{_{\!}}\bullet{_{\!}}\bullet ;\hbox{\tenbf G}) = $}
[\hbox{{{\tenrm Def.}}\ {\tenrm p.\ \!7}}] =
\doubleH$_{_{1}}\!((\bullet,\{\emptyset\})\ \!{\ast}_{_{\!^\cup}}\
\! (\bullet\bullet\bullet,\emptyset) ;\hbox{\tenbf G})$
=\ [Eq.\ 3\ {\tenrm p.$\ \!$10}]=
\noindent$={\rlap{$_{_{_{i+j=0}}}$} {\
\raise2pt\hbox{$\bigoplus$}}}\
\!$\doubleH$_{_{i}}(\bullet,\{\emptyset\};\hbox{\tenbf
R})\otimes\!\!_{_{\hbox{\fivebf R}}}
\mdoubleH_{_{\!j}}\!(\bullet\bullet\bullet,\emptyset;\hbox{\tenbf
G})\!=\ \!
$\doubleH$_{_{0}}(\bullet,\{\emptyset\};\hbox{\tenbf
R})\otimes\!\!_{_{\hbox{\fivebf R}}}
\mdoubleH_{_{0}}\!(\bullet\bullet\bullet,\emptyset;\hbox{\tenbf
G}) \!= $
\hbox{[Lemma\nobreak\ p.\ $\!$6]=}
$=\hbox{\tenbf R}\otimes\!\!_{_{\hbox{\fivebf R}}}(\hbox{\tenbf
G}\oplus\hbox{\tenbf G})$
$ =\hbox{\tenbf G}\oplus\hbox{\tenbf G}. $
\indent\indent
$(\ _{\!^{\!}}\mdoubleH_{_{\!i\!}}(\raise1pt\hbox{\eightsy
{\char"09}};\hbox{\tenbf G})\!=\!0$ if $i\!\ne\! 1$.)
\endexample

\head I: General Topological Properties for Realizations of
Simplicial Complexes
\endhead

\subhead I:1 Realizations and Local Homology Groups Related to
Simplicial Products and Joins\endsubhead Only fairly recent the
first examples of non-triangulable topological (i.e.
$0$-differentiable) manifolds have been successfully constructed,
cp.
\cite{23} %
\S5,
and since all differentiable manifolds are triangulable,
cf. \cite{20} %
p.\ 103 Th.\ 10.6, essentially all spaces within mathematical
physics are triangulable.
The realization$_{_{\!^o}}\!$, p.\ 5, of any simplicial
complex$_{_{\!^o}}\!$ is a CW-complex$_{_{\!^o}}\!$.
Our CW-complexes$_{_{\!^o\!}}\!$ will have, as the so called
``relative CW-complexes" defined in \cite{9} %
p.$\ _{\!}$326, a $(-1)$-cell $\{\wp\}$ ({\it ``an ideal cell"}),
but their topology is that of the spaces in the category
$\hbox{\tensy D}_{_{\!\!^\wp}}\!$ defined\nobreak\ in\nobreak\ p.\
3.
CW-complexes are compactly generated, perfectly normal spaces,
\cite{9} %
pp.\ $22,112,242$,\break
that are locally\ contractible in a strong sense and
(hereditarily) paracompact, cf. \cite{9} %
pp.\ 28\hbox{-}29\ Th.\ \!1.3.2\ \!\raise1pt\hbox{\eightrm+}\ Th.\
1.3.5\ (Ex.\ 1\ p.\ 33).
A topological space has the homotopy type of a CW-complex
\underbar{\tenrm iff}
it has the type of the realization of a simplicial complex, which
it has \underbar{\tenrm iff} it has the type of an ANR,
{\tenrm cf.\ \cite{9} %
p.\ 226\ Th.\ 5.2.1.}

The {\bf k}-iffikation {\bf k}({\teni X}) of {\teni X} is {\teni
X} with its topology enlarged to the weak topology w.r.t.
its compact subspaces, cf. \cite{4}. %
Put
${X}\bar{\times}\ \!Y\!:=\hbox{\tenbf k}({X}{{{\times}}}\ \!Y)$.
If $X$ and $Y$\nobreak\ are CW-complexes, this is a proper
topology-enlargement only if none of the two underlying complexes
are locally\ finite and at least one is uncountable.
Let ${X}\bar{\ast}\ \!Y\!$ be the quotient space w.r.t.
$\!${\hbox{\teni {\char"70}}} $\!$: $\!$({\teni X}%
$\!\ _{\!}\rlap{\lower1.0pt\hbox{$\bar{\phantom{-}
_{\!}}$}}{_{^{\!}}
\hbox{\tensy {\char"02}}}$%
{\teni Y}%
){\tensy {\char"02}}{\tenbf I}%
\hbox{\tensy {\char"21}}%
{\teni X}{\hbox{\tensy {\char"0E}}}%
{\teni Y} from\nobreak\ p.\ $\!$7.
See \cite{16} %
p.\ 214 for relevant distinctions.
Now, simplicial ${\times},{\ast}$ ``commute" with realization by
turning into $\bar{\times},\bar{\ast}$ respectively.

Unlike $\ast$ defined in p.\ 7, $\bar{\ast}$ is actually
associative for arbitrary topological spaces.

\definition{Definition}
\hbox{(cf. \cite{7} %
Def.\ 8.8 p. 67.)}
Given ordered simplicial complexes $\Delta{\!^{^{_{\prime}}}}\!$
and $\Delta{\!^{^{_{\prime\prime}}}}\!$ i.e. the vertex sets
$V\!\!\!_{{^{\Delta^{^{_{\!{{\!}}\prime}}}}}}\!\!$ and
$V\!\!\!_{{^{\Delta^{^{_{\!\!\prime{_{\!}}\prime}}}}}}\!\!$ are
partially ordered so that each simplex becomes linearely ordered
resp.
{\it The Ordered Simplicial Cartesian Product}
$\Delta{\!^{^{_{\prime}}}}\!
{\raise1.0pt\hbox{\ninesy{\char"02}}}
\Delta{\!^{^{_{\prime\prime}}}}\!\!$ of
${\Delta{\!^{^{_{\prime}}}}}\!$ and
$\Delta{\!^{^{_{\prime\prime}}}}\!\!$ $($\hbox{\tenrm
triangulates} $|_{\!}\Delta{\!^{^{_{\prime}}}}{\!}|
\bar{\raise0.5pt\hbox{\ninesy{\char"02}}}|_{\!}\Delta{\!^{^{_{\prime\prime}}}}\!|$
\hbox{\tenrm and}$)$ is defined through
$\!V\!\!\!{_{\!}}_{{^{{\Delta^{^{_{\!\prime}}}}\!{\!}\times{\!}
\Delta^{^{_{\!\!\prime{_{\!}}\prime}}} }}}\!\!:= $
$\{_{\!}(v\!\raise1pt\hbox{$_{_{^{i}}}$}\!^{_{\!}\prime}_{\!},_{_{\!}}
v\!\raise1pt\hbox{$_{_{^{j}}}$}\!^{_{\!}\prime\prime})_{\!}\}
\!\!=\!\!V\!\!\!_{{^{\Delta^{^{_{\!\prime}}}}}}\!\!{\raise1.0pt\hbox{\ninesy{\char"02}}}
V\!\!\!_{{^{\Delta^{^{_{\!\!\prime{_{\!}}\prime}}}}}}\!\! .$
Put 
$w\!_{_{^{i,j}}}\!_{\!}{_{\!}}:=\!(v\!\raise1pt\hbox{$_{_{^{i}}}$}\!^{_{\!}\prime}_{\!},_{_{\!}}
v\!\raise1pt\hbox{$_{_{^{j}}}$}\!^{_{\!}\prime\prime}).$ Now,
simplices in
$\Delta{\!^{^{_{\prime}}}}\!{\raise1.0pt\hbox{\ninesy{\char"02}}}\Delta{\!^{^{_{\prime\prime}}}}\!_{\!}$
are sets $\{ w_{i_0,j_0},$ $w_{i_1,j_1},\!.., $
$w_{i_k,j_k}\}_{_{\!}},\!$
$\hbox{\tensl with}\ w\!_{_{^{{i_{\!s\!}},{j_{\!s}}}}}\!\!\!
\neq\!w_{{i_{s+1}},{j_{s+1}}}\ \!\hbox{\tensl and}\
v_{i_0}^\prime\!\!\le\!\!v_{i_1}^\prime\!\!\le\!\!..\!\! \le\!\!
v_{i_k}^\prime$
$(v_{j_0}^{\prime\prime}\!\!\le\!\!v_{j_1}^{\prime\prime}\!\!
\le\!\!..\!\!\le\!\!v_{j_k}^{\prime\prime})$ where\
$v_{i_0}^\prime, v_{i_1}^\prime,.., v_{i_k}^\prime
(v_{j_0}^{\prime\prime},
v_{j_1}^{\prime\prime},..,v_{j_k}^{\prime\prime})$ is a sequence
of vertices, with repetitions possible,
constituting a simplex in $\Delta^\prime\
(\Delta^{\prime\prime})$.
\enddefinition

\proclaim{Lemma}
\hbox{\tenrm(cp. \cite{7} %
p.\ 68.)} %
$\{\eta:=(|p_{_{^{\!1}}}\!|,|p_{_{^{\!2}}}\!|),\Sigma_1\times
\Sigma_2\}$, $
p_{{{i}}}\!:\Sigma_{_{^{\!1}}}\!\!\times\!\Sigma_{_{^{\!2}}}
\!\!\rightarrow\! \Sigma_{{{i}}} $ the simplicial projection,
triangulates $|\Sigma_1|\bar{\times}|\Sigma_2|$. If $L_1$ and
$L_2$ are subcomplexes of $\Sigma_1$ and $\Sigma_2$, then $\eta$
carries $|L_1\times L_2|$ onto $|L_1|\bar{\times}|L_2|$.
Furthermore, this triangulation has the property that, for each
vertex $B$ of $\Sigma_2$, say, the correspondence
$x\rightarrow(x,B)\ \hbox{\tensl is\ a}$ simplicial map of
$\Sigma_1$ into $\Sigma_1\times \Sigma_2$. Similarly
for joins,\ \hbox{\tenrm cp.\ \cite{28} %
p.\ 99}.
\endproclaim

\demo{Proof}
{\bf($\times$)}. The simplicial projections $
p_{{{i}}}\!:\Sigma_{_{^{\!1}}}\!\times\Sigma_{_{^{\!2}}}
\!\!\rightarrow\! \Sigma_{{{i}}} $ gives realized continuous maps
$|p_{i}|, i\!=\!1,2.$ $ \eta:=
\!(|p_{_{^{\!1}}}\!|,|p_{_{^{\!2}}}\!|)\!:\!
|\Sigma_{_{^{\!1}}}\!\!\times\!\Sigma_{_{^{\!2}}}\!|
\!\!\rightarrow\!
|\Sigma_{_{^{\!1}}}\!|\bar\times|\Sigma_{_{^{\!2}}}\!| $ is
bijective and continuous cf. \cite{2} %
2.5.6 p.\ $\!$32\ $\!$\raise1pt\hbox{\eightrm+}\ $\!$Ex.\ $\!$12,
14 p. 106-7.
$\tau\!\!{_{_{^{|\!\Sigma_{_{^{\!1\!\!}}}\times\Sigma_{_{^{\!2\!\!}}}|}}}}\!\!$
$(\tau\!\!{_{_{^{|\!\Sigma_{_{^{\!1\!\!}}}|
\!\bar\times\!|\!\Sigma_{_{^{\!2\!\!}}}|}}}}\!)\!$ is the weak
topology w.r.t. $\!$the compact subspaces
$\{|\Gamma\!_{_{^{\!1}}}\!\times\Gamma\!_{_{^{\!2}}}\!|\}
_{_{^{\!\Gamma\!_{{^{i}}}\!\subset\Sigma_{{^{_{\!}i}}}}}}\!\! $ ($
\{\widetilde=\
\!|\Gamma\!_{_{^{\!1}}}\!|\times|\Gamma\!_{_{^{\!2}}}\!|\}
_{_{^{\!\Gamma\!_{{^{i}}}\!\subset\Sigma_{{^{_{\!}i}}}}}}\!),\!\hbox{\tenrm\
cf.}$
\noindent
\cite{9} %
p.\ 246 Prop.\ A.2.1.
$\! (|p_{_{^{\!1}}}\!|,|p_{_{^{\!2}}}\!|)(
|\Gamma\!_{_{^{\!1}}}\!\times\Gamma\!_{_{^{\!2}}}\!| \cap A) \!=\!
\nobreak
(|\Gamma\!_{_{^{\!1}}}\!|^{^{_{\!}}}\times^{^{_{\!}}}|\Gamma\!_{_{^{\!2}}}\!|)
\cap(|p_{_{^{\!1}}}\!|,|p_{_{^{\!2}}}\!|)(A) \hbox{\tenrm\ i.e.}\
\! (|p_{_{^{\!1}}}\!|,
|p_{_{^{\!2}}}\!|)^{^{_{\!{\raise0.5pt\hbox{\fivebf-} \!1}}}}\!\!$
{\tenrm\ is\ continuous}.$\ \triangleright$\break
\noindent
({$\ast$}).
As for $\times_{_{_{\!}}}$, cp. \cite{29} %
(3.3) p.\ 59, with
$\Sigma_{_{^{\!1}}}\!\ast \Sigma_{_{^{\!2}}}\!:=
\{\sigma\!{_{_{^1 \!\!}}}\cup \sigma\!{_{_{^2\!\!}}}\ \vert\
\sigma\!{_{_{i \!\!}}}\in \Sigma_{_{^{\!i}}}\!\ (i=1,2)\}$,
using;
$$ \eta
:|\Sigma_{_{^{\!1}}}\!\ast \Sigma_{_{^{\!2}}}\!|
= |\{\sigma\!{_{_{^1 \!\!}}}\cup \sigma\!{_{_{^2\!\!}}}\ \vert\
\sigma\!{_{_{i \!\!}}}\in \Sigma_{_{^{\!i}}}\!\
(i=1,2)\}|{\rlap{$\hookrightarrow$}{\ \rightarrow}}
$$
$$
{\rlap{$\hookrightarrow$}{\ \rightarrow}}\ |\Sigma_{_{^{\!1}}}\!|\
\!\bar{\ast}\ \!|\Sigma_{_{^{\!2}}}\!|\ ;\
t_{_{^{\!1}}}^{^{_{_{\ }}}}\!\hbox{\tenrm
v}_{_{^{\!1}}}^{^{_{_{\prime}}}}\! \
\!\raise1pt\hbox{\eightrm+}\!\dots\!\raise1pt\hbox{\eightrm+}\ \!
t_{_{^{\!\hbox{\fiverm q}}}}^{^{_{_{\ }}}}\! \hbox{\tenrm
v}_{_{^{\!\hbox{\fiverm q}}}}^{^{_{_{\prime}}}}\! \
\raise1pt\hbox{\eightrm+}\ t\raise2pt\hbox{$
_{_{^{\!\raise0pt\hbox{\fiverm q{\fivebf+}1} }}}^{^{_{_{\ }}}}\!
$} \hbox{\tenrm
v}^{^{_{_{\prime\prime}}}}\!\!\!\!\!\raise2pt\hbox{$
_{_{^{\!\raise0pt\hbox{\fiverm q{\fivebf+}1} }}}^{^{_{_{\ }}}}\!
$} \
\!\raise1pt\hbox{\eightrm+}\!\dots\!\raise1pt\hbox{\eightrm+}\ \!
t\raise2pt\hbox{$ _{_{^{\!\raise0pt\hbox{\fiverm q{\fivebf+}r}
}}}^{^{_{_{\ }}}}\! $} \hbox{\tenrm
v}^{^{_{_{\prime\prime}}}}\!\!\!\!\!\raise2pt\hbox{$
_{_{^{\!\raise0pt\hbox{\fiverm q{\fivebf+}r} }}}^{^{_{_{\ }}}}\!
$}\mapsto$$
$$\hskip1.0cm \mapsto\big(\ \ \ \!\Sigma\!\!\!\!\!\!\!\!_{_{_{{1\le i\le \hbox{\fiverm
q}}}}} \!\!\!\! t\raise2pt\hbox{$_{_{^{\raise0pt\hbox{\fivei i}
}}}^{^{_{_{\ }}}}\! $} \big)( {{ t\raise2pt\hbox{$
_{_{^{\raise0pt\hbox{\fiverm 1} }}}^{^{_{_{\ }}}}\! $}} \over { \
\Sigma t\raise2pt\hbox{$_{_{_{\!\raise0pt\hbox{\fivei i} }}}$} }}
\hbox{\tenrm v}^{^{_{_{\prime}}}}\!\!\!\raise2pt\hbox{$
_{_{^{\!\raise0pt\hbox{\fiverm 1} }}}^{^{_{_{\ }}}}\! $} \
\!\raise1pt\hbox{\eightrm+}\!\dots\!\raise1pt\hbox{\eightrm+}\ \!
{{ t\raise2pt\hbox{$ _{_{^{\raise0pt\hbox{\fiverm q} }}}^{^{_{_{\
}}}}\! $} } \over { \ \Sigma
t\raise2pt\hbox{$_{_{_{\!\raise0pt\hbox{\fivei i} }}}$} }}
\hbox{\tenrm v}^{^{_{_{\prime}}}}\!\!\!\raise2pt\hbox{$
_{_{^{\!\raise0pt\hbox{\fiverm q} }}}^{^{_{_{\ }}}}\! $}) \
\raise1pt\hbox{\eightrm+}\big(\ \ \ \ \
\Sigma\!\!\!\!\!\!\!\!\!\!\!\!\!\! _{_{_{{\hbox{\fiverm
q{\fivebf+}}1\le j\le \hbox{\fiverm q{\fivebf+}r}}}}}
\!\!\!\!\!_{\!} t\raise2pt\hbox{$_{_{^{_{\!}\raise0pt\hbox{\fivei
j} }}}^{^{_{_{\ }}}}\! $} \ \big)( {{ t\raise2pt\hbox{$
_{_{^{\!\raise0pt\hbox{\fiverm q{\fivebf+}1} }}}^{^{_{_{\ }}}}\!
$} } \over { \ \Sigma
t\raise2pt\hbox{$_{_{_{\!\raise0pt\hbox{{\fiverm }{\fivebf}{\fivei
j} } }}}$} }} \hbox{\tenrm
v}^{^{_{_{\prime\prime}}}}\!\!\!\!\!\raise2pt\hbox{$
_{_{^{\!\raise0pt\hbox{\fiverm q{\fivebf+}1} }}}^{^{_{_{\ }}}}\!
$} \
\!\raise1pt\hbox{\eightrm+}\!\dots\!\raise1pt\hbox{\eightrm+}\ \!
{{ t\raise2pt\hbox{$ _{_{^{\!\raise0pt\hbox{\fiverm q{\fivebf+}r}
}}}^{^{_{_{\ }}}}\! $} } \over { \ \Sigma
t\raise2pt\hbox{$_{_{_{\!\raise0pt\hbox{{\fiverm }{\fivebf}{\fivei
j} } }}}$} }} \hbox{\tenrm
v}^{^{_{_{\prime\prime}}}}\!\!\!\!\!\raise2pt\hbox{$
_{_{^{\!\raise0pt\hbox{\fiverm q{\fivebf+}r} }}}^{^{_{_{\ }}}}\!
$}).
\qed
$$
\enddemo
\normalbaselines

So; $\eta\!:\!(|\Sigma_{_{\!1}}\hbox{\lower1pt\hbox{\mynabla{7}}}
\Sigma_{_{\!2}}\!|,
|\Sigma_{_{\!1}}\hbox{\lower1pt\hbox{\mynabla{7}}}
\Sigma_{_{\!2}}\!|\setminus\!_{_{^o}}{_{\!}}
{(\!\widetilde{\alpha_{_{\!1}}\!,\alpha_{_{\!2}}\!})}\})
\rlap{$\longrightarrow$}{\raise4pt\hbox{$\ \simeq$}}\ \
(|\Sigma_{_{\!1}}\!|\ _{\!}{{\mybarnabla{7}}} |\Sigma_{_{\!2}}\!|,
|\Sigma_{_{\!1}}\!|\ _{\!}{{\mybarnabla{7}}} |\Sigma_{_{\!2}}\!|
\setminus\!_{_{^o}}{_{\!}}
\{(\alpha_{_{\!1}}\!,\alpha_{_{\!2}}\!)\})$
is a homeomorphism if
\hbox{$\eta(\!\widetilde{\alpha_{_{\!1}}\!,\alpha_{_{\!2}}\!})\!=\!
(\!\alpha_{_{\!1}}\!,\alpha_{_{\!2}}\!)$}
and
it's easilly seen that
\hbox{{\bf k}($
|{\Sigma}\!^{^{_{\hbox{\fivebf\ \!\!}}}}|_{_{^{\!\
\!}}}\!{{{\ast}}}\ \!\! |\Delta|)
\!=\!
|{\Sigma}\!^{^{_{\hbox{\fivebf\ \!\!}}}}{{{\ast}}} \Delta|
\!=\!
\hbox{\tenbf k}(|{\Sigma}^{^{_{\hbox{\fivebf\ \!\!}}}}\!|_{_{^{\!\
\!}}}{\rlap{{\lower2.5pt\hbox
{\vbox{\moveright0.1pt\hbox{$^{^{\land}}$}}}}}{{\ast}}}\ \!
|\Delta|) $,}
since these spaces have the same topology on their compact
subsets.
Moreover, if
${\Delta^{^{_{\!{{\!}}\prime}}}}\!\!\subset\!{\Sigma}^{^{_{\prime}}}\!$,
${\Delta^{^{_{\!{{\!}}\prime\prime}}}}\!\!\subset\!{\Sigma}^{^{_{\prime\prime}}}$
then
\hbox{${\vert{\Delta^{^{_{\!{{\!}}\prime}}}}_{\!}\vert}
\bar{\ast}\ \!\!\vert\Delta{\!^{^{_{\prime\prime}}}}\!\!\vert $
is a subspace of
${\vert{\Sigma}^{^{_{\prime}}}\!\vert} \bar{\ast}\
\!\!\vert{\Sigma}^{^{_{\prime\prime}}}\!\!\vert.$}

\example{Example}
Let $\hbox{\tenbf Z}^{^{_{\!\hbox{\fivebf +}}}}\!$ be the positive
integers i.e. 1,2,..., regarded as a simplicial complex and let
$|\hbox{\tenbf Z}^{^{_{\!\hbox{\fivebf +}}}}\!|$ be its
realization according to p.\ 5.
Imbed $|\hbox{\tenbf Z}^{^{_{\!\hbox{\fivebf
+}}}}\!|_{_{^{\!d\!}}}\simeq|\hbox{\tenbf Z}^{^{_{\!\hbox{\fivebf
+}}}}\!|$ into the positive $\hbox{\tenbf R}_{_{^{\hbox{\fivebf
1}}}}\!\!\!^{^{_{\hbox{\fivebf }}}}$-axis
in the real plane
${\hbox{\tenbf R}_{_{^{\hbox{\fivebf \!\
}}}}\!\!\!^{^{_{\hbox{\fivebf 2}}}}\!\!:={\hbox{\tenbf
R}_{_{^{\hbox{\fivebf 1}}}}\!\!\!^{^{_{\hbox{\fivebf
}}}}\!\times\!\hbox{\tenbf R}_{_{^{\hbox{\fivebf
2}}}}\!\!\!^{^{_{\hbox{\fivebf \ }}}}}}$
and connect, with a straight line, each of its points with
$(0,1)\!\in\! \hbox{\tenbf R}_{_{^{\hbox{\fivebf \!\
}}}}\!\!\!^{^{_{\hbox{\fivebf 2}}}}$ and denote the result
$\hbox{\tenbf Z}^{^{_{\hbox{\!\fivebf +}}}}\!{{{ \circ}}}\
\!\bullet$
and let $(\cdot)^\tau\!$ denote the topology of $\cdot$, then,
from right to left, we have natural one-one maps implying;
\smallskip

\nointerlineskip
$$
(|\hbox{\tenbf Z}^{^{_{\hbox{\!\fivebf +}}}}\!{{{ \ast}}}
\bullet|_{_{^{\!d}}})^\tau\!
\ \!\mmysubsetneqq\ \!
(\hbox{\tenbf Z}^{^{_{\hbox{\!\fivebf +}}}}\!{{{ \circ}}}\ \!
\bullet)^\tau\!
\ \!\mmysubsetneqq\ \!
(|\hbox{\tenbf Z}^{^{_{\!\hbox{\fivebf
+}}}}\!|_{_{^{\!d\!}}}{\rlap{{\lower2.5pt\hbox
{\vbox{\moveright0.1pt\hbox{$^{^{\land}}$}}}}}{{ \ast}}}\ \!
\bullet)^\tau\!
\ \!\subseteqq\ \!
(|\hbox{\tenbf Z}^{^{_{\!\hbox{\fivebf +}}}}\!\!{{{ \ast}}}
\bullet|)^\tau\!
\ \!=\ \!
(|\hbox{\tenbf Z}^{^{_{\!\hbox{\fivebf +}}}}\!|_{_{^{\!d\!}}}{{{
\ast}}}\ \!\bullet)^\tau\!.
$$

\nointerlineskip
The first ``\mysubsetneqq" \ follows from a simple check of metric
topologies.
The equality follows since the topology of both these
non-metrizable {\bf k}-spaces, i.e. also of the first, is
determined as the identification topology w.r.t. the canonical
quotient map from
$ |\hbox{\tenbf Z}\!^{^{_{\hbox{\fivebf +}}}}\!|{{{\times}}}\ \!
\hbox{\tenbf I}$ where {\bf I} is the unit interval.
Let $V_{_{^{{ \!n}}}}\!$ be all points on the segment from the
cone point $w_{_{^{{ \!0}}}}\!:=(0,1)$ within $1/n$ from $w_{_{^{{
\!0}}}}\!$, and $V\!\!:=\!\bigcup_{_{{_{\!\!\!\!\!n}}}}\!V_{_{^{{
\!n}}}}\!$.
$V$ isn't open in
$
(\hbox{\tenbf Z}{\raise5pt\hbox{\!\fivebf +}}\!{{{ \circ}}}
\bullet)^\tau\!
$
but it is open in the final topology
$
(|\hbox{\tenbf Z}{\raise5pt\hbox{\!\fivebf +}}\!\!{{{ \ast}}}
\bullet|)^\tau\!
$
cp. \cite{6} %
p.\ $\!$127.
The last paragraph prior to this example gives the rest since {\bf
k}-iffications preserves the subset-relations w.r.t. topologies,
implying that the ``$\ \!\subseteqq\ \!$" is an equality
\underbar{iff}\ \ \
$|\hbox{\tenbf Z}^{^{_{\!\hbox{\fivebf
+}}}}\!|_{_{^{\!d\!}}}{\rlap{{\lower2.5pt\hbox
{\vbox{\moveright0.1pt\hbox{$^{^{\land}}$}}}}}{{ \ast}}}\ \!
\bullet $\ \ \
is a \hbox{\tenbf k}-space.
\endexample

$\dim (\Sigma\times \Delta)=\dim \Sigma+ \dim \Delta$
\ and \
$\dim(\Sigma\ast \Delta)=\dim \Sigma+ \dim \Delta +1$.

If $\alpha_i\!\in \hbox{\tenrm
Int}\sigma\!_{i}\!\subset\!|\Sigma_{i}|,$
${(\widetilde{\alpha_{_{^{\!1}}},\alpha_{_{^{\!2}}}})}
:=\eta^{{_{-1}}}\!(\alpha_{_{\!1}}\!,\alpha_{_{\!2\!}})
\!\in\hbox{\tenrm Int}\sigma
\subset|\Sigma_{_{^{\!1}}}\times\Sigma_{_{^{\!2}}}|$ and
$c_{\sigma}\!:=\!\dim\sigma_{_{\!1}}\!+
\dim\sigma_{_{\!2}}\!-\dim\sigma\  \hbox{\tenrm then;}$
$c_{\sigma}\!\ge\!0\ \hbox{\tenrm and}\ [c_{\sigma}\!\!=\!0 $
\underbar{iff} $\sigma$ is a maximal simplex in
$\bar\sigma_{1}\times\bar\sigma_{2}\!\subset\!\Sigma_{1}\times\Sigma_{2}].\
$

\proclaim{Corollary} {\rm(to Th.\ 6)}
{Let $\hbox{\tenbf G}$, $\hbox{\tenbf G}^{\prime}$ be arbitrary
modules over a  {\bf PID} {\bf R} such that
$\hbox{\tenrm Tor}_1^{\hbox{\fivebf R}}(\hbox{\tenbf
G},\hbox{\tenbf G}^\prime)=0$, then, for any
$\underline{\emptyset_o\ne\sigma}\in\Sigma_{1}\times\Sigma_{2}$
with $\eta\big(\hbox{\tenrm Int}(\sigma_{})\big)\!\subset
\hbox{\tenrm Int}(\sigma_{1})\times\hbox{\tenrm Int}(\sigma_{2});$
}

\smallskip%
\leftline{\indent $
{\underline{\underline{\mdoubleH_{i+c_{\!\sigma}+1} (\hbox{\tenrm
Lk}\!\!\!\!\!\!
 \lower1.1pt\hbox{$_{_{\Sigma_{1}\!\times\Sigma_{2}}}$}
\!\!\!\! \sigma; \hbox{\tenbf G}\otimes_{_{\hbox{\fivebf
R}}}{\hbox{\tenbf G}}^{\prime})}}}
\ \! {{{_{\hbox{\fivebf R}}}}\atop{{\raise2pt\hbox{$\cong$}} }} \
\!
{\rlap{$_{_{_{{{{{p+q=i}\atop{p,q\ge -1}}}}}}}$} {\ \
\raise2pt\hbox{$\ \bigoplus$}}}\ [\mdoubleH_{p}
           (\hbox{\tenrm Lk}\!\!
\lower1.1pt\hbox{${_{_{\Sigma\!_{_{1}}}}}$}\!\!\!\sigma_{1};\hbox{\tenbf
G}) \otimes_{_{\hbox{\fivebf R}}}
        \mdoubleH_{q}
(\hbox{\tenrm Lk}\!\!
\lower1.1pt\hbox{${_{_{\Sigma\!_{_{2}}}}}$}\!\!\!\sigma_{2});\hbox{\tenbf
G}^{\prime})]\oplus $ }

\smallskip
\rightline
{${\rlap{$\oplus$}{\ {\rlap{$_{_{_{_{_{{{p+q=i-1}\atop{p,q\ge
-1}}}}}}}$} {\ \ \raise2pt\hbox{$\ \bigoplus$}}}}} \hbox{\tenrm
Tor}_1^{\hbox{\fivebf R}} \bigl(\mdoubleH_{p}
           (\hbox{\tenrm Lk}\!\!
\lower1.1pt\hbox{${_{_{\Sigma\!_{_{1}}}}}$}\!\!\!\sigma_{1};\hbox{\tenbf
G}),\
 \mdoubleH_{q}
      (\hbox{\tenrm Lk}\!\!
\lower1.1pt\hbox{${_{_{\Sigma\!_{_{2}}}}}$}\!\!\!\sigma_{2};\hbox{\tenbf
G}^{\prime})\bigr)
\ \! {{{_{\hbox{\fivebf R}}}}\atop{{\raise2pt\hbox{$\cong$}} }} \
\!
{\underline{\underline{\mdoubleH_{i+1}
        (\hbox{\tenrm Lk}\!\!\!\!\!\!
 \lower1.1pt\hbox{$_{_{\Sigma_{1}\!\ast\Sigma_{2}}}$}
\!\!\!\! (\sigma_1
\lower0.5pt\hbox{\eightsy {\char"5B}} 
\sigma_2); \hbox{\tenbf G}\otimes_{_{\hbox{\fivebf
R}}}\hbox{\tenbf G}^{\prime})}}}.$}

\nobreak
\smallskip
\indent
So, if $\emptyset\!_{_{^{o}}}\!\ne\!\sigma$ and $c_{\sigma}\!=\!0$
then $\mdoubleH_{i} (\hbox{\tenrm Lk}\!\!\!\!\!\!
 \lower1.1pt\hbox{$_{_{\Sigma_{1}\!\times\Sigma_{2}}}$}
\!\!\!\! \sigma;\hbox{\tenbf G})
\ \! {{{_{\hbox{\fivebf R}}}}\atop{{\raise2pt\hbox{$\cong$}} }} \
\!
\mdoubleH_{i} (\hbox{\tenrm Lk}\!\!\!\!\!\!
 \lower1.1pt\hbox{$_{_{\Sigma_{1}\!\ast\Sigma_{2}}}$}\!\!(\sigma_{_{\!1}}\!
\lower0.5pt\hbox{\eightsy {\char"5B}} 
\sigma_{_{\!2}})
                ;\hbox{\tenbf G})$ \ \ and

\smallskip
\noindent $\mdoubleH_{_{0}}
 (\hbox{\tenrm Lk}\!\!\!\!\!\!
 \lower1.1pt\hbox{$_{_{\Sigma_{1}\!\times\Sigma_{2}}}$}
\!\!\!\! \sigma; \hbox{\tenbf G}\otimes\hbox{\tenbf G}^{\prime})
\ {{{_{\lower0pt\hbox{\fivebf R}}}}\atop{\raise2pt\hbox{$\cong$}
}}\
{\mdoubleH_{_{0}}(\hbox{\tenrm Lk}\!\!
 \lower1.1pt\hbox{${_{_{\Sigma\!_{_{1}}}}}$}\!\!\!\sigma\!_{_{1}}; \hbox{\tenbf G})
\otimes\mdoubleH\!\!_{_{-\!1}}\! (\hbox{\tenrm Lk}\!\!
 \lower1.1pt\hbox{${_{_{\Sigma\!_{_{2}}}}}$}\!\!\!\sigma\!_{_{2}};
\hbox{\tenbf G}^{\prime}) \oplus \mdoubleH\!\!_{_{-\!1}}\!
(\hbox{\tenrm Lk}\!\!
\lower1.1pt\hbox{${_{_{\Sigma\!_{_{1}}}}}$}\!\!\!\sigma\!_{_{1}};
 \hbox{\tenbf G})}
\otimes
        \mdoubleH_{_{0}}(\hbox{\tenrm Lk}\!\!
 \lower1.1pt\hbox{${_{_{\Sigma\!_{_{2}}}}}$}\!\!\!\sigma\!_{_{2}};
       \hbox{\tenbf G}^{\prime}).$
\endproclaim

\demo{Proof}
Note that
$\sigma\!\ne\emptyset_o\!\Rightarrow\!\sigma_{_{\!\!j}}\!\!\ne\emptyset_o,\
j\!=\!1,2.$
The isomorphisms of the underlined modules are, by $\hbox{\tenrm
Proposition\ 1\ p.\ 11},$
Theorem\ 6{\bf i} p.\ 11 in simplicial disguise, and holds even
without the {\bf PID}-assumption.
Prop.\ 1 p.\ 30 and Theorem\ 4\ p.\ 10 gives the second isomorpism
even for $\sigma\!_{_{1}}\!\!=\!\emptyset\!_{_{^{o}}}\
\hbox{\tenrm and}\ \!\!/\ \!\!\hbox{\tenrm or}\
\sigma\!_{_{2}}\!\!=\!\emptyset\!_{_{^{o\!}}}.$
\qed
\enddemo

\subhead {\rm I:}2
Simplicial Connectedness Properties under Products and Joins
\endsubhead

\normalbaselines

\smallskip%
Any $\Sigma$ is representable as {$\Sigma=
\bigcup\!\!\!\!\!\!\!\!\!\!_{_{_{_{_{\sigma ^m \in
\Sigma}}}}}\!\!\! \overline {\sigma ^m},$ where $\sigma^m$\
\hbox{\tenrm denotes\ ``maximal\ simplex"$\!$.}}

\definition{Definition 1}
Two maximal faces $\sigma,\tau\!\in\!\Sigma$ are {\it strongly
connected} if they can be connected by a finite sequence $\sigma=
{{\delta\!_{_{0}}}},..,{\delta\!_{_{i}}},..,{{\delta\!_{_{^{q}}}}}\!
= \tau $ of maximal faces with
$\#({\delta\!_{_{i}}}\!\cap{\delta\!_{_{i+1}}}\!)=
\hbox{max}\!\!\!\!\!\!\!\!\!\!\!\! _{_{_{{{0\le j\le q}}}}}\!
\#{\delta\!_{_{^{j}}}}\!-_{\!}1$ for consecutives.
Strong connectedness imposes an equivalence relation among the
maximal faces, the equivalence classes of which defines the {\it
maximal strongly connected components} of $\Sigma,$ {\tenrm cp.\
\cite{2} %
p.\ 419ff.}
$\Sigma$ is said to be {\it strongly connected} if each pair of of
its  maximal simplices are strongly connected.\break
\hbox{\spaceskip1.8pt A $submaximal\ face$  has exactly one vertex
less then some maximal face  containing it.}
\enddefinition

\remark{Note}
Strongly connected complexes are pure, i.e.
${_{^{\!\!}}}\sigma {_{^{\!\!\ }}} \!\!\in\!\Sigma{_{^{\!\ }}}\
\!\!\hbox{\tenrm maximal}
\Rightarrow\dim\!\sigma\nobreak=\nobreak\dim\!\Sigma.$
\endremark

\definition{Definition 2}
(cp. \cite{2} %
p.\ 419ff.)
$ \Delta{\raise1.5pt\hbox{\eightmsbm
\char"72}}\Delta\!^{^{\!_{o}}} \!:=\{\delta\in\Delta\ |\
\delta\!\not\in\!\Delta\!^{^{\!_{o}}} \} $ is {\it connected as a
poset} $($\hbox{\tenrm partially ordered set}$)$ w.r.t.
$\!$simplex inclusion if for every pair
$\sigma,\tau\!\in\!\Delta{\raise1.5pt\hbox{\eightmsbm
\char"72}}\Delta\!^{^{\!_{o}}}$ there is a chain
$\sigma\!=\!\sigma\!_{_{^{0}}}\!,\sigma\!_{_{^{1}}}\!,
...,\sigma\!_{_{^{k}}}\!\!=\!\tau $ where
$\sigma\!_{_{^{i}}}\!\!\in\! \Delta{\raise1.5pt\hbox{\eightmsbm
\char"72}}\Delta\!^{^{\!_{o}}} $ and
$\sigma\!_{_{^{i}}}\!\!\subseteq\!\sigma\!_{_{^{\!i+1}}}\!\!$ or
$\sigma\!_{_{^{i}}}\!\!\supseteq\!\sigma\!_{_{^{\!i+1\!}}}$.
\enddefinition

\remark{Note}
$\!$(cf. \cite{10} %
p.\ 162.) $\!\Delta_{\!}{\raise1.5pt\hbox{\eightmsbm
\char"72}}\Delta\!^{^{\!_{o}}}\!,$
             $\{\emptyset_{_{^{\!o}}}\!\}\mmysubsetneqq\!\Delta\!^{^{\!_{o}}}\!$
             is connected as a poset \underbar{iff}
             $|\Delta^{}|\ \!{{{\setminus}\!_{_{^{o}}}}}|\Delta\!^{^{\!_{o}}}|$
             is pathwise connected.
             When $\Delta\!^{^{\!_{o}}}=\{\emptyset_{_{^{\!o}}}\!\},$
             then the notion of connectedness as a
             poset\nobreak\ is\nobreak\ equivalent to
             the usual one for $\Delta^{}$.
             $|\Delta^{}|$ is connected \underbar{iff}
             $|\Delta\!^{^{\!_{(\!1\!)}}}\!|$
             is.
             $\Delta\!^{^{\!_{(\!p\!)}}}\!\!:=
             \{\sigma\in\Delta\ \!\vert\ \!\#\sigma\le p+1\}$.
\endremark

\proclaim {Lemma 1}
{\rm(\cite{10} %
p.\ 163)}
                    { $\Delta{\raise1.5pt\hbox{\eightmsbm \char"72}}\Delta\!^{^{\!_{o}}}$
                     is connected as a poset
                     \underbar{\tenrm iff}
                    to each pair of maximal $\hbox{\tenrm simplices}\
              \sigma,\tau\in\Delta{\raise1.5pt\hbox{\eightmsbm \char"72}}\Delta\!^{^{\!_{o}}}$
                there is a chain in
                $\Delta{\raise1.5pt\hbox{\eightmsbm \char"72}}\Delta\!^{^{\!_{o}}},$
                $\sigma\!=\!\sigma\!_{_{^{0}}}\!\supseteq\!\sigma\!_{_{^{1}}}\!
                 \subseteq\!\sigma\!_{_{^{2}}}\!\supseteq\!...\!
                    \subseteq\!\sigma\!_{_{^{2m}}}\!\!=\!\tau\!$,
               where the $\sigma_{_{{\!\!2i}}}\!$s are maximal faces and
                $\sigma\!_{_{^{2i}}}\!{\raise0.5pt\hbox{\eightmsbm \char"72}}
                 \sigma\!_{_{^{2i+1}}}\!$ and
               $\sigma\!_{_{^{2i+2}}}\!
{\raise0.5pt\hbox{\eightmsbm \char"72}}%
\sigma\!_{_{^{2i+1}}}\!$
                    are situated in different components of
                    $\hbox{\tenrm Lk}\!_{_{^{\Delta}}}\!\!\sigma\!_{_{^{2i+1}}}
                   (i\!=\!0,1,..., m\!-\!1).$}
\qed
\endproclaim

\smallskip
Lemma\ 1 gives Lemma\ 2 found in \cite{1} %
p.\ $_{\!}$1856;

\proclaim {Lemma 2}
{Let $\Sigma$ be a finite$\ \!$-$\ \!$dimensional simplicial
           complex, and assume that \hbox{\tenrm Lk}$_{_{\!\Sigma}}\sigma$ is
connected for all $\sigma\!\in\!\Sigma$, i.e. inkluding\ \
$\emptyset_{_{^{\!o}}}\!\!\in\Sigma$, such that $\dim$\hbox{\tenrm
Lk}$_{_{^{\!\Sigma}}}\!\sigma\!\geq\!1.$ Then $\Sigma$ is pure and
strongly connected.} \qed
\endproclaim

Lemma\ 3 (4) is related to the defining properties for quasi-
(pseudeo-)manifolds.

\smallskip%
Read {\sevensy \char"72} in Lemma\ 3 as ``$\times\!$" or all
through as ``$\ast$" when it's trivially true if any
$\Sigma_{i}\!\!=\!\{\emptyset_{\!o}\}$ and else for $\ast$,
$\Sigma_{_{i}}$ are assumed to be \underbar{connected or
$0$-dimensional}.%
$\ \!$``codim$\lower3.5pt\hbox{$^{\hbox{$\sigma\ge 2$}}$}$"
($\Rightarrow$
$\dim\hbox{\tenrm Lk}_{_{^{\!\Sigma}}}\!\sigma\!\geq\!1$)
                  means that a maximal
                  simplex,  $\lower3.5pt\hbox{$^{\hbox{$\tau$}}$}$
say, containing  $\lower3.5pt\hbox{$^{\hbox{$\sigma$}}$}$ always
fulfills
$\lower3.5pt\hbox{$^{\hbox{$\dim\tau\!\ge\dim\sigma\!+\!2$}}$}$.
$\hbox{\tenbf G}\!_{_{1\!}},\hbox{\tenbf G}\!_{_{2}}\!$ are
\hbox{\tenbf A}-modules.
For definitions of Int$\sigma$ and $\bar\sigma$ see p.\ 30.

\proclaim {Lemma 3}
If $\dim\Sigma_{i}\geq0$ and $v_i:=\dim \sigma_{i}\ (i=\ ,1,2)$
then \hbox{\tenbf D}$_1$-{\rm{\tenbf D} are all equivalent;}

\smallskip
\noindent\hbox{\tenbf D}$_1)$\ $\mdoubleH_{0} (\hbox{\tenrm
Lk}_{(\Sigma_{1}\nabla\Sigma_{2})}\sigma; \hbox{\tenbf
G}_1\otimes\hbox{\tenbf G}_2)=0$
                       for $\emptyset_o \neq \sigma\in
                   \Sigma_{1}{^{_{_{\hbox{\sevensy \char"72} }}}}\Sigma_{2}$,
                    whenever codim$\sigma\ge 2.$

\smallskip
\noindent\hbox{\tenbf D$_1^{\prime})$}\
            $\mdoubleH_{0}
                      (\hbox{\tenrm Lk}_{\Sigma_{i}}\sigma_{i};\hbox{\tenbf G}_i)=0$
                           for $\emptyset_o \neq \sigma_{i}\in \Sigma_{i}$,
                         whenever codim$\sigma_{i}\ge2\ (i=1,2)$.

\smallskip
\noindent\hbox{\tenbf D$_2)$} $\mdoubleH_{v+1}
{(\Sigma_{1}\!\lower1pt\hbox{$^{_{_{\hbox{\sevensy
\char"72}}}}$}\! \Sigma_{2}}, \hbox{\tenrm
cost}_{_{\!\Sigma_{1}\!\!\nabla\!\Sigma_{2}}}\!\!\sigma;
\hbox{\tenbf G}_1\otimes\hbox{\tenbf G}_2)=0$
                       for $\emptyset_o \neq \sigma\in
                   \Sigma_{1}\lower1pt\hbox{$^{_{_{\hbox{\sevensy \char"72}}}}$}\!\Sigma_{2}$,
                   if
                   codim$\sigma\ge 2.$

\smallskip
\noindent\hbox{\tenbf D$_2^{\prime})$} $\mdoubleH_{v_i+1}
                    (\Sigma_{i},\hbox{\tenrm
                     cost}\lower2pt\hbox{\sixrm {\char"06}}_{_{\!i}}\!\sigma\!_{_{^{i}}}\!;\hbox{\tenbf G}_i)=0$
                           for $\emptyset_o \neq \sigma_{i}\in \Sigma_{i}$,
                           whenever codim$\sigma_{i}\ge2\ (i=1,2)$.

\smallskip
\noindent\hbox{\tenbf D$_3)$}\ \ $\!\mdoubleH_{v+1}
                       (|\Sigma_{1}{^{_{_{\hbox{\sevensy \char"72}}}}}\!\Sigma_{2}|,
                         |\Sigma_{1}{^{_{_{\hbox{\sevensy \char"72}}}}}_{\!}\Sigma_{2}|
                         \setminus_{_{\!^o}}{_{\!}}\alpha;
                        \hbox{\tenbf G}_1\!\otimes\hbox{\tenbf G}_2)=0$
                     for all $\alpha_0\ne\alpha_{\!}\in\hbox{\tenrm Int}(\sigma)$
                           if codim$\ \!\sigma\!_{\!}\ge2.$

\smallskip
\noindent\hbox{\tenbf D$_3^{\prime})$}\
                  $\mdoubleH_{v_i+1}
                      (|\Sigma_{i}|,|\Sigma_{i}|\setminus\!_{_{\!^o}}{_{\!}}
                      \alpha_i;\hbox{\tenbf G}_i)=0$
                            for $\alpha_0\ne \alpha_i\in \hbox{\tenrm Int}(\sigma_{i})$,
                            if codim$\sigma_{i}\ge2\ (i=1,2)$.

\smallskip
\noindent \hbox{\tenbf D$)$} $\mdoubleH_{v+1}
                   (\ \!|\Sigma_{1}|\lower1pt\hbox{$^{_{_{\hbox{\sevensy \char"72}}}}$}\!|\Sigma_{2}|,
                    |\Sigma_{1}|\lower1pt\hbox{$^{_{_{\hbox{\sevensy \char"72}}}}$}\!|\Sigma_{2}|
                    \setminus\!_{_{\!^o}}{_{\!}}
                    (\alpha_1,\alpha_2);
                   \hbox{\tenbf G}\!_{_{1}}\!\!\otimes\!\hbox{\tenbf G}\!_{_{2}})\!\!=\!0\!$
              for all $\alpha_0\ne{(\widetilde{\alpha_1,\alpha_2})}\!\in\!
             \hbox{\tenrm Int}(\sigma_{ })\break
             \subset|\Sigma_{1\!}
             \lower1pt\hbox{$^{_{_{\hbox{\sevensy \char"72}}}}$}\!\Sigma_{2}|$
\noindent
if codim$\sigma\ge2$, where
$\eta\!:|\Sigma_{1}\!\lower1pt\hbox{$^{_{_{\hbox{\sevensy
\char"72}}}}$}\!\Sigma_{2}|\ \!
\rlap{$\longrightarrow$}{\raise4pt\hbox{$\ \simeq$}}\ \ \ \!
|\Sigma_{1}|{{\bar{\lower1pt\hbox{$^{_{_{\hbox{\sevensy
\char"72}}}}$}\!}}}|\Sigma_{2}|$ and
${\eta(\widetilde{\alpha_1,\alpha_2})}\!=\!(\alpha_1,\alpha_2),$\break
{\rm({\bf k}-iffikations never effect the homology modules.)}
\endproclaim

\demo{Proof}
By the homogenity of the interior of
$\vert\bar\sigma\!_{_{^{1}}}\!\times\bar\sigma\!_{_{^{2}}}\!\vert$
we only need to deal with simplices $\sigma$ fulfilling
$c_{_{^{\!\sigma}}}\!\!\!=\!0$. Proposition\ 1 p.\ 11
+ p.\ 15 top, implies that all non-primed resp. $\!$primed items
are equivalent among themselves. The above connectedness
conditions are not coefficient sensitive, so suppose
$\hbox{\tenbf G}_i\!:=\!\hbox{\tenbf k}$, a field.
\hbox{\tenbf D$\!_1$} $\Leftrightarrow$ \hbox{\tenbf D$\!_1^{\ \!
\prime}\ $} by Cor. p.\ 15.
For joins of finite complexes, this is done explicitly in \cite{10} %
p.\ 172.
\qed
\enddemo

\proclaim {Lemma 4}
{\rm(\cite{8} %
p.\ 81 gives a proof, valid for any finite-dimensional
complexes.)}

\smallskip
\noindent
\hbox{{\tenbf A}$)$} If $d_i\!:=\dim\Sigma_{i}\geq0$
then$;$ {$\Sigma_{1} \!\times\! \Sigma_{2}$
is pure        
$\Longleftrightarrow$ $\Sigma_{1}$ and $\Sigma_{2}\ \hbox{\tenrm
are\ both\ pure.}$}

\smallskip
\noindent\hbox{\tenbf B}$)$ If $\dim\sigma_{i} ^m\!\geq1$ for each
maximal simplex $\sigma_{i}^m \!\in\! \Sigma_{i}$ then$;$

\smallskip
Any submaximal face in $\Sigma_{_{\!1}}
 \!\times\!
\Sigma_{_{\!2}}\!$ lies in at most $($exactly$)$ two maximal faces
$\Longleftrightarrow $
Any submaximal face in $\Sigma_{_{\!i}}\!$ lies in at most
$($exactly$)$ two maximal faces of $\Sigma_{_{\!i}},\ i=1,2.$

\smallskip
\noindent\hbox{\tenbf C}$)$ If $d_i>0$
then;
$\Sigma_{_{\!1}} \!\times\! \Sigma_{_{\!2}}$ strongly connected
$\Longleftrightarrow$ $\Sigma_{_{\!1}}$, $\Sigma_{_{\!2}}$ both
strongly connected. \qed
\endproclaim

\remark{Note} Lemma\ 4 is true also for $\ast$ with exactly the
same reading but now with no other restriction than that
$\Sigma_{_ {\!i}}\!\ne\emptyset$ and this includes in particular
item {\bf B}.
\endremark

%


\head II: Concepts Related to Combinatorics and Commutative
Algebra
\endhead

\subhead
{{\rm II:1} Definition of Stanley-Reisner rings}
\endsubhead

\medskip
$\!${\teni Stanley}-{\teni Reisner} (St-Re) {\teni ring} {\teni
theory} is a basic tool within combinatorics, where it supports
the use of commutative algebra.
The definition of $m_\delta$, in
${\hbox{\tenbf I}}_{_{^{\!\Delta}}}\!\!=
\hbox{\tenbf(}\{m_\delta \ |\ \delta \!\not\in\!\Delta
\}\hbox{\tenbf)}$,
in existing literature is so vague that it allows You to state
nothing but
${\hbox{\tenbf A}} \hbox{\tenbf[$\!\!{_{_{\!}}}{_{_{\!}}}$[}
\{\emptyset_{_{\!^o}}\!\}
\hbox{\tenbf]$\!\!{_{_{\!}}}{_{_{\!}}}$]}=\ \!?\!=
          {\hbox{\tenbf A}}\hbox{\tenbf[$\!\!{_{_{\!}}}{_{_{\!}}}$[}\
\!\emptyset\ \!\hbox{\tenbf]$\!\!{_{_{\!}}}{_{_{\!}}}$]}$.
Our definition of $m_\delta$ below rectifies this and we conclude
that
${\hbox{\tenbf A}} \hbox{\tenbf[$\!\!{_{_{\!}}}{_{_{\!}}}$[}
\{\emptyset_{_{\!^o}}\!\}
\hbox{\tenbf]$\!\!{_{_{\!}}}{_{_{\!}}}$]}\cong \hbox{\tenbf A} \ne
\!0\!=
           $``The
           trivial ring"$= {\hbox{\tenbf A}}\hbox{\tenbf[$\!\!{_{_{\!}}}{_{_{\!}}}$[}\
\!\emptyset\ \!\hbox{\tenbf]$\!\!{_{_{\!}}}{_{_{\!}}}$]}$.

\definition{Definition} {A subset $s\!\!\subset\!\!\hbox{\tenbf W}\!\!\supset\!\!
V_\Delta$ is said to be a \hbox{\teni non}-\hbox{\teni simplex}
(w.r.t.\ \hbox{\tenbf W}) of a simplicial complex $\Delta$,
denoted $s \rlap{\raise4pt\hbox{$\ n$}}{\notin}\Delta,\ {if}\
s\!\not\in\!\Delta$ but
 $\dot s\!=\!{(\bar s)}^{(dim\ \!s)-1}\!\!\!\!\subset\!\Delta\ ($
i.e.\ the $(\dim s-1)$-dimensional skeleton of $\bar s$,
consisting of all proper subsets of $s$, is a subcomplex of
$\Delta).$ For a simplex $\delta = \{ \hbox{\tenrm v}_{i_1},\dots,
\hbox{\tenrm v}_{i_k} \}$ we define $m_{\delta}$ to be the
squarefree monic monomial
 $m_{\delta} :=1_{_{\hbox{\fivebf A}}}\!\!\cdot \hbox{\tenrm v}
 _{i_1}\!\!\cdot\!\dots\!
\cdot \hbox{\tenrm v}_{i_k}\in \hbox{\tenbf A}[\hbox{\tenbf W}]$
where $ \hbox{\tenbf A}[\hbox{\tenbf W}]$ is the \hbox{graded
polynomial algebra} on the variable set $\hbox{\tenbf W}$ over the
commutative ring $\hbox{\tenbf A}$ with unit $1\!_{_{\hbox{\fivebf
A}}}.$
\hbox{\tenrm So},\ $m_{\emptyset\!_o}\!\!=1\!_{_{\hbox{\fivebf
A}}}.
$

Let $\hbox{\tenbf A}
{\hbox{\tenbf[}\!\!{_{_{\!}}}{_{_{\!}}}\hbox{\tenbf[} } \Delta
{\hbox{\tenbf]}\!\!{_{_{\!}}}{_{_{\!}}}\hbox{\tenbf]}} :=
\hbox{\tenbf A}[\hbox{\tenbf W}]/{\hbox{\tenbf I}_{\Delta}}$ where
                     {\bf I}$_{\Delta}$
is the ideal generated by $\{m_\delta \ |\ \delta
\rlap{\raise4pt\hbox{$\ n$}}{_{\!}\!_{\!}\not\in}\ \!\Delta \}$.
$\hbox{\tenbf A}
\hbox{\tenbf[}\!\!{_{_{\!}}}{_{_{\!}}}\hbox{\tenbf[} \Delta
\hbox{\tenbf]$\!\!{_{_{\!}}}{_{_{\!}}}$]}$ is called the ``face
ring" or ``Stanley-Reisner $($St-Re$)$ ring" of $\Delta$ over
$\hbox{\tenbf A}$. Frequently $\hbox{\tenbf A}\!=\hbox{\tenbf k},$
a field.}
\enddefinition

\remark{Note}{\bf i.}
Let $\overline{\hbox{\tenbf P}}$ be the set of finite subsets of
the set $\hbox{\tenbf P}$ then;
$\hbox{\tenbf A}
\hbox{\tenbf[}\!\!{_{_{\!}}}{_{_{\!}}}\hbox{\tenbf[} \ \!
\overline{\hbox{\tenbf P}}\ \!
\hbox{\tenbf]$\!\!{_{_{\!}}}{_{_{\!}}}$]}=\hbox{\tenbf A[{\bf
P}]}.$
$ \overline{\hbox{\tenbf P}}$ is known as ``{\it the full
simplicial complex on $\hbox{\tenbf P}$}" and the natural numbers
{\bf N} gives
$\overline{\hbox{\tenbf N}}$ as ``{\it the infinite simplex}$\
\!$". \  $\Delta=\bigcap\!\!\!\!\!\!\!\!\!
_{_{_{_{\lower4.5pt\hbox{$^{s}$}
\rlap{\raise4pt\hbox{$\lower7.5pt\hbox{$^{^{\ n}}$}$}}{\not\in}\
\!\! {\lower4.5pt\hbox{$^{\Delta}$}}}}}} \!\!\hbox{\tenrm
cost}\!\!\!_ {_{_{\overline{\hbox{\fivebf
W}}}}}\!\!{{{s}}\!_{_{^{\ }}}}\!$. See pp.\ 30-1 for
$\!$definitions.

\noindent {\bf ii.}
$\hbox{\tenbf A}
\hbox{\tenbf[}\!\!{_{_{\!}}}{_{_{\!}}}\hbox{\tenbf[} \Delta
\hbox{\tenbf]$\!\!{_{_{\!}}}{_{_{\!}}}$]}
\ \!\cong\ \!
\lower0.0pt\hbox{\sixrm${ {{\lower0.0pt\hbox{${\hbox{\eightbf
A}}$}}[{\lower0pt\hbox{\eighti V}}\!\!_{_{\Delta}}] \over
(\{\hbox{\seveni m}\!
\lower2.5pt\hbox{\fivei {\char"0E}} 
\in{\lower0.0pt\hbox{${\hbox{\eightbf A}}$}}
[{\lower0pt\hbox{\eighti V}}\!\!_{_{\Delta}}]
\lower0.0pt\hbox{\sevensy {\char"6A}}         
\ \lower0.5pt\hbox{\seveni {\char"0E}}\ \!    
\rlap{\raise4pt\hbox{$\lower7.5pt\hbox{$^{^{\
n}}$}$}}{^{\!}\not\in}\ \!\!
{\lower0.0pt\hbox{\eightrm {\char"01}}        
} \})}}$} \raise3pt\hbox{\ \ },$ if $\Delta\neq \emptyset$,
           ${\{{\emptyset}\!_o\!\}}.$
           So, the choice of the universe {\bf W}
           isn't all that critical.
If  $\Delta=\emptyset$, then
           the set of non-simplices equals $\{\emptyset\}$,
since
           $\emptyset_o\!\not\in\emptyset$, and
           ${\overline{\emptyset}\!_o\!\!}
           ^{_{((\dim\emptyset_{\!o})-1)}}\!\!\!\! =\!
           \overline{\{\emptyset_o\}}^{_{(-2)}}\!\!=\emptyset
           \subset\emptyset$, implying; $\hbox{\tenbf A}[\ \!\emptyset\ \!]\!=\!0\!=
           $``The
           trivial ring", since\nobreak\ $m_{\emptyset\!_{_{^o}}}\!\!=
           \!1_{_{\!\hbox{\fivebf A}}}$.

Since $\emptyset \in \Delta$ for
           every simplicial complex $\Delta\neq \emptyset$,
           ${\{\hbox{\tenrm v}\}}$ is a non-simplex of
           $\Delta$\nobreak\ for\nobreak\ every
           $\hbox{\tenrm v}\!\in\!\hbox{\tenbf W}\setminus \hbox{\teni V}\!_{\Delta}$,
           i.e.\ [$\hbox{\tenrm v}\!\not\in\!\hbox{\teni V}\!_{\Delta}\!\neq\!\emptyset
           {\rlap{{\rlap{$]$}{$\ \Leftrightarrow$}}}{\ \ \ \ [}}\{\hbox{\tenrm v}\}\
            \rlap{\raise4pt\hbox{$\ n$}}_{{\!}}{\not\in}\ \Delta
           \neq \emptyset].$
           So
$\hbox{\tenbf A}
\hbox{\tenbf[}\!\!{_{_{\!}}}{_{_{\!}}}\hbox{\tenbf[}
\{\emptyset_{_{\!^o}}\!\}
\hbox{\tenbf]$\!\!{_{_{\!}}}{_{_{\!}}}$]}= \hbox{\tenbf A}$ since
$\{{\delta}\
            \rlap{\raise4pt\hbox{$\ n$}}_{{\!}}{\not\in}\ \Delta\}= \hbox{\tenbf W}.$

\medskip
\noindent {\bf iii.}
$\hbox{\tenbf k}
{\hbox{\tenbf[}\!\!{_{_{\!}}}{_{_{\!}}}\hbox{\tenbf[}}
\Delta_{_{^{_{\!}1\!}}}\!\!\ast\!\Delta_{_{^{_{\!}2\!}}}
{\hbox{\tenbf]$\!\!{_{_{\!}}}{_{_{\!}}}$]}\!}
$
$\ \!\cong\ \!$
$\hbox{\tenbf k}
{\hbox{\tenbf[}\!\!{_{_{\!}}}{_{_{\!}}}\hbox{\tenbf[}}
\Delta_{_{^{_{\!}1\!}}}
{\hbox{\tenbf]$\!\!{_{_{\!}}}{_{_{\!}}}$]}}
\ \!{\otimes}\ \!
\hbox{\tenbf k}
{\hbox{\tenbf[}\!\!{_{_{\!}}}{_{_{\!}}}\hbox{\tenbf[}}
\Delta_{_{^{_{\!}2\!}}}
{\hbox{\tenbf]$\!\!{_{_{\!}}}{_{_{\!}}}$]}}
$
$
\ \hbox{\tenrm with}\
\hbox{\tenbf I}
{\lower2.0pt\hbox{\sevenrm {\char"01}$_{_{1}}\!$ {\sevensy
{\char"03}}\sevenrm {\char"01}$_{_{2}}$
}}
\!=
\hbox{\tenrm (} \{m_{\delta}\ \!\big\vert\ \!\![
{\delta}\ \rlap{\raise4pt\hbox{$\ n$}}_{{\!}}{\not\in}\
\Delta_{_{^{_{\!}1\!}}}\!\!\lor
{\delta}\ \rlap{\raise4pt\hbox{$\ n$}}_{{\!}}{\not\in}\
\Delta_{_{^{_{\!}2\!}}}\!]
\land[\delta\!\notin\!{\Delta_{_{^{_{\!}1\!}}}\!
\!\ast\!{\Delta_{_{^{_{\!}2\!}}}\!}}]\}\hbox{\tenbf )} $
by \cite{8} %
Example 1\ p.\ $\!$70.

\medskip
\noindent {\bf iv.}
If $\Delta_i\neq\emptyset\ i=1,2$, it's well known that

{\bf a.}
\hbox{\tenbf I}$\!\!$
{\lower2.0pt\hbox{\sevenrm {\char"01}$_{_{1}}\!$        
{\sevensy {\char"5B}}\sevenrm {\char"01}$_{_{2}}$
}}
{$\!=\!
\hbox{\tenbf I}\!
{\lower2.5pt\hbox{\sevenrm {\char"01}$_{_{^1}}\!$       
}}
\!\cap
\hbox{\tenbf I}\!
{\lower2.5pt\hbox{\sevenrm {\char"01}$_{_{^2}}$
}}
\!=
(\{m=\hbox{\tenrm Lcm}(m_{\delta_1},m_{\delta_2})\ \!\big\vert\ \!
\delta_i{\rlap{\raise4pt\hbox{$\ n$}}{\notin}\ \!\Delta_i},\
\!i=1,2 \})$}

\indent {\bf b.} \hbox{\tenbf I}$\!\!$
{\lower2.0pt\hbox{\sevenrm {\char"01}$_{_{1}}\!$     
{\sevensy {\char"5C}}\sevenrm {\char"01}$_{_{2}}$
}}
$\!=\!
\hbox{\tenbf I}\!
{\lower2.5pt\hbox{\sevenrm {\char"01}$_{_{^1}}\!$     
}}
\!\!+
\hbox{\tenbf I}\!
{\lower2.5pt\hbox{\sevenrm {\char"01}$_{_{^2}}$
}}
= (\{m_{{\delta}}\ \big\vert\ \delta{\rlap{\raise4pt\hbox{$\
n$}}{\notin}\ \!\Delta_1}\lor \delta{\rlap{\raise4pt\hbox{$\
n$}}{\notin}\ \!\Delta_2}\} )$ in {\bf A}[{\bf W}]

\smallskip
$
\hbox{\tenbf I}\!
{\lower2.5pt\hbox{\sevenrm {\char"01}$_{_{^1}}\!$    
}}
\!\cap
\hbox{\tenbf I}\!
{\lower2.5pt\hbox{\sevenrm {\char"01}$_{_{^2}}$
}}
$
and
\hbox{\tenbf I}\!
{\lower2.5pt\hbox{\sevenrm {\char"01}$_{_{^1}}\!$      
}}
\!+
\hbox{\tenbf I}\!
{\lower2.5pt\hbox{\sevenrm {\char"01}$_{_{^2}}$
}}
are generated by a set (no restrictions on its cardinality) of
squarefree monomials, if both {\bf I}$_{_{\Delta_1}}$ and {\bf
I}$_{_{\Delta_2}}$ are.
These squarefree monomially generated ideals form a distributive
sublattice $(\hbox{\tensy J}\ \!^{\circ};\cap,+,\hbox{\tenbf
A[W]}_+)$, of the ordinary lattice structure on the set of ideals
of the polynomial ring {\bf A[W]}, with a counterpart, with
reversed lattice order, called the squarefree monomial rings
\underbar{with unit}, denoted ({\bf
A}$^{\!\circ}${\bf[W]}$;\cap,+,\hbox{\tenbf A})$. {\bf A[W]}$_+$
is the unique homogeneous maximal ideal and zero element. We can
use the ordinary subset structure to define a distributive lattice
structure on ${\Sigma}_{\hbox{\fivebf W}}^\circ:=$ ``The set of
\underbar{non-empty} simplicial complexes over {\bf W}", with
$\{\emptyset\}$ as zero element and denoted
(${\Sigma}_{\hbox{\fivebf W}}^\circ;{\cup},\cap,\{\emptyset\})$.
The {\it Weyman/Fr\"oberg/Schwartau Construction} eliminates the
above
{squarefree}-demand, cf. \cite{27} %
p. 107.
\endremark

\proclaim{Proposition}
$\!{\!}$The Stanley-Reisner Ring Assignment Functor defines a
monomorphism on distributive lattices from
$({\Sigma}_{\hbox{\fivebf W}}^\circ,{\cup},\cap,\{\emptyset\})$ to
$(\hbox{\tenbf A}^{\!\circ}\hbox{\tenbf[W]},\cap,+,\hbox{\tenbf
A})$, which is an isomorphism for finite {\bf W}.
\qed
\endproclaim

%


\subhead
{{\rm II:2} Buchsbaum, Cohen-Macaulay and 2-Cohen-Macaulay
Complexes}
\endsubhead

\medskip
\medskip
We're now in position to give combinatorial/algebraic counterparts
of the weak homology manifolds defined in \S3.4 p.\ 12.
\ Prop.\ 1 and Th.\ 8 below, together with Th.\ 11 p\ 25, show how
simplicial homology manifolds can be inductively generated.

Combinatorialists call a finite simplicial complex a Buchsbaum
({Bbm}$_{\!_{\hbox{\fivebf k}}}\!$) (Cohen-Macaulay
({CM}$_{\!_{\hbox{\fivebf k}}}\!$)) complex if its Stanley-Reisner
Ring is a Buchsbaum (Cohen-Macaulay) ring.
We won't use the ring theoretic definitions of Bbm or C-M and
therefore we won't write them out.
Instead we'll use some homology theoretical characterizations
found in
\tenrm \cite{26} %
pp.$\ \!$73, 60-1  resp. 94 %
to deduce, through Prop.\ 1 p.\ 11, the following consistent
definitions for arbitrary modules and topological spaces.

\definition{Definition}
{$\!X$ is ``{Bbm}$_{_{\!\hbox{\fivebf G}}}\!\!\!"$
($\!$``{CM}$_{_{\!\hbox{\fivebf G}}}\!\!"$,
2-$``${CM}$_{_{\!\hbox{\fivebf G}}}\!\!"$) if} $\!X$ is an
$n$-{whm}$_{_{\!\hbox{\fivebf G}}}\!$
($n$-{jwhm}$_{_{\!\hbox{\fivebf G}}}$,
$n$-whsp$_{_{\!\hbox{\fivebf G}\!}}).\!$
\indent
A simplicial complex {\tensy $\Sigma$} is
``$\ \!${\tenrm Bbm}$_{_{\!\hbox{\fivebf G}}}\!\!",$
``{CM}$_{_{\!\hbox{\fivebf G}}}\!\!"$
resp. 2-$``${CM}$_{_{\!\hbox{\fivebf G}}}\!\!"$ if {\tensy
j$\Sigma$j} is.
In particular,
$\Delta$ is $2$-``{\rm CM}$\!_{_{\hbox{\fivebf G}}}\!\!"\!$ \
\underbar{iff}\ \
$\Delta$ is ``{\rm CM}$\!_{_{\hbox{\fivebf G}}}\!\!"\!$
and\ \
$ \mdoubleH\!\!\!\!_{_{n\!_{_{}}\!-\!1}}\! (\hbox{\tenrm
cost}\!_{_{\Delta}}\!{{\delta\!_{_{\ }}}}\!;\hbox{\tenbf G})=0,\
\! \forall\ \delta\!\in\!\Delta, $
\hbox{cp. \cite{26} Prop.\ 3.7 p.\ 94.}

The {\teni n} in
``$n$-{\tenrm Bbm}" (``$n$-{\tenrm CM}" resp.
2-``{$n$-CM}$_{_{\!\hbox{\fivebf G}}}\!\!"$) is deleted since any
interior point $\alpha$ of a realization of a maxi-dimensional
simplex gives;
$ \mdoubleH\!\!\!\!\!\raise0.5pt\hbox{{$_{_{_{\dim\Delta}}}$}}
\!\!(|\Delta|,|\Delta|\setminus_{_{\!o}\!}\alpha;\hbox{\tenbf G})
=\hbox{\tenbf G}. $
$\sigma\in\Sigma$ is maxi-dimensional if $\dim\sigma = \dim\Sigma$
and ${``\!\setminus_{_{\!o}\!}\!\!"}\!$ is defined in p.\ 4.
\enddefinition

So we're simply renaming
$n$-{whm}$_{_{\!\hbox{\fivebf G}}},$
$n$-{jwhm}$_{_{\!\hbox{\fivebf G}}}$\nobreak\ and
$n$-whsp${_{_{\!\hbox{\fivebf G}}}}_{\!}$
to
$_{\!}$``{Bbm}$_{_{\!\hbox{\fivebf G}}}\!\!\!"_{\!}$,
$\!``${CM}$_{_{\!\hbox{\fivebf G}}}\!\!"\!$\noindent\ resp.

\noindent 2-$``${CM}$_{_{\!\hbox{\fivebf G}}}\!\!"$,
where the quotation marks indicate that we're not limited to
compact spaces nor to just {\bf Z} or {\bf k} as coefficient
modules.
{\bf N.B.}, the definition p.\ 31 of
\rlap{\raise1.5pt\hbox{\tent{\char"5E}}}{\vbox{\moveleft-1.2pt\hbox{\tensy{\char"6A}}}}
$\!\!$%
{\tenrm{\char"01}}$\!_{\!}${{\raise4.5pt\hbox{\fivei{\char"7D}}}}%
$\!$({\teni X})%
\rlap{\raise1.5pt\hbox{\tent{\char"5E}}}{\vbox{\moveleft-1.2pt\hbox{\tensy{\char"6A}}}},
through which each topological space $_{\!}X_{\!}$ can be provided
with a Stanley-Reisner\nobreak\ ring, which, w.r.t.
``{Bbm}$_{_{\!\hbox{\fivebf G}}}\!\!\!"$-,
$\!``${CM}$_{_{\!\hbox{\fivebf G}}}\!\!"$-
and
2-$``${CM}$_{_{\!\hbox{\fivebf G}}}\!\!"$-ness is
\hbox{triangulation invariant.}

\proclaim{Proposition 1}
{The following conditions are equivalent}:\ \ {\rm (We assume}
$\dim\Delta=n.${\rm )}

\indent\indent\indent{\bf a.} $\Delta$ is $\hbox{\tenrm
``Bbm}\!_{_{\hbox{\fivebf G}}}\!\!"\ \!\!$,

\indent\indent\indent{\bf b.} {\rm(Schenzel)} { $\Delta$ is pure
and $\hbox{\tenrm Lk}\!\!_{_{_{\Delta}}}\!\!{{\delta\!_{_{ }}}} \
\hbox{\tenrm is}\ \!\!\hbox{\tenrm ``CM}\!_{_{\hbox{\fivebf
G}}}\!\!"\ \forall\
\!\emptyset\!_{_{^{o}}}\!\!\ne\!{\delta}\!\in\!{\Delta}$, }

\indent\indent\indent{\bf c.} {\rm(Reisner)} \  { $\Delta$ is pure
and $ \hbox{\tenrm Lk}\!\!_{_{_{\Delta}}}\!\!{{\hbox{\tenrm
v}\!_{_{ }}}}\ \hbox{\tenrm is}\ \hbox{\tenrm
``CM}\!_{_{\hbox{\fivebf G}}}\!\!"\ \forall\ \!\hbox{\tenrm
v}\!\in\!\hbox{\teni V}_{_{\!\!\!\Delta}} $. }
\endproclaim

\demo{Proof}
\hbox{\spaceskip2.3pt Use Proposition\ 1 p.\ 11 and Lemma\ 2 p.\
16 and then use Eq. {\bf I} p.\ 30.}
\qed
\enddemo

\example{Example}
When limited to compact polytopes  and {\bf k} as koefficient
module, we add, from \cite{26} %
p.\ 73, %
the following Buchsbaum-equivalence using {\it local\ cohomology};

\smallskip
{\bf d.} {\rm(Schenzel)} \  $\Delta$ {\it is Buchsbaum}
\underbar{\it iff}\ \
$\dim{
_{{\lower2.1pt\hbox{\tenbf--}\!\!\!}}\!_{\!} \hbox{\tenrm {H}} }
^{^{_{i}}}\!\!
\raise0.8pt\hbox{$_{_{\!{\hbox{\fivebf k}
\hbox{\fivebf[}\Sigma\hbox{\fivebf]}_{\!\hbox{\fivebf+}\!\!}}} }$}
(\hbox{\tenbf k} \hbox{\tenbf[$\!\!{_{_{\!}}}$[} \Sigma
\hbox{\tenbf]$\!\!{_{_{\!}}}{_{_{\!}}}$]})
\le \infty
\ if\
0\le i < \dim\hbox{\tenbf k} \hbox{\tenbf[$\!\!{_{_{\!}}}$[}
\Sigma \hbox{\tenbf]$\!\!{_{_{\!}}}{_{_{\!}}}$]}),
$
in which case
${
_{{\lower2.1pt\hbox{\tenbf--}\!\!\!}}\!_{\!} \hbox{\tenrm {H}} }
^{^{_{i}}}\!\!
\raise0.8pt\hbox{$_{_{\!{\hbox{\fivebf k}
\hbox{\fivebf[}\Sigma\hbox{\fivebf]}_{\!\hbox{\fivebf+}\!\!}}} }$}
(\hbox{\tenbf k} \hbox{\tenbf[$\!\!{_{_{\!}}}$[} \Sigma
\hbox{\tenbf]$\!\!{_{_{\!}}}{_{_{\!}}}$]}) \ \!\widetilde=\ \!
\widetilde{\mdoubleH}_{_{i\!-\!1}}\!(|\Sigma|;\hbox{\tenbf k}), $
cf.\ \cite{27} %
p.\ 144 for proof.
Here, ``$\dim$" is {\it Krull} {\it dimension}, which for
Stanley-Reisner Rings is simply 1+ the simplicial dimension.

For
$\Gamma\!\!_{_{^{1}}}\!, \Gamma\!\!_{_{^{2}}}\!$ finite and {CM}$\
\!\! _{_{\hbox{\fivebf k}}}\!\!$
we get the following K\"unneth formula for ring theoretical local
cohomology; (``$\cdot_+$" indicates the unique homogeneous maximal
ideal of ``$\cdot$".)

\medskip
\item{\bf 1.}
$\ \
{{_{{\lower-0.5pt\hbox{%
{
\vbox{\moveright0.0cm\hbox{\vrule width 0.09 true in height 0.05pt
depth 0.015 true cm}}}
}
}%
}}%
\vbox{\moveleft3.65mm\hbox{\hbox{\tenrm {H}}}}}\!^{^{_{\!q}}}
%
%
\vbox{\moveleft0.7cm\hbox{%
\lower3.5pt\hbox{$_{_{\! \hbox{\sixbf(} \hbox{\sixbf k}
\raise1pt\hbox{%
{\hbox{\fivebf[\vbox{\moveleft1.9pt\hbox{[}}}}%
}%
{\hbox{\sixrm{\char"00}}}\!\!_{_{^{1}}}\!\!   
\raise1pt\hbox{%
{\hbox{\fivebf]%
\vbox{\moveleft1.9pt\hbox{]}}}}%
}
\otimes \hbox{\sixbf k}
\raise1pt\hbox{%
{\hbox{\fivebf[\vbox{\moveleft1.9pt\hbox{[}}}}%
}%
{\hbox{\sixrm{\char"00}}}\!\!_{_{^{2}}}\!\!       
\raise1pt\hbox{%
{\hbox{\fivebf]\vbox{\moveleft1.9pt\hbox{]}}}}%
}
\hbox{\sixbf)} _{\!\hbox{\fivebf+}\!\!} }}$}
}}
$%
$ \!\!\! \!\!\! \!\!\! \!\!\! \!\!\! \!\!
( \hbox{\tenbf k} {
\raise1pt\hbox{%
{\hbox{\eightbf[%
\vbox{\moveleft2.45pt\hbox{[}}%
}}%
}
\Gamma\!\!_{_{^{1}}}\!
\raise1pt\hbox{%
{\hbox{\eightbf]%
\vbox{\moveleft2.45pt\hbox{]}}%
}}%
}%
} \otimes \hbox{\tenbf k} {
\raise1pt\hbox{%
{\hbox{\eightbf[%
\vbox{\moveleft2.45pt\hbox{[}}%
}}%
}
\Gamma\!\!_{_{^{2}}}\!
\raise1pt\hbox{%
{\hbox{\eightbf]%
\vbox{\moveleft2.45pt\hbox{]}}%
}}%
}
} ) \ \!\widetilde=\! $
$ \Big[ {{
\hbox{\sevenbf k}
\raise1pt\hbox{%
{\hbox{\sixbf[\vbox{\moveleft1.9pt\hbox{[}}}}%
}%
{\hbox{\sevenrm{\char"00}}}\!\!_{_{^{1}}}\!        
\raise1pt\hbox{%
{\hbox{\sixbf]%
\vbox{\moveleft1.9pt\hbox{]}}}}%
}
\otimes
\hbox{\sevenbf k}
\raise1pt\hbox{%
{\hbox{\sixbf[\vbox{\moveleft1.9pt\hbox{[}}}}%
}%
{\hbox{\sevenrm{\char"00}}}\!\!_{_{^{2}}}\!         
\raise1pt\hbox{%
{\hbox{\sixbf]\vbox{\moveleft1.9pt\hbox{]}}}}%
}
\ \!\tilde=\ \!
\hbox{\sevenbf k}
\raise1pt\hbox{%
{\hbox{\sixbf[\vbox{\moveleft1.9pt\hbox{[}}}}%
}%
\Gamma\!\!_{_{^{1}}}\!\ast \Gamma\!\!_{_{^{2}}}\!%
\raise1pt\hbox{%
{\hbox{\sixbf]\vbox{\moveleft1.9pt\hbox{]}}}}%
},\
\hbox{\sevenrm cf.\ p.$\ \!$17$\ \!$} } \atop { \hbox{\sevenrm
\raise1pt\hbox{\fivebf+} Eq. 3\ p.$\ \!$10
 and  Corollary\ {\sevenrm i} p.\ 12%
} }} \Big]
\ \!\widetilde=\ \ \!\!
\raise0pt\hbox {$
{\rlap{$\!\!_{_{_{{{i+j=}{\hbox{\fiverm q}}}}}}$} {\
\raise2pt\hbox{$\bigoplus$}}} $}\
%
%
{{_{{\lower-0.5pt\hbox{%
{
\vbox{\moveright0.0cm\hbox{\vrule width 0.09 true in height 0.05pt
depth 0.015 true cm}}}
}
}%
}}%
\vbox{\moveleft3.65mm\hbox{\hbox{\tenrm {H}}}}}\!^{^{_{\!i}}}
\vbox{\moveleft0.4cm\hbox{%
%
%
\lower3.5pt\hbox{$_{_{\!
\hbox{\sixbf k}%
\raise1pt\hbox{%
{\hbox{\fivebf[\vbox{\moveleft1.9pt\hbox{[}}}}%
}%
{\hbox{\sixrm{\char"00}}}\!\!_{_{^{1}}}\!\!        
\raise1pt\hbox{%
{\hbox{\fivebf]\vbox{\moveleft1.9pt\hbox{]}}}}%
}%
_{\!\hbox{\fivebf+}\!\!} }}$}
%
}}
$%
$ \!\!\! \!\!\!
( \hbox{\tenbf k}
{
\raise1pt\hbox{\eightbf[$\!{_{\!}}{{_{_{_{\!}}}}}$[}
\Gamma\!\!_{_{^{1}}}\!
\raise1pt\hbox{\eightbf]$\!{_{\!}}{_{_{_{\!}}}}$]}
} )
\!\otimes\!
%
%
{{_{{\lower-0.5pt\hbox{%
{
\vbox{\moveright0.0cm\hbox{\vrule width 0.09 true in height 0.05pt
depth 0.015 true cm}}}
}
}%
}}%
\vbox{\moveleft3.65mm\hbox{\hbox{\tenrm {H}}}}}\!^{^{_{\!j}}}
%
%
\vbox{\moveleft0.4cm\hbox{%
%
%
\lower3.5pt\hbox{$_{_{\!
\hbox{\sixbf k}%
\raise1pt\hbox{%
{\hbox{\fivebf[\vbox{\moveleft1.9pt\hbox{[}}}}%
}%
{\hbox{\sixrm{\char"00}}}\!\!_{_{^{2}}}\!\!      
\raise1pt\hbox{%
{\hbox{\fivebf]\vbox{\moveleft1.9pt\hbox{]}}}}%
}%
_{\!\hbox{\fivebf+}\!\!} }}$}
%
%
}}
$%
$ \!\!\! \!\!\!
( \hbox{\tenbf k}
{
\raise1pt\hbox{\eightbf[$\!{_{\!}}{{_{_{_{\!}}}}}$[}
\Gamma\!\!_{_{^{2}}}\!
\raise1pt\hbox{\eightbf]$\!{_{\!}}{_{_{_{\!}}}}$]}
} ).
$

\medskip
\item{\bf 2.}
%
Put:
$ \beta_{_{\!\hbox{\fivebf G}}}\!(X):= \hbox{\tenrm inf}\{j\
\!\vert\ \!\exists x;\ \! x\!\in\! X \land \ \!
\mdoubleH_{_{\!j}}\! (X,X\setminus_{_{\!o}}\!\! x;\hbox{\tenbf
G})\ne0\}. $
For a finite $\Delta$, $\beta_{_{\!\hbox{\fivebf k}}}\!(\Delta)$
is related to the concept ``depth of the ring
$ \hbox{\tenbf k}\hbox{\tenbf[$\!\!{_{_{\!}}}{_{_{\!}}}$[} \Delta
\hbox{\tenbf]$\!\!{_{_{\!}}}{_{_{\!}}}$]}" $
and ``C-M-ness of $\hbox{\tenbf k}
\hbox{\tenbf[$\!\!{_{_{\!}}}{_{_{\!}}}$[} \Delta
\hbox{\tenbf]$\!\!{_{_{\!}}}{_{_{\!}}}$]}" $ through
$\beta_{_{\!\hbox{\fivebf k}}}\!(\Delta)=\hbox{\tenrm depth} (
\hbox{\tenbf k}\hbox{\tenbf[$\!\!{_{_{\!}}}{_{_{\!}}}$[} \Delta
\hbox{\tenbf]$\!\!{_{_{\!}}}{_{_{\!}}}$]} )-\!1,$
in \cite{3} %
\hbox{Ex.\ 5.1.23 and
\cite{26} %
p.\ $\!$142 Ex.\ $\!$34. See also \cite{22}.} %
\endexample

\proclaim{Proposition 2}\ {\bf a.}
$[{\Delta}\ \hbox{\tenrm is ``CM}\!_{_{\hbox{\fivebf G}}}\!\!"]
\!\Longleftrightarrow\! [\mdoubleH_{_{i}}({\Delta},\!\hbox{\tenrm
cost}\!_{_{\Delta}}\!{{\delta\!_{_{ }}}}; \hbox{\tenbf G})=0\
\forall\ \!\delta\!\in\!\Delta \ \hbox{\tenrm and}\ \forall\
\!i\!\leq\!n\!-\!1].$

\indent\indent{\bf b.}
$[{\Delta}\ \hbox{\tenrm is 2-``CM}_{_{\!\hbox{\fivebf G}}}\!\!"]
\!\Longleftrightarrow\! [\mdoubleH_{_{i}}(\hbox{\tenrm
cost}\!_{_{\Delta}}\!{{\delta\!_{_{ }}}}; \hbox{\tenbf G})=0\
\forall\ \!\delta\!\in\overline{\hbox{\teni$\!$
V}}_{_{\!\!\!\Delta}} \ \hbox{\tenrm and}\ \forall\
\!i\!\leq\!n\!-\!1].$
\endproclaim

\demo{Proof}
Use the \hbox{\tenbf LHS} w.r.t.
$\!({\Delta}, \hbox{\tenrm cost}\!_{_{\Delta}}\!{{\delta\!_{_{
}}}})$, Prop.\ 1 p.\ 11 and the fact that $\hbox{\tenrm
cost}\!_{_{\Delta}}\!{{\emptyset\!_{_{o}}}}\!\!=\!\emptyset$
respectively
$\hbox{\tenrm cost}\!_{_{\Delta}}\!{{\delta}}\!\!=\!\Delta$ if
$\delta\not\in\Delta$ .\nobreak\qed
\enddemo

\goodbreak%

\noindent Put
${ \Delta\!^{^{_{(_{\!}p_{\!})}}}\!\!:= \{\delta\in\Delta\
\!\vert\ \!\#\delta\le p+1\}},$
${\Delta\!^{^{_{p\!}}}\!\!:=
\Delta\!\!^{^{_{(_{\!}p_{\!})}}}\!\!\setminus
\Delta\!\!^{^{_{(_{\!}p\!-\!1_{\!})}}} }\!\!\!,\
$
${\Delta\!^{^{_{\prime\!}}}\!:=
\Delta\!^{^{_{(\!n{_{_{^{_{\!}}}}}\!-\!1_{\!})}}}}\!\!\!$,
$ n\!:=\!\dim\Delta.$
So,
$ {
\Delta\!^{^{_{(\!n{_{_{^{_{}}}}}\!\!_{\!})}}}\!\!=\!\Delta
}.$

\proclaim{Theorem 8.\ a}
$\Delta$ is ``CM$_{_{\!^{\hbox{\fivebf G}}}}\!\!$" \underbar{iff}\
\ $\Delta\!^{^{_\prime}}$ is 2-$\!$``CM$_{_{^{\hbox{\fivebf
G}}}}\!\!\!$" and
$\mdoubleH\!\!\!\!_{_{^{n-1}}}\!\!(\Delta,\hbox{\tenrm
cost}\!_{_{^{\!\Delta}}}\!\!\delta ;\hbox{\tenbf G})=0$
$\forall\ \delta\!\in\! \Delta$.

\hskip1.5cm
{\bf b.}
$\Delta$ is 2-$\!$``CM$_{_{\!^{\hbox{\fivebf G}}}}\!\!$"
\underbar{iff}\ \ $\Delta\!^{^{_\prime}}$ is
2-$\!$``CM$_{_{^{\hbox{\fivebf G}}}}\!\!\!$" and
$\mdoubleH\!\!\!\!_{_{^{n-1}}}\!\!(\hbox{\tenrm
cost}\!_{_{^{\!\Delta}}}\!\!\delta ;\hbox{\tenbf G})=0$
$\forall\ \!\delta\!\in\overline{\hbox{\teni$\!$
V}}_{_{\!\!\!\Delta}}$.
\endproclaim

\demo{Proof}
Proposition\ 2 together with the fact that adding or deleting
$n$-simplices does not effect homology groups of degree $\leq
n-2$. See Proposition\ 3.{\tenbf e} p.\ 31.
\qed
\enddemo

Th.\ 8 is partially deducible from
{\rm \cite{13} %
p.\ 359-60}. %

\proclaim{\bf Lemma 1}
$
\Delta\ \hbox{\tenrm ``CM}_{_{\!\!\hbox{\fivebf
G}}}\!\!"\Longrightarrow
\cases
(\hbox{\tenbf a})\ \mdoubleH\!_{_{i}}(\!\hbox{\tenrm
cost}\!_{_{\Delta}}\!{{\delta\!_{_{ }}}}; \hbox{\tenbf G})=0\
\forall\ \delta\ \in\!\Delta & if\ i\leq n-2 \cr
\hbox{\tenrm and} \cr
(\hbox{\tenbf b})\ \mdoubleH\!_{_{i}}( \!\hbox{\tenrm
cost}\!\!\!\!\! _{_{\hbox{\sevenrm
cost}\!\!_{_{\Delta}}\!\!\!{{{\hbox{\sixi
{\char"0E}}}\!_{_{\hbox{\fiverm 2}}} } }}}\!
\!\!{{\delta\!_{_{1}}}}\!; \hbox{\tenbf G})=0\ \forall\
\delta\!_{_{1}},\delta\!_{_{2}}\in\!\Delta & if\ i\leq n-3. \cr
\endcases
$
\endproclaim

\medskip
\noindent{\bf Proof.} {\bf (a)} Use Proposition\ 1 p.\ 11, the
definition of ``CM"$\!$-ness and the {\bf LHS} w.r.t.\ $\!$
$({\Delta},\hbox{\tenrm cost}\!_{_{\Delta}}\!{{\delta\!_{_{
}}}})\!$,\ {which\ reads;}

\smallskip
\noindent
$\lower4.5pt\hbox{$^{
\dots\
\!{\buildrel \hbox{\eighti{\char"0C} $_{_{^{1\ast}}}$} \over
\longrightarrow }
{
\medoubleH\!_{_{n\!\!}} (\hbox{\eightrm{\char"01}}              
,\hbox{\eightrm cost}\!_{_{\!\Delta\!}}\!{{\delta\!_{_{
}}}};\hbox{\eightbf G})}
\ {\buildrel \hbox{\eighti{\char"0E}                            
$_{_{^{1\ast}}}$} \over \longrightarrow}\
{ \medoubleH\!\!\!\!_{_{n\!-\!1\!\!}} (\hbox{\eightrm
cost}\!_{_{\!\Delta\!}}\!{{\delta\!_{_{ }}}};\hbox{\eightbf G})}
\ \!{\buildrel \hbox{\eighti{\char"0B}                          
$_{_{^{1\ast}}}$} \over \longrightarrow } \ \!
{\medoubleH\!\!\!\!_{_{n\!-\!1\!\!}} (\hbox{\eightrm{\char"01}} 
;\hbox{\eightbf G})}
\ \!{\buildrel \hbox{\eighti{\char"0C}                          
$_{_{^{1\ast}}}$} \over \longrightarrow } \ \!
\medoubleH\!\!\!\!_{_{n\!-\!1\!\!}} (\hbox{\eightrm{\char"01}}  
,\hbox{\eightrm cost}\!_{_{\!\Delta\!}}\!{{\delta\!_{_{
}}}};\hbox{\eightbf G}) \ \!
\!{\buildrel \hbox{\eighti{\char"0E}                            
$_{_{^{1\ast}}}$} \over \longrightarrow }
{ \medoubleH\!\!\!\!_{_{n\!-\!2\!\!}} (\hbox{\eightrm
cost}\!_{_{\!\Delta\!}}\!{{\delta\!_{_{ }}}};\hbox{\eightbf G}) }
\
\!{\buildrel \hbox{\eighti{\char"0B}                            
$_{_{^{1\ast}}}$} \over \longrightarrow }\dots
}$}$%
\hfill$\triangleright$

\medskip
\noindent$\!${\bf (b)} Apply Prop.\ 3. p.\ 31 {\bf a}+{\bf b} to
the {\bf M-$\!$V$\!$s} w.r.t.
$\! (\!{\hbox{\tenrm cost}\!_{_{{\!}\Delta}}\!\!{{\delta\!_{_{1}}
} }}\!, {\hbox{\tenrm cost}\!_{_{{\!}\Delta}}\!\!{{\delta\!_{_{2}}
} }}\!),$ {\rm then\ use\ {\bf a}, i.e;}

\smallskip
\noindent
$\lower4.5pt\hbox{$^{
\dots\ =\ \dots { \buildrel \beta_{_1\ast} \over \longrightarrow
}\ { \medoubleH_{_{\!n}}\! (\hbox{\eightrm
cost}\!_{_{\Delta}}\!(\!{{\delta\!_{_{1}}}}\!\cup
{{\delta\!_{_{2}}}}\!);\hbox{\eightbf G}) } \ {\buildrel
\delta_{_1\ast} \over \longrightarrow }\
\medoubleH\!\!\!\!_{_{n\!-\!1}}\! (\!\hbox{\eightrm
cost}\!\!\!\!\! _{_{\hbox{\sixrm
cost}\!\!_{_{\Delta}}\!\!\!{{\delta\!_{_{2}} } }}}\!
\!\!{{\delta\!_{_{1}}}}\!;\hbox{\eightbf G}) \ {\buildrel
\alpha_{_1\ast} \over \longrightarrow }\
\medoubleH\!\!\!\!_{_{n\!-\!1}}\! (\hbox{\eightrm
cost}\!_{_{\Delta}}\!{{\delta\!_{_{1}}}};\hbox{\eightbf G}) \
\oplus\ \medoubleH\!\!\!\!_{_{n\!-\!1}}\! (\hbox{\eightrm
cost}\!_{_{\Delta}}\!{{\delta\!_{_{2}}}};\hbox{\eightbf G}) \
{\buildrel \beta_{_1\ast} \over \longrightarrow}}$}$

\medskip
$\lower4.5pt\hbox{$^{
{\buildrel \beta_{_1\ast} \over \longrightarrow }\ {
\medoubleH\!\!\!\!_{_{n\!-\!1}}\! (\hbox{\eightrm
cost}\!_{_{\Delta}}\!(\!{{\delta\!_{_{1}}}}\!\cup{{\delta\!_{_{2}}}}\!);
\hbox{\eightbf G}) } \ {\buildrel \delta_{_1\ast} \over
\longrightarrow }\ \medoubleH\!\!\!\!_{_{n\!-\!2}}\!
(\!\hbox{\eightrm cost}\!\!\!\!\! _{_{\hbox{\sixrm
cost}\!\!_{_{\Delta}}\!\!\!{{\delta\!_{_{2}} } }}}\!
\!\!{{\delta\!_{_{1}}}}\!;\hbox{\eightbf G}) \ {\buildrel
\alpha_{_1\ast} \over \longrightarrow }\
\medoubleH\!\!\!\!_{_{n\!-\!2}}\! (\hbox{\eightrm
cost}\!_{_{\Delta}}\!{{\delta\!_{_{1}}}};\hbox{\eightbf G}) \
\oplus\ \medoubleH\!\!\!\!_{_{n\!-\!2}}\! (\hbox{\eightrm
cost}\!_{_{\Delta}}\!{{\delta\!_{_{2}}}};\hbox{\eightbf G}) \
{\buildrel \beta_{_1\ast} \over \longrightarrow} \ \dots }$}$
\hfill$\square$

\medskip
\noindent\hbox{\tenbf Observation.}
To turn the implication in lemma\ 1 into an equivalence we just
have to add
$\mdoubleH\!_{_{i}}\! (\Delta;\!\hbox{\tenbf G})\!=\!0\
\hbox{\tenrm for}\ i\leq n\!-\!1 $, giving us the equivalence in;

\medskip
\noindent $\Delta\ \!\hbox{\tenrm ``CM}\!_{_{\hbox{\fivebf
G}}}\!\!\!" \!\!\Longleftrightarrow\!\!
\cases
\!\!\hbox{\tenbf i)}\ \ \ \mdoubleH_{_{i}}\! (\Delta;\hbox{\tenbf
G})\!=\!0 & \!\!\!\!\!i\!\leq\!n\!-\!1 \cr
\!\!\hbox{\tenbf ii)}\ \ \mdoubleH_{_{i}}(\!\hbox{\tenrm
cost}\!_{_{\Delta}}\!{{\delta\!_{_{ }}}}; \hbox{\tenbf G})=0\
\forall\ \!\delta\!\!\in\!\Delta & \!\!\!\!\!i\!\leq\!n\!-\!2 \cr
\!\!\hbox{\tenbf iii)}\ \mdoubleH_{_{i}}( \!\hbox{\tenrm
cost}\!\!\!\!\!\!\!\! _{\lower0.8pt\hbox{$ {_{\hbox{\sevenrm
cost}\!\!_{_{\Delta}}\!\!{{\delta\!_{_{2}} } }}} $}}\!
\!{{\delta\!_{_{1}}}}\!;\hbox{\tenbf G})\! =\!0\ \!\forall\
\!\delta\!_{_{1}}\!,\delta\!_{_{2}}\!\!\in\!\Delta &
\!\!\!\!\!i\!\leq\!n\!-\!3 \cr
\endcases
$
\hfill$(\hbox{\tenbf 1})$

\medskip
$\mdoubleH\!\!\!\!\!_{_{^{n-1}}}\!\!(\hbox{\tenrm
cost}\!_{_{^{\!\Delta}}}\!\!\delta ;\hbox{\tenbf G})\!=0\ \forall\
\delta\!\in\! \Delta$
allows one more step in the proof of Lemma\ 1b, i.e.;

\medskip
\noindent $\Delta\
\!\hbox{\tenrm ``CM}\!_{_{\hbox{\fivebf G}}}\!\!\!"\!$
\&
$\mdoubleH\!\!\!\!\!_{_{^{n-1}}}\!\!(\hbox{\tenrm
cost}\!_{_{^{\!\Delta}}}\!\!\delta ;\hbox{\tenbf G})\!=0\ \forall\
\!\delta_{\!\!}\in_{\!\!}\Delta
\!\!\Longleftrightarrow\!\!
\cases
\!\!\hbox{\tenbf i)}\ \ \ \!\mdoubleH_{_{i}}\!
(\Delta;\hbox{\tenbf G})\!=\!0 & \!\!\!\!\!i\!\leq\!n\!-\!1 \cr
\!\!\!\hbox{\tenbf ii)}\ \ \mdoubleH_{_{i}}\!(\!\hbox{\tenrm
cost}\!_{_{\Delta}}\!{{\delta\!_{_{ }}}}; \hbox{\tenbf G})=0\
\forall\ \!\delta\!\!\in\!\Delta & \!\!\!\!\!i\!\leq\!n\!-\!1 \cr
\!\!\!\hbox{\tenbf iii)}\ \mdoubleH_{_{i}}\!( \!\hbox{\tenrm
cost}\!\!\!\!\!\!\!\!\! _{\lower0.8pt\hbox{$ {_{\hbox{\sevenrm
cost}\!\!_{_{\Delta}}\!\!{{\delta\!_{_{2}} } }}} $}}\!
\!{{\delta\!_{_{1}}}}\!;\hbox{\tenbf G})\! =\!0\ \!\forall\
\!\delta\!_{_{1}}\!,\delta\!_{_{2}}\!\!\in\!\!\Delta &
\!\!\!\!\!i\!\leq\!n\!-\!2\cr
\endcases
$
\hfill$(\hbox{\tenbf 2})$

\smallskip
Note that {\bf a:} i) in Eq.\ 2 is a consequence of ii) and iii)
by the {\bf M-$\!$V$\!$s} above and that {\bf b:} the l.h.s. is by
definition equivalent to $\Delta$ being 2-``CM".

\smallskip
Since, [${\Delta}\ \!\!\hbox{\tenrm ``CM}\!_{_{\hbox{\fivebf
G}}}\!\!"] \!\Longleftrightarrow\!
[\mdoubleH_{_{i}}({\Delta},\!\hbox{\tenrm
cost}\!_{_{\Delta}}\!{{\delta\!_{_{ }}}}; \hbox{\tenbf G})=0\
\forall\ \!\delta\!\in\!\Delta \ \hbox{\tenrm and}\ \forall\
\!i\!\leq\!n\!-\!1] $, item \hbox{\tenbf iii} in Eq.\ 1 is, by the
\hbox{\tenbf LHS} w.r.t.
$\!({\Delta},\!\hbox{\tenrm cost}\!_{_{\Delta}}\!{{\delta\!_{_{
}}}})$,
totally superfluous as far as the equivalence is concerned but
never the less it becomes quite useful when substituting
$\hbox{\tenrm cost}\!_{_{\Delta}}\!\!{{\delta\!_{_{}}}}$
for every occurrence of $\Delta$ and using that $ \hbox{\tenrm
cost}\!\!\!\!\!\!\! _{\lower0.8pt\hbox{$ {_{\hbox{\sixrm
cost}\!\!_{_{\Delta}}\!\!{{\delta\!_{_{\ }}}} }} $}}
\!\!{{\delta\!_{_{1}}}} = \hbox{\tenrm cost}\!
_{_{{\Delta}}}\!\!{{\delta\!_{_{\ }}}} $ for $ {{\delta\!_{_{1}}
}}\!\!\notin \!\hbox{\tenrm cost}\!
_{_{{\Delta}}}\!\!{{\delta\!_{_{\ }}}}\!\! $, we get;

\medskip
$\hbox{\tenrm cost}\!_{_{\Delta}}\!{{\delta\!_{_{}}}} \
\!\hbox{\tenrm ``CM}\!_{_{\hbox{\fivebf G}}}\!\!\!"\! \ \ \forall\
\!\delta\!\!\in\!\Delta \Longleftrightarrow
\cases
\hbox{\tenbf (i)}\ \ \ \mdoubleH_{_{i}}\! (\hbox{\tenrm
cost}\!_{_{\Delta}}\!{{\delta\!_{_{}}}} ;\hbox{\tenbf G})\!=\!0 &
\!\!\!i\!\leq\!n\!_{_{\delta}}\!-\!1\ \ \!\forall\
\!\delta\!\!\in\!\Delta  \cr
\hbox{\tenbf (ii)}\ \ \mdoubleH_{_{i}}(\!\hbox{\tenrm
cost}\!\!\!\!\!\!\! _{\lower0.8pt\hbox{$ {_{\hbox{\sixrm
cost}\!\!_{_{\Delta}}\!\!{{\delta\!_{_{\ }}}} }} $}}
\!\!{{\delta\!_{_{1}}}}; \hbox{\tenbf G})=0\ \ &
\!\!\!i\leq\!n\!_{_{\delta}}\!-\!2\ \forall\
\delta,\delta\!_{_{1}}\!\!\in\! {\Delta}.  \cr
\endcases
$\hfill(\hbox{\tenbf 3})

\remark{Remark 1}
The \hbox{\tenbf LHS} w.r.t. $\!({\Delta},\hbox{\tenrm
cost}\!_{_{\Delta}}\!{{\delta_{_{ }}}}),$ Corollary\ iii p.\ 12
and Note\ 1 p.\ 25 gives;\hfill\break
\noindent
If ${{\Delta\!}}$ is $2$-$\!$``CM$\!_{_{\hbox{\fivebf G}}}\!\!"$
then Note\ 1 p.\ 13 plus Prop.\ 1 p.\ 11 implies that $
\mdoubleH_{_{{\!\hbox{\seveni n}}_{^{_{\!\!\Delta}}}}}\!
\!\!(\!{\Delta};\!\hbox{\tenbf G})\!\not=\!0_{\!}$,\nobreak\ i.e.
Hip$(|\Delta|)\!,_{\!}$ def. ${\!}$p.\ $_{\!}$12, is empty if it's
a subcomplex, which it is for manifolds, cf. p.\nobreak\ 25.
\endremark

\goodbreak

\remark{Remark 2}
$\!$Set $\bullet\ \!\bullet\!\!\mthickline\!\bullet:=
\{{\emptyset}\!_{_{^{o}}}\!,\{\hbox{\tenrm
v}\!_{_{^{1}}}\!\},\{\hbox{\tenrm v}\!_{_{^{2}}}\!\},
\{\hbox{\tenrm v}\!_{_{^{3}}}\!\},\{\hbox{\tenrm
v}\!_{_{^{2}}},\hbox{\tenrm v}\!_{_{^{3}}}\!\} \},$
$\!n\!_{_{\varphi}}\!\!\!:=\! \dim\hbox{\tenrm
cost}_{_{\!\Delta}}\!\varphi \ \hbox{\tenrm and}\nobreak \
n\!_{_{\Delta}}\!\!\!:=\dim{{\Delta}}.$

Note that it's always true that;
$n\!_{_{\Delta}}\!\!-\!1\ \!\le
n_{_{\!\raise1pt\hbox{${_{{{\tau}}}}$}} } \!\le
n\!_{_{\raise1pt\hbox{${_{{{\delta}}}}$}} } \!\le n\!_{_{\Delta}}
$
if $\tau\!\subset\!\delta$.

\smallskip
Now,
$ {
[n\!_{_{\hbox{\fiverm v}}}\!\! =\!n\!_{_{\Delta}}
\hbox{\tenrm and}\ \hbox{\tenrm
cost}_{_{\!\!\Delta}}\!\!\hbox{\tenrm v} \ \hbox{\tenrm pure}\
\forall\ \!\hbox{\tenrm v}\!_{_{}}\! \in\!\hbox{\tenrm
V}\!\!_{_{\Delta}}] \!\Longleftrightarrow\! [n\!_{_{{\delta}}}\!\!
=\!n\!_{_{\Delta}}
\hbox{\tenrm and}\
\hbox{\tenrm cost}\!_{_{\Delta}}\!{{\delta\!_{_{ }}}}\
\hbox{\tenrm pure}\ \forall\ \!\emptyset\!_{_{^{o}}} \!\!\ne\!
\delta\!_{_{}}\!\in\!\Delta] \!\Longrightarrow
}$\hfill\break
${ \!\! \indent \Longrightarrow [\ ^{_{\!}}\!\! \hbox{\tenrm
cost}_{_{\!\!\Delta}}\!\!\delta \
(=\!\cup\!\!\!\!\!\!{_{_{_{\hbox{\fiverm v}\in\delta}}}}
\!\hbox{\tenrm cost}\!_{_{\Delta}}\!\!{\hbox{\tenrm v}\!_{_{^{\
}}}}) \
\hbox{\tenrm pure}\ \forall\ \!\emptyset\!_{_{^{o}}} \! \ne
\delta\!_{_{}}\in\Delta \not=\bullet\
\!\bullet\!\!\mthickline\!\bullet ] \Longrightarrow [{\Delta} \
\hbox{\tenrm is\ pure}] \Longrightarrow
}$ \hfill\break
${
\Longrightarrow
\hbox{\tenbf
\big[$\!\!{_{_{\!}}}{_{_{\!}}}$\big[$\!\!{_{_{\!}}}{_{_{\!}}}$\big[%
}
[\hbox{\tenrm v}\!\in\!\hbox{\tenrm V}\!_{_{\!\Delta}}
}$
is a cone point (p.\ 22)]
${ \Longleftrightarrow
[n\!_{_{\hbox{\fiverm v}}}\!\! =\!n\!_{_{\Delta}}\!\!-\!1]\ \!
\Longleftrightarrow [\hbox{\tenrm
v}\!\in\!{{\delta\!_{_{m}}}}\!\!\!\in\!\Delta\ \hbox{\tenrm if}\
\dim{{{{\delta\!_{_{m}}}}}}\!\!=\!n\!_{_{\Delta}}]
\hbox{\tenbf\big]$\!\!{_{_{\!}}}{_{_{\!}}}$\big]$\!\!{_{_{\!}}}{_{_{\!}}}$\big]}.
} \! $

Since, $
{
n\!_{_{\varphi}}:=\dim\hbox{\tenrm cost}_{_{\!\Delta}}\!\varphi\!
=n\!_{_{\Delta}}\!\!\!-\!1 \!\Longleftrightarrow
\emptyset\!_{_{^{o}}} \!\ne \varphi \!\subset\!
{{\delta\!_{_{m}}}}\!\!\!\in\!\Delta\  \forall\
\!\delta\!_{_{m}}\!\!\in\Delta}$ \hbox{\tenrm with}\
${\dim{{{{\delta\!_{_{m}}}}}}\!\! = n\!_{_{\Delta}}{_{\!}},}$ we
\hfill\break
conclude that;
If $\Delta$ is pure then $\varphi$ consists of nothing but cone
points, cf. p.\ 22.

So,
$ [\Delta\ \hbox{\tenrm pure\ and}\ n\!_{_{{\delta}}}\!
\!=\!\dim{{\!\Delta}}\ \forall\ \delta\in\Delta]
\Longleftrightarrow [\Delta\ \hbox{\tenrm pure}\ \hbox{\tenrm and\
has\ no\ cone\ points}].$
\break
\indent{\bf(}$[n\!_{_{\hbox{\fiverm v}}}\!
=n\!_{_{\Delta}}\!\!-\!1\ \! \forall\ \hbox{\tenrm
v}\!\in\hbox{\teni V}_{_{\!\!\!\Delta}}] \Longleftrightarrow
\hbox{\teni V}_{_{\!\!\!\Delta}}\ \hbox{\tenrm is\ finite\ and}\
\Delta=\overline{\!\hbox{\teni V}}_{_{\!\!\!\Delta}} \
(:=\!\hbox{\tenrm the\ full\ complex\ w.r.t.}\ \hbox{\teni
V}_{_{\!\!\!\Delta}}).${\bf)}
\endremark

\proclaim{Theorem 9}
The following two conditions are equivalent to
``$\Delta$ is $2$-$\!``\hbox{\tenrm CM}\!_{_{\hbox{\fivebf
G}}}\!\!"";$ \hfill\break
\noindent{\bf a.}\indent {$
{ \hbox{\tenbf i.}\ \hbox{\tenrm
cost}\!_{_{\Delta}}\!{{\delta\!_{_{}}}}\ is
\ ``\hbox{\tenrm CM}\!_{_{\hbox{\fivebf G}}}\!\!"\!,\ \!\forall\
\delta\!_{_{}}\in\!\Delta
\ and\
\hbox{\tenbf ii.}\ \ \!n\!_{_{\delta}} \!:=\! \dim\hbox{\tenrm
cost}\!_{_{\Delta}}\!{{\delta\!_{_{}}}}=\dim\Delta
=:n\!_{_{\Delta}}\ \!\forall\ \emptyset\!_{_{^{o}}} \!\! \ne\!
\delta\!_{_{}}\in\!\Delta.
} $}

\smallskip
\noindent{\bf b.}\indent
{$
{\hbox{\tenbf i.}\ \hbox{\tenrm
cost}\!_{_{\Delta}}\!{{\hbox{\tenrm v}\!_{_{}}}}\ \hbox{\tensl
is}\ ``\hbox{\tenrm CM}\!_{_{\hbox{\fivebf G}}}\!\!"\!,\ \
\forall\ \hbox{\tenrm v}\!_{_{}}\in\!\hbox{\tenrm V}\!_{_{\Delta}}
\indent and \indent
\hbox{\tenbf ii.} \ \{\bullet\
\!\bullet\!\!\mthickline\!\bullet\}\ne{\Delta}\ \hbox{\tensl has\
no\ cone\ points}.
}$}
\endproclaim

\demo{Proof}
If there is no dimension collapse in Eq.\ 3 it is equivalent to
Eq.\ 2. \qed
\enddemo

Our next corollary was, for $\hbox{\tenbf G}=\hbox{\tenbf k},$
originally ring theoretically proven by T. Hibi. We'll essentially
keep Hibi's formulation, though using that
$ {\Delta\!^{^{\!_{o}}}\!} \!:=\!
\Delta{\raise1.5pt\hbox{\eightmsbm \char"72}}           
\{\tau\!\in\!\Delta\ \!\vert\ \!\tau\!\supset\!\delta_{_{\!i}} \
\hbox{\tenrm for\ some}\ i\!\in\!\hbox{\tenbf I}\}
=
\bigcap\!\!\!\!\!\!\!_{_{_{_{_{i\in\
_{\!\!\!}\lower0.8pt\hbox{\sevenbf I}}}}}} \hbox{\tenrm
cost}\!_{_{\Delta}}\!\!{{\delta}\!_{_{^{i}}}}\!. $

\proclaim{Corollary}
{\rm(\cite{12} %
Corollary\ p.\ 95-6)}
Let $\Delta$ be a pure simplicial complex of dimension $n$ and
$\{{\delta}\!_{_{^{i}}}\!\}\!_{_{^{i\in\lower0.5pt\hbox{\eightbf
I}}}}$, a finite set of faces in $\Delta$ satisfying
$\delta_{_{\!i}}\raise0.5pt\hbox{\sevensy {\char"5B}}          
\ \!\delta_{_{\!j}}\!\! \notin\!\Delta\ \forall\ \!i\!\ne\!j$.
Set,
${\Delta\!^{^{\!_{o}}}\!}:=
\!\bigcap\!\!\!\!\!\!\!_{_{_{_{_{i\in\lower0.2pt\hbox{\fivebf
I}}}}}} \hbox{\tenrm
cost}\!_{_{\Delta}}\!\!{{\delta}\!_{_{^{i}}}}$.

\noindent
$\hbox{\tenbf a)}$  If $\Delta$ is $\hbox{\tenrm
``CM}\!_{_{\hbox{\fivebf G}}}\!\!"\!$ and
$\dim\Delta\!^{^{_{\!o}}}\!<\!n$, then
$\dim\Delta\!^{^{_{\!o}}}\!=\!n\!-\!1$ and $\Delta\!^{^{_{\!o}}}$
is $\hbox{\tenrm ``CM}\!_{_{\hbox{\fivebf G}}}\!\!"\!$.

\medskip
\noindent
$\hbox{\tenbf b)}$ If ${\bar{\hbox{\tenrm
st}}}\!_{_{\Delta}}\!{{\delta\!_{_{i}}}} \ \hbox{\tenrm is}\
\hbox{\tenrm ``CM}\!_{_{\hbox{\fivebf G}}}\!\!"\ \forall\
{i\!\in\!\lower0.0pt\hbox{\tenbf I}}$ and $\Delta\!^{^{_{\!o}}}$
is $\hbox{\tenrm ``CM}\!_{_{\hbox{\fivebf G}}}\!\!"\!$ of
dimension $n$, then $\Delta$ is also
$\hbox{\tenrm ``CM}\!_{_{\hbox{\fivebf G}}}\!\!"\!$.
\endproclaim

\demo{Proof}
\noindent\hbox{\tenbf a.} \noindent $
[\delta_{_{\!i}}\!\raise0.5pt\hbox{\sevensy {\char"5B}}   
\delta_{_{\!j}}\!\!\notin\!\Delta\ \! \forall\
\!i\!\ne\!j\!\in\!\lower0.2pt\hbox{\tenbf I}]
\!\Longleftrightarrow\!\! [\hbox{\tenrm
cost}\!_{_{\Delta}}\!\!{{\delta}\!_{_{^{i}}}}
\raise0.0pt\hbox{\ninesy {\char"5B}}                      
\big( \bigcap\!\!\!\!\!\!\!_{_{_{_{_{j\not={i}}}}}} \hbox{\tenrm
cost}\!_{_{\Delta}}\!\!{{\delta}\!_{_{^{j}}}}\!\big)\!=\!\Delta]
\Longrightarrow\!
\big[\!\!{_{_{\!}}}{_{_{\!}}}\big[\!\!{_{_{\!}}}{_{_{\!}}}\big[
{\hbox{\tenrm cost}\!\!_{_{\Delta^{^{{\!{_{\!\prime}}}}}}}
\!\!{{\delta}\!_{_{^{i}}}}}  \raise0.0pt\hbox{\ninesy {\char"5B}}
\big( \bigcap\!\!\!\!\!\!\!_{_{_{_{_{j<{i}}}}}} \hbox{\tenrm
cost}\!\!_{_{\Delta^{^{{\!{_{\!\prime}}}}}}}
\!\!{{\delta}\!_{_{^{j}}}}\! \big)
\!=\! \big(\hbox{\tenrm
cost}\!_{_{\Delta}}\!\!{{\delta}\!_{_{^{i}}}}
\!\raise0.0pt\hbox{\ninesy {\char"5B}} \big(
\bigcap\!\!\!\!\!\!\!_{_{_{_{_{j<{i}}}}}} \!\hbox{\tenrm
cost}\!_{_{\Delta}}\!\!{{\delta}\!_{_{^{j}}}}\big) \big)
\raise0.0pt\hbox{\ninesy {\char"5C}}
\Delta\!^{^{_{\!\hbox{$_{\prime}$}}}}
\!=\!\Delta\!^{^{_{\!\hbox{$_{\prime}$}}}}
\big]\!\!{_{_{\!}}}{_{_{\!}}}\big]\!\!{_{_{\!}}}{_{_{\!}}}\big]
\Longrightarrow $ $ \dim{\Delta\!^{^{\!_{o}}}}\!\!=\!n\!-\!1
\Longrightarrow {\Delta\!^{^{\!_{o}}}} \!\!=\!
{\Delta\!^{^{\!_{o}}}}\! \raise0.0pt\hbox{\ninesy {\char"5C}}
\Delta\!^{^{_{\!\hbox{$_{\prime}$}}}} \!=\!
\big(\bigcap\!\!\!\!\!\!\!_{_{_{_{_{i\in\lower0.8pt\hbox{\sevenbf
I}}}}}} \hbox{\tenrm
cost}\!_{_{\Delta}}\!\!{{\delta}\!_{_{^{i}}}}\!\big)
\raise0.0pt\hbox{\ninesy {\char"5C}}
\Delta\!^{^{_{\!\hbox{$_{\prime}$}}}} \!\!=\break
=\! \bigcap\!\!\!\!\!\!\!_{_{_{_{_{i\in\lower0.8pt\hbox{\sevenbf
I}}}}}} \!{\hbox{\tenrm
cost}\!\!_{_{\Delta^{^{{\!{_{\!\prime}}}}}}}
\!\!{{\delta}\!_{_{^{i}}}}}\!. $
\noindent By Th.\ 8 we know that
$\Delta\!^{^{_{\!\hbox{$_{\prime}$}}}}$ is 2-$\!``\hbox{\tenrm
CM}\!_{_{\hbox{\fivebf G}}}\!\!"\!$, implying that $ {\hbox{\tenrm
cost}\!\!_{_{\Delta^{^{{\!{_{\!\prime}}}}}\!}}
\!{{\delta}}}\!_{_{^{i}}}\ \hbox{\tenrm is}\ ``\hbox{\tenrm
CM}\!_{_{\hbox{\fivebf G}}}\!\!"\!\ \forall\
\!{i\!\in\lower0.2pt\hbox{\tenbf I}}. $
Induction using the
\hbox{\tenbf M-$\!$Vs} w.r.t. $ (\hbox{\tenrm
cost}\!_{_{\!\Delta^{^{{\!{_{\!\prime}}}}}}}
\!\!{{\delta}\!_{_{^{i}}}},
\bigcap\!\!\!\!\!\!\!_{_{_{_{_{j<{i}}}}}} \hbox{\tenrm
cost}\!\!_{_{\Delta^{^{{\!{_{\!\prime}}}}}}}
\!\!{{\delta}\!_{_{^{j}}}}\!) $ gives $
\mdoubleH_{_{\!i}}\!(\Delta\!^{^{_{\!o}}};\hbox{\tenbf G})\!=\!0\
\forall\ \! i\!<\!n\!_{_{\Delta}}\!\!\!-\!\!1. $ For links, use
{\rm Prop.}$\ \!2\ \!\hbox{\tenbf a}\ \!\!+\ \!\!\hbox{\tenbf b}$,
p.\ 30.
E.g. $ {\hbox{\tenrm Lk}\!\!\raise0.5pt\hbox{$_{_{{
\Delta\!\!^{^{_{\circ}}}}}}$}\!\!\delta} \!=\! \!{\hbox{\tenrm
Lk}\!\!\!\!\! \raise-0.5pt\hbox{$_{_{{
\Delta\!\!^{^{_{\circ}}}\!\!\cap
{\Delta^{^{{\!{_{\prime}}}}}}}}}$}\!\!\!\!\!\!\delta} =
{\hbox{\tenrm Lk}\!\!\!\!\!\! \lower0.5pt\hbox{$ {_{_{
\!\!\!\!\!\!\!\!(\ \cap\!\!\!_{\!}\cap\!\!\!\!\cap
\!\!\!\!\!\!\!_{_{_{{{i\in\lower0.2pt\hbox{\fivebf I}}}}}}
\!\hbox{\sevenrm
cost}\!_{_{\!\Delta}}\!\!{{\!\delta}\!_{_{^{i}}}}\!\!) }}} $}
\raise0pt\hbox{$_{_{{ \!\cap
{\Delta^{^{_{\!{\prime}}}}}}}}$}\!\!\!\!\!\!\!\!\!\!\!\!\delta} \
\ =\! $ \hbox{\tenrm [Prop.}$\ \!$2$\ \!\!$\hbox{\tenbf a} p.\
30]= $ \bigcap\!\!\!\!\!\!\!_{_{_{_{_{i\in{\lower1pt\hbox{\sevenbf
I}}}}}}} {\hbox{\tenrm Lk}\!\!\!\!\!\raise0.75pt\hbox{$_{_{_{
\hbox{\sevenrm
cost}\!\!_{_{\!\Delta^{^{{\!\!{_{\prime}}}}}}}\!{{\!\delta}\!_{_{^{i}}}}
}}}$}\!\!\!\!\!\!\delta}\ \! $ where $ {\hbox{\tenrm
cost}\!\!\!_{_{\Delta^{^{{\!{_{\prime}}}}}\!}}
\!{{\delta}}}\!_{_{^{i}}}\ \hbox{\tenrm and\ so}\ {\hbox{\tenrm
Lk}\!\!\!\!\!\raise0.75pt\hbox{$_{_{_{ \hbox{\sevenrm
cost}\!\!_{_{\!\Delta^{^{{\!\!{_{\prime}}}}}}}\!{{\!\delta}\!_{_{^{i}}}}
}}}$}\!\!\!\!\!\!\delta}\ \ \hbox{\tenrm is}\ ``\hbox{\tenrm
CM}\!_{_{\hbox{\fivebf G}}}\!\!"\!,\ \forall\
\!{i\!\in\lower0.2pt\hbox{\tenbf I}}. \hfill\triangleright$

\noindent\hbox{\tenbf b.} $ \Delta\!=\!{\Delta\!^{^{\!_{o}}}}\!
\bigcup\!\!\!\!\!\!\!_{_{_{_{_{i\in\lower0.2pt\hbox{\fivebf
I}}}}}} \overline{\hbox{\tenrm
{st}}}\!_{_{\Delta}}\!\delta\!_{_{^{j}}} $. $\
\!\overline{\hbox{\tenrm {st}}}\!_{_{\Delta}}\!\delta \cap
{\Delta\!^{^{\!_{o}}}} = \Big[{{_{
\hbox{\sevenrm  Eq. {\sevenbf II}+{\sevenbf III}\ p\ 30-31
}}} \atop {\hbox{\sevenrm plus that}\
{\delta_{_{\!i}}\!\cup\delta_{_{\!j}}\!\notin\Delta}}}\Big] =
\overline{\hbox{\tenrm {st}}}\!_{_{\Delta}}\!\delta \cap
{\hbox{\tenrm cost}\!_{_{\Delta}}\!\!{{\delta}}} =
\dot{\delta}\ast {\hbox{\tenrm Lk}\!_{_{\Delta}}\!\!{{\delta}}}. $
{So, by Cor.\ $\!$i p.\ 12,} $ \ \!\overline{\hbox{\tenrm
{st}}}\!_{_{\Delta}}\!\delta \cap {\Delta\!^{^{\!_{o}}}}\
\hbox{\tenrm is} $
\noindent $\hbox{\tenrm ``CM}\!_{_{\hbox{\fivebf G}}}\!\!"$
\underbar{iff} $\ \!\overline{\hbox{\tenrm
{st}}}\!_{_{\Delta}}\!\delta$ is.
Induction, using the \hbox{\tenbf M-$\!$Vs}$\ \!$ w.r.t.
$(\overline{\hbox{\tenrm {st}}}\!_{_{\Delta}}\!\delta\!_{_{^{i}}},
{\Delta\!^{^{\!_{o}}}}\! \bigcup\!\!\!\!\!\!\!_{_{_{_{_{j<{i}}}}}}
\overline{\hbox{\tenrm {st}}}\!_{_{\Delta}}\!\delta\!_{_{^{j}}}),$
gives $ \mdoubleH_{_{\!i}}\!(\Delta\!^{^{_{\ }}}\!;\hbox{\tenbf
G})\!=\!0\ \ \! \forall\ \!i\!<\!n\!_{_{\Delta}}.$
$\hbox{\tenrm End\ as\ in\ \hbox{\tenbf a}.}\hfill \qed $
\enddemo

\definition{Definition}
$ \Delta\setminus[\{\hbox{\tenrm v}\!_{_{^{1}}},\dots,\hbox{\tenrm
v}\!_{_{^{p}}}\}]:= \{ \delta\in\Delta\vert\ \!
\delta\cap\{\hbox{\tenrm v}\!_{_{^{1}}},\dots,\hbox{\tenrm
v}\!_{_{^{p}}}\}=\emptyset \}. $
$(\Delta\!\!\setminus\![\{\hbox{\tenrm v}\}] \!=\! \hbox{\tenrm
cost}\!_{_{\Delta}}\!\!{\hbox{\tenrm v}} $.)
\enddefinition

Permutations and partitions within $\{\hbox{\tenrm
v}\!_{_{^{1}}},\dots,\hbox{\tenrm v}\!_{_{^{p}}}\}$ doesn't effect
the result, i.e;

\proclaim{Lemma 2}
$ \Delta\!\!\setminus\![\{\hbox{\tenrm
v}\!_{_{^{1}}}\!,\dots,\!\hbox{\tenrm v}\!\!_{_{{p+q}}}\}] \!=\!
(\Delta\!\!\setminus\! [\{\hbox{\tenrm
v}^{{_{\prime}}}\!\!\!_{_{^{1}}}\!,\dots,\!\hbox{\tenrm
v}^{{_{\prime}}}\!\!\!_{_{^{p}}}\}]) \!\!\setminus\!
[\{\hbox{\tenrm
v}^{{_{\prime\prime}}}\!\!\!\!_{_{^{1}}}\!,\dots,\! \hbox{\tenrm
v}^{{_{\prime\prime}}}\!\!\!\!_{_{^{q}}}\}]$\ \ {and}

\smallskip
$\Delta\!\setminus[\{\hbox{\tenrm
v}\!_{_{^{1}}}\!,\dots,\!\hbox{\tenrm v}\!_{_{^{p}}}\}]
$
$
=
 \bigcap\!\!\!\!\!\!\!\!_{_{_{_{i=1,p}}}}
\!\!\!\hbox{\tenrm cost}\!_{_{\Delta}}\!\!{{\hbox{\tenrm
v}}\!_{_{^{i}}}}\!
=
\hbox{\tenrm cost} \!\!\!\!\!\!\!\!_{_{\lower1pt\hbox{\sixrm
cost}\!_{_{ {^{^{\hbox{\tenbf .}}}\!_{\!}{\hbox{\tenbf
.}}_{_{\!\hbox{\tenbf .}}}} }}\!\!\!\!\!{{\hbox{\sevenrm
v}}\!_{_{^{2}}}}}}\!\! \!{{\hbox{\tenrm v}}\!_{_{^{1}}}} $
\hfill$\square$

\nointerlineskip $
\hskip5.0cm
{{_{_{{\hbox{\fiverm cost}\!\!_{_{\Delta}}\!\!\!{{\hbox{\sixrm
v}\!_{_{p}}}}}}}}} \indent\ $
\endproclaim

\definition{Alternative Definition}
For $k\in\hbox{\tenbf N}$,
$\Delta$ is $k$-``\hbox{\tenrm CM}$\!_{_{\hbox{\fivebf
G}}}\!\!"\!$ if for every subset $T\subset\hbox{\tenrm
V}_{_{\!\!\!\Delta}}$ such that $\#T=k-1$, we have:
\indent{\tenbf i.} $\Delta\setminus[T]$ is ``\hbox{\tenrm
CM}$\!_{_{\hbox{\fivebf G}}}\!\!"\!,$\
\indent \hbox{\tenbf ii.}\ \
$\dim\Delta\setminus[T]=\dim\Delta=:\!n\!_{_{\Delta}}\!\!=:\!n. $
\enddefinition

Changing ``$\#T\!=k-1$" to ``$\#T\!<k$" doesn't 
alter the extension of the definition.
(Iterate in Th.\ 9{\bf b} mutatis mutandis.)
So, for $\hbox{\tenbf G}\!=\!\hbox{\tenbf k}$, a field, it's
equivalent to Kenneth Baclawski's original definition in Europ. J.
Combinatorics {\bf 3} (1982)\nobreak\ p.\nobreak\ 295.

%


\subhead
{{\rm II:3} \ \ \ {Segre Products}}
\endsubhead
The St-Re ring for a simplicial product is a {\it Segre product}.

\medskip
\noindent{\bf Definition.}
The {\it Segre product} of the graded ${\hbox{\tenbf A}}$-algebras
$R_1$ and $R_2$, denoted
$R=\sigma\!_{_{\hbox{\fivebf A}}}(R_1,R_2)$ or
$R=\sigma(R_1,R_2)$, is defined through;
$[R {\rlap{{\rlap{$]_{_p}$}{\ \ $=$}}}{\ \ \ \ \ [}} R_1]_{_p}{
\otimes_{_{\hbox{\fivebf A}}}}[R_2]_{_p}\!,\  \forall
p\!\in\!\hbox{\tenbf N}.$

\example{Example 1}
The trivial Segre product, $R_1\hbox{\tenbf 0}R_2$, is equipped
with the trivial prod-
\indent uct,
i.e. every {product} of elements, both of which lacks ring term,
equals 0.

\smallskip
\item{\bf 2.}
The ``canonical" $\!$Segre product, $R_1\underline\otimes R_2$, is
equipped with a product induced by extending $(\!$ linearly and
distributively$)$ the
{\it componentwise multiplication}
on simple homogeneous \hbox{elements: If $_{\!}m_1^\prime
\!_{\!}\otimes_{\!} m_1^{\prime\prime}$}
$\in \bigl[R_1\underline\otimes R_2\bigr]_\alpha$ and $m_2^\prime
\otimes m_2^{\prime\prime} \in \bigl[R_1\underline\otimes
R_2\bigr]_\beta$ then $(m_1^\prime \otimes m_1^{\prime\prime})
(m_2^\prime \otimes m_2^{\prime\prime}):= m_1^\prime m_2^\prime
\otimes m_1^{\prime\prime} m_2^{\prime\prime}\in
 \bigl[R_1\underline\otimes R_2\bigr]_{\alpha+\beta}$.

\smallskip
\item{\bf 3.}
{The ``canonical" generator-order sensitive Segre product,
$R_1{\bar\otimes} R_2$, of two graded ${\hbox{\tenbf k}}$-standard
algebras $R_1$ and $R_2$ presupposes the existence of a uniquely
defined partially ordered minimal set of generators for $R_1$
$(R_2)$ in $[R_1]_{_1}\ ([R_2]_{_1})$ and is equipped with a
product induced by extending $($distributively and linearly$)$ the
following operation defined on simple homogeneous elements, each
of which, now are presumed to be written, in product form, as an
increasing chain of the specified linearly ordered generators: If
$(m_{11} \otimes m_{21}) \in \bigl[R_1{\bar\otimes}
R_2\bigr]_\alpha$ and $(m_{12} \otimes m_{22}) \in \bigl[R_1{
\bar\otimes} R_2\bigr]_\beta$ then $(m_{11} \otimes m_{21})
(m_{12} \otimes m_{22}):= ((m_{11} m_{12} \otimes m_{21}
m_{22})\in \bigl[R_1{\bar\otimes} R_2\bigr]_{\alpha+\beta}$ if by
``pairwise" permutations, $(m_{11} m_{12}, m_{21} m_{22})$ can be
made into a chain in the product ordering, and 0 otherwise. Here,
$(x,y)$ is a pair in $(m_{11} m_{12}, m_{21} m_{22})$ if $x$
occupy the same position as $y$ counting from left to right in
$m_{11} m_{12}$ and $m_{21} m_{22}$ respectively. }
\endexample

\remark{Note 1}
(\cite{27} %
{\rm p.}\ 39-40$)$ Every Segre product of $R_1$ and $R_2$ is
module-isomorphic by definition and
{so}, they all have the same Hilbert series.
$\!${The {\it Hilbert series} of a graded $\hbox{\tenbf
k}$-algebra $R\!=\!\bigoplus\!\!_{_{^{i\ge 0}}}\!R_i\!$ is
{$\hbox{\tenbf Hilb}_{R}(t)\!:=\!\sum_{i\ge 0}(\hbox{\tenrm
H}(R,i))t^i\!:=\! \sum_{i\ge 0}(\hbox{\tenrm dim}_{\hbox{\fivebf
k}}R_i)t^i.$\

\vskip1.5pt
Here, ``$\dim$" is {\it Krull} {\it dimension} \rm and ``H" stands
for the ``{\it The Hilbert function}".}}

\smallskip
\noindent \noindent {\it {2}.}
If $R_1$ and $R_2$ are graded algebras finitely
generated (over $\hbox{\tenbf k}$) by $x_1,\dots,x_n\in[R_1]_{_1},
y_1,\dots,y_m\in[R_2]_{_1}$, resp., then $R_1{\underline\otimes}
R_2$ and $R_1{\bar\otimes} R_2$ are generated by $(x_1\otimes
y_1),\dots,(x_n\otimes y_m)$, and $\dim R_{_1}{\bar\otimes}
R_{_2}=\dim R_{_1}
{\underline\otimes}
R_{_2} =\dim R_{_1} + \dim R_{_2}-1.$

\smallskip
\noindent {\it {3}.}
The generator-order sensitive case covers all cases above.
In the theory of Hodge Algebras and in particularly in its
specialization to {\it Algebras with Straightening Laws}  (ASLs),
the generator-order is the
\noindent
main issue, cf. \cite{3} %
\S\ 7.1\
and \cite{14} %
p.\nobreak\ 123\nobreak\ ff.
\endremark

\medskip
\cite{8} %
p.\ 72 Lemma gives a reduced (Gr\"obner) basis, $C^\prime\cup D$,
for ``{\bf I}" $\!$in $\hbox{\tenbf k}
{\hbox{\tenbf[$\!\!{_{_{\!}}}{_{_{\!}}}$[}}
\Delta_{_{^{_{\!^{\!}}1\!}}}\!\times\Delta_{_{^{_{\!}2\!}}}
{\hbox{\tenbf]$\!\!{_{_{\!}}}{_{_{\!}}}$]}\!}
$ $\widetilde=\!\raise1pt\hbox{ $\hbox{\tenbf
k}[V\!\!\!_{_{^{\Delta_1}}}\!\!\!\times\!
V\!\!\!_{_{^{\Delta_2}}}\!]$} / \lower1pt\hbox{$\hbox{\tenbf I}$}\
\hbox{\tenrm with}$
$C^\prime\!\!:=\!\{ w_{_{^{\!\lambda,\mu}}}w\!_{_{^{\nu,\xi}}}|\
\lambda\!<\!\nu\land \mu\!>\!\xi\},
w_{\!_{^{\lambda,\mu}}}\!\!:= \!(v_{\!_{^{\lambda}}}\!,
v\!_{_{^{\mu}}}\!), v_{\!_{^{\lambda}}}\!\!\in\!
V\!\!\!_{_{^{\Delta_{^{\!}1}}}}\!,
v\!_{_{^{\mu}}}\!\!\in\!V\!\!\!_{_{^{\Delta_2}}}\!\!$ where the
subindices reflect the assumed linear ordering on the factor
simplices and with $\overline{p_i}$ as the projection down onto
the $i$:th factor;

\medskip
\noindent
{$\raise10.5pt\hbox{$D:=$}\raise3pt\hbox{ $^{
\Bigl\{\hbox{\eightbf w}\ \!=\
\!w_{_{^{{\!\lambda_{_{^{\!1}}}}\!\!,{{\mu_{_{^{\!1}}}}}}}} \!\!\!
\cdot\dots\cdot
w_{_{^{{\!\lambda_{_{^{\!k}}}}\!\!,{{\mu_{_{^{\!k}}}}}}}}
\big\vert
\ \Bigl[\bigl[\bigl[\{\overline{p_1}(\hbox{\eightbf w})\} \rlap{\
\raise4pt\hbox{\seveni n}}{\notin} {\Delta_1}\bigr]\ \!\land\ \!
\bigl[\{\overline{p_2}(\hbox{\eightbf w})\} \in{\Delta_2} \bigr] \
\!\land\ \! {{\lambda_1<\dots<\lambda_k} \brack {\mu_1\le\dots\le
\mu_k}}\bigr]\lor
\bigl[\bigl[ \{\overline{p_1}(\hbox{\eightbf w})\} \in{\Delta_1}
\bigr]\ \!\land
}$}$

\hfill
{$\raise0pt\hbox {$^{%
\land
\ \! \bigl[\{\overline{p_2}(\hbox{\eightbf w})\} \rlap{\
\raise4pt\hbox{\seveni n}}{\notin} {\Delta_2}\bigr] \ \!\land\
\!{{\lambda_1\le\dots\le \lambda_k} \brack {\mu_1<\dots<
\mu_k}}\bigr] \ \!\lor\ \!
\bigl[\bigl[\{\overline{p_1}(\hbox{\eightbf w})\} \rlap{\
\raise4pt\hbox{\seveni n}}{\notin} {\Delta_1}\bigr] \ \!\land\ \!
\bigl[\{\overline{p_2}(\hbox{\eightbf w})\} \rlap{\
\raise4pt\hbox{\seveni n}}{\notin} {\Delta_2}\bigr] \ \!\land\ \!
{{\lambda_1<\dots<\lambda_k}\brack{\mu_1<\dots<\mu_k}}\bigr]\Bigr]\Bigr\}.
}$}$}}

\smallskip
$
C^\prime\cup D\!=\!\{m_\delta\ \!\big\vert\ \! \delta\
{{\rlap{\raise4pt\hbox{$_{_{^{\!}}}n$}}\!\!{\not\in}}}\ \!
{\Delta_{_{^{_{\!}1\!}}}\!
\!\times\!{\Delta_{_{^{_{\!}2\!}}}}\}}$
and the identification
$ v_{\!_{^{\lambda}}}\!\otimes v\!_{_{^{\mu}}}\!\leftrightarrow
(v_{\!_{^{\lambda\!}}}, v\!_{_{^{\mu}}}\!)$
gives, see \cite{8} %
Theorem 1\ p.\ $\!$71, the following graded ${\hbox{\tenbf
k}}$-algebra isomorphism of degree zero;
$\hbox{\tenbf k} {\hbox{\tenbf[$\!\!{_{_{\!}}}{_{_{\!}}}$[}}
\Delta_{_{^{_{\!}1\!}}}\!\times\Delta_{_{^{_{\!}2\!}}}
{\hbox{\tenbf]$\!\!{_{_{\!}}}{_{_{\!}}}$]}\!}
$
$\ \!\widetilde=\!$
$\hbox{\tenbf k} {\hbox{\tenbf[$\!\!{_{_{\!}}}{_{_{\!}}}$[}}
\Delta_{_{^{_{\!}1\!}}}
{\hbox{\tenbf]$\!\!{_{_{\!}}}{_{_{\!}}}$]}}
\ \!{\bar\otimes}\ \!
\hbox{\tenbf k} {\hbox{\tenbf[$\!\!{_{_{\!}}}{_{_{\!}}}$[}}
\Delta_{_{^{_{\!}2\!}}}
{\hbox{\tenbf]$\!\!{_{_{\!}}}{_{_{\!}}}$]}},$
which, in the Hodge Algebra terminology,
is {\it the} \underbar{\it discrete} {\it algebra with the same
data} as
$\hbox{\tenbf k} {\hbox{\tenbf[$\!\!{_{_{\!}}}{_{_{\!}}}$[}}
\Delta_{_{^{_{\!}1\!}}}
{\hbox{\tenbf]$\!\!{_{_{\!}}}{_{_{\!}}}$]}}
\ \!{\underline\otimes}\ \!
\hbox{\tenbf k} {\hbox{\tenbf[$\!\!{_{_{\!}}}{_{_{\!}}}$[}}
\Delta_{_{^{_{\!}2\!}}}
{\hbox{\tenbf]$\!\!{_{_{\!}}}{_{_{\!}}}$]}}, $
cf.\ \cite{3}\ %
\S\ 7.1.\
If the discrete algebra is ``C-M" or Gorenstein (Definition p.\
22), so is the original by
\cite{3} 
Corollary\ 7.1.6.                
Any finitely generated graded {\bf k}-algebra has a Hodge
Algebra structure, see \hbox{\cite{14} %
p.\ 145.}

%


\subhead
{\rm II:4}\ \ \ {
Gorenstein Complexes}
\endsubhead
\normalbaselines

\definition{Definition 1}
$\hbox{\tenrm v}\!\in\! V_{_{^{\!\Sigma}}}\!$
is a {\it cone point} if $\hbox{\tenrm v}$\ is\ a\ vertex in every
maximal simplex in $\Sigma$.
\enddefinition

\definition{Definition 2}
\hbox{\tenrm (cp. \cite{26} %
Prop 5.1)}
Let $\Sigma$ be an arbitrary $($finite$)$ complex and put;
                   $\Gamma:=\hbox{\tenrm core}\Sigma
                   :=
                  \{\sigma\in\Sigma\ |\ \sigma\ \hbox{contains\ no\
                  cone\ points}\}.
                  $
Then;
$\emptyset$ isn't Gorenstein$_{\!_{\hbox{\fivebf G}}}\!$,\nobreak\
while

\noindent$\emptyset\not=\Sigma$ is Gorenstein$_{\!_{\hbox{\fivebf
G}}}\!$ (Gor$_{\!_{\hbox{\fivebf G}}}\!$) if \
{$
\mdoubleH_{i}(|\Gamma|,|\Gamma|\setminus_{_{\!o}\!} \alpha;
              \hbox{\tenbf G})=
\cases
0 & if\ i\ne \dim\Gamma
\cr
\hbox{\tenbf G} & if\ i= \dim\Gamma \cr
\endcases
\ \ \forall\ \!\alpha\in |\Gamma|.
$}
\enddefinition

\remark{Note}
$\Sigma$ Gorenstein$_{\!_{\hbox{\fivebf G}}}\!\!$
$\Longleftrightarrow \Sigma$ finite and $|\Gamma|$ is a
homology$_{_{^{\!\hbox{\fivebf G}}}}\!$ \hbox{\spaceskip1.81pt
sphere as defined in page 13.}
$\delta\!_{_{^{\Sigma}}}\!\!:=\{v\!\in\! V\!_{_{^{\Sigma}}}\!\!\
|\ \!v\ \hbox{\tenrm is\ a\ cone\ point}\}\!\in\!\Sigma.$
Now;
v is a cone point \underbar{iff} $\
\overline{{\hbox{\tenrm{st}}}}_{_{^{\!\Sigma\!}}}\!\hbox{\tenrm
v}\!=\!\Sigma$ and so,
$\ \overline{{\hbox{\tenrm{st}}}}_{_{^{\!\Sigma\!}}}\!
{{\delta}}_{_{^{\!\Sigma}}}\!
\!=\!
\Sigma=$ (core$\Sigma)\ast {\bar{\delta}}_{_{^{\!\Sigma}}}\!$ and
core$\Sigma$
$= \hbox{\tenrm Lk}_{_{^{\!\Sigma}}}\! \delta_{_{^{\!\Sigma}}}\!
:= \{\tau\!\in\! \Sigma| [\delta_{_{^{\!\Sigma}}}\!\cap \tau
=\emptyset]\land [\delta_{_{^{\!\Sigma}}}\!\cup \tau\in
\Sigma]\}.$
\endremark

\proclaim{Proposition 1}
{\rm (Cf. \cite{8} %
p.\ 77. %
)} {\sl $\!\Sigma\!_{_{^{1}}}
 \!\ast\Sigma \!_{_{^{2}}} $ Gorenstein$\!_{_{^{\hbox{\fivebf G}}}}
\!\!\!\Longleftrightarrow\! \Sigma\!_{_{^{1}}},\Sigma\!_{_{^{2}}}$
both Gorenstein$\!_{_{^{\hbox{\fivebf G}}}}\!.
\!\!\!\!\!$}%
\qed
\endproclaim

Gorensteinness is, unlike
``{Bbm}$_{_{\!\hbox{\fivebf G}}}\!\!\!"$-,
$\!``${CM}$_{_{\!\hbox{\fivebf G}}}\!\!"\!$-
and
2-$``${CM}$_{_{\!\hbox{\fivebf G}}}\!\!"$-ness,
triangulation-sensi-tive and in\noindent\ particular, the
Gorensteinness for products is sensitive to the partial orders,
assumed in the $\hbox{\tenrm definition,\ given\ to}$ the vertex
sets of the factors.
\hbox{See p.\ 21 for}
$\{m_\delta\ \!\big\vert\ \!
\delta{\rlap{\vbox{\moveright4pt\hbox{\raise4pt\hbox{$n$}}}}{\not\in}
}\ \!
{\Delta_{_{^{_{\!}1\!}}}\!
\!\times\!{\Delta_{_{^{_{\!}2\!}}}\!}\}}.\
$
In \cite{8} %
p.\ 80, the product is represented in the form of matrices, one
for each pair $(\delta_{_{^{\!1}}},\delta_{_{^{\!2}}})$ of maximal
simplices $\delta\!_{_{^{i}}}\!\in\!\Delta_{_{\!i}},\ i=1,2.$
It is then easily seen that a cone point must occupie the upper
left corner in each matrix or the lower right corner in each
matrix. So a product ($\dim\Delta_{_{i}}\ge1$) can never have more
than two cone points. For Gorensteinness to be preserved under
product the factors must have at least one cone point to preserve
even $``\hbox{\tenrm CM}_{\!_{\hbox{\fivebf G}}}\!\!"$-ness, by
Corollary\ iii\nobreak\ p.$\ \!$12.

Bd(core$(\Delta_{_{1}}\!\!\times\! \Delta_{_{2}}))=\emptyset$
demands each $\Delta_{_{i}}\!$ to have as many cone points as
$\Delta_{_{1}}\!\!\times\! \Delta_{_{2}}\!.\ $So;

\proclaim{Proposition 2}
{\rm (Cf. \cite{8} %
p.\ 83ff. for proof.)} Let $\Delta_{_{1}}\!,\ \!\Delta_{_{2}}$ be
two arbitrary finite simplicial complexes
 with $\dim\!\Delta_{_{i}}\!\geq\!1,(i\!=\!1,2)$ and
                    a linear order defined on their
                    vertex sets $V_{\Delta_1}, V_{\Delta_2}$ respectively,
then;

$\Delta_{_{1}}\!\!\times\!\Delta_{_{2}}$
Gor$_{\!_{\hbox{\fivebf G}}}\!$
$\!\Longleftrightarrow\!\Delta_{_{1}},\ \!\Delta_{_{2}}$
both Gor$_{\!_{\hbox{\fivebf G}}}\!$
\indent\hbox{\tenbf and}\indent  condition \hbox{\tenbf I}
\hbox{\tenbf or} \hbox{\tenbf II}\nobreak\ holds, where;

\smallskip
\noindent {\bf (I)} \ $\Delta_{_{1}},\ \Delta_{_{2}}$
has exactly one cone point each - either both minimal {\tenbf or}
both maximal.
\hfill\break
{\bf (II)}\ $\Delta_{_{\!i}}\!$\ 
has exactly two cone points,
one minimal \hbox{\tenbf and} the other maximal,\ $i\!=\!1,2.\!$
\qed
\endproclaim

\example{Example}
Gorensteinness is character sensitive!
Let $\Gamma:=\hbox{\rm core}\Sigma=\Sigma$ be a 3-dimensional
$\hbox{\tenrm Gor}\!_{_{^{\hbox{\fivebf k}}}}\!$ complex where
${\hbox{\tenbf k}}$ is the prime field $\hbox{\tenbf
Z}_{_{^{\!\hbox{\fivebf p}}}}\!$ of characteristic {\bf p}.
This implies, in particular, that $\Gamma$ is a
homology$\!_{_{^{{\hbox{\fivebf Z}}}}}$ $3$-manifold. Put
$\mdoubleH_{i}:=\mdoubleH_{i}(\Gamma;\hbox{\tenbf Z})$, then;
$ \mdoubleH_{_{^{_{\!}0}}}\!=0, $
$ \mdoubleH_{_{^{3}}}=\hbox{\tenbf Z}\ \hbox{\tenrm and}\
\mdoubleH_{_{^{2}}} $ has no torsion by Lemma\ $\!$1.{\bf i} p.\
24.
Poincare' duality \hbox{\tenrm and} \cite{25} %
p.\ 244\ Corollary\ 4 gives $
\mdoubleH_{_{^{1}}}\!\!=\mdoubleH_{_{^{2}}}\!\oplus\ \!
\hbox{\tenbf t}\mdoubleH_{_{^{1}}},\ \hbox{\tenrm where}\
\hbox{\tenbf t}\circ := \hbox{\tenrm the\
torsion\hbox{-}submodule\ of}\ \circ. $ So a $ \Sigma=\Gamma\
\hbox{\tenrm with\ a\ pure\ torsion}\ \mdoubleH_{1}=\hbox{\tenbf
Z}_{{\hbox{\fivebf p}}}, $ say, gives us an example  of a
presumptive character sensitive Gorenstein complex. Examples of
such orientable compact combinatorial manifolds without boundary
is given by the projective space of dimension 3, $ {\hbox{\tenbf
P}^3} $ and the lens space $\hbox{\tenbf L}(n,k) $ where $
\mdoubleH_{_{1}}\!({\hbox{\tenbf P}^3};\hbox{\tenbf Z})=
\hbox{\tenbf Z}_{_{{\hbox{\fivebf 2}}}} $ and $
\mdoubleH_{_{1}}\!({{\hbox{\tenbf L}}}(n,k);\hbox{\tenbf Z})=
\hbox{\tenbf Z}_{_{{\hbox{\fivebf n}}}}. $ So, $ {\hbox{\tenbf
P}^3}\ast\bullet $ \big( $ \hbox{\tenbf L}(n,k)\ast\bullet $ \big)
is $ \hbox{\tenrm Gor}\!{\lower2pt\hbox{\fivebf k}} $ for $
\hbox{\tenrm char}\hbox{\tenbf k}\ne2\ (\hbox{\tenrm
char}\hbox{\tenbf k}\ne n) $ while it is not even Buchsbaum for
$\hbox{\tenrm
char}\hbox{\tenbf k}=2\ (\hbox{\tenrm char}\hbox{\tenbf k}=n) $, cf. \cite{21} %
p. 231-243 for details on $ {\hbox{\tenbf P}^3}\ \hbox{\tenrm
and}\ \hbox{\tenbf L}(n,k). $
Cf. \cite{26} %
Prop.\ 5.1 p.\ 65 or \cite{8} %
p.\ 75 for Gorenstein equivalences.
A $\hbox{\tenrm Gorenstein}_{_{^{\hbox{\fivebf k}}}}\Delta$ isn't
in general shellable, since if so, $\Delta$ would be
CM$_{_{^{\!\hbox{\fivebf Z}}}}$ but $ \hbox{\tenbf L}(n,k) $ and $
\hbox{\tenbf P}^3\ \hbox{isn't.} $
Indeed, in 1958 M.E. Rudin published {\it An Unshellable
Triangulation of the Tetrahedron}.

Other examples are given by Jeff Weeks' computer program $\!
$``SnapPea" hosted at
http://thames.northnet.org/weeks/index/SnapPea.html,
e.g.
$_{\!}\mdoubleH_{_{1}}\!({\Sigma}\!\!\!\!{{{\lower3.5pt\hbox{\fivebf
fig8}}}}\!(5,1);\hbox{\tenbf Z})\!=\!\nobreak \hbox{\tenbf
Z}_{\hbox{\fivebf 5}}$
for the old tutorial example of SnapPea
$ {\Sigma}\!\!\!\!{{{\lower3.5pt\hbox{\fivebf fig8}}}}\!(5,1),$
i.e. the Dehn surgery filling w.r.t. $\!$a figure eight complement
with diffeomorphism kernel generated by (5,1).
$\!{\Sigma}\!\!\!\!{{{\lower3.5pt\hbox{\fivebf fig8}}}}\!(5,1)
\ast\bullet\ \hbox{\tenrm is}$
Gor$\!_{\hbox{\fivebf k}}$ if char\hbox{\tenbf k}$\ne5$ but not
even {\tenrm Bbm}$_{\hbox{\fivebf k}}$ if char$\hbox{\tenbf
k}=\!5$,
\hbox{cf. \cite{24} %
Ch. $\!$9 %
for more on surgery.}
\endexample

%


\head
{III:\ \ \ Simplicial Manifolds}
\endhead

\subhead
{{III:}1 \ {Definitions}}
\endsubhead

\normalbaselines

\medskip
We will make extensive use of Proposition\ 1 p.\ 11 without
explicit notification.

\definition{Definition 1}
An $n$-dimensional pseudomanifold is a locally finite $n$-complex
${\Sigma}$ such that;\hfill \break
($\alpha$) ${\Sigma}$ is pure,\ \
i.e. the maximal simplices in ${\Sigma}$ are all
$n$-dimensional.\hfill \break
{\bf($\beta$)} Every $(n-1)$-simplex of ${\Sigma}$ is the face of
at most two $n$-simplices of ${\Sigma}$.\hfill \break
{\bf($\gamma$)} If $s$ and $s'$ are $n$-simplices in ${\Sigma}$,
there is a finite sequence $s=s_0,s_1,\ldots s_m=s'$ of
$n$-simplices in ${\Sigma}$ such that $s_i\cap s_{i+1}$ is an
$(n-1)$-simplex for $0\le i<m$.

The {\it boundary}, {\rm Bd}${\Sigma}$, of an $n$-dimensional
pseudomanifold ${\Sigma}$, is the subcomplex generated by those
$(n-1)$-simplices which are faces of exactly one $n$-simplex in
${\Sigma}$.
\enddefinition

\definition{Definition 2}
$\Sigma={\bullet}{\bullet}$ is a {\it quasi-$0$-manifold}. Else,
$\Sigma$ is a {\it quasi-$n$-manifold} if it's an $n$-dimensional,
locally finite complex fulfilling; \hfill\break
                   \hbox{\tenbf($\alpha$)}
\hbox{\tenbf $\Sigma$} is pure.\indent
\hbox{\tenrm($\alpha$ is redundant since it is a consequence of
$\gamma$ by Lemma\ 2 p.\ 16.)}\hfill\break
                   \hbox{\tenbf($\beta$)} Every $(n-1)$-simplex of \hbox{\tenbf $\Sigma$}
                           is the face of at most two
                            $n$-simplices of $\Sigma$.\hfill\break
\hbox{\tenbf($\gamma$)} \hbox{\tenrm Lk}$_{_{\Sigma}}\sigma$ is
connected i.e. $\mdoubleH _{0}(\hbox{\tenrm
Lk}_{_{\Sigma}}\sigma;\hbox{\tenbf G})=0$ for all
                                         $\sigma\in\Sigma$,     
                                        s.a. dim$\sigma<n-1$.

\smallskip
The {\it boundary with respect to \hbox{\tenbf G}}, denoted
$\hbox{\tenrm Bd}_{_{\!\hbox{\fivebf G}}}\!\Sigma,$
\footnote{For a classical manifold
{\eighti X}
\rlap{\raise0.3pt\hbox{\sevenrm {\char"2F}}}{$\!=$}$\!$
{\eightsy {\char"0F}}
\underbar{s.a. {\eightbf Bd}{\eighti X}={\eightsy{\char"3B}}};
{\eightbf Bd}{\eightsy{\char"46}}
$\!\!\!$\lower1.5pt\hbox{\fivei{\char"7D}}({\eighti X}) $\!$=$\!$
{\eightsy{\char"3B}} if {\eighti X} is compact
and orientable and {\eightbf Bd}{\eightsy{\char"46}}$
\!\!$\lower1.5pt\hbox{\fivei{\char"7D}}({\eighti X}) $\!$=$\!$
{\eightsy{\char"66}}$_{\!}${\eighti{\char"7D}}$_{\!}${\eightsy{\char"67}}
else.%
}
of a quasi-$n$-manifold $\Sigma$, is the set of simplices
{$\hbox{\tenrm Bd}_{_{\!\hbox{\fivebf G}}}\Sigma
\!:=\!\{\sigma\!\in\! \Sigma\ |\ \mdoubleH_{n}(\Sigma,\hbox{\tenrm
cost}_{_{\Sigma}}\sigma;\hbox{\tenbf
G})=0\}$}, where \hbox{\tenbf G} is a unital module  over a ${%
commutative\ ring\ _{\!}\hbox{\tenbf A}.}\!$
$(\beta$ in Def. 1-2 $\Rightarrow \bullet\ \hbox{\tenrm and}\
\bullet_{\!}\bullet$ {\rm are the only 0-manifolds}.$)$
\enddefinition

\remark{Note 1}
$\!\Sigma$ is (locally) finite $\!\Longleftrightarrow\!$
$\vert\Sigma\vert$ is (locally) compact.
By Th.\ 5 p.\ 10; if $X$ is a homology$_{_{\!\!\hbox{\fivebf
R}}}\!$ $n$-manifold ($n$-hm$_{_{\!\!\hbox{\fivebf R}}}\!$) then
$X_{\!}$ is a $n$-hm$_{_{\!\!\hbox{\fivebf G}}}\!$
for any \hbox{\tenbf R}-$_{\!}\hbox{\tenbf PID}\ \hbox{\tenrm
module}\ \hbox{\tenbf G}$.\break
\tenrm \cite{25}\ p.\ 207-8 \hbox{\tenbf +} p.\ 277-8 treats the
classical standard case
$_{\!}\hbox{\tenbf R}\!=\!\hbox{\tenbf
G}\!_{\!}=\!_{\!}\hbox{\tenbf Z}\!_{\!}:=\!_{\!}\{\hbox{\tenrm the
integers}\}_{\!}$, which in our exposition more represents a
particularly straightforward
\hbox{extreme case.}
A simplicial complex $\Sigma$ is called a hm$_{_{\!\!\hbox{\fivebf
G}}}$ if $\ \!|\Sigma|$ is, -- now $n\!=\dim\Sigma.$

From a purely technical point of view we really don't need the
``locally finiteness"-assumption, as is seen from Corollary p.\
12.
\endremark

\definition{Definition 3}
\hbox{\tenrm (Let ``manifold" stand for pseudo-, quasi- or
homology manifold.)}

\noindent A compact $n$-manifold, $\hbox{\tensy S}$, is {\it
orientable}$_{_{{^{^{\!\!\hbox{\fivebf G}}}}}}\!\!\ { if}\ \!$
$\mdoubleH\!_{_{{n}}}\!(\hbox{\tensy S}\!,\hbox{\tenrm
Bd}\hbox{\tensy S};\!\hbox{\tenbf G})\tilde =\hbox{\tenbf G}.$
An $n$-manifold is {\it orientable}$_{_{{^{^{\!\!\hbox{\fivebf
G}}}}}}\!\!$ if all its compact $n$-submanifolds are orientable
$\!-\!$ else, {\it non-orientable}$_{_{{^{^{\!\hbox{\fivebf
G}}}}}}$.
Orientability is left undefined for $\emptyset$.
\enddefinition

\definition{Definition 4}
$\{\hbox{\tenbf B}_{_{\!\!{\hbox{\fivebf
G}\!,_{\!}j}}}\!\!\!\!\!^{^{_{\Sigma\!\!}}}\
\!\}\!_{_{_{{{j\in\lower0.6pt\hbox{\sevenbf I}}}}}}\!\!\!$
{is {\it the set of strongly connected boundary components of}}
$\Sigma$
{if}
$\{\hbox{\tenbf B}_{_{\!\!{\hbox{\fivebf
G}\!,_{\!}j}}}\!\!\!\!\!^{^{_{\Sigma\!\!}}}\
\!\}\!_{_{_{{{j\in\lower0.6pt\hbox{\sevenbf I}}}}}}\!\!\!$
{is the maximal strongly connected components of} $ \hbox{\tenrm
Bd}\!\!\raise0.6pt\hbox{$_{_{_{\hbox{\fivebf G
}}}}$}\!\!\!\!{{\Sigma}},$ from Definition\ 1 p.\ 15.
$(\Rightarrow\hbox{\tenbf B}_{_{\!\!{\hbox{\fivebf
G}\!,_{\!}j}}}\!\!\!\!\!^{^{_{\Sigma\!\!}}}\ \!$ pure
and if $\sigma$ is a \underbar{maximal} simplex in $\hbox{\tenbf
B}_{_{{\!{j}}}}\!$, then;
$\hbox{\tenrm Lk}\!\!\!\!\! \raise0.5pt\hbox{$_{_{_{\hbox{\fiverm
Bd}\!_{_{\!\hbox{\fivebf G}}}\!\!\!\Sigma}}}$}\!\!\!
\sigma\!_{_{^{\!{\ }}}}\!\!=\! \hbox{\tenrm Lk}\!\!
\raise0.5pt\hbox{$_{_{_{\hbox{\fivebf B}_{_{{\!{j}}}}\!}}}$}\!\!\!
\sigma\!_{_{^{\!{\ }}}} \!\!=
\!\{{\emptyset}_{_{^{\!o}}}\!_{\!}\}) $.
\enddefinition

\remark{Note 2}
$\!\emptyset, \{\emptyset_o\!\}\!,$ and $0$-dimensional complexes
with either one, ${\bullet},$ or two, ${\bullet}{\bullet},$
vertices are the only manifolds in dimensions $\le\!0$, and the
$|$1-manifolds$|$ are finite/infinite $1$-circles and (half)lines,
while [${\Sigma}$ is a quasi-2-manifold]$\ \!\Longleftrightarrow\!
[{\Sigma}\ \hbox{\tenrm is\ a\ homology}\!_{_{^{\hbox{\fivebf
Z}}}}\!$ 2-manifold].$\ \!$ Def. $\hbox{\tenbf 1}.\gamma\nobreak\
\hbox{\tenrm is}$ paraphrased by $\!``{\Sigma}\ \hbox{\tenrm is}$
{\it strongly connected}$\ \!"$ and $\bullet\bullet$-complexes,
though strongly connected, $\hbox{\tenrm are\ the\ only\nobreak\
non\hbox{-}connected\nobreak\ manifolds.}$%
\indent Note also that 
$\hbox{\tenbf S}^{^{_{_{^{\!}}\hbox{\raise0.4pt\hbox{\fivebf
-}$_{\!}$\fiverm 1}}}}\!\!\! :=\{\emptyset_o\!\} $ is the boundary
of the $0$-ball, $\bullet$, the double of which is the $0$-sphere,
$\bullet\bullet$. \noindent Both the $(-1)$-sphere
$\{\emptyset_o\!\}$ and the $0$-sphere $\bullet\bullet$ has, as
preferred, empty boundary.

$\bullet$ is the only compact orientable manifold with its
boundary equal to
$\{\wp_{\!}\}\ \!(\{_{\!}\emptyset\!_{_{^{o}}}\!_{\!}\})$.
\endremark

%


\subhead
III:2 \ \ \ Auxiliaries
\endsubhead

\normalbaselines

\proclaim{Lemma 1}
{ For a finite n-pseudomanifold $\Sigma;\ \!$

\smallskip
\item {\bf i.} {\rm $($cf. \cite{25} %
p.\ 206 Ex.\ {\bf E}2.$)$}
\smallskip
$\mdoubleH_{n}\!(\Sigma,\hbox{\tenrm Bd}\Sigma_{_{}};\hbox{\tenbf
Z}) \!=\!\hbox{\tenbf Z}$ and
$\mdoubleH_{n-1}\!(\Sigma,\hbox{\tenrm
Bd}\Sigma_{_{}};\hbox{\tenbf Z})$ has no torsion,
\noindent
or
\noindent
$\mdoubleH_{n}\!(\Sigma,\hbox{\tenrm Bd}\Sigma_{_{}};\hbox{\tenbf
Z}) \!=\!0$ and the torsion submodule of
$\mdoubleH_{n-1}\!(\Sigma,\hbox{\tenrm
Bd}\Sigma_{_{}};\hbox{\tenbf Z})$ is isomorphic to $\hbox{\tenbf
Z}_{_{^{_{\!}2}}}.$

\smallskip
\item {\bf ii.}
$\Sigma$ non-orientable$_{_{\!\hbox{\fivebf G}}}\!\!$
$\Longleftrightarrow\!$ $\Sigma$
non-orientable$_{_{\!\hbox{\fivebf Z\!\!}}}$ and $\hbox{\tenrm
Tor}\!_{_{^{1}}}\!\!\!^{^{_{\hbox{\fivebf
Z}}}\!}\!\bigl(\hbox{\tenbf Z}_{_{\!2}}\!,\!\hbox{\tenbf G}
\!\bigr)\!\!\ne\nobreak\!\nobreak\hbox{\tenbf G}.\!$

\indent
So in particular:
Manifolds\ are\ orientable w.r.t. $\!\hbox{\tenbf Z}_{_{{ 2}}}.$ }
\smallskip
\noindent
\item {\bf iii.}
$\mdoubleH^{{n}}\!(\Sigma,\hbox{\tenrm
Bd}\Sigma_{_{}};_{\!}\hbox{\tenbf Z}) \!=\!\hbox{\tenbf Z}\
(\hbox{\tenbf Z}_{_{^{_{\!}2}}}\!)$ when $\Sigma$ is $({
non}\hbox{-}\!){ orientable}.$
\endproclaim

\demo{Proof}
\noindent
By conditions $\alpha$ and $\beta$ in Def. $\!$1, a possible
relative n-cycle in ${C}_{{n}}^{{o}}\!(\Sigma,\hbox{\tenrm
Bd}\Sigma_{_{}};\hbox{\tenbf Z}_{_{^{\!\hbox{\fivebf m}}}})$ must
include all oriented n-simplices all of which with coefficients of
one and the same value.
When the boundary function is applied to such a possible relative
n-cycle the result is an $(n-1)$-chain that includs all oriented
$(n-1)$-simplices, not supported by the boundary, all of which
with coefficients 0 or
${\rlap{\raise2pt\hbox{\fiverm+}}{_{\!\!}\lower2pt\hbox{\
-}}}2c\in \hbox{\tenbf Z}_{_{^{\!\hbox{\fivebf m}}}}.$ So,
$\mdoubleH_{{n}}\!(\Sigma,\hbox{\tenrm Bd}\Sigma_{_{}};{{\tenbf
Z}_{_{^{_{\!}2}}}}) \!=\!\hbox{\tenbf Z}_{_{^{_{\!}2}}}\ \!-\!$
allways, and $\mdoubleH_{{n}}\!(\Sigma,\hbox{\tenrm
Bd}\Sigma_{_{}};\hbox{\tenbf Z}) \!=\!\hbox{\tenbf Z}\ (0)$
$\Longleftrightarrow$ $\mdoubleH_{{n}}\!(\Sigma,\hbox{\tenrm
Bd}\Sigma_{_{}};\hbox{\tenbf Z}_{_{^{_{\!}\!\hbox{\fivebf m}}}})
\!=\!\hbox{\tenbf Z}_{_{^{_{\!}\!\hbox{\fivebf m}}}}\ (0)$ if
$\hbox{\tenbf m}\ne2.$ The Universal Coefficient Theorem (=Th.\ 5\
p.\ 10)\nobreak\ now\nobreak\ gives;

\medskip
$
\hbox{\tenbf Z}_{_{^{2}}}=
\mdoubleH_{_{n\!}} (\Sigma,\hbox{\tenrm
Bd}\Sigma_{_{}};\hbox{\tenbf Z}_{_{^{ 2}}})
\ \! {{{_{\hbox{\fivebf Z}}}}\atop{{\raise2pt\hbox{$\cong$}} }} \
\!
\mdoubleH_{_{n\!}} (\Sigma,\hbox{\tenrm Bd}\Sigma_{_{}};
\hbox{\tenbf Z}\otimes\!_{_{\hbox{\fivebf Z}}}\!\hbox{\tenbf
Z}_{_{^{2}}})
\ \! {{{_{\hbox{\fivebf Z}}}}\atop{{\raise2pt\hbox{$\cong$}} }} \
\!
$

\smallskip
\hfill
{$ \ \! {{{_{\hbox{\fivebf Z}}}}\atop{{\raise2pt\hbox{$\cong$}} }}
\ \!
\big[\mdoubleH_{_{n\!}}(\Sigma,\hbox{\tenrm
Bd}\Sigma_{_{}};\hbox{\tenbf Z}) \otimes\!_{_{\hbox{\fivebf
Z}}}\!\hbox{\tenbf Z}_{_{^{2}}}\big]
\oplus
\hbox{\tenrm Tor}^{\!\hbox{\fivebf Z}}_1\!
\big(\mdoubleH_{_{{n\!-\!1}}}\! (\Sigma,\hbox{\tenrm
Bd}\Sigma_{_{}};\hbox{\tenbf Z}), {\hbox{\tenbf Z}}_{_{^{2}}})
$}
\indent and,

\medskip
$
\mdoubleH_{_{n\!}} (\Sigma,\hbox{\tenrm
Bd}\Sigma_{_{}};\hbox{\tenbf Z}_{_{^{\hbox{\fivebf m}}}})
\ \! {{{_{\hbox{\fivebf Z}}}}\atop{{\raise2pt\hbox{$\cong$}} }} \
\!
$
{$ \mdoubleH_{_{n\!}} (\Sigma,\hbox{\tenrm Bd}\Sigma_{_{}};
\hbox{\tenbf Z}\otimes\!_{_{\hbox{\fivebf Z}}}\!\hbox{\tenbf
Z}_{_{^{\hbox{\fivebf m}}}})
\ \! {{{_{\hbox{\fivebf Z}}}}\atop{{\raise2pt\hbox{$\cong$}} }} \
\!
$}

\smallskip
\hfill
{$ \ \! {{{_{\hbox{\fivebf Z}}}}\atop{{\raise2pt\hbox{$\cong$}} }}
\ \!
\big[\mdoubleH_{_{n\!}}(\Sigma,\hbox{\tenrm
Bd}\Sigma_{_{}};\hbox{\tenbf Z}) \otimes\!_{_{\hbox{\fivebf
Z}}}\!\hbox{\tenbf Z}_{_{^{\hbox{\fivebf m}}}}\big]
\oplus
\hbox{\tenrm Tor}^{\!\hbox{\fivebf Z}}_1\!
\big(\mdoubleH_{_{{n\!-\!1}}}\! (\Sigma,\hbox{\tenrm
Bd}\Sigma_{_{}};\hbox{\tenbf Z}), {\hbox{\tenbf
Z}}_{_{^{\hbox{\fivebf m}}}}), \indent\indent
$} \noindent

\smallskip
\noindent where the last homology module in each formula,
by \cite{25} %
p.\ 225 Cor.\ 11, can be substituted by its torsion submodule.
Since $\Sigma$ is finite, $\mdoubleH_{_{{n\!-\!1}}}\!
(\Sigma,\hbox{\tenrm Bd}\Sigma_{_{}};\hbox{\tenbf Z})
=\hbox{\tenbf c}_{_{^{_{\!}1}\!}}\oplus \hbox{\tenbf
c}_{_{^{_{\!}2}\!}}\oplus...\oplus \hbox{\tenbf
c}_{_{^{\hbox{\tenrm s}}\!}} $ by The { Structure Theorem} for
Finitely Generated Modules over {\bf PID}s,
cf. \cite{25} %
p.\ 9.
Now, a simple check, using \cite{25} %
p.\ 221 Example\ 4, gives {\bf i}, which gives {\bf iii}
by \cite{25} %
p.\ 244 Corollary.
Theorem\ 5 p.\ 10 and {\bf i} implies \hbox{\tenbf ii}.
\qed
\enddemo

Proposition\ 1 p.\ $\!$18 together with Proposition\ 1 p.\ 11
gives the next Lemma.

\proclaim{Lemma 2.i}
$\!\!\Sigma$ is a $n$-hm$_{_{^{\!\hbox{\fivebf G}}}}\!$
$_{\!}{{_{\Longrightarrow}}\atop{^{\
\!\not\!\!\!\Longleftarrow}}}\!$  $\Sigma$ is a
quasi-$n$-manifold$\ _{\!}{{_{\Longrightarrow}}\atop{^{\
\!\not\!\!\!\Longleftarrow}}}\Sigma$ is an $n$-pseudomanifold.

\noindent {\bf ii.}
$\!\Sigma$ is a $n$-hm$_{_{^{\!\hbox{\fivebf G}}}}\!$
\underbar{iff} it's a $\!$``Bbm$_{_{\!^{\hbox{\fivebf
G}\!\!\!}}}"\!$ pseudomanifold
$ {and}\
\mdoubleH\!\!\!\!\!\!\lower1pt\hbox{$_{_{^{\!{n\lower1pt\hbox{-}\hbox{\fivebf\#}\sigma\!}}}}$}
\!\!\!\!\raise1pt\hbox{\eightbf(}\hbox{\tenrm
Lk}\!_{_{_{{{\Sigma}}}}}\!\!\!\sigma;\!\hbox{\tenbf
G}\raise1pt\hbox{\eightbf)}\!=_{\!}0\ _{\!} {or}\ \! \hbox{\tenbf
G}\
\!\forall\ \!\sigma{\!}\not=\!\emptyset_{_{^{\!o}}}\!.
\square_{\!}$
\endproclaim

The ``only if"$\!$-part of Th.\ 10 was given
for finite quasi-manifolds in {\rm \cite{10} %
p.\ 166}. Def.\ 1 p.\ 15 makes perfect sense even for
non-simplicial posets like $\Gamma {\raise1.5pt\hbox{\eightmsbm
\char"72}}                                             
\Delta\!^{^{\!_{\ }}}$ (Def. 2 p.\ 16) which allow us to say that
$\Gamma{\raise1.5pt\hbox{\eightmsbm \char"72}}\Delta\!^{^{\!_{\
}}}$
is or is not {\it strongly connected} ({\it as a poset}) depending
on whether $\Gamma{\raise1.5pt\hbox{\eightmsbm
\char"72}}\Delta\!^{^{\!_{\ }}}$ fulfills Def.\ 1 p.\ 15 or not.
Now, for quasi-manifolds $\Gamma{\raise1.5pt\hbox{\eightmsbm
\char"72}}\Delta\!^{^{\!_{\ }}}$ {\it connected as a poset} (Def.
2 p.\ 16) is equivalent to $\Gamma{\raise1.5pt\hbox{\eightmsbm
\char"72}}\Delta\!^{^{\!_{\ }}}$ {\it strongly connected} which is
a simple consequence of Lemma\ 1 p.\ 16 and the
definition of quasi-manifolds, cf. \cite{10} %
p.\ 165 Lemma\ 4.
$\Sigma{\raise1.5pt\hbox{\eightmsbm \char"72}}\hbox{\tenrm
cost}_{_{\!\Sigma}}\!\sigma\!^{^{\!_{\ }}}$ is connected as a
poset for any simplicial complex $\Sigma$ and any
$\sigma\!\ne\!\emptyset$ and $\Sigma{\raise1.5pt\hbox{\eightmsbm
\char"72}}\hbox{\tenrm Bd}_{\!}{{\Sigma}}$ is strongly connected
for any pseudomanifold $\Sigma$. The ``if"-part of Theorem\ 10 can
fail for an infinite $\Sigma$.

\proclaim{Theorem 10}
If $\hbox{\tenbf G}$ is a module over a commutative ring {\bf A}
with unit, $\Sigma$ a finite $n$-pseudomanifold, and
$\Delta \raise2pt\hbox{\mysubsetneqq}  \Gamma
\raise2pt\hbox{\mysubsetneqq} \Sigma$ and
$\dim\Gamma\!=\!\dim\Sigma$
{\rm(Injectivity otherwise trivial!)}
then;
$\Sigma{\raise1.5pt\hbox{\eightmsbm \char"72}}\Delta\!^{^{\!_{\
}}}$ is strongly connected
\underbar{iff} \
$\mdoubleH_{n}(\Sigma,\Delta;\hbox{\tenbf G})\longrightarrow
 \mdoubleH_{n}(\Sigma,\Gamma;\hbox{\tenbf G})\ {is\ an\
 injection}.$
\endproclaim

\demo{Proof}
Each strongly connected $n$-component $\Gamma\!_{^{_{i}}}$ of
$\Gamma$ is an $n$-pseudomanifold with
(Bd$\Gamma\!_{^{_{i}}}\!)\!^{^{_{n-1}}}\!\!\!\ne\!\emptyset$
since $\Gamma \raise2pt\hbox{\mysubsetneqq} \Sigma$.
Now; Bd$\Gamma\!_{^{_{i}}}\!$ is all imbedded in
$\Delta_{^{_{i\!}}}\!\!:=_{\!}\Gamma\!_{^{_{i\!}}}\!\cap\!\Delta$
\underbar{iff}
$\mdoubleH_{n}\!(\Gamma\!_{^{_{i}}}\!,\Delta_{^{_{i}}}\!;\hbox{\tenbf
G})\!\ne\!0$,
i.e. \underbar{iff} every sequence from
$\gamma\!\in\!(\Gamma\!_{^{_{i}}}\!{\raise1.5pt\hbox{\eightmsbm
\char"72}}\Delta_{^{_{i}}}\!\!)\!^{^{_{n}}}\!$
to
$\sigma\!\in\!(\Sigma{\raise1.5pt\hbox{\eightmsbm
\char"72}}\Gamma)\!^{^{_{n}}}\!$
connects via\nobreak\ $\Delta_{^{_{i}}}
\!\!\!\!^{^{_{_{\!}n_{\!}-1}}}\!\!.$
Now; the relative {\bf LHS} w.r.t. $\!(\Sigma,\Gamma,\Delta)$
gives our claim. %
(Cp. the Jordan Curve Theorem.)$\ \!\square$
\enddemo

\goodbreak

$\Gamma\!=\hbox{\tenrm cost}\!_{_{\Sigma}}\tau\ \hbox{\tenrm and}\
\Delta\!=\hbox{\tenrm
cost}\!_{_{\Sigma\!}}\emptyset\!_{_{^{o}}}\!\!=\emptyset$ resp.
$\hbox{\tenrm cost}\!\!\!\! {\lower1.0pt\hbox{${
_{_{{\!\!\hbox{\fiverm Bd}\!_{_{_{\ }}}\!\!\!{{{\Sigma}} }}}}\!
}$}} \!\tau$ gives {\bf b} in the next Corollary.

\proclaim{Corollary 1}
If $\Sigma$ is an $n$-pseudomanifold and
$\sigma\ \!\raise2pt\hbox{\mysubsetneqq} \tau\!\in\!\Sigma$, then;
\smallskip
\noindent $\hbox{\tenbf a.}\ \ \mdoubleH_{n}(\Sigma,\hbox{\tenrm
cost}\!\!_{_{_{^{\Sigma}}}}\!\!\sigma;\hbox{\tenbf G})
\longrightarrow \mdoubleH_{n}(\Sigma,\hbox{\tenrm
cost}\!\!_{_{_{^{\Sigma}}}}\!\!\tau;\hbox{\tenbf G})\ {
injective}\ \underline{if_{\!}f} \
\Sigma{\raise1.5pt\hbox{\eightmsbm \char"72}}
\hbox{\tenrm cost}\!\!_{_{_{^{\Sigma}}}}\!\!\sigma$ strongly
connected.

\noindent {\bf b.}
$\mdoubleH_{n}(|\Sigma|\setminus_o\alpha;\hbox{\tenbf G})
\ \raise3pt\hbox{${{{\lower3pt\hbox{\fivebf A}}}}\atop{\cong}$}\
\mdoubleH_{n}\!(\hbox{\tenrm cost}{\lower1.0pt\hbox{${
_{_{{\!\!\hbox{\tenrm }\!\!_{_{_{\ }}}\!\!\!{{{\Sigma}} }}}}\!
}$}}\tau;\hbox{\tenbf G}) = \mdoubleH_{n}\!(\hbox{\tenrm
cost}{\lower1.0pt\hbox{${ _{_{{\!\!\hbox{\tenrm }\!\!_{_{_{\
}}}\!\!\!{{{\Sigma}} }}}}\! }$}}\tau,\hbox{\tenrm cost}\!\!\!\!
{\lower1.0pt\hbox{${ _{_{{\!\!\hbox{\fiverm Bd}\!\!_{_{_{\
}}}\!\!\!{{{\Sigma}} }}}}\! }$}} \!\tau;\hbox{\tenbf G}) =0$, if
$\alpha\in\hbox{\tenrm Int}(\tau)$ for any $ \tau\in \Sigma.$
{\rm (Cp. the proof of Proposition\ 1 p.\ 11.)}
\qed
\endproclaim

\remark{Note 1} Corollary\ $\!$1.{\bf a} implies that the boundary
of a any manifold is a subcomplex.\hfill\break
\indent\ %
{\bf b} implies that any simplicial manifold is ``ordinary"$\!\!$,
def.\ p.\ 12, and that
$\hbox{\tenrm Bd}_{_{\!{\hbox{\fivebf G}}}} \!\Sigma\not=\emptyset
\Longleftrightarrow
\mdoubleH_{n}(\Sigma;\hbox{\tenbf G})=0.$%
\endremark

\proclaim{Corollary 2}
{\bf i.}
$_{\!}$If $\Sigma$ is an n-pseudomanifold with $\#\hbox{\tenbf
I}\!\geq\!2$ then;
\smallskip
{\hfill$\mdoubleH\!_{_{{{n\!\!}}}}
({\Sigma},\cup\!_{_{\!}}\!\!_{_{\!}}\!\!_{_{\!}}\cup
\!\!\!\!\!\!\!_{_{_{{{j\not=i}}}}}\hbox{\tenbf B}_{_{\!\!{j\!}}}
;\!\hbox{\tenbf G}) \!=\!0 $
\indent and \indent
$ \mdoubleH\!_{_{{{n\!\!}}}}({\Sigma},\hbox{\tenbf
B}_{_{\!{i}}\!};\!\hbox{\tenbf G}) \!=\!0 $.\hfill}

\smallskip
{\bf ii.}
Both $ \mdoubleH\!_{_{{{_{\!}n\!\!}}}}(\Sigma;_{\!}\hbox{\tenbf
Z})$ and
$\mdoubleH\!_{_{{{_{\!}n\!\!}}}}(\Sigma,\hbox{\tenrm
Bd}\Sigma;_{\!}\hbox{\tenbf Z})_{\!}$ equals $0$ or {\bf Z}.
\endproclaim

\demo{Proof}
{\bf i.}
$ {\cup\!{_{\!}}\!_{_{\!}}\!\!_{{\!}}\cup
\!\!\!\!\!\!\!_{_{_{{{j\not=i}}}}}} \! \hbox{\tenbf B}_{_{\!j}}
\raise2pt\hbox{\mysubsetneqq} \hbox{\tenrm
cost}_{_{_{\!{\Sigma}}}}\!\!\sigma\ \hbox{\tenrm for\ some}\
\hbox{\tenbf B}_{_{\!{i}}}\!\hbox{-}\hbox{\tenrm maxidimensional}\
\sigma\!\in\!\hbox{\tenbf B}_{_{\!{i}}}$ { and \ vice\ versa.}
\hfill$\triangleright$

\noindent {\bf ii.}
$\!\![\dim\!\tau\!=\!n] \Rightarrow\!
[[\mdoubleH\!_{_{{{n\!\!}}}}(\Sigma,\hbox{\tenrm
cost}\!_{_{\Sigma}}\!\tau;\hbox{\tenbf G})\!=\! \hbox{\tenbf G}]
\land [\hbox{\tenrm Bd}\Sigma\subset\hbox{\tenrm
cost}\!_{_{\Sigma}}\!\tau]]\ \hbox{\tenrm and} \
\!\Sigma\!\setminus\! \hbox{\tenrm Bd}\Sigma$ strongly connected.
Th.\ $_{\!}$10 and the {\bf LHS} gives the injections \hbox{$
\mdoubleH\!_{_{{{_{\!}n\!\!}}}}(\Sigma;_{\!}\hbox{\tenbf
G})\!\hookrightarrow\!
\mdoubleH\!_{_{{{_{\!}n\!\!}}}}(\Sigma,\hbox{\tenrm
Bd}\Sigma;_{\!}\hbox{\tenbf G})\!\hookrightarrow\!\hbox{\tenbf
G}.$}\qed
\enddemo

\remark{Note 2} $\!\!\Sigma$ Gorenstein $\!\Rightarrow\!$ $\Sigma$
finite.
$_{\!}[\Sigma_{_{^{\!{1\!}}}}\!\ast_{\!}\Sigma_{_{^{\!{2\!}}}}\!_{\!}\not\subset\!\Sigma_{_{^{\!{i\!}}}}$
locally finite]
$\!\Leftrightarrow\!$
[$\Sigma_{_{^{\!{1\!}}}}, \!\Sigma_{_{^{\!{2\!}}}}\!\!\neq\emptyset,\{\emptyset_{_{^{\!{o\!}}}}\!\}$%
\nobreak\ $\!$\hbox{both\nobreak\ finite}].$\!\!$
\endremark

\remark{Note 3}
$\!$For a finite n-manifold $\Sigma$ and a n-submanifold $\Delta$,
put $\hbox{\tensy U}:=\vert\Sigma\vert\setminus\vert\Delta\vert$,
implying that $\vert\hbox{\tenrm Bd}_{_{\!{\hbox{\fivebf
G}}}}\!\Delta\vert\cup\hbox{\tensy U}$
\noindent is the polytope of a subcomplex, $\Gamma$, of $\Sigma$
i.e. $\vert\Gamma\vert\!=\! \vert\hbox{\tenrm
Bd}_{_{\!{\hbox{\fivebf G}}}}\!\Delta\vert\cup\hbox{\tensy U}$,
and $\hbox{\tenrm Bd}_{_{\!{\hbox{\fivebf
G}}}}\!\Sigma\subset\Gamma$,
cp. \cite{21} %
p.\ 427-429. %
Consistency  of Definition\ 3 follows by excision in simplicial
$\mdoubleH$omology since;
$\mdoubleH\!_{_{{{n\!\!}}}}(\Sigma, \hbox{\tenrm
Bd}_{_{\!{\hbox{\fivebf G}}}} \!\Sigma;\hbox{\tenbf
G})\!\hookrightarrow \mdoubleH\!_{_{{{n\!\!}}}}
(\Sigma,\Gamma;\hbox{\tenbf G})
\cong
\mdoubleH\!_{_{{{n\!\!}}}}(\Sigma\setminus\hbox{\tensy
U},\Gamma\setminus\hbox{\tensy U};\hbox{\tenbf G})
\!=
\mdoubleH\!_{_{{{n\!\!}}}}(\Delta,\hbox{\tenrm
Bd}_{_{\!{\hbox{\fivebf G}}}}\!\Delta;\hbox{\tenbf G}).$
E.g., any $\hbox{\tenrm Lk}{\!_{_{\Sigma}}}\!\sigma$ as well as
its (cf. Th.\ 7) iterated cone $\overline{\hbox{\tenrm
st}}{\!_{_{\Sigma}}}\!\sigma\!$, are orientable quasi manifolds
if\nobreak\ $\Sigma_{_{^{\!{\hbox{\fiverm q}}}}}$ {\rm is}.
Moreover, $\delta\!\subset\!\sigma\!\Rightarrow\!
\overline{\hbox{\tenrm st}}{\!_{_{\Sigma}}}\!\delta$
non-orientable  if $\overline{\hbox{\tenrm
st}}{\!_{_{\Sigma}}}\!\sigma$ {\rm is}.
\endremark

\proclaim{Theorem 11.i.a}
\hbox{\rm(cp. \cite{10} %
p.\ 168, \cite{11} %
p.\ 32.)} {$\Sigma$ is a quasi-manifold \underbar{iff}\
$\Sigma\!=\! \bullet \bullet$ or $\Sigma$ is connected and
$\hbox{\tenrm Lk}\!_{_{\Sigma}}\!\sigma$ is a finite
quasi-manifold for all\ \ $\emptyset\neq\sigma\in\Sigma$.}

\noindent
{\bf 1.b.} $\Sigma$ is a homology$_{_{\!\hbox{\fivebf G}}}\!\ n
$-manifold \underbar{\hbox{\tenrm iff}}\
$\Sigma\!=\!\bullet{\!}\bullet\ {or}$
$\mdoubleH_{_{_{\!\!{0\!}}}} \!({{{{{\Sigma}}}}};\!\hbox{\tenbf
G})\!=\!0$\
and $\hbox{\tenrm Lk}\!_{_{\Sigma}}\!\sigma$ is a finite
\hbox{\tenrm ``CM}$\!_{_{\hbox{\fivebf G}}}\!"
\!$-homology$_{_{\!\hbox{\fivebf G}}}\!$
$(n$-$\#\sigma)$-manifold$\ \ \forall\ \
\emptyset\!_{_{^{o}}}\!\neq\!\sigma\!\in\!\Sigma$.

\smallskip
\noindent {\bf ii.}\ \ $\Sigma$ {quasi-manifold} $\Longrightarrow$
$
\hbox{\tenrm Bd}\!
\! \raise0.5pt\hbox{$_{_{_{\hbox{\fivebf G}}}}$}\!\!
_{\!}(\hbox{\tenrm
Lk}\!\!\raise0.5pt\hbox{$_{_{_{\Sigma}}}$}\!\!\sigma\!)\! =
\hbox{\tenrm Lk}\!\!\!\!\! \raise0.5pt\hbox{$_{_{_{\hbox{\fivebf
Bd}\!_{_{\!\hbox{\fivebf G}}}\!\!\!\Sigma}}}$}\!\!\! \sigma
\ {if}\ \sigma\in\hbox{\tenrm Bd}\!\!
\raise0.5pt\hbox{$_{_{_{\hbox{\fivebf G}}}}$}\!\!\!\Sigma \ {and}\
\emptyset\ {else}.
$

\smallskip
\noindent {\bf iii.} {$\Sigma$ is a quasi-manifold\
$\Longleftrightarrow$ $\hbox{\tenrm Lk}\!_{_{\Sigma}}\!\sigma$ is
a pseudomanifold} $\forall\  \sigma\in\Sigma.$
\endproclaim

\demo{Proof. {\bf i}}
A simple check confirms all our claims for $\dim\Sigma\leq1$, cf.
Note\ 2 p.\ 23. %

So, assume $\dim\Sigma\geq2$ and note that
\hbox{
$\!\sigma\!\in\!\hbox{\tenrm Bd}\!\!
\raise0.5pt\hbox{$_{_{_{\hbox{\fivebf G}}}}$}\!\!\Sigma $
\ \underbar{iff}\ \
$\hbox{\tenrm Bd}\!
\! \raise0.5pt\hbox{$_{_{_{\hbox{\fivebf G}}}}$}\!\!
(_{\!}\hbox{\tenrm Lk}
\!\!\raise0.5pt\hbox{$_{_{_{\Sigma}}}$}\!\!\sigma\!)\!\not=\!\emptyset$.}

\smallskip
\noindent {\bf i.a.} $\!(\!\Leftarrow\!$)
That $\hbox{\tenrm Lk}\!_{_{\Sigma}}\!\sigma$, with
$\dim\hbox{\tenrm Lk}\!_{_{\Sigma}}\!\sigma\!=\!0$, is a
quasi-0-manifold implies definition condition 2$\beta\
\hbox{\tenrm p}.\ 23$ and since the other ``links" are all
\underbar{connected} condition 2$\gamma$ follows.\nobreak\
$\triangleright$

\smallskip
\noindent $(\Rightarrow$) Definition condition 2$\beta$ p.\ 23 %
implies that 0-dimensional links are $\bullet$ or $\bullet
\bullet$ while Eq.\ {\bf I} p.\ 30 gives the necessary
connectedness of `links of links'${{\!}}$, cp.  Lemma 2
p.\nobreak\ 16.\nobreak\ $\triangleright$

\smallskip
\noindent{\bf i.b.}
Lemma\ $\!{_{\!}}$2.\ $\!${\bf ii} above plus Proposition\ 1 p.\
18 and Eq.\ {\bf I} p.\ 30
\hfill$\triangleright$

\medskip
\noindent{\bf ii.} Pureness is a local property, i.e. $\Sigma$\
{pure}\ $\Longrightarrow$ \hbox{\tenrm Lk}$\!_{_{\Sigma}}\!\sigma$
pure. Put $n:=\dim{\Sigma}$. Now;

\smallskip
\noindent
$\epsilon\!\in\!\hbox{\tenrm
Bd\!\!\raise0.5pt\hbox{$_{_{_{\hbox{\fivebf G}}}}$}\!\!
(Lk}\!_{_{\Sigma}}\!\sigma) \!\Leftrightarrow\! 0\!=\!
\mdoubleH\!\!\!\!\!\!\!\!\!\!_{_{_{\!{n\lower1pt\hbox{-}\hbox{\fivebf\#}\sigma
\lower1pt\hbox{-}\hbox{\fivebf\#}\epsilon}}}} \!\!(\hbox{\tenrm
Lk}\!\!\!\!_{_{_{\hbox{\fivebf
Lk}\!_{_{\Sigma}}\!\!\sigma}}}\!\!\!\! \epsilon;\hbox{\tenbf G})
\!=\! \big[{{\hbox{\spaceskip2.7pt\sevenrm Eq.\ \hbox{\sevenbf I}
p.\ 30
}} \atop{ \epsilon\in\hbox{\sevenrm
Lk}\!_{_{\Sigma}}\!\sigma}}\big] \!=\!
\mdoubleH\!\!\!\!\!\!\!\!\!_{_{_{\!{n\lower1pt\hbox{-}
\hbox{\fivebf\#}(\sigma\cup\epsilon)}}}} \!\!\hbox{\tenbf
(}\hbox{\tenrm Lk}\!\!_{_{_{{{\Sigma}}}}}\!\!
(\sigma\cup\epsilon);\hbox{\tenbf G}\hbox{\tenbf )}\ \hbox{\tenrm
and}\ \epsilon\!\in\!\hbox{\tenrm Lk}\!_{_{\Sigma}}\!\sigma.$
So;

\noindent
$\epsilon\!\in\!\hbox{\tenrm
Bd\!\!\raise0.5pt\hbox{$_{_{_{\hbox{\fivebf G}}}}$}\!\!
(Lk}\!_{_{\Sigma}}\!\sigma)
\!\Leftrightarrow\!
[\sigma\cup\epsilon \!\in\!\hbox{\tenrm
Bd}\!\!\raise0.5pt\hbox{$_{_{_{\hbox{\fivebf
G}}}}$}\!\!{{\Sigma}}\ _{\!}\hbox{\tenrm and}\
\epsilon\!\in\!\hbox{\tenrm Lk}\!_{_{\Sigma}}\!\sigma]
\!\Leftrightarrow\!
[\sigma_{\!}\cup_{\!}\epsilon \!\in\!\hbox{\tenrm
Bd}\!\!\raise0.5pt\hbox{$_{_{_{\hbox{\fivebf
G}}}}$}\!\!{{\Sigma}}\ \hbox{\tenrm and}\
\sigma\cap\epsilon_{\!}=_{\!}\emptyset]
\!\Leftrightarrow\!
[\epsilon\!\in\!\hbox{\tenrm Lk}\!\!\!\!\!
\raise0.5pt\hbox{$_{_{_{\hbox{\fivebf Bd}\Sigma}}}$}
\!\!\!\sigma].$

\smallskip
\noindent {\bf iii.} ($\Rightarrow$)
Lemma\ 2.\ $\!${\bf i} and {\bf i.{\bf a}} above.
{$(\Leftarrow)$} All links are connected, except for
$\bullet\bullet$.
\qed
\enddemo

\proclaim{Corollary 1}
For any quasi-$n$-manifold $\Sigma$ except infinite $1$-circles;

\medskip
\centerline{$ \dim\hbox{\tenbf B}_{_{{\!{j}}}}\!\!\ge\!_{\!}
n_{_{{\!{\ }}}}\!\!\!-\!2 \ \Longrightarrow\ \dim\hbox{\tenbf
B}_{_{{\!{j}}}}\!\!=\! _{\!}n_{_{{\!{\ }}}}\!\!\!-\!1 . $
}%
\endproclaim

\demo{Proof}
Check $n \!\le\! 1$. Now; assume $n \!\ge\! 2$. If
$\dim\sigma\!=\dim\hbox{\tenbf B}_{_{{\!{j}}}}\!\!= n_{_{{\!{\
}}}}\!\!-2 $ and
$\sigma\!_{_{^{\!{\ }}}}\!\!\in\! \hbox{\tenbf B}_{_{{\!{j}}}}\!$
then;
$ \hbox{\tenrm Lk}\!\!\!\!\! \raise0.5pt\hbox{$_{_{_{\hbox{\fivebf
Bd}\!_{_{\!\hbox{\fivebf G}}}\!\!\!\Sigma}}}$}\!\!\!
\sigma\!_{_{^{\!{\ }}}}\!\!=\! \hbox{\tenrm Lk}\!\!
\raise0.5pt\hbox{$_{_{_{\hbox{\fivebf B}_{_{{\!{j}}}}\!}}}$}\!\!\!
\sigma\!_{_{^{\!{\ }}}} \!\!= \!\{{\emptyset}_{_{^{\!o}}}\!_{\!}\}
$.

\noindent By Th.\ 11{\bf ii};
$ \hbox{\tenrm Lk}\!\!\!\!\! \raise0.5pt\hbox{$_{_{_{\hbox{\fivebf
Bd}\!_{_{\!\hbox{\fivebf G}}}\!\!\!\Sigma}}}$}\!\!\!
\sigma\!_{_{^{\!{\ }}}} = \hbox{\tenrm Bd}\!_{_{\hbox{\fivebf
G}}}\!(\hbox{\tenrm Lk}\!_{_{\Sigma}}\sigma\!_{_{^{\!{ }}}})
=
\big[{{\hbox{\eightrm{Lk}}\!_{_{\Sigma}}\!\sigma\!_{_{^{\!{\
}}}}\hbox{\eightrm is, by Th. 11\hbox{\eightbf i},  a finite
quasi\lower0.6pt\hbox{-}}}\atop {\hbox{\eightrm
1\lower0.6pt\hbox{-}manifold i.e. (a circle or) a line.}}}\big]
=
(\emptyset\ \hbox{\tenrm or})\ \bullet\bullet. $\ \ Contradiction!
\qed
\enddemo

Denote $\Sigma$ by
$\Sigma\!\!\raise0.9pt\hbox{$_{_{_{\!\hbox{\fiverm ps}}}}$}\!$,
$\!\Sigma_{_{^{\!\hbox{\fiverm q}}}}\!$ and
$\Sigma_{_{^{{_{\!}}\hbox{\fiverm h}\!}}}$ when it's assumed to be
a $\hbox{\tenrm pseudo\ \!\hbox{-,}}$ quasi- resp. a homology
manifold.
Note also that;  $\!\sigma\in\hbox{\tenrm Bd}\!\!
\raise0.5pt\hbox{$_{_{_{\hbox{\fivebf G}}}}$}\!\!\!\Sigma
_{_{^{\!\hbox{\fiverm q}}}}$%
$\Longleftrightarrow$
$
\hbox{\tenrm Bd}\!
\! \raise0.5pt\hbox{$_{_{_{\hbox{\fivebf G}}}}$}\!\!
_{\!}(\hbox{\tenrm
Lk}\!\!\raise0.5pt\hbox{$_{_{_{\Sigma}}}$}\!\!\sigma\!)\! =
\hbox{\tenrm Lk}\!\!\!\!\! \raise0.5pt\hbox{$_{_{_{\hbox{\fivebf
Bd}\!_{_{\!\hbox{\fivebf G}}}\!\!\!\Sigma}}}$}\!\!\!
\sigma\not=\emptyset$
by Th.\ 11.ii.

\proclaim{Corollary 2.i}
Each boundary component
$\hbox{\tenbf B}_{_{\!{j\!}}}$, $j\!_{\!}\in\! \hbox{\tenbf I}$,
of
$\Sigma_{_{^{\!\hbox{\fiverm q}}}}\!$ is a pseudomanifold.

\smallskip
\noindent \hbox{\tenbf ii.}
If ${{\Sigma}}_{_{^{\!\hbox{\fiverm q}}}} \!$ is finite with
$ \hbox{\tenrm Bd}\!\raise0.8pt\hbox{$_{_{_{\!\hbox{\fivebf
G}}}}$}\!\!{{\Sigma}}_{_{^{\!\hbox{\fiverm q}}}} \!\!=\!
{\cup\!{_{\!}}\!_{_{\!}}\!\!_{{\!}}\cup \!\!\!\!\!\!\!_{_{_{{{j\in
\lower0.5pt\hbox{\sixbf I}}}}}}} \!\hbox{\tenbf B}_{_{\!\!^{j}}}
\!$
and
$-1\!\le\!\dim\hbox{\tenbf B}_{_{{\!{i}\!}}}\!
<\!\dim\Sigma\!-\!1$ for some $i\!\in\! \hbox{\tenbf I}$ then
$\Sigma_{_{^{\!\hbox{\fiverm q}}}}\!$
is nonorientable$_{_{\!\hbox{\fivebf G}}}\!$.

\smallskip
\noindent \hbox{\tenbf iii.} For any
orientable$_{_{\!\hbox{\fivebf G}}}\!\!$ quasi-$n$-manifold
$\Sigma_{_{^{\!\hbox{\fiverm q}\!}}}\!$ each boundary component
$\hbox{\tenbf B}_{_{\!{i}\!}}\!:=\!\hbox{\tenbf
B}\!\raise0.6pt\hbox{$_{_{\!{_{\hbox{\fivebf
G}\!,_{\!}i}}}}$}\!\!\!\!\!^{^{_{\Sigma\!}}}\ \!\!
\nsubseteqq\!
\{\emptyset{_{_{^{\!o\!}}}}\!\}
$
is an orientable $(n-1)$-pseudomanifold without boundary.

\smallskip
\noindent \hbox{\tenbf iv.}
$\mdoubleH\!\!{_{_{n\lower1pt\hbox{-}i}}}\! (\Sigma,\hbox{\tenrm
Bd}\Sigma\!\!_{_{\hbox{\fiverm ps}}}\!;\hbox{\tenbf
G}^{\!\prime})\!\!=\! \mdoubleH\!\!{_{_{n\lower1pt\hbox{-}i}}}\!
(\Sigma,\hbox{\tenrm Bd}\!\raise0.8pt\hbox{$_{_{_{\!\hbox{\fivebf
G}}}}$}\!\!\Sigma\!_{_{\hbox{\fiverm q}}};\hbox{\tenbf
G}^{\!\prime})\!\!=\! \mdoubleH\!\!{_{_{n\lower1pt\hbox{-}i}}}\!
(\Sigma,\hbox{\tenrm Bd}\!\raise0.8pt\hbox{$_{_{_{\!\hbox{\fivebf
G}}}}$}\!\!\Sigma_{_{\!\hbox{\fiverm h}}};\hbox{\tenbf
G}^{\!\prime}),\  i=0,1$
even if $\hbox{\tenbf G}\!\neq\!\hbox{\tenbf G}^{\!\prime}.$

\noindent
Orientability is independent of $\ \!\hbox{\tenbf
G}^{\!\prime}\!$, as long as $ \hbox{\tenrm Tor}_1^{\hbox{\fivebf
Z}}\bigl(\hbox{\tenbf Z}_{_{\!2}},\!\hbox{\tenbf G}^{\prime}
\bigr)\!\ne\hbox{\tenbf G}^{\prime} $
\hbox{\tenrm (Lemma\ $\!{_{\!}}$1.\ $\!$\hbox{\tenbf ii}
p.\ 24).}
Moreover,
$\hbox{\tenrm Bd}\!\!\lower3.4pt\hbox{{{{\fivebf
G}}}}\!_{\!}\Sigma\!_{_{\hbox{\sixrm q}}}
\!\!=
\hbox{\tenrm Bd}\!\!\lower3.4pt\hbox{{{{\fivebf
G}}}}\!_{\!}\Sigma_{_{\hbox{\sixrm h}}}\!
$
always, while
$(\hbox{\tenrm Bd}\Sigma\!\!_{_{\hbox{\sixrm
ps}}})^{^{{\!{n\!\lower1.0pt\hbox{-}\!\hbox{\fiverm 1}}}}}
\!\!\!\!=
(\hbox{\tenrm Bd}\!\!\lower3.4pt\hbox{{{{\fivebf
G}}}}\!_{\!}\Sigma\!_{_{\hbox{\sixrm
q}}})^{^{{\!{n\!\lower1.0pt\hbox{-}\!\hbox{\fiverm 1}}}}} and
$
$(\hbox{\tenrm Bd}\Sigma\!\!_{_{\hbox{\sixrm
ps}}})^{^{{\!{n\!\lower1.0pt\hbox{-}\hbox{\fiverm 2}}}}}
\!\!\!\!=
(\hbox{\tenrm Bd}\!\!\lower3.4pt\hbox{{{{\fivebf
G}}}}\!_{\!}\Sigma\!_{_{\hbox{\sixrm
q}}})^{^{{\!{n\!\lower1.0pt\hbox{-}\hbox{\fiverm 2}}}}}
$
except for infinite $1$-circles in which case
$(\hbox{\tenrm Bd}\Sigma\!\!_{_{\hbox{\sixrm
ps}}})^{^{{\!{n\!\lower1.0pt\hbox{-}\hbox{\fiverm 2}}}}}
\!\!=\emptyset\not=\{\emptyset_{_{\hbox{\sixrm o}}}\}=
(\hbox{\tenrm Bd}\!\!\lower3.4pt\hbox{{{{\fivebf
G}}}}\!_{\!}\Sigma\!_{_{\hbox{\sixrm
q}}})^{^{{\!{n\!\lower1.0pt\hbox{-}\hbox{\fiverm 2}}}}}$.

\medskip%
\noindent
\ $\!$\hbox{\tenbf v.}\
$ \hbox{\tenrm Tor}\!\lower1.0pt\hbox{$_{^{_{\hbox{\fivebf
1}}}}$}\!\!\!{^{^{_{\hbox{\fivebf Z}}}}\!}\!(\hbox{\tenbf
Z}_{_{^{\!2}}}\!,\!\hbox{\tenbf G})\!=0$
$\!\Longrightarrow\!$
$\hbox{\tenrm Bd}\!\raise0.6pt\hbox{$_{_{_{\!\hbox{\fivebf
G}}}}$}\!\!\!\Sigma_{_{^{\!{\hbox{\fiverm q}}}}} \!= \hbox{\tenrm
Bd}\!\raise0.5pt\hbox{$_{_{_{\!\hbox{\fivebf
Z}}}}$}\!\!\Sigma_{_{^{\!{\hbox{\fiverm q}}}}}.$

\smallskip
\noindent
{\bf vi.}
{$_{\!} \hbox{\tenrm Bd}\Sigma\!\!_{_{\hbox{\fiverm ps}}} \!\!=
\!\hbox{\tenrm Bd}\!\raise0.5pt\hbox{$_{_{_{\!\hbox{\fivebf
Z}_{^{^{\!2}}}}}}$}\!\!\!\!\Sigma_{_{^{\hbox{\fiverm q}}}}
\!\!\subseteq \!\hbox{\tenrm
Bd}\!\raise0.8pt\hbox{$_{_{_{\!\hbox{\fivebf
G}}}}$}\!\!\Sigma_{_{^{\!q}}} \!\!\subseteq \!\hbox{\tenrm
Bd}\!\raise0.5pt\hbox{$_{_{_{\!\hbox{\fivebf
Z}}}}$}\!\!\Sigma_{_{^{\hbox{\fiverm q}}}}\! $ with}
{equality if
$ \hbox{\tenrm Bd}\!\!\raise0.6pt\hbox{$_{_{_{\hbox{\fivebf Z
}}}}$}\!\!{{\Sigma}}\!=\!\emptyset$
or
$ \dim\!\hbox{\tenbf B}\!\raise0.6pt\hbox{$_{_{\!{_{\hbox{\fivebf
Z}\!,_{\!}j}}}}$}\!\!\!\!\!^{^{_{\Sigma\!\!}}}\ \!\! =\!n\!-\!1\
\forall j\!\in\!\hbox{\tenbf I}, $
except if $\Sigma$ is infinite and
$ \hbox{\tenrm Bd}\hbox{$_{_{\!^{\hbox{\fivebf
Z}_{_{^{\!2}}}}}}$}\!\!\!\!\Sigma\!_{_{^{q}}}\!
\!=\!\{\emptyset_{_{^{\!o}}}\!_{\!}\} \!\neq\!\emptyset
\!=\!\hbox{\tenrm Bd}\Sigma\!\!_{_{\hbox{\fiverm ps}}} {{\!}}$}
{\rm (by Lemma\ $\!{_{\!}}$1.\ $\!$\hbox{\tenbf
i}\raise1pt\hbox{\eightrm+}\hbox{\tenbf ii} p.\ $\!$ plus Th.\ 5
p.\ 10
since
$\emptyset\!_{_{^{o}}}\!\in\hbox{\tenrm
Bd}\!\!\lower3.4pt\hbox{{{{\fivebf
G}}}}\!_{\!}\Sigma\!\lower2.4pt\hbox{{{{\fiverm
q}}}}{\!}\neq\emptyset$
if $\Sigma\!\lower2.4pt\hbox{{{{\fiverm q}}}}_{\!} $ is
infinite.).}
{\rm {If $\Sigma$}$_{_{^{\hbox{\fiverm h}}}}\!\!$
{\hbox{\tenrm Gorenstein}}\hbox{$_{_{\!^{\hbox{\fivebf
Z}_{_{^{\!2}}}}}}\!\!\!$}
{then};
$\hbox{\tenrm Bd}\hbox{$_{_{\!^{\hbox{\fivebf
Z}_{_{^{\!2}}}}}}$}\!\!\!\!\Sigma\!_{_{^{\hbox{\fiverm h}}}}\!\!$
=
$\hbox{\tenrm Bd}\Sigma\!_{_{^{\hbox{\fiverm ps}}}}\!. $
}%
\endproclaim

\normalbaselines
\demo{Proof}
\noindent \hbox{\tenbf i.} The claim is 
true if $\dim \Sigma\le 1$ and \hbox{assume it's true for
dimensions $\le_{\!} n_{\!}-_{\!}1$.}
$\alpha$ and $\gamma$ are true by definition of $\hbox{\tenbf
B}_{_{\!{i}}}\!$ so only $\beta$ remains.
If $\dim\hbox{\tenbf B}_{_{\!{i}}}\!\!=\!m$ and $\sigma\!\in\!
\hbox{\tenbf B}\raise0.5pt\hbox{$_{_{\!^{i}}}$}\!$
with $\dim \sigma \!=\! m_{_{\!}}-_{\!}1$,
then
$\hbox{\tenrm Bd}\!\!_{_{_{\hbox{\fivebf G}}}}\!\! ({\hbox{\tenrm
Lk}\!_{_{_{\!\Sigma\!\!}}}\!\sigma})\!_{_{_{\!\hbox{\fiverm
q}}}}\!
\!=\!
{{{\hbox{\tenrm
Lk}_{_{{\!\!\!\!\!\!\!\lower1.3pt\hbox{$_{^{{\hbox{\sixrm
Bd}\!\!_{_{_{\hbox{\fivebf G}}}}\!\!\!\Sigma}}}$}}
}}\!\!\!{\sigma\!^{_{^{\hbox{\fivebf \ }\! }}}}}}}\! \!=\!
{{{\hbox{\tenrm Lk}_{_{{\!\!\!\lower1.3pt\hbox{$_{^{{\hbox{\fiverm
B}_{^{_{\!i}}}}}}$}} }}\!\!\!\sigma}}} $
where the r.h.s. is zero dimensional and so, strongly connected,
implying, by
the induction assumption, that the sole component on the l.h.s. is
a \hbox{$0$-pseudomanifold i.e. $\bullet$ or $\bullet\bullet$.}
\hfill$\triangleright$

\smallskip%
\noindent \normalbaselines \hbox{\tenbf ii.}
$\dim\hbox{\tenbf B}_{_{\!{i}}}\!\!<\!{n}_{_{{\!{\Sigma}}}}\!-\!2$
by Cor.$\ \!1$ 
giving the 2:nd equality
$\hbox{\tenrm and\ Th.}\ 10\ \!\hbox{\tenbf a}\ \hbox{\tenrm
gives\ the}\!$
arrow in;

\bigskip
$ \mdoubleH\!_{_{{{n}_{_{{\!{\ }}}}\!\!\!\!}}}
({\Sigma},\hbox{\tenrm Bd}\!
\raise0.8pt\hbox{$_{_{_{\!\hbox{\fivebf
G}}}}$}\!\!{{\Sigma}};\hbox{\tenbf G}) \!=
\mdoubleH\!_{_{{{n}_{_{{\!{\ }}}}\!\!\!\!}}}
({\Sigma},\hbox{\tenbf
B}_{_{\!{i}}}\!\cup({\cup\!\!_{_{\!}}\!\!_{_{\!}}\cup}
\!\!\!\!\!_{_{_{{{j\not=i}}}}}\hbox{\tenbf
B}_{_{\!{j}}});\hbox{\tenbf G}) \!= \mdoubleH\!_{_{{{n}_{_{{\!{\
}}}}\!\!\!\!}}} ({\Sigma},{\cup\!\!_{_{\!}}\!\!_{_{\!}}\cup}
\!\!\!\!\!_{_{_{{{j\not=i}}}}}\hbox{\tenbf B}_{_{\!{j}}}
;\hbox{\tenbf G})
\hookrightarrow
$

\bigskip
\hfill\hfill
$ \hookrightarrow
\mdoubleH\!_{_{{{n}_{_{{\!{\ }}}}\!\!\!\!}}}
({\Sigma},\hbox{\tenrm
cost}_{_{_{\!{\Sigma}}}}\!\!\!\sigma;\hbox{\tenbf G})\!=\!0$
\ \
if\ \ $\sigma_{{\!}}\in_{{\!}}\hbox{\tenrm
Bd}\!\!_{_{_{\hbox{\fivebf G}}}}\!\!\Sigma\
\!^{_{_{\setminus}}}{{\cup_{_{\!}}\!_{^{\!}}\!_{{\!}}\!_{^{\!}}\cup}}
\!\!\!\!_{_{\!}}\!_{_{_{{{j\not=i}}}}}\hbox{\tenbf
B}_{_{\!\!{j\!}}}\!^{^{\!_{\ }}}.$ \hfill$\triangleright$

\medskip
\noindent \hbox{\tenbf iii.} We can w.l.o.g. assume that
$\Sigma_{_{^{\hbox{\fiverm q}}}}$ is finite.
If $n_{_{\!}}-_{\!}2\!=\!\dim \!\sigma$ and
$\sigma\!_{\!}\in\!_{\!}
\hbox{\tenbf B}\raise0.5pt\hbox{$_{_{\!^{j}}}$}\!
$ then
$
\dim_{\!}\hbox{\tenrm Lk} \!\raise0.5pt\hbox{$_{_{_{\Sigma}}}$}\!
\!\!\sigma\!_{_{^{\!{\ }}}}\!\!=\!1
$
and so,\
$
\hbox{\tenrm Lk}\!\! \raise0.5pt\hbox{$_{_{_{\hbox{\fivebf
B}_{_{{\!{j}}}}\!}}}$}\!\!\! \sigma\!_{_{^{\!{\ }}}}\!
$
=
$ \hbox{\tenrm Lk}\!\!\!\!\! \raise0.5pt\hbox{$_{_{_{\hbox{\fivebf
Bd}\!_{_{\!\hbox{\fivebf G}}}\!\!\!\Sigma}}}$}\!\!\!
\sigma\!_{_{^{\!{\ }}}}
$
=
$\hbox{\tenrm Bd}\!\!\raise0.5pt\hbox{$_{_{_{\hbox{\fivebf
G}}}}$}\!\!\!$
$
\!(\hbox{\tenrm Lk}\!\raise0.5pt\hbox{$_{_{_{\Sigma}}}$}
\!\!\sigma\!_{_{^{\!{\ }}}}\!_{\!})
=
$
${\bullet\bullet}$\ \ gives \
$\hbox{\tenrm Bd}_{\!}(_{\!}\hbox{\tenbf
B}\raise0.5pt\hbox{$_{_{\!^{j}}}$}\!)_{_{^{\hbox{\fiverm
ps}}}}\!=\ \! \emptyset.$\ \

Now, from the boundary component definition we learn that,
$\dim\hbox{\tenbf B}_{_{\!{i}}}\!= {n}_{_{{\!{\ }}}}\!\!-\!1\
\forall\ i\in\hbox{\tenbf I}\Longrightarrow\dim\hbox{\tenbf
B}_{_{\!{j}}}\cap \hbox{\tenbf B}_{_{\!{i}}}\!\!<\!n-2\ \forall\
j\!\neq\! i,$
which gives us;

\medskip
$ \dots\longrightarrow
\!\!\!\!\!\underbrace{ \mdoubleH\!_{_{{{n}_{_{{\!{\ }}}}\!\!}}}\!
(\hbox{\tenrm Bd}\!\raise0.8pt\hbox{$_{_{_{\!\hbox{\fivebf
G}}}}$}\!\!{{\Sigma}}, {\cup\!\!\!\!_{_{\!}}\cup}
{\!\!\!\!\!_{_{_{{{j\not=i}}}}}}\!\hbox{\tenbf B}_{_{\!{j}}}
;\!\hbox{\tenbf G})}
_{=\ \!0\ \!\hbox{\eightrm for\ \!dimensional\ \!reasons.}}
\vbox{\moveleft0.2cm\hbox{$\longrightarrow$}}
\underbrace{ \mdoubleH\!_{_{{{n}_{_{{\!{\ }}}}\!\!}}}\!
({{\Sigma}}, {\cup\!\!\!\!_{_{\!}}\cup}
{\!\!\!\!\!_{_{_{{{j\not=i}}}}}}\!\hbox{\tenbf B}_{_{\!{j}}}
;\!\hbox{\tenbf G})}
_{=\ \!0\ \! \hbox{\eightrm by\ \!Cor.\ \!\hbox{\eightrm 2}\ p.\
\!25.} } \longrightarrow $

\bigskip
\hfill\hfill
$
\longrightarrow
\underbrace{ \mdoubleH\!_{_{{{n}_{_{{\!{\
}}}}\!\!}}}\!({{\Sigma}}, \hbox{\tenrm
Bd}\!\raise0.8pt\hbox{$_{_{_{\!\hbox{\fivebf G}}}}$}\!\!{{\Sigma}}
;\hbox{\tenbf G})}
_{=\hbox{\eightbf G}}
\longrightarrow\!
\!\!\!\!\!\!\!\!\!\!\underbrace{ \mdoubleH\!\!_{_{^{
\!{n}_{_{{\!{\ }}}}\!\!\lower1pt\hbox{\tenbf-}1}}}\!( \hbox{\tenrm
Bd}\!\raise0.8pt\hbox{$_{_{_{\!\hbox{\fivebf
G}}}}$}\!\!{{\Sigma}}, {\cup\!\!\!\!_{_{\!}}\cup}
{\!\!\!\!\!_{_{_{{{j\not=i}}}}}}\!\hbox{\tenbf B}_{_{\!{j}}}
;\hbox{\tenbf G})}
_{=\ \!\mminidoubleH_{\!\!n\hbox{\fivebf-}1\!}(\hbox{\eightbf
B}_{_{\!{i}}} ;\hbox{\eightbf G})\ \!\hbox{\eightrm for\
\!dimensional\ \!reasons.}}
\!\vbox{\moveleft0.8cm\hbox{$\longrightarrow\cdots. $}}
\hfill \triangleright
$

\noindent \normalbaselines
\hbox{\tenbf iv.}
\hbox{$_{\!}$Corollary\ 1 \ \ and \ \ [$
\sigma\!\!\in\!\Sigma^{^{{\!{n\lower1.0pt\hbox{-}\!\hbox{\fiverm
1} \lower0.7pt\hbox{}}}}} \!_{\!}\cap\hbox{\tenrm Bd}\Sigma $ \ \
{\underbar{iff}} \ \
${{{\hbox{\tenrm
Lk}\raise0.5pt\hbox{$\!\!_{_{_{{\Sigma}_{_{^{\!\hbox{\tenrm }}}}}
}}$}\!\!{\sigma^{_{^{\hbox{\fivebf }}}}\!_{\!}}}}}=\bullet]$.
}\hfill$\triangleright$

\noindent \hbox{\tenbf v.}
$\hbox{\tenrm Bd}\raise0.5pt\hbox{$\!_{_{_{\hbox{\fivebf
G}}}}$}\!\!\!\Sigma\!_{_{^{\hbox{\fiverm q}}}} \!\!\not\ni_{\!}
\sigma\!\in\! \hbox{\tenrm
Bd}\raise0.5pt\hbox{$\!_{_{_{\!{\hbox{\fivebf
Z}}}\!}}$}\!\!\Sigma\!_{_{^{\hbox{\fiverm q}}}}\! $
$ \Leftrightarrow\ \!
\mdoubleH \vbox{\moveleft0.40cm\hbox{$
_{_{_{\!{n \lower1.0pt\hbox{-} \hbox{\fivebf\#}\sigma}}}}
$}}
\!\!\!\!
(({\hbox{\tenrm
Lk}\!_{_{_{\!\Sigma\!\!}}}\!\sigma}\raise1.5pt\hbox{\eightbf
{\char"29}}\!_{_{^{\hbox{\fiverm q}}}}\! ;_{\!}\hbox{\tenbf G})
\!\not= 0 {\!}={\!}
\mdoubleH
\vbox{\moveleft0.40cm\hbox{$
_{_{_{\!{n \lower1.0pt\hbox{-} \hbox{\fivebf\#}\sigma}}}}
$}}
\!\!\!\!
(({\hbox{\tenrm Lk}\!_{_{_{\!\Sigma\!\!}}}\!\sigma}
\raise1.5pt\hbox{\eightbf {\char"29}}\!_{_{^{\hbox{\fiverm q}}}}\!
;_{\!}\hbox{\tenbf Z})
\Leftrightarrow
$%
$ \hbox{\tenrm Bd}\!\!\raise0.5pt\hbox{$_{_{_{\hbox{\fivebf
G}}}}$}\!\!\!$
$
\!(\hbox{\tenrm Lk}\!\raise0.5pt\hbox{$_{_{_{\Sigma}}}$}
\!\!\sigma\!_{_{^{\!{\ }}}}\!_{\!})\!_{_{^{\hbox{\fiverm q}}}}\!
=\emptyset\not=\
$%
\hbox{\tenrm Bd$\!\!$}%
\raise0.5pt\hbox{$_{_{_{\hbox{\fivebf Z}}}}$}%
$\!
\!(\hbox{\tenrm Lk}\!\raise0.5pt\hbox{$_{_{_{\Sigma}}}$}
\!\!\sigma\!_{_{^{\!{\ }}}}\!_{\!})\!_{_{^{\hbox{\fiverm q}}}}\!.
$

$\hbox{\tenbf iv}\Longrightarrow$
$
n-\hbox{\eightrm\#}_{\!}\sigma-3\ge m:=\dim{\hbox{\tenrm Bd$\!\!$}%
\raise0.5pt\hbox{$_{_{_{\hbox{\fivebf Z}}}}$}%
\!
\!(\hbox{\tenrm Lk}\!\raise0.5pt\hbox{$_{_{_{\Sigma}}}$}
\!\!\sigma\!_{_{^{\!{\ }}}}\!_{\!})\!_{_{^{\hbox{\fiverm q}}}}\!
}$
$\Longrightarrow$

\medskip
$\Longrightarrow$
$
0\neq \mdoubleH
\vbox{\moveleft0.40cm\hbox{$
_{_{_{\!{n \lower1.0pt\hbox{-} \hbox{\fivebf\#}\sigma}}}}
$}}
\!\!\!\!
(({\hbox{\tenrm
Lk}\!_{_{_{\!\Sigma\!\!}}}\sigma}\raise1.5pt\hbox{\eightbf
{\char"29}}\!_{_{^{\hbox{\fiverm q}}}}\! ;_{\!}\hbox{\tenbf G})
\!=
$
$
\underbrace{ \mdoubleH \vbox{\moveleft0.40cm\hbox{$
_{_{_{\!{n \lower1.0pt\hbox{-} \hbox{\fivebf\#}\sigma}}}}
$}}
\!\!\!\!
(({\hbox{\tenrm
Lk}\!_{_{_{\!\Sigma\!\!}}}\sigma}\raise1.5pt\hbox{\eightbf
{\char"29}}\!_{_{^{\hbox{\fiverm q}}}}\! ;_{\!}\hbox{\tenbf Z})
}
_{=\ \!0\ \!
\hbox{\eightrm by\ \!\hbox{\eightrm assumption}.}}%
\!\!\! \otimes\hbox{\tenbf G}
\ \oplus
$
$
\hbox{\tenrm Tor}\!\lower1.0pt\hbox{$_{^{_{\hbox{\fivebf
1}}}}$}\!\!\!{^{^{_{\hbox{\fivebf Z}}}}\!}\!
(\phantom{\ }
\mdoubleH \vbox{\moveleft0.40cm\hbox{$
_{_{_{\!{n \lower1.0pt\hbox{-}
\hbox{\fivebf\#}\sigma}\lower1.0pt\hbox{-}1}}}
$}}
\!\!\!\!
(({\hbox{\tenrm
Lk}\!_{_{_{\!\Sigma\!\!}}}\sigma}\raise1.5pt\hbox{\eightbf
{\char"29}}\!_{_{^{\hbox{\fiverm q}}}}\! ;_{\!}\hbox{\tenbf
Z}),\hbox{\tenbf G})
= $

=$\big[{{\hbox{\eightrm dimensional}}\atop{{\hbox{\eightrm
reasons.}}}}\big]$
$=
\hbox{\tenrm Tor}\!\lower1.0pt\hbox{$_{^{_{\hbox{\fivebf
1}}}}$}\!\!\!{^{^{_{\hbox{\fivebf Z}}}}\!}\!
(\phantom{\ }
\mdoubleH \vbox{\moveleft0.40cm\hbox{$
_{_{_{\!{n \lower1.0pt\hbox{-}
\hbox{\fivebf\#}\sigma}\lower1.0pt\hbox{-}1}}}
$}}
\!\!\!\!
(({\hbox{\tenrm
Lk}\!_{_{_{\!\Sigma\!\!}}}\sigma}\raise1.5pt\hbox{\eightbf
{\char"29}}\!_{_{^{\hbox{\fiverm q}}}}\!,
\hbox{\tenrm Bd$\!\!$}%
\raise0.5pt\hbox{$_{_{_{\hbox{\fivebf Z}}}}$}%
\!
\!(\hbox{\tenrm Lk}\!\raise0.5pt\hbox{$_{_{_{\Sigma}}}$}
\!\!\sigma\!_{_{^{\!{\ }}}}\!_{\!})\!_{_{^{\hbox{\fiverm q}}}}\!
;_{\!}\hbox{\tenbf Z}),\hbox{\tenbf G})
= 0$. Contradiction! - since the torsion module of \
$
\mdoubleH \vbox{\moveleft0.40cm\hbox{$
_{_{_{\!{n \lower1.0pt\hbox{-}
\hbox{\fivebf\#}\sigma}\lower1.0pt\hbox{-}1}}}
$}}
\!\!\!\!
(({\hbox{\tenrm
Lk}\!_{_{_{\!\Sigma\!\!}}}\sigma}\raise1.5pt\hbox{\eightbf
{\char"29}}\!_{_{^{\hbox{\fiverm q}}}}\!,
\hbox{\tenrm Bd$\!\!$}%
\raise0.5pt\hbox{$_{_{_{\hbox{\fivebf Z}}}}$}%
\!
\!(\hbox{\tenrm Lk}\!\raise0.5pt\hbox{$_{_{_{\Sigma}}}$}
\!\!\sigma\!_{_{^{\!{\ }}}}\!_{\!})\!_{_{^{\hbox{\fiverm q}}}}\!
;_{\!}\hbox{\tenbf Z}) $ is either 0 or homomorphic to
$\hbox{\tenbf Z}_{_{^{\!2}}}$
by Lemma\ $\!$1.\ $\!$\hbox{\tenbf i} p.\ 24. %
Now, use that only the torsion sub-modules matters in the torsion
product by \cite{25} %
Corollary\ 11 p.\ 225.
\hfill$\triangleright$

\medskip
\noindent
{\bf vi.}
\ \ \hbox{\tenbf iii} \ and \ \hbox{$\hbox{\tenrm
Bd}\!\raise0.6pt\hbox{$_{_{_{\!\!\hbox{\fivebf
G}}}}$}\!\!\Sigma_{_{^{\!\hbox{\fiverm q}}}} \subseteq
\hbox{\tenrm Bd}\!\raise0.5pt\hbox{$_{_{_{\!\hbox{\fivebf
Z}}}}$}\!\!\Sigma_{_{^{\!\hbox{\fiverm q}}}}$, by Corollary\ 2 p.\
25 plus Theorem\ 5 p.\ 10.}
\qed
\enddemo

\normalbaselines

\subhead
III:3 \ \ \ Products and Joins of Simplicial Manifolds
\endsubhead

\medskip
Let in the next theorem, when
\lower1.0pt\hbox{$^{_{_{\hbox{\sevensy \char"72} }}}\!$},  
all through, is interpreted as $\times$, the word ``manifold(s)"
in {\bf 12.1} temporarily excludes $\emptyset,\{\emptyset_o\}$ and
$\bullet\bullet$.

When
\lower1.0pt\hbox{$^{_{_{\hbox{\sevensy \char"72} }}}\!$},  
all through, is interpreted as $\ast$ let the word ``manifold(s)"
in {\bf 12.1} stand\nobreak\ for
\underbar{finite}``pseudo-manifold(s)" ($\!$``quasi-manifold(s)"),
\hbox{\tenrm cf.\ \cite{10} %
4.2\ pp.\ 171\hbox{-}2.}
We conclude, w.r.t. joins, that
\hbox{\tenrm Th.\ 12  is trivial if}\
$\Sigma_{_{\!1}}\!$ or $\!\Sigma_{_{\!2}}\!\!
=\!\{\emptyset\!_{_{^{o\!}}}\}
$ and else, $\Sigma_{_{\!1}}\!, \Sigma_{_{\!2}}\!$ must be
\underbar{finite}
since otherwise, theire join isn't locally \underbar{finite}.
$\epsilon=0 /1$ if
$\lower1.0pt\hbox{$^{_{_{\hbox{\sevensy \char"72} }}}$}\!=\times /
\ast.$

\proclaim{Theorem 12}
If $\hbox{\tenbf G}\!_{_{}}$ is an {\bf A}-module, {\bf A}
commutative with unit,
and
$V\!\!\!_{_{\Sigma_{^{i\!}}}}\!\ne\emptyset$ then$;$

\smallskip
\noindent {{\bf 12.1.}\
$\Sigma_{_{\!1}}\!
\lower1.0pt\hbox{$^{_{_{\hbox{\sevensy \char"72} }}}$}\!   
\Sigma_{_{\!2}}\!$ is a $(n_{_{\!1}}\!+n_{_{\!2}}\!+\epsilon)$%
$
\hbox{-manifold}\ \! \Longleftrightarrow \Sigma_{i}$ is a
$n_{_{\!i\!}}$-manifold.

\smallskip
\noindent {\bf 12.2.} $\!\hbox{\tenrm Bd} (\bullet\times \Sigma)=
\bullet\times (\hbox{\tenrm Bd} \Sigma)$. {Else;} $\hbox{\tenrm
Bd} (\Sigma_{_{\!1}}\!
\lower1.0pt\hbox{$^{_{_{\hbox{\sevensy \char"72} }}}$}\!  
\Sigma_{_{\!2}}\!)= ((\hbox{\tenrm Bd} \Sigma_{_{\!1}}\!)
\lower1.0pt\hbox{$^{_{_{\hbox{\sevensy \char"72} }}}$}\!  
\Sigma_{_{\!2}}\!)\cup (\Sigma_{_{\!1}}\!
\lower1.0pt\hbox{$^{_{_{\hbox{\sevensy \char"72} }}}$}\!  
(\hbox{\tenrm Bd} \Sigma_{_{\!2}}\!)).\
$

\smallskip
\noindent {\bf 12.3.}\
If\ any\ side\ of\ $ 12.1$\ holds$;$
$\Sigma_{_{\!1}}
\lower1.0pt\hbox{$^{_{_{\hbox{\sevensy \char"72} }}}$}\!  
\ \!\Sigma_{_{\!2}}\!$ is\ orientable$_{_{\!\hbox{\fivebf G}}}$%
$\Longleftrightarrow$
$\Sigma_{_{\!1}},\Sigma_{_{\!2}}\!$ are both
orientable$_{_{\!\hbox{\fivebf G}}}.$ }
\endproclaim

\demo{Proof}
{\bf (12.1)}\
[{\rm Pseudomanifolds}] Lemma\ 4 p.\ 16.\
$\triangleright$\
[{\rm Quasi-${\sim}$}] Lemma\ 3$\!$ + $\!$4 p.\ 16.\
$\!\triangleright$

\smallskip
The rest of this proof could be substituted for a reference to the
proof of Theorem\ 7 p.\ 13, Proposition\ 1 p.\ 11 and Corollary\
2.{\tenbf iv} p.\ 26.

\smallskip
\noindent{\bf (12.2)} \noindent[{Quasi-manifolds}] 
Put $n:=\dim{\Sigma_{1}\!\times\!\Sigma_{2}}=
\dim{\Sigma_{1}\!+\!\dim\Sigma_{2}}= n_{_1}+n_{_2}$.
The invariance of local $\mdoubleH$omology within 
$\hbox{\tenrm Int}\sigma_{_{\!1}}\!\times \hbox{\tenrm
Int}\sigma_{_{\!2}}$
implies, through Prop\ 1 p.\ 11, %
that w.l.o.g. we'll only  need to study simplices with
$c_{_{^{\!\sigma}}}\!\!=0$ (Def. p.\ 15).
Put $v\!:=\!\dim\sigma\!=\!
\dim{\sigma_{_{^{_{\!}1\!\!}}}}\!+\dim\sigma_{_{^{_{\!}2\!}}}\!=:\nobreak
v_{_{^{_{\!}1\!\!}}}\!+v_{_{^{_{\!}2\!}}}.$
We need to prove that; {$\sigma\!\in\!\hbox{\tenrm
Bd}\!_{_{\hbox{\fivebf
G}}}{\Sigma_{^{_{^{\!}1}}}\!\!\times\!\Sigma_{^{_{^{\!}2}}}}\!$
$\Longleftrightarrow\!$ $\sigma\!_{^{_{1}}}\!\in\!\hbox{\tenrm
Bd}\!_{_{\hbox{\fivebf G}}}{\Sigma_{^{_{\!1}}}}$ or
$\sigma_{^{_{\!2}}}\!\!\in\!\hbox{\tenrm Bd}\!_{_{\hbox{\fivebf
G}}} {\Sigma_{^{_{^{\!}2}}}}_{\!}$.} $
\!$This follows\nobreak\ from\nobreak\ Lemma\ 2.{\bf ii} p.\ 25,
Th.\ 11 and Corollary p.\ 15 %
with $\hbox{\tenbf G}^{\prime}:=\hbox{\tenbf Z}$ which, after
simplification through Note 1 p.\nobreak\ 3, %
gives;
\noindent
$\diamondsuit$
$ \indent
\mdoubleH\!\!\!\!{_{_{n\lower0pt\hbox{\sevenrm-}v\lower0pt\hbox{\sevenrm-}1}}}\!
(\hbox{\tenrm
Lk}\!\!\!\!\!\!\!_{_{_{\Sigma_{1}\!\times\!\Sigma_{2}}}}\!\!\!\!\!\sigma;
\hbox{\tenbf G})
{\lower1pt\hbox{${^{{{_{\hbox{\fivebf
Z}}}}\atop{\hbox{$\cong$}}}}$}}
$
\hfill$
\mdoubleH\!\!\!\!
{_{_{n\!_{_1}\!\!\lower0pt\hbox{\sevenrm-}v\!_{_1}\!\!
\lower0pt\hbox{\sevenrm-}1}}}\!
                 (\hbox{\tenrm Lk}\!\!_{_{_{{\Sigma_{1}}}}}\!\!\!\sigma\!_{_1};\hbox{\tenbf Z})
\otimes\!\!{_{_{\hbox{\fivebf Z}}}}\ \! \mdoubleH\!\!\!\!
{_{_{n\!_{_2}\!\!\lower0pt\hbox{\sevenrm-}v\!_{_2}\!\!
\lower0pt\hbox{\sevenrm-}1}}}\!
           (\hbox{\tenrm Lk}\!\!_{_{_{{\Sigma_{2}}}}}\!\!\!\sigma\!_{_2};\hbox{\tenbf G})
\oplus \hbox{\tenrm Tor}_1^{\hbox{\fivebf Z}}\bigl(
\mdoubleH\!\!\!\!
{_{_{n\!_{_1}\!\!\lower0pt\hbox{\sevenrm-}v\!_{_1}\!\!
\lower0pt\hbox{\sevenrm-}2}}}\!
                  (\hbox{\tenrm Lk}\!\!_{_{_{{\Sigma_{1}}}}}\!\!\!\sigma\!_{_1};\hbox{\tenbf Z})
,\mdoubleH\!\!\!\!
{_{_{n\!_{_2}\!\!\lower0pt\hbox{\sevenrm-}v\!_{_2}\!\!
\lower0pt\hbox{\sevenrm-}1}}}\!
           (\hbox{\tenrm Lk}\!\!_{_{_{{\Sigma_{2}}}}}\!\!\!\sigma\!_{_2};\hbox{\tenbf G})\bigr).
\ \ \ {\triangleright}
$

\smallskip
\noindent A similar reasoning holds also for pseudomanifolds with
$\sigma$ restricted to the submaxidimensional simplices. Use
Corollary p.\ 15 %
also for joins.
\hfill$\triangleright$

\smallskip
\noindent{\bf (12.3)}
\noindent $
(\Sigma_{_{^{_{\!}1}}}\!\lower1.0pt\hbox{$^{_{_{\hbox{\sevensy
\char"72} }}}$}\!\Sigma_{_{^{_{\!}2\!}}}, \hbox{\tenrm
Bd}_{_{\!\hbox{\fivebf
G}}}{\!}_{\!}(\Sigma_{_{^{_{\!}1}}}\!\lower1.0pt\hbox{$^{_{_{\hbox{\sevensy
\char"72} }}}$}\!\Sigma_{_{^{_{\!}2\!}}}))
_{\!}=_{\!}
[{\hbox{\tenrm Th.\ 12.2}}]
_{\!}=_{\!}
(\Sigma_{_{^{_{\!}1}}}\!\lower1.0pt\hbox{$^{_{_{\hbox{\sevensy
\char"72} }}}$}\!\Sigma_{_{^{_{\!}2\!}}},
\Sigma_{_{^{_{\!}1}}}\!\lower1.0pt\hbox{$^{_{_{\hbox{\sevensy
\char"72} }}}$}\!\hbox{\tenrm Bd}_{_{\!\hbox{\fivebf
G}}}{\!}_{\!}\Sigma_{_{^{_{\!}2\!}}}\cup \hbox{\tenrm
Bd}_{_{\!\hbox{\fivebf
G}}}{\!}_{\!}\Sigma_{_{^{_{\!}1}}}\!\lower1.0pt\hbox{$^{_{_{\hbox{\sevensy
\char"72} }}}$}\! \Sigma_{_{^{_{\!}2\!}}})
_{\!}=_{\!}
$

$
_{\!}=_{\!}
[{\hbox{\tenrm Def.}
\ {\!}\hbox{\tenrm p.\ 5}}]%
_{\!}=_{\!}
(\Sigma_{_{^{_{\!}1\!}}},\hbox{\tenrm Bd}_{_{\!\hbox{\fivebf
G}}}{\!}_{\!}\Sigma_{_{^{_{\!}1\!}}})
\lower1.0pt\hbox{$^{_{_{\hbox{\sevensy \char"72} }}}$}\!
(\Sigma_{_{^{_{\!}2\!}}},\hbox{\tenrm Bd}_{_{\!\hbox{\fivebf
G}}}{\!}_{\!}\Sigma_{_{^{_{\!}2\!}}}).
$

By Cor. $\!$2.iv p.\ 26 %
we can w.l.o.g.$\ $confine our study to pseudomanifolds, and
\nobreak{choose} the coefficient module to be, say, a field
$\hbox{\tenbf k}$ (char$\hbox{\tenbf k}\!\!\neq\!\!2$) or {\bf Z}.
Since any finite maxi-dimensional submanifold, i.e. a submanifold
of maximal dimension, in
$\Sigma_{_{^{_{\!}1\!}}}\!\!\times\!\!\Sigma_{_{^{_{\!}2\!}}}$

\noindent ($\Sigma_{_{^{_{\!}1\!}}}\!\ast\Sigma_{_{^{_{\!}2\!}}},$
cp.
\cite{29} %
(3.3) p.\ $_{\!}$59) %
can be embedded in the product ({\rm join}) of two \hbox{finite
maxi-}

\noindent dimensional submanifolds,
we confine, w.l.o.g., our attention to \hbox{finite maxi-dimen-}

\noindent sional submanifolds ${\tensy S}\!_{_{1}}$, ${\tensy
S}\!_{_{2}}$ of $\Sigma_{_{^{_{\!}1\!}}}$ resp.
$\Sigma_{_{^{_{\!}2\!}}}$. Now, use
\hbox{$\hbox{\tenrm Eq.\ 1\ p.\ 8}$ %
{\bf(}\hbox{\tenrm Eq.\ 3\ p.\ 10}{\bf)}.} %
\qed
\enddemo

%


\example{Example 1}
For a triangulation $\Gamma$ of a two-dimensional cylinder
$\hbox{\tenrm Bd}\!\lower1.3pt\hbox{$_{_{\!\hbox{\fivebf Z}}}$}\!
\Gamma\!_{_{^{\!\hbox{\fiverm h}}}}\!\!=\! \hbox{\tenrm two\
circles}$.
\indent
$\hbox{\tenrm Bd}\!\lower1.3pt\hbox{$_{_{\!\hbox{\fivebf
Z}}}$}\!\bullet =\{\emptyset_o\!\}.$ By Th.\ 12
$\hbox{\tenrm Bd}\!\lower1.3pt\hbox{$_{_{\!\hbox{\fivebf Z}}}$}\!
(\Gamma\!\ast\bullet)_{_{\!\hbox{\fiverm q}}}
\!\!=\!
\Gamma\cup(\{\hbox{\tenrm two\ circles}\}\ast\bullet) $.
So, $\hbox{\tenrm Bd}_{_{\hbox{\fivebf
Z}}}\!(\Gamma\ast\bullet)_{_{\!\hbox{\fiverm q}}}$,
{\bf R}$^{^{_{3}}}$-
\indent
realizable as a pinched torus, is a 2-pseudomanifold but not a
quasi-manifold.

\smallskip
\item{\bf 2.}
``The boundary w.r.t. {\bf Z} of the cone of the M\"obius band" =
$\hbox{\tenrm Bd}\!\lower1.3pt\hbox{$_{_{\!\hbox{\fivebf Z}}}$}\!
({\hbox{\tensy M}}\!\ast\bullet)_{_{\!\hbox{\fiverm q}}}
\!\!=\!
(\hbox{\tensy M}\ast\{\emptyset_o\!\})\cup(\{\hbox{\tenrm a\
circle}\}\ast\bullet)
\!\!=\!
{\hbox{\tensy M}}\cup \{\hbox{\tenrm a\ $2$-disk}\}$
which is a well-known representation of the real projective plane
{\tensy P}\raise4pt\hbox{{\fiverm 2}}$_{\!}$.
So, %
$\hbox{\tenrm Bd}\!\lower1.3pt\hbox{$_{_{\!\hbox{\fivebf Z}}}$}\!
({\hbox{\tensy M}}\!\ast\bullet)_{_{\!\hbox{\fiverm q}}}
\!\!=\!
\hbox{\tensy P}\raise4pt\hbox{{\fiverm 2}}_{\!}
$
is a homology\lower1.0pt\hbox{$\!\!_{_{\hbox{\fivebf
Z}_{_{\!\hbox{\fivebf p}}}}}\!\!\!$}
$2$-manifold  \underbar{with} boundary
$\{\emptyset_o\!\}\neq\emptyset$ if  $\hbox{\tenbf p}\neq2$. \ \
$\!{{\hbox{\bf Z}_{_{\!\hbox{\fivebf p}}}}}\!\!:=\!$ The
prime-number field modulo {\bf p}. \hfill\break
\indent
$
\hbox{\tensy P}^{_2}\#{\hbox{\tenbf S}}^{_{^{2}}}\!\!
\!=\!
\hbox{\tensy P}^{_2}_{\!}
\!=\!
{\hbox{\tensy M}}\
\raise1pt\hbox{$\!\cup\!\!_{_{_{\!\!\!\hbox{\fiverm Bd}}}}$}
\!\{\hbox{\tenrm a\ $2$-disk}\}
$
confirms the obvious, i.e. that the $n$-sphere is the unit element
w.r.t. the connected sum of two $n$-manifolds, cf. \cite{21} %
p.\ 38ff   + Ex.\ 3 p. \ 366.
``\ \raise1pt\hbox{$\!\cup\!\!_{_{_{\!\!\!\!\hbox{\fiverm Bd}}}}
$}"
is ``union through homeomorphic identification of\nobreak\
boundaries".
\indent
$\hbox{\tenrm Bd}\!\lower1.3pt\hbox{$_{_{\!\hbox{\fivebf Z}}}$}\!
({\hbox{\tensy M}}\!\ast\!{\hbox{\tenbf
S}}^{_{^{1}}})_{_{\!\hbox{\fiverm q}}}
\!\!=\!
(\hbox{\tensy M}\ast\emptyset)\cup({\hbox{\tenbf
S}}_{\!}^{_{^{1}}}\ast{\hbox{\tenbf S}}^{_{^{1}}})
\!\!=\!
{\hbox{\tenbf S}}_{\!}^{_{^{3}}}.
$
So,  $\hbox{\tenrm Bd}\!\lower1.3pt\hbox{$_{_{\!\hbox{\fivebf
Z}}}$}\!({\hbox{\tensy M}}_{_{\!1}}{\!}\ast_{\!} \hbox{\tensy
M}_{_{\!2}}\!)_{_{\!\hbox{\fiverm q}}}\!{\!}
={\!}({\hbox{\tenbf S}}_{\!}^{_1}\!\ast {\hbox{\tensy
M}}_{_{\!2}}{\!})_{_{\!\hbox{\fiverm q}}}\!\ {\!}\cup\ {\!}
(_{\!}\hbox{\tensy M}_{_{\!1}}\!\ast {\hbox{\tenbf
S}}_{\!}^{_1})_{_{\!\hbox{\fiverm q}}}\!
$
is a quasi-$4$-manifold without boundary, represented as the
connected sum of two copies of $ ({\hbox{\tensy
M}}\!\ast\!{\hbox{\tenbf S}}^{_{^{1}}}\!)_{_{\!\hbox{\fiverm q}}}.
$

\smallskip
\item{\bf 3.}
Let {\tensy P}\raise4pt\hbox{{\fiverm 2}}$_{\!}$
({\tensy P}\raise4pt\hbox{{\fiverm 4}}$_{\!}$) be a triangulation
of the projective plane (projective space with $\dim\hbox{\tensy
P}_{^{\!}}\raise4pt\hbox{{\fiverm 4}}_{\!}\!\!= \!4$) implying
$\hbox{\tenrm Bd}\!\lower1.3pt\hbox{$_{_{\!\hbox{\fivebf Z}}}$}\!
\hbox{\tensy P}\raise4pt\hbox{{\fiverm
2}}_{\!}\!\!_{_{\!\hbox{\fiverm h}}}\!\!=\! \hbox{\tenrm
Bd}\!\lower1.3pt\hbox{$_{_{\!\hbox{\fivebf Z}}}$}\!\hbox{\tensy
P}\raise4pt\hbox{{\fiverm 4}}_{\!}\!_{_{\!\hbox{\fiverm
h}}}\!\!=\! \{\emptyset_o\!\}.$
So, by $\hbox{\tenrm Th.\ 12,}$ $\hbox{\tenrm
Bd}\!\!\lower1.3pt\hbox{$_{_{\hbox{\fivebf Z}_{{\!\hbox{\fivebf
p}}}}}$}\!\!\!\! (\hbox{\tensy P}\raise4pt\hbox{{\fiverm
2}}_{\!}\!\ast\hbox{\tensy P}\raise4pt\hbox{{\fiverm
4}}_{\!})_{_{\!\hbox{\fiverm h}}} \!\!=\!\hbox{\tensy
P}\raise4pt\hbox{{\fiverm 4}}\cup\hbox{\tensy
P}\raise4pt\hbox{{\fiverm 2}}_{\!}\!$,
$\hbox{\tenbf p}\!\neq\!2$
$(\hbox{\tenrm Bd}\!\!\lower1.3pt\hbox{$_{_{\hbox{\fivebf
Z}_{{\!\hbox{\fivebf p}}}}}$}\!\!\!\!
(\hbox{\tensy P}\raise4pt\hbox{{\fiverm 2}}_{\!}\!\ast\hbox{\tensy
P}\raise4pt\hbox{{\fiverm 4}}_{\!})_{_{\!\hbox{\fiverm q}}}
\!=\!\emptyset\ \hbox{\tenrm if}\ \hbox{\tenbf p}\!=\!2)$,
and $\dim\hbox{\tenrm Bd}\!\!\lower1.3pt\hbox{$_{_{\hbox{\fivebf
Z}_{{\!\hbox{\fivebf p}}}}}$}\!\!\!\! (\hbox{\tensy
P}\raise4pt\hbox{{\fiverm 2}}_{\!}\!\ast\hbox{\tensy
P}\raise4pt\hbox{{\fiverm 4}}_{\!})_{_{\!\hbox{\fiverm h}}}
\!\!=\!4$
while
$\dim(\hbox{\tensy P}\raise4pt\hbox{{\fiverm
2}}_{\!}\!\ast\hbox{\tensy P}\raise4pt\hbox{{\fiverm
4}}_{\!})_{_{\!\hbox{\fiverm h}}} \!\!=\!7$, cp. Corollary\ 1 p.\
26. \
$\hbox{\tenrm Bd}_{_{\hbox{\fivebf Z}}}\!(\hbox{\tensy P}
\raise4pt\hbox{{\fiverm
2}}_{\!}\!\ast\bullet\bullet)_{_{\!\hbox{\fiverm
q}}}\!\!=\!{\bullet\bullet}$
 \ \ and  \ \
$\hbox{\tenrm Bd}_{_{\hbox{\fivebf Z}}}\!(\hbox{\tensy
P}\raise4pt\hbox{{\fiverm
2}}_{\!}\!\ast\bullet)_{_{\!\hbox{\fiverm q}}}\!
{=\!\hbox{\tensy P}\raise4pt\hbox{{\fiverm 2}}_{\!}\cup\bullet}$\
.

\smallskip
\item{\bf 4.}
$ \hbox{\tenbf E}^{_{^m}}\!\!:= $ the $m$-unit ball. With
$n:=p+q,\ p,q\ge0,$
$\hbox{\tenbf S}^{_{^n}}\!\! =$ Bd$\hbox{\tenbf E}^{^{_{n+1}}}\!
\simeq$ Bd$(\hbox{\tenbf E}^{^{_p}}\!\!\ast \hbox{\tenbf
E}^{^{_q}}\! ) \ \!\!\simeq $\break
$\simeq\hbox{\tenbf E}^{^{_{p}}}\!\!\!\ast \hbox{\tenbf
S}^{^{_{q-1}}}\!\!
\cup \hbox{\tenbf S}^{^{_{p-1}}}\!\! \ast\hbox{\tenbf
E}^{^{_{q}}}\!\simeq\! $
Bd$ (\hbox{\tenbf E}^{^{_{p+1}}}\!\!\!\times\! \hbox{\tenbf
E}^{^{_{q}}}\!)
\simeq$
$ \hbox{\tenbf E}^{^{_{p+1}}}\!\!\times\! \hbox{\tenbf
S}^{^{_{q-1}}}\!\cup \hbox{\tenbf S}^{^{_{p}}}\!\!\!\times\!
\hbox{\tenbf E}^{^{_{q}}}\!$
by Th.\ 12 and Lemma p.\nobreak\ 14. %
Cp. \cite{17} %
p.\ 198 Ex.\ 16 on surgery. %
$\!\hbox{\tenbf S}^{^{_{n+1}}}\!\!\!\simeq
\hbox{\tenbf S}^{^{_{p}}}\!\!\ast \hbox{\tenbf S}^{^{_{q}}}\!\!$
and
$ \hbox{\tenbf E}^{^{_{n+1}}}\!\!\!\!\simeq
\hbox{\tenbf E}^{^{_{p}}}\!\!\!_{\!}\ast_{\!} \hbox{\tenbf
S}^{^{_{q}}}\!$ also hold.

\smallskip
\noindent
\item{\bf 5.}
$_{\!}\hbox{\tensy P}\raise4pt\hbox{{\fiverm 2}}_{\!}$ w.r.t.
$\!\hbox{\tenbf G}_{\!}:=_{\!}\hbox{\tenbf Z}_{_{\!\hbox{\fivebf
2}}}\!_{\!}\oplus_{\!}\hbox{\tenbf Z}_{_{\!\hbox{\fivebf 3}}}\!$
is a non-orientable homology$\!_{_{\hbox{\fivebf G}}}$ 2-manifold
without boundary.

\smallskip
See also \cite{21} %
p.\ 376 for some non-intuitive manifold examples.
Also \cite{27} %
pp.\ 123-131 gives insights on different aspects of different
kinds of simplicial manifolds.
\endexample

\proclaim{Proposition}
If ${{\Sigma}}_{_{^{\!\hbox{\fiverm q}}}} \!$ is finite and
$\ -1\!\le\!\dim\hbox{\tenbf B}_{_{{\!{\hbox{\fivebf G},i}\!}}}\!
<\!\dim\Sigma\!-\!1$
then $\hbox{\tenrm Lk}_{_{\!\Sigma}}\!{\delta}$ is
non-orientable$_{_{^{\!\hbox{\fivebf G}}}}\!\!$
\hbox{for all $ {\delta}\!\in\!\hbox{\tenbf
B}_{_{{\!{\hbox{\fivebf G},i}\!}}}.$}
{\rm(Note that Cor.\ 2.\ {\bf ii} p.\ 26 is the special case
$\hbox{\tenrm
Lk}_{_{\!\Sigma}}\!{\emptyset}_{_{^{\!o}}}\!=\Sigma$.)}
\endproclaim
\normalbaselines

\demo{Proof}
Use the proof of Cor.\ 2.{\bf ii} p.\ 26 %
plus the end of Note\ 3 p.\ 25. %
\qed
\enddemo

\subhead
III:4 \ \ \ Simplicial Homology$_{_{\!\hbox{\fivebf G}}}\!\!$
Manifolds and Their Boundaries
\endsubhead

\normalbaselines

\medskip
In this section we'll work mainly with finite simplicial complexes
and though we're still working with arbitrary coefficient modules
we'll delete those annoying quotation marks surrounding
\hbox{``}{\rm CM}$\!_{_{\hbox{\fivebf G}}}\!\!"$. The coefficient
module plays, through the St-R ring functor, a more delicate role
in commutative ring theory than it does here in our
$\mdoubleH$omology theory, so when it isn't a Cohen-Macaulay ring
we can not be sure that a {\rm CM} complex gives rise to a {\rm
CM} St-R ring.

\proclaim{Lemma 1}
\hbox{\tenbf i.}
$\Sigma$ is a
homology$\!_{_{\!\hbox{\fivebf G}}}\!$ manifold \underbar{iff}
{\hbox{\tenbf [}}{\hbox{\tenbf [}}%
$\Sigma\!=\!\bullet\bullet$ or $\Sigma$ is connected and
$\hbox{\tenrm Lk}\!_{_{\Sigma}}\!\!\hbox{\tenbf v}$ is a finite
\noindent {\rm CM}$\!_{_{\hbox{\fivebf G}}}\!$
{\rm pseudo} manifold$\ \!$ for all vertices\ \ $\hbox{\tenbf
v}\in\!{V}_{_{^{\!\Sigma}}}.$
{\hbox{\tenbf ]}}
\underbar{\hbox{\tenbf and}}
{\hbox{\tenbf [}}
$ \hbox{\tenrm Bd}_{_{\!\hbox{\fivebf
G}}}\!\Sigma_{_{^{\hbox{\fiverm q}}}} \!\!=\! \hbox{\tenrm
Bd}_{_{\!\hbox{\fivebf Z}}}\!\Sigma_{_{^{\hbox{\fiverm q}}}}\! $
or else$;$
$[[\hbox{\tenrm Bd}_{{_{{\!\hbox{\fivebf
Z}}}}_{_{^{\!\hbox{\fivebf 2}}}}}\!\!\!\Sigma_{_{^{\hbox{\fiverm
q}}}}\!\!= \hbox{\tenrm Bd}_{_{\!\hbox{\fivebf
G}}}\!\Sigma_{_{^{\hbox{\fiverm
q}}}}\!\!=\emptyset\not=\!\{\emptyset_{_{^{\!o}}}\!\}\!=
\hbox{\tenrm Bd}_{_{\!\hbox{\fivebf
Z}}}\!\Sigma_{_{^{\hbox{\fiverm q}}}}]$
or
$[\ \!\exists\  \{\emptyset{_{_{^{\!o\!}}}}\!\}
\nsupseteqq
\hbox{\tenbf B}\!\raise0.6pt\hbox{$_{_{\!{_{\hbox{\fivebf
Z}\!,_{\!}j}}}}$}\!\!\!\!\!^{^{_{\Sigma\!\!}}}\ \!\!
\subset
\Sigma^{^{{\!{n\!\lower1.0pt\hbox{-}\hbox{\fiverm 3}}}}}\!\!
 $
and
$ \hbox{\tenrm Tor}_1^{\hbox{\fivebf Z}}(\hbox{\tenbf
Z}_{_{^{\!2}}}\!,\!\hbox{\tenbf G})\!=\!\hbox{\tenbf G}]]$%
{\hbox{\tenbf ]}}{\hbox{\tenbf ]}}.

\noindent
\hbox{\tenbf ii.}
\nobreak If $\Delta$ is a {\rm CM}$\!_{_{\hbox{\fivebf Z}}}\!$
homology$_{_{\!\hbox{\fivebf Z}}}\!$
manifold then
$\hbox{\tenrm Bd}\!_{_{\hbox{\fivebf
G}}}\!\!\Delta_{_{^{\!\hbox{\fiverm h}}}}\!\!=\! \hbox{\tenrm
Bd}_{_{\hbox{\fivebf Z}}}\!\Delta_{_{^{\!\hbox{\fiverm h}}}}$ for
any module {\bf G}, and so, for a homology$_{_{\!\hbox{\fivebf
Z}}}\!$ manifold $\Sigma$,
$\hbox{\tenrm Bd}_{_{\hbox{\fivebf
G}}}\!\Sigma=\emptyset,\{\emptyset_{_{^{\!o}}}\!\}$ or
$\dim\hbox{\tenbf B}_{_{{\!{_{\!}i}\!}}}=(n\!-\!1)$
for each boundary component.
\endproclaim

\demo{Proof} {\bf i.}
Prop.\ 1 p.\ 18, %
Lemma\ $\!$2.{\bf ii} p.\ 24 %
and the proof of Corollary\ 2.{\bf v} p.\ 27.\break
{\bf ii.} Use {\rm Theorem\ 5 p.\ 10} %
to prove the boundary equality and then,
{\rm put $\hbox{\tenbf G}_{\!}=_{\!}\hbox{\tenbf Z}_{_{^{\!2}}}$
in Corollary\ 2{\bf iii} p.\ 26.}
\qed
\enddemo

\proclaim{Lemma 2}
{For a finite $\Sigma_{_{^{\!\hbox{\fiverm q}}}}\!$, $
\delta\!_{_{n\!_{_{\Sigma}}}}\!\!\!\!\!:\!
\mdoubleH\!\!\!_{_{n\!_{_{\Sigma}}}}\!\! ( {\Sigma},\hbox{\tenrm
Bd}\!\raise0.8pt\hbox{$_{_{_{\!\hbox{\fivebf
G}}}}$}\!\!{{\Sigma}}; \hbox{\tenbf G}) \!\rightarrow
\mdoubleH\!\!\!\!_{_{{{n\!_{_{\Sigma}}\!\!-\!1\!\!}}}} (
{\Sigma},\hbox{\tenbf B}_{_{\!{i}}}\cap
({\cup\!{_{\!}}\!_{_{\!}}\!\!_{{\!}}\cup
\!\!\!\!_{_{\!}}\!\!{_{_{_{{j\not=i}}}}}} \!\!\hbox{\tenbf
B}_{_{\!\!{j}}}); \hbox{\tenbf G}) $ in the relative {\bf
M-$\!$Vs}$_{_{\!^{o}}}$ \noindent w.r.t. $ \{
({\Sigma},\hbox{\tenbf B}_{_{\!{i}}}),
({\Sigma},{\cup\!_{\!}\!_{_{\!}}\!_{_{\!}}\cup}
\!\!\!\!_{\!}\!_{_{_{{{j\not=i}}}}}\!\!\hbox{\tenbf
B}{_{_{\!{j}}}}{_{\!}}) \}$
is injective if $\#\hbox{\tenbf I}\geq2$ in Definition\ 4
p.\nobreak\ 23.

\goodbreak

So, $\! [\mdoubleH\!\!\!\!_{_{{{n\!_{_{\Sigma}}\!\!-\!1\!\!}}}}
({\Sigma},\{\emptyset\};\hbox{\tenbf G})\!=\!0]
\!\Longrightarrow\! $ $[ \mdoubleH\!_{_{n\!_{_{\Sigma}}}}\!\! (
{\Sigma},\hbox{\tenrm
Bd}\!\!\raise0.6pt\hbox{$_{_{_{\!\hbox{\fivebf
G}}}}$}\!\!{{\Sigma}}; \hbox{\tenbf G}) \!=\!0 $ or $\hbox{\tenrm
Bd}\!\!\raise0.6pt\hbox{$_{_{_{\!\hbox{\fivebf
G}}}}$}\!\!{{\Sigma}}\ {is\ strongly\ connected}
],$ %
} %
e.g., if $\Sigma\!\ne\!\bullet, \bullet\bullet$ is a {\rm
CM}$\!_{_{\hbox{\fivebf G}}}\!$ {quasi-n-manifold}.
If $\Sigma$ is a finite {\rm CM}$\!_{_{\hbox{\fivebf Z}}}\!$
quasi-n-manifold then $\hbox{\tenrm
Bd}\!\raise0.8pt\hbox{$_{_{_{\!\hbox{\fivebf
Z}}}}$}\!\!{{\Sigma}}= \emptyset$ or it's\ strongly\ connected and
$\dim(\hbox{\tenrm Bd}\!\raise0.8pt\hbox{$_{_{_{\!\hbox{\fivebf
Z}}}}$}\!\!{{\Sigma}})= n_{_{\!\Sigma}}\!\!-1$ $
($%
\hbox{\rm by Lemma\ 1}$)$.
\endproclaim

\demo{Proof}
Use the relative {\bf M-$\!$V$_{\!}$s}$_{_{\!^{o}}}\!$ w.r.t. $\!
\{ ({\Sigma},\hbox{\tenbf B}_{_{\!{i}}}),
({\Sigma},{\cup\!\!\!_{_{\!}}\!_{_{\!}}\cup}
\!\!\!\!\!_{_{_{{{j\not=i}}}}}\hbox{\tenbf B}_{_{\!{j}}}) \}, $ $
\dim(\hbox{\tenbf B}{_{_{\!{i}}}}_{\!}\cap\
{\cup\!\!\!_{_{\!}}\!_{_{\!}}\cup}
\!\!\!\!\!_{_{_{{{j\not=i}}}}}\!\hbox{\tenbf B}_{_{\!{j}}})\le
n_{_{\!\Sigma}}\!-3$ and Corollary\ 2.{\bf i}\ p.\ 25. %
\qed
\enddemo

\proclaim{Theorem 13.i}
If $\Sigma$ is a finite orientable$_{_{^{\!\!\hbox{\fivebf
G}}}}\!\!$ {\rm CM}$\!_{_{\hbox{\fivebf G}}}\!$
homology$_{_{^{\!\!\hbox{\fivebf G}}}}$ $n$-manifold with boundary
then, $\hbox{\tenrm Bd}_{_{{\!\hbox{\fivebf
G}}}}\!\Sigma$\nobreak\ is an orientable$_{_{^{\!\!\hbox{\fivebf
G}}}}\!\!$
homology$_{_{^{\!\!\hbox{\fivebf G}}}}$ $(n-1)$-manifold without
boundary.\break
\ \ \ {\bf ii.} Moreover; $\hbox{\tenrm Bd}_{{_{{\!\hbox{\fivebf
G}}}}}\!\!\Sigma$ is Gorenstein$_{{_{{\!\hbox{\fivebf G}}}}}\!.$
\endproclaim

\demo{Proof}{\bf i.}
Induction over the dimension, using Th.\ 11$\ \!${\bf i},{\bf ii},
once the connectedness
of the\nobreak\ boundary is established through Lemma\ 2, while
orientability$_{_{^{\!\!\hbox{\fivebf G}}}}\!$ resp.
 $\hbox{\tenrm Bd}_{_{{\!\hbox{\fivebf
G}}}}\!(\hbox{\tenrm Bd}_{_{{\!\hbox{\fivebf
G}}}}\!\Sigma)\!=\!\emptyset$ follows from
Cor.\ 2.iii-iv p.\ 26. \hfill$\triangleright$

\noindent{\bf ii.} $\Sigma\ast(\bullet\bullet)$ is a finite
orientable$_{_{^{\!\!\hbox{\fivebf G}}}}
\!$
{\rm CM}$\!_{_{\hbox{\fivebf G}}}\!$
homology$_{_{^{\!\!\hbox{\fivebf G}}}}\!\!$ $(n\!+\!1)$-manifold
with boundary by Th.\ 12 + Cor.\ i p.\ 12.
$\ \!\hbox{\tenrm Bd}_{_{{\!\hbox{\fivebf
G}}}}\!(\Sigma\ast\bullet\bullet) \!=\! [{{\hbox{\eightrm Th.\
12.}\ \!\hbox{\eightbf 2}}\atop{ \hbox{\eightrm or Th.\ 7.}\
\!\hbox{\eightbf 2}}}] \!=\! \big(\Sigma\ast\hbox{\tenrm
Bd}_{_{{\!\hbox{\fivebf G}}}}\!(\bullet\bullet)\big) \cup
\big((\hbox{\tenrm Bd}_{_{{\!\hbox{\fivebf
G}}}}\!\Sigma)\ast(\bullet\bullet)\big) \!=\! \Sigma\ast\emptyset
\cup (\hbox{\tenrm Bd}_{_{{\!\hbox{\fivebf
G}}}}\!\Sigma)\ast(\bullet\bullet) \!=\! \big(\hbox{\tenrm
Bd}_{_{{\!\hbox{\fivebf G}}}}\!\Sigma\big)\ast(\bullet\bullet)$
where the l.h.s.\ is an orientable$_{_{^{\!\!\hbox{\fivebf
G}}\!\!}}$ homology$_{_{{\!\!_{\!}\hbox{\fivebf G}}\!\!}}$
$n$-manifold without boundary by the first part.
So, $\hbox{\tenrm Bd}\!\!\raise0.7pt\hbox{$_{_{_{\!\hbox{\fivebf
G}}}}$}\!\!\Sigma$ is an orientable$_{_{^{\!\!\hbox{\fivebf
G}}}}\!$
{\rm CM}$\!_{_{\hbox{\fivebf G}}}\!$
homology$_{_{^{\!\!\hbox{\fivebf G}}}}$ $\!(n-1)$-manifold
without\ boundary by Th.$\ 12$\ \ $\!$i.e.
it's Gorenstein$\!_{_{\hbox{\fivebf G}}}\!.$
\qed
\enddemo

\remark{Note}
$\emptyset\!\neq\!{{\Delta}}$ is a 2-{\rm CM}$\!_{_{\hbox{\fivebf
G}}}\!$ homology$\!_{_{\!\hbox{\fivebf G}}}\!$ manifold $
\!\Longleftrightarrow\! {{\Delta}} \!=\! \hbox{\tenrm
core}{\Delta}$ is Gorenstein$\!_{_{\hbox{\fivebf G}}}\!\!$ $
\Longleftrightarrow\break
\Longleftrightarrow\!
 {{\Delta}}$ is a {\rm
homology}$\!_{_{\!\hbox{\fivebf G}}}$ sphere.
\endremark

\proclaim{Corollary 1}
{\rm(Cp. \cite{17} %
p.\ 190.)} %
If $\Sigma$ is a finite orientable$_{_{^{\!\!\hbox{\fivebf G}}}}$
homology$_{_{^{\!\!\hbox{\fivebf G}}}}n$-manifold with boundary,
so is $\hbox{\tenbf2}\Sigma$ except that $\hbox{\tenrm
Bd}_{_{{\!\hbox{\fivebf G}}}}\!(\hbox{\tenbf2}\Sigma)=\emptyset.$
$\hbox{\tenbf2}\Sigma=\!$``the double of $\Sigma"
\!:=\Sigma\ \!$%
\raise1pt\hbox{$\cup\!\!_{_{_{\!\!\!\!\hbox{\fiverm Bd}}}}$}%
\rlap{\raise2.0pt\hbox{\tent{\char"7F}}}{\tenrm$_{\!}${\char"06}}
where
\rlap{\raise2.0pt\hbox{\tent{\char"7F}}}{\tenrm$_{\!}${\char"06}}
is a disjoint mirrored copy of $\Sigma$
and
``\ \raise1pt\hbox{$\!\cup\!\!_{_{_{\!\!\!\!\hbox{\fiverm
Bd}}}}\!\!$}"
is ``the union through identification of the boundary
vertices$"\!$. If $\Sigma$ is $\hbox{\tenrm
CM}_{_{{\!\hbox{\fivebf G}}}}\!$ \hbox{then
\hbox{$\hbox{\tenbf2}\Sigma$ is
$\hbox{\tenrm2}\hbox{-}\hbox{\tenrm CM}_{_{{\!\hbox{\fivebf
G}}}}\!.$ }}
\endproclaim

\demo{Proof}
Apply the {\bf M-$_{\!}$Vs} to the pair $(\hbox{\tenrm
Lk}\raise0.0pt\hbox{$_{_{{\!{\Sigma^{{_{\
}}}}}}}$}\!\!\!\!{\hbox{\tenrm v}}, \hbox{\tenrm Lk}\!
\raise0.0pt\hbox{$_{_{{\!{\Sigma^{^{_{\prime}}}}}}}$}\!\!{\hbox{\tenrm
v}})$
using Prop.\ 2. {\bf a} p.\ 30 %
and then to
$(\Sigma,$ %
\rlap{\raise2.0pt\hbox{\tent{\char"7F}}}{\tenrm$_{\!}${\char"06}}%
)
for the $\hbox{\tenrm CM}_{_{{\!\hbox{\fivebf G}}}}\!$ case,
or even simpler, use \cite{14} %
p.\ 57 (23.6) Lemma, %
where also the non-relative augmental {\bf M-$_{\!}$Vs} is used.
\qed
\enddemo

\proclaim{Theorem 14}
If $\Sigma$ is a finite {\rm CM}$\!_{_{\hbox{\fivebf
Z}}}\!$-homology$_{_{\!\hbox{\fivebf Z}}}\!\!$ $n$-manifold, then
$\Sigma$ is orientable$_{_{{\!\hbox{\fivebf Z}}}}.$
\endproclaim

\demo{Proof}
$\Sigma$  finite {\rm CM}$_{_{\!\hbox{\fivebf Z}}}\! $ $
\Longleftrightarrow \Sigma\ ${\rm CM}${{_{_{\hbox{\fivebf
Z}}}}_{_{{\!\hbox{\fivebf p}}}}}\! $ for all prime fields $
{{{{\hbox{\tenbf Z}}}}_{_{{\!\hbox{\fivebf p}}}}}\!\! $ of
characteristic ${{\hbox{\tenbf p}}},$ by (M.A. Reisner, 1976)
induction over $\dim \Sigma$,
Theorem 5 p. 10 %
and the {Structure Theorem for Finitely Generated Modules over
{\bf PID}s}.
So, $\Sigma$ is a finite {\rm CM}$_{{_{{\!\hbox{\fivebf
Z}}}}_{_{^{\!\hbox{\fivebf p}}}}}\!\!\!$
homology$_{{_{{\!\hbox{\fivebf Z}}}}_{_{^{\!\hbox{\fivebf
p}}}}}\!\!\!$ $n$-manifold for any prime number $\hbox{\tenbf p}$
by Lemma\ 1.{\bf i}, since $\hbox{\tenrm
Bd}\raise0.0pt\hbox{$_{_{{\!\hbox{\fivebf Z}}}}$}\!\Sigma
=
\hbox{\tenrm Bd}_{{_{{\!\hbox{\fivebf Z}}}}_{_{^{\!\hbox{\fivebf
p}}}}}\!\!\!\Sigma$ by Lemma\ 1.{\bf ii} above.
In particular, $%
\hbox{\tenrm Bd}_{{_{{\!\hbox{\fivebf Z}}}}_{_{^{\!\hbox{\fivebf
2}}}}}\!\!\!\Sigma$ is Gorenstein$_{{_{{\!\hbox{\fivebf
Z}}}}_{_{^{\!\hbox{\fivebf 2}}}}}\!\!$ by Lemma\ 1.{\bf ii} p.\ 24
and Theorem\ 13 above.
If $\hbox{\tenrm Bd}_{_{{\!\hbox{\fivebf Z}}}}\!\Sigma
=
\emptyset$
then $\Sigma$ is orientable$_{_{{\!\hbox{\fivebf Z}}}}$.
Now, if $\hbox{\tenrm Bd}_{{_{{\!\hbox{\fivebf
Z}}}}_{_{^{\!\hbox{\fivebf 2}}}}}\!\!\!\Sigma\!\ne\!\emptyset$
then $\dim\hbox{\tenrm Bd}_{{_{{\!\hbox{\fivebf
Z}}}}_{_{^{\!\hbox{\fivebf 2}}}}}\!\!\!
\Sigma\!=\!n\!-\!1$ by Cor. 2{\rm iii+iv} p.\ 26 %
and, in particular, $\hbox{\tenrm Bd}_{_{{\!\hbox{\fivebf
Z}}}}\!\Sigma
=
\hbox{\tenrm Bd}_{{_{{\!\hbox{\fivebf Z}}}}_{_{^{\!\hbox{\fivebf
2}}}}}\!\!\!\Sigma$ is a quasi-$(n\!-\!1)$-manifold.
\normalbaselines

$\hbox{\tenrm Bd}_{{_{{\!\hbox{\fivebf Z}}}}_{_{^{\!\hbox{\fivebf
2}}}}} (\hbox{\tenrm Bd}_{{_{{\!\hbox{\fivebf
Z}}}}_{_{^{\!\hbox{\fivebf 2}}}}}\!\!\!\Sigma)= \emptyset$ since
$\hbox{\tenrm Bd}_{{_{{\!\hbox{\fivebf Z}}}}_{_{^{\!\hbox{\fivebf
2}}}}}\!\!\!\Sigma$ is Gorenstein$_{{_{{\!\hbox{\fivebf
Z}}}}_{_{^{\!\hbox{\fivebf 2}}}}}\!\!$ and so, $\dim\hbox{\tenrm
Bd}_{{_{{\!\hbox{\fivebf Z}}}}} (\hbox{\tenrm
Bd}_{{_{{\!\hbox{\fivebf Z}}}}}\!\Sigma)\le n-4 $ by Cor.\ 1 p.\
26.
So if $\hbox{\tenrm Bd}_{{_{{\!\hbox{\fivebf Z}}}}} (\hbox{\tenrm
Bd}_{{_{{\!\hbox{\fivebf Z}}}}}\!\Sigma)\neq\emptyset$ then, by
Cor.\ 2{\bf ii}\ p.\ 26,
$ \hbox{\tenrm Bd}_{{_{{\!\hbox{\fivebf Z}}}}}\!\Sigma$ is
nonorientable$_{{_{{\!\hbox{\fivebf Z}}}}}$ i.e. $
\mdoubleH\!\!\!_{_{^{\!n-1}}}\! (\hbox{\tenrm
Bd}_{{_{{\!\hbox{\fivebf Z}}}}}\!\Sigma, \hbox{\tenrm
Bd}_{{_{{\!\hbox{\fivebf Z}}}}}(\hbox{\tenrm
Bd}_{{_{{\!\hbox{\fivebf Z}}}}}\!\Sigma);\hbox{\tenbf Z})
$%
$\ \!= \ \!$%
$\mdoubleH\!\!\!_{_{^{\!n-1}}}\! (\hbox{\tenrm
Bd}_{{_{{\!\hbox{\fivebf Z}}}}}\!\Sigma;\hbox{\tenbf Z})\!=\!0 $
and the torsion submodule of $ \mdoubleH\!\!\!_{_{^{\!n-2}}}\!
(\hbox{\tenrm Bd}_{{_{{\!\hbox{\fivebf Z}}}}}\!\Sigma,
\hbox{\tenrm Bd}_{{_{{\!\hbox{\fivebf Z}}}}}(\hbox{\tenrm
Bd}_{{_{{\!\hbox{\fivebf Z}}}}}\!\Sigma);\hbox{\tenbf Z}) =
\Big[{\hbox{\eightrm For dimen-}\atop \hbox{\eightrm sional
reasons.}}\Big]$
$ =\!\mdoubleH\!\!\!\!_{_{^{\!n-2}}}\! (\hbox{\tenrm
Bd}_{{_{{\!\hbox{\fivebf Z}}}}}\!\Sigma;\hbox{\tenbf Z}) $ is
isomorphic to $\hbox{\tenbf
Z}_{_{^{\!\hbox{\fivebf 2}}}}\!$ by Lemma.\ $\!${\bf i} p.\ 24. %
In particular $ \mdoubleH\!\!\!\!_{_{^{\!n-2}}}\! (\hbox{\tenrm
Bd}_{{_{{\!\hbox{\fivebf Z}}}}}\!\Sigma;\hbox{\tenbf Z}) \otimes
\hbox{\tenbf Z}_{_{^{\!\hbox{\fivebf 2}}}}\!\ne\!0$ implying, by
Th.\ 5 p.\ 10, that
$ \mdoubleH\!\!\!\!_{_{^{\!n-2}}}\! (\hbox{\tenrm
Bd}{_{_{{\!\hbox{\fivebf Z}}}}}_{_{^{\!\hbox{\fivebf
2}}}}\!\!\!\Sigma;\hbox{\tenbf Z}_{_{^{\!\hbox{\fivebf 2}}}}) \!=
\mdoubleH\!\!\!\!_{_{^{\!n-2}}}\! (\hbox{\tenrm
Bd}{_{_{{\!\hbox{\fivebf Z}}}}}\!\Sigma;\hbox{\tenbf
Z}_{_{^{\!\hbox{\fivebf 2}}}}) \!=
\mdoubleH\!\!\!\!_{_{^{\!n-2}}}\!(\hbox{\tenrm
Bd}{_{_{{\!\hbox{\fivebf Z}}}}}\!\Sigma;\hbox{\tenbf Z}) \otimes
\hbox{\tenbf Z}{_{_{^{\!\hbox{\fivebf 2}}}}}\! \oplus \hbox{\tenrm
Tor}{^{^{_{\!\hbox{\fivebf Z}}}}}\!\!\!{_{_{^{\!\hbox{\fivebf
1}}}}} \big(\mdoubleH\!\!\!_{_{^{\!n-3}}}\! (\hbox{\tenrm
Bd}{_{_{{\!\hbox{\fivebf Z}}}}}\!\Sigma;\hbox{\tenbf
Z}),\hbox{\tenbf Z}{_{_{^{\!\hbox{\fivebf 2}}}}}\!
\big)\nobreak\!\ne\nobreak0 $
contradicting the Gorenstein$_{{_{{\!\hbox{\fivebf
Z}}}}_{_{^{\!\hbox{\fivebf 2}}}}}\!\!$-ness of $\hbox{\tenrm
Bd}\!\raise0.5pt\hbox{$_{_{_{\hbox{\fivebf
Z}_{^{\!2}}}}}$}\!\!\!_{\!}\Sigma.$
\qed
\enddemo

By Proposition\ 1 p.\ 18 we now get.
\proclaim{Corollary 2}
Each simplicial homology$_{_{\!\hbox{\fivebf Z}}}\!\!$
$n$-manifold $\Sigma$ is locally orientable.
\qed
\endproclaim

%


\head
Appendix: Simplicial Calculus and Simplicial Sets
\endhead

The complex$_{_{^o}}\!\!$ of all subsets of a simplex$_{_{^o}}\!$
$\sigma$ is denoted ${\bar{\sigma}}$, while the {\it boundary
of}\nobreak\ $\sigma$, ${\dot{\sigma}}$, is the set of all proper
subsets.
(${ \dot{\sigma}}\!:=\{\tau\ |\ \tau{{\lower2pt\hbox {$\subset$}}
\atop{\ne}} \sigma \}={\bar{\sigma}}{\raise0.5pt\hbox{\eightmsbm
\char"72}} \{\sigma\}$,
 ${ \bar{\emptyset}}\!_{_{^o}}\!\!=\!\{\emptyset\!_{_{^o}}\!\}$ and
${ \dot{\emptyset}}\!_{_{^o}}\!\!=\!\emptyset$.)

\noindent ``The {\it closed star} of $\sigma$ w.r.t. ${\Sigma}"=$
$\overline{\hbox{\tenrm {st}}}_{_{{\!\Sigma}}}\! \sigma\! :=
\{\tau\in \Sigma|\ \!\sigma\cup \tau\! \in \Sigma\}. $

\smallskip
\noindent ``The {\it open star} of $\sigma$ w.r.t. $\Sigma"\!=\!$
${\hbox{\tenrm {st}}}_{_{{\!\Sigma}}}\!\sigma\!:=\!
\{\alpha\in|\Sigma|\ |\ [\hbox{\tenrm v}\in\sigma] \Longrightarrow
[\alpha(\hbox{\tenrm v})\neq0] \}.$
So,
$\alpha_0\not\in{\hbox{\tenrm {st}}}_{_{{\!\Sigma}}}\!\sigma$
except for
${\hbox{\tenrm {st}}}_{_{{\Sigma}}}\!
\emptyset\!_{_{^o}}\!=|\Sigma|.$
$\bigl({\hbox{\tenrm {st}}}_{_{{\Sigma}}}\! (\sigma\!_{_{ }})\!
=\!\{{\rlap{$\alpha$}{\ \ \in}}|\Sigma|\big\vert{\rlap{$\alpha$}
{\ \ \in}} {\rlap{$_{_{_{_{_{v\in\sigma\!\!_{_{ }}}}}}}$}{\
\raise0pt\hbox{$\bigcap$}}}\ {{\hbox{\tenrm {st}}}_{_{
\!{\Sigma\!_{_{\ }}}}}\!\!\{v\}} \}.\bigr)$

\noindent ``The {\it closure} of $\sigma$ w.r.t. $\Sigma"\!=\!$
$|\sigma|\!:=\!\{\alpha\in|\Sigma|\ |\ { [\alpha(\hbox{\tenrm
v})\!\neq0] \Longrightarrow [\hbox{\tenrm v}\!\in\! \sigma] \}.}$
So, $|\emptyset\!_{_{^o}}\!|\!:= \{\alpha_0\}.$

\smallskip
\noindent ``The {\it interior} of $\sigma$ w.r.t.
${\Sigma}"\!\!=\! \hbox{\tensy h}\sigma\hbox{\tensy i} \!=\!
\hbox{\tenrm Int} (\sigma) \!:=\! \{\alpha\in|\Sigma|\ |\
[\hbox{\tenrm v}\in\nobreak\sigma ] \Longleftrightarrow
[\alpha(\hbox{\tenrm v})\neq0] \}.$
So,
$\hbox{\tenrm Int}_{_{\!^o}}\!(\sigma\!_{_{\!^{\ \!\!}}}) $ is an
open subspace of $|\Sigma|$ \underbar{iff} $\sigma$ is a maximal
simplex in $\Sigma$ and
$\hbox{\tenrm Int}(\emptyset\!_{_{^o}}\!)\!:= \{\alpha_0\}$.

\smallskip
\noindent
The {\it barycenter} $\hat\sigma$ of\ $\sigma$ is the
$\alpha\!\in\!\hbox{Int}(\sigma)\!\in\!\vert\Sigma\vert$
fulfilling
\hbox{$\hbox{\tenrm v}\!\in\! \sigma
\!\Rightarrow\!
\alpha({\sevenrm v})\!=\!{1 \over{\#\sigma}}\!$ while
$\hat{\emptyset}\!_{_{^o}}\!\!:= _{\!}\alpha_0.$}

\smallskip%
\noindent
$\!$``The {\it link} of $\sigma$ w.r.t. $\Sigma"=$ $\hbox{\tenrm
{Lk}}_{_{\!_{\Sigma}}}\!\!\!\sigma\!:= \{\tau\!\in \Sigma|
[\sigma\cap \tau =\emptyset]\land [\sigma\cup \tau \in \Sigma]\}.$
So, $\hbox{\tenrm
Lk}\!\!_{_{_{\Sigma}}}\!\emptyset\!_{_{^o}}\!=\!\Sigma,$
$\sigma\!\in\!\hbox{\tenrm
Lk}\!\!_{_{_{\Sigma}}}\!\tau\Longleftrightarrow\tau\!\in\!\hbox{\tenrm
Lk}\!\!_{_{_{\Sigma}}}\!\sigma$
and
$\hbox{\tenrm
Lk}\!\!_{_{_{\Sigma}}}\!\sigma\!=\!\emptyset\!_{_{}}$
\underbar{iff} $\sigma\!\not\in\!\Sigma,$
while
$\hbox{\tenrm
Lk}\!\!_{_{_{\Sigma}}}\!\tau\!=\{\emptyset\!_{_{^o}}\!\} $
\underbar{iff} $\tau\!\in\!{\Sigma}$ is maximal.

\smallskip
\noindent
$ \Sigma_{1}\!\ast\! \Sigma_{2}:= \{\sigma\!_{1}\cup
\sigma\!_{2}\vert \sigma\!_{i}\in \Sigma_{i}\ (i=1,2)\}.\ \ \
\hbox{\tenrm In\ particular}\ \Sigma\ast\! \{\emptyset_o\!\}=
\{\emptyset_o\}\ast\! \Sigma = \Sigma. $

\proclaim{Proposition 1}
{{If} $\ V_{\Sigma_{1}}\cap V_{\Sigma_{2}}=\emptyset \ $ {then},

\smallskip
$
[\tau\in \Sigma_{1}\ast_o \Sigma_{2}]
\Longleftrightarrow
[\exists !\ \sigma_{i}\in \Sigma_{i}$ {so\ that} $\tau
=\sigma_{1}\cup \sigma_{2}],$
{\rm (Direct from definition.)}\

\noindent {and}

$ \hbox{\tenrm Lk}\!\!\!\!\!\!
 \lower1.1pt\hbox{$_{_{\Sigma_{1}\!\ast\Sigma_{2}}}$}\!\!(\sigma_{_{\!1}}\!
\cup\sigma_{_{\!2}})
=
(\hbox{\tenrm Lk}\!\!
\lower1.1pt\hbox{${_{_{\Sigma\!_{_{1}}}}}$}\!\!\!\sigma\!_{_{1}})\ast
(\hbox{\tenrm Lk}\!\!
\lower1.1pt\hbox{${_{_{\Sigma\!_{_{2}}}}}$}\!\!\!\sigma\!_{_{2}}).
$} \hfill $(${\rm Proved by bracket juggling}.$)$ \indent
\qed
\endproclaim

Any link is an iterated link of vertices and

\smallskip
\noindent
$\hbox{\tenrm Lk}
\!\!\!\!_{\lower2.4pt\hbox {$_{{{{\hbox{\sevenrm
Lk}}}}\!\!_{_{_{\Sigma}}}\!\!\sigma}$}}
\!\!\!\!\!\tau
\!=\!
\{\emptyset\!_{_{^o}}\!\}
\ast \hbox{\tenrm Lk}
\!\!\!\!_{\lower2.4pt\hbox {$_{{{{\hbox{\sevenrm
Lk}}}}\!\!_{_{_{\Sigma}}}\!\!\sigma}$}}
\!\!\!\!\tau
\!=\! \hbox{\tenrm Lk}
\!\!_{\lower2.4pt\hbox {$_{{{\bar{\hbox{\seveni {\char"1B}}}}}}$}}
\sigma\!\ast
\hbox{\tenrm Lk}
\!\!\!\!_{\lower2.4pt\hbox {$_{{{{\hbox{\sevenrm
Lk}}}}\!\!_{_{_{\Sigma}}}\!\!\sigma}$}}
\!\!\!\!\tau
\!=\!\!
\big[\!{\hbox{\eightrm Prop.\ 1}\atop{\hbox{\eightrm above}}}
\!\big]\!\!=\! \hbox{\tenrm Lk}\!\!
\!\!\!\!_{\lower2.4pt\hbox {$_{{{{\bar{\hbox{\seveni
{\char"1B}}}\ast\hbox{\sevenrm
Lk}}}}\!\!_{_{_{\Sigma}}}\!\!\sigma}$}}
\!\!\!\!\!\!(\sigma
\hbox{\sevensy {\char"5B}}
\tau\!)
\!=\!
\hbox{\tenrm Lk}\!\!\!\!_{\hbox{$_{_{{\overline{\hbox{\sevenrm
st}}}\!\!_{_{_{\Sigma}}}\!\!\hbox{\seveni {\char"1B}}}}$}}
\!\!\!(\sigma
\hbox{\sevensy {\char"5B}}
\tau)
\!=\!
\hbox{\tenrm Lk}\!\!\!\!\!\!_{\hbox{$_{_{{
\overline{\hbox{\sevenrm st}}}\!\!_{_{_{\Sigma}}}\!\!
(\hbox{\seveni {\char"1B}}
\hbox{\fivesy {\char"5B}}
\hbox{\seveni {\char"1C}}
)
}}$}} \!\!\!\!\!\!\!\!\!\!\hbox{\sevenrm(}\sigma
\hbox{\sevensy {\char"5B}}
\tau\hbox{\sevenrm)}
\!=\!
\hbox{\tenrm Lk}\!_{_{\Sigma}} \!(\sigma
\hbox{\sevensy {\char"5B}}
\tau).$

\smallskip
\noindent %
So, $\tau\!\notin\!\hbox{\tenrm
Lk}\!_{_{\Sigma}}\!\sigma\!\Rightarrow\!$
$\hbox{\tenrm Lk} \!\!\!\lower1.1pt\hbox{$\!_{_{\hbox{\fiverm
Lk}\!_{_{\Sigma}}\!\sigma}}$} \!\!\!\!\!\tau
\!=\!
\emptyset$\ while $\tau\!\in\!\hbox{\tenrm
Lk}\!\!_{_{_{\Sigma}}}\!\sigma\Rightarrow $
$\hbox{\tenrm Lk}
\!\!\!\!_{\lower2.4pt\hbox {$_{{{{\hbox{\sevenrm
Lk}}}}\!\!_{_{_{\Sigma}}}\!\!\sigma}$}}
\!\!\!\!\!\tau
\!=\!
\hbox{\tenrm Lk}\!\!_{_{_{\Sigma}}} \!(\sigma\!\cup\!\tau)
\!=\!
\hbox{\tenrm Lk}
\!\!\!\!_{\lower2.4pt\hbox {$_{{{{\hbox{\sevenrm
Lk}}}}\!\!_{_{_{\Sigma}}}\!\!\tau}$}}
\!\!\!\!\!\sigma
\ \big(\!\!\subset\! \hbox{\tenrm Lk}\!_{_{\Sigma}}\!\sigma
\!\cap\! \hbox{\tenrm Lk}\!_{_{\Sigma}}\!\tau\big).$ \hfill({\bf
I})

\smallskip
\noindent
$\tau\!\in\!\hbox{\tenrm Lk}\!\!_{_{_{\Sigma}}}\!\sigma
\Longrightarrow
\overline{\hbox{\tenrm st}}
\!\!\!\!_{\lower2.4pt\hbox {$_{{{{\hbox{\sevenrm
Lk}}}}\!\!_{_{_{\Sigma}}}\!\!\tau}$}}
\!\!\!\!\sigma
\!=\!
\bar\sigma\ast
\hbox{\tenrm Lk}
\!\!\!\!_{\lower2.4pt\hbox {$_{{{{\hbox{\sevenrm
Lk}}}}\!\!_{_{_{\Sigma}}}\!\!\tau}$}}
\!\!\!\!\sigma
\!=\!
[\hbox{{\eightrm Eq.$\ \!${I}}}]
\!=\!
\bar\sigma
\ast \hbox{\tenrm Lk}
\!\!\!\!_{\lower2.4pt\hbox {$_{{{{\hbox{\sevenrm
Lk}}}}\!\!_{_{_{\Sigma}}}\!\!\sigma}$}}
\!\!\!\!\tau
\!=\! \hbox{\tenrm Lk}
\!\!_{\lower2.4pt\hbox {$_{{{\bar{\hbox{\seveni {\char"1B}}}}}}$}}
\emptyset\!_{_{^o}} \!\ast \hbox{\tenrm Lk}
\!\!\!\!_{\lower2.4pt\hbox {$_{{{{\hbox{\sevenrm
Lk}}}}\!\!_{_{_{\Sigma}}}\!\!\sigma}$}}
\!\!\!\!\tau
\!=\!\!
\big[\!{\hbox{\eightrm Prop.\ 1}\atop{\hbox{\eightrm above}}}
\!\big]\!
\!=\!
\hbox{\tenrm Lk}\!\!
\!\!\!\!_{\lower2.4pt\hbox {$_{{{{\bar{\hbox{\seveni
{\char"1B}}}\ast\hbox{\sevenrm
Lk}}}}\!\!_{_{_{\Sigma}}}\!\!\sigma}$}}
\!\!\!\!\!\!
\tau
\!=\!
\hbox{\tenrm Lk}
\!\!\!_{\lower2.4pt\hbox {$_{{{\overline{\hbox{\sevenrm
st}}}}\!\!_{_{_{\Sigma}}}\!\!\sigma}$}}\!\!\!
\tau.$

\smallskip
\noindent %
Put; $n:=\dim\Delta,$
${
\Delta\!^{^{_{(_{\!}p_{\!})}}}\!\!\!:= \{\delta\!\in\!\Delta\
\!\vert\ \!\#\delta\le p+1\},\ \!
\Delta\!^{^{_{(\!n{_{_{^{_{}}}}}\!\!_{\!})}}}\!\!\!=\!\Delta,\ \!
\Delta\!^{^{_{\prime\!}}}\!:=
\Delta\!^{^{_{(\!n{_{_{^{_{\!}}}}}\!-\!1_{\!})}}}\!
}\!$\
\hbox{\tenrm and}\
$\Delta\!^{^{_{p\!}}}\!\!:= \Delta\!\!^{^{_{(_{\!}p_{\!})}}}\!
{\raise1.5pt\hbox{\eightmsbm \char"72}}            
\Delta\!\!^{^{_{(_{\!}p\!-\!1_{\!})}}}\!\!.$

\smallskip
\noindent %
$\Gamma\subset\Sigma$ is {\it full} in $\Sigma$ if for all
$\sigma\!\in\!\Sigma$; $\sigma\subset V_{_{^\Gamma}}
\Longrightarrow
\sigma\in \Gamma.$

\proclaim{Proposition 2}
{\bf a.} {\bf i)}\ $\!{\hbox{\tenrm Lk}\!\!\!\!\!\!\!
\lower1.5pt\hbox{$_{_{{ \Delta\!_{_{{1\!\!}}}\cup
{\Delta_{_{{\!{{2}}}}}}}}}$}\!\!\!\!\delta} =\ \!{\hbox{\tenrm
Lk}\! \lower0.5pt\hbox{$ _{_{{ {\Delta_{_{{\!{{1}}}}}}}}} $}
\!\!\!\!\delta} \cup {\hbox{\tenrm Lk}\! \lower0.5pt\hbox{$ _{_{{
{\Delta_{_{{\!{{2}}}}}}}}} $} \!\!\!\!\delta}\ ,\ \
\hbox{\tenbf ii)}$\ ${\hbox{\tenrm Lk}\!\!\!\!\!\!\!
\lower1.5pt\hbox{$_{_{{ \Delta\!_{_{{1\!\!}}}\cap
{\Delta_{_{{\!{{2}}}}}}}}}$}\!\!\!\!\delta} =\  _{\!}
\!{\hbox{\tenrm Lk}\! \lower0.5pt\hbox{$ _{_{{
{\Delta_{_{{\!{{1}}}}}}}}} $} \!\!\!\!\delta} \cap {\hbox{\tenrm
Lk}\! \lower0.5pt\hbox{$ _{_{{ {\Delta_{_{{\!{{2}}}}}}}}} $}
\!\!\!\!\delta},$ \ \

\noindent
\hbox{\tenbf iii)} %
$\big({{\Delta\!^{^{\!_{\ \!}}}}\!_{_{{1}}}\!{{\ast}}}
{{\Delta\!^{^{\!_{\ \!}}}\!_{_{{2}}}}} \big)
\!^{^{{\hbox{$_{\prime}$}\!}}}\!\!=_{\!}
({{\Delta\!^{^{_{{\hbox{$_{\prime}$}\!}}}}\!_{_{1}}}}\!{{\ast}}
\Delta\!^{^{\!_{\ \!}}}\!_{_{{2}}}) \cup (\Delta\!^{^{\!_{\
\!}}}\!_{_{{1}}}\!{{\ast}}
{{\Delta\!^{^{_{{\hbox{$_{\prime}$}\!}}}}\!_{_{2}}}}), $
\noindent
{\bf iv)} $\Delta\ast(\ \!
{\cap\!\!\!\!\!_{_{_{{{i\in\hbox{\sevenbf I}}}}}}}\!
{\Delta\!^{^{\!_{\ \!}}}}\!_{_{{i}}}\!)\!=
{\cap\!\!\!\!\!_{_{_{{{i\in\hbox{\sevenbf I}}}}}}}
(\Delta\ast{\Delta\!^{^{\!_{\ \!}}}}\!_{_{{i}}}) $,\
$ \Delta\ast(\ \! {\cup\!\!\!\!\!_{_{_{{{i\in\hbox{\sevenbf
I}}}}}}}\! {\Delta\!^{^{\!_{\ \!}}}}\!_{_{{i}}})\! =
{\cup\!\!\!\!\!_{_{_{{{i\in\hbox{\sevenbf I}}}}}}}\!\!
(\Delta\ast{\Delta\!^{^{\!_{\ \!}}}}\!_{_{{i}}}).$

{\rm (iv holds also for topological spaces under the
{\rlap{{\lower2.5pt\hbox
{\vbox{\moveright0.17pt\hbox{$^{^{\land}}$}}}}}{{$\ast$}}}-join.)}

\noindent{\bf b.}\ $\Delta\ {pure}\ \!\Longleftrightarrow\!
\big({\hbox{\tenrm Lk}\!_{_{\Delta}}\!\!{{\delta}}}\big)
\!^{^{{\hbox{$_{\prime}$}\!}}} \!= {\hbox{\tenrm
Lk}\!\!_{_{\Delta\!^{^{_{\prime\!\!}}}}}{{\delta}}}\ \ \forall\
\emptyset\!_{_{^{o}}}\!\!\ne\! \delta\!_{_{}}\in\!\Delta.$

\smallskip
\noindent{\bf c.}\ $ [{\Gamma} \cap \hbox{\tenrm
Lk}\!_{_{\Delta}}\!\delta = \hbox{\tenrm
Lk}_{_{\!_{\Gamma}}}\!\!\delta $ $\ \! \forall\ \!
\delta\in{\Gamma}] \!\Longleftrightarrow\! [{\Gamma} $ is full in
$\Delta]$ $ \!\Longleftrightarrow\! [{\Gamma} \cap
\overline{\hbox{\tenrm st}}\!_{_{\Delta}}\!\delta =
\overline{\hbox{\tenrm st}}_{_{\!_{\Gamma}}}\!\!\delta $ $\ \!
\forall\ \! \delta\in{\Gamma}]. $
\hfill {\rm({\bf II})}
\endproclaim

\demo{Proof}
{\bf a.} Associativity and distributivity of the logical
connectives. \hfill$\triangleright$

\noindent{\bf b.}\ $(\Longrightarrow)$ $ \Delta\ \hbox{\tenrm
pure}\ \Longrightarrow {\Delta\!^{^{\!_{\prime}}}\!}\ \hbox{\tenrm
pure},\ \hbox{\tenrm so;} $ \noindent ${\hbox{\tenrm
Lk}\!_{_{\Delta\!^{^{\!_{\prime}}}\!}}\!{{\delta}}} \!=\! \{
\tau\!\in\!\Delta {|{_{\!}} _{\!\!}|} \
\!\!\tau\!\cap\!{\delta}\!=\!\emptyset \land\
\!\!\tau\!\cup\!{\delta}\!\in\!\Delta\!^{^{\!_{\prime}}}\! \}
\!=\! \{ \tau\!\in\!\Delta{|{_{\!}} _{\!\!}|}\
\!\!\tau\!\cap\!{\delta}\!=\!\emptyset \land\
\!\!\tau\!\cup\!{\delta}\!\in\!\Delta\ \land\
\!\#(\tau\!\cup\!{\delta})\le n\!_{_{^{\ \!}}} \} \!= $ $ \{
\tau\!\in\!{\hbox{\tenrm Lk}\!_{_{\Delta}}\!\!{{\delta}}}\ \!
{|{_{\!}} _{\!\!}|}\ \!\#(\!\tau\!\cup\!{\delta})\!\le\!
n\!_{_{^{\ \!}}} \!\} \!=\! \{ \tau\!\in\!{\hbox{\tenrm
Lk}\!_{_{\Delta}}\!\!{{\delta}}}\ \! {|{_{\!}} _{\!\!}|}\
\!\#\tau\!\le\! n\!_{_{^{\ }}}\!\!\!-\!\#{\delta} \} \!=\! \{
\tau\!\in\!{\hbox{\tenrm Lk}\!_{_{\Delta}}\!\!{{\delta}}}\ \!
{|{_{\!}} _{\!\!}|}\ \!\#\tau\!\le\! \dim{{{{\hbox{\tenrm
Lk}\!_{_{\Delta}}\!\!{{\delta}}}}}}
\} \!=\! \big({\hbox{\tenrm Lk}\!_{_{\Delta}}\!\!{{\delta}}}\big)
\!^{^{{\hbox{$_{\prime}$}\!}}}\!.$

\noindent $(\Longleftarrow)$ If $\Delta$ non-pure, then $\exists\
\!\delta\!_{_{^{\hbox{\fiverm m}}}}\!\!\!
\in\!{\Delta}\!^{^{{\hbox{$_{\prime}$}\!}}}$ maximal in both
${\Delta}\!^{^{{\hbox{$_{\prime}$}\!}}}$ and ${\Delta}$ i.e.,

\hfill
$_{\!}\big({\hbox{\tenrm
Lk}\!_{_{\Delta}}\!\!{{\delta}\!_{_{^{\hbox{\fiverm m}}}}\!}}\big)
\!^{^{{\hbox{$_{\prime}$}\!}}} \!=\!
\big(\{\emptyset\!_{_{o}}\}\big)\!^{^{{\hbox{$_{\prime}$}\!}}}
\!=\! \emptyset \!\not=\! \{\emptyset\!_{_{o}}\} \!=\!
{\hbox{\tenrm
Lk}\!_{_{\Delta\!^{^{\!_{\prime}}}\!}}\!{{\delta}\!_{_{^{\hbox{\fiverm
m}}}}\!}}. \ \hfill\triangleright$

\noindent{\bf c.} $ \hbox{\tenrm Lk}\!_{_{_{\Gamma}\!\!}}\sigma
\!= \{\tau\vert\ \tau\cap\sigma=\emptyset\ \&\
\tau\cup\sigma\in\Gamma\} = \!\big[ {{\hbox{\eightrm true}\
\forall\ \!\!\sigma\in\Gamma} \atop {\underline{\hbox{\eightrm
iff}}\ \Gamma\ \hbox{\eightrm full}}} \big]\! =
\{\tau\in\Gamma\vert\ \tau\cap\sigma =\emptyset\ \&\
\tau\cup\sigma\!\in\!\Delta\}=\break
= \Gamma\cap\hbox{\tenrm Lk}\!_{_{\Delta\!\!}}\sigma. $
\indent
$ {\Gamma} \cap \overline{\hbox{\tenrm st}}\!_{_{\Delta}}\!\delta
= {\Gamma} \cap \{\tau\in\Gamma\vert\ \tau\cup\sigma\in\Delta\}
\cup \{\tau\in\Delta\vert\ \tau\notin\Gamma\ \&\
\tau\cup\sigma\in\Delta\} =\break =\{\tau\in\Gamma\vert\
\tau\cup\sigma\in\Delta\}$
$ = \big[{{\forall \sigma\in\ \!\Gamma\ \underline{\hbox{\eightrm
iff}}\ {\Gamma}\ \hbox{\eightrm is\ a}}\atop{\hbox{\eightrm full\
subcomplex.}}}\big] = \{\tau\in\Gamma\vert\
\tau\cup\sigma\in\Gamma\} = \overline{\hbox{\tenrm
st}}_{_{\!_{\Gamma}}}\!\!\delta.
$ \qed
\enddemo

\goodbreak

\noindent
The {\it contrastar} of\ $\sigma\!\in\!\Sigma\!=\hbox{\tenrm
cost}_{_{\!{\Sigma}}}\!\sigma\!:=\! \{\tau\!\in\!\Sigma|\
\tau\!\not\supseteq\sigma\}.$
$\hbox{\tenrm
cost}_{_{{\!\Sigma}}}\!\emptyset\!_{_{^o}}\!\!=\!\emptyset\!$
\ and\ %
${\hbox{\tenrm cost}_{_{{\!\Sigma}}}\!\sigma\!=\!\Sigma\
\underline{\hbox{\tenrm iff}}\ \sigma\!\not\in\!\Sigma}.$

\proclaim{Proposition 3}
Without assumptions whether
${{\delta\!_{_{1}}}}\!,{{\delta\!_{_{2}}}} \!\!\in\!\Delta$ or
not, the following holds;

\noindent
{\bf a.}\ $ \hbox{\tenrm
cost}\!_{_{\Delta}}\!({{\delta\!_{_{1}}}}\!\!\cup\!{{\delta\!_{_{2}}}})
= \hbox{\tenrm cost}\!_{_{\Delta}}\!{{\delta\!_{_{1}}}} \!\cup
\hbox{\tenrm cost}\!_{_{\Delta}}\!{{\delta\!_{_{2}}}}$ \ \ \
and
\ \ \ $\delta\!=\!\{\hbox{\tenrm v}\!_{_{^{1}}},...,\hbox{\tenrm
v}\!_{_{^{p}}}\} \Longrightarrow \hbox{\tenrm
cost}\!_{_{\Delta}}\!{{\delta\!_{_{ }}}}
=\bigcup\!\!\!\!\!\!\!_{_{_{_{i=1,p}}}} \!\!\!\!\hbox{\tenrm
cost}\!_{_{\Delta}}\!\!{\hbox{\tenrm v}\!_{_{^{i}}}}. $

\noindent {\bf b.}\ $ \hbox{\tenrm cost}\!\!\!\!\!\!\!
_{\lower0.8pt\hbox{$ {_{\hbox{\eightrm
cost}\!\!_{_{\Delta}}\!\!{{\delta\!_{_{1}} }}  } }$}}\!
\!{{\delta\!_{_{2}}}}\! = \hbox{\tenrm
cost}\!_{_{\Delta}}\!{{\delta\!_{_{1}}}} \!\cap \hbox{\tenrm
cost}\!_{_{\Delta}}\!{{\delta\!_{_{2}}}}\! = \hbox{\tenrm
cost}\!\!\!\!\!\!\! _{\lower0.8pt\hbox{$ {_{\hbox{\eightrm
cost}\!\!_{_{\Delta}}\!\!\!{{\delta\!_{_{2}} }} } }$}}\!
\!{{\delta\!_{_{1}}}} \ \ \hbox{\tenrm and}\ \
\bigcap\!\!\!\!\!\!\!\!_{_{_{_{i=1,q}}}} \!\!\!\hbox{\tenrm
cost}\!_{_{\Delta}}\!\!{{\delta}\!_{_{^{i}}}}\! = \hbox{\tenrm
cost} \!\!\!\!\!\!\!\!\! _{\lower0.8pt\hbox{$ {_{\hbox{\eightrm
cost}\!_{_{ {^{^{\hbox{\tenbf .}}}\!{\hbox{\tenbf
.}}_{_{\hbox{\tenbf .}}}}
}}\!\!\!\!\!\!\!{{\delta}\!_{_{^{2}}}}}}\!\! $}}
{_{\!}}{{\delta}\!_{_{^{1}}}}. $

\nointerlineskip
$ \hfill\hfill\hfill\hfill\indent\indent {{ _{\lower0.8pt\hbox{$
{_{{\hbox{\fiverm cost}\!\!_{_{\Delta}}\!\!{{{
\hbox{\fivei{\char"0E}}
}\!_{_{\hbox{\fivei q}}}}}}}} $}} }} \hfill \indent\ \ \
$

\nointerlineskip \noindent {\bf c.} $ \
\cases
\hbox{\tenbf i.}\ \ {\delta}\!\notin\!{\Delta} \Longleftrightarrow
\hbox{\tenbf[} \hbox{\tenrm Lk}\!\!\!\!\!\! {\lower0.4pt\hbox{${
_{_{\hbox{\eightrm cost}\!\!_{_{\Delta}}\!\!\!{\hbox{\eightbf v}}
}}\!\!\!\!{{\delta\!_{_{ }}}} }$}}
=
\hbox{\tenrm cost}\!\!\!\!\!
\lower1pt\hbox{$_{_{{\!\!\hbox{\eightrm Lk}\!\!_{_{\Delta}}
}}}\!\!\lower3pt\hbox{$_{\delta\!}$} $}\! {\hbox{\tenrm v}}
=
\emptyset \hbox{\tenbf]}. \cr
\hbox{\tenbf ii.}\ {\delta}\!\in\!{\Delta} \Longleftrightarrow
\hbox{\tenbf[} \hbox{\tenrm v}\notin{\delta} \Longleftrightarrow
{\delta}\in\hbox{\tenrm cost}\!_{_{\Delta}}\!\!{\hbox{\tenrm v}}
\Longleftrightarrow \hbox{\tenrm Lk}\!\!\!\!\!\!
{\lower0.4pt\hbox{${ _{_{\hbox{\eightrm
cost}\!\!_{_{\Delta}}\!\!\!{\hbox{\eightbf v}}
}}\!\!\!\!{{\delta\!_{_{ }}}} }$}}
=
\hbox{\tenrm cost}\!\!\!\!\!
\lower1pt\hbox{$_{_{{\!\!\hbox{\eightrm Lk}\!\!_{_{\Delta}}
}}}\!\!\lower3pt\hbox{$_{\delta\!}$} $}\! {\hbox{\tenrm v}}
\!\supset\!
\{\emptyset\}\ (\neq\!\!\emptyset) \hbox{\tenbf]}.
\cr
\phantom{\hbox{\tenbf (3)}}
\big(\!\!\!\big( {\delta}\!\in\!{\Delta} \Longleftrightarrow
\hbox{\tenbf[} \hbox{\tenrm v}\in{\delta} \Longleftrightarrow
{\delta}\notin\hbox{\tenrm cost}\!_{_{\Delta}}\!\!{\hbox{\tenrm
v}} \Longleftrightarrow \hbox{\tenrm Lk}\!\!\!\!\!\!
{\lower0.4pt\hbox{${ _{_{\hbox{\eightrm
cost}\!\!_{_{\Delta}}\!\!\!{\hbox{\eightbf v}}
}}\!\!\!\!{{\delta\!_{_{ }}}} }$}} = \emptyset \neq \{\emptyset\}
\subset \hbox{\tenrm Lk}\!_{_{\Delta}}\!\!{{\delta\!_{_{ }} }} =
\hbox{\tenrm cost}\!\!\!\! \lower1pt\hbox{$_{_{{\!\!\hbox{\eightrm
Lk}\!\!_{_{\Delta}}\!\!{{{\hbox{\eighti {\char"0E}}} }} } }}$}\!
{\hbox{\tenrm v}} \hbox{\tenbf]}. \big)\!\!\!\big) \cr
\endcases
$

\smallskip
\noindent
{\bf d.}\ $ {\delta},\tau\!\in\!{\Delta} \!\Longrightarrow\!
\hbox{\tenbf[} {\delta}\in{\hbox{\tenrm
cost}\!_{_{\Delta}}\!\!{{\tau}}} \!\Longleftrightarrow\!
\emptyset\not=\hbox{\tenrm Lk}\!\!\!\!\!\!\! {\lower0.5pt\hbox{${
_{_{\!\hbox{\eightrm
cost}\!\!_{_{\Delta}}\!\!\!{\lower2.4pt\hbox{${^{_{{{\hbox{\eighti
{\char"1C}}}}}}}$}} }}\! }$}} \!{{\delta\!_{_{ }}}} = \hbox{\tenrm
Lk}\!\!\!\!\!\!\!\!\!\! {\lower1.4pt\hbox{$_{_{\!\hbox{\eightrm
cost}\!\!_{_{\Delta}}\!\!\!\!\! {\lower2.4pt\hbox{
${^{_{{{\hbox{{\sixrm (}{\eighti {\char"1C}}}}\
\!\!{{{{\setminus}}}}\ {\!}\!{\hbox{{\eighti {\char"0E}}{\sixrm
)}}}}}}}$} } }}$}}\! \!\!\!\!\!{{\delta\!_{_{ }}}}
= \hbox{\tenrm cost}\!\!\!\!
{\lower1.0pt\hbox{${ _{_{{\!\!\hbox{\eightrm
Lk}\!\!_{_{\Delta}}\!\!\!{{{{\hbox{\eighti {\char"0E}}}}} }}}}\!
}$}}
(\tau\lower1pt\hbox{$\ _{\!}\!^{_{_{\setminus}}}\
_{\!}\!\delta$})\hbox{\tenbf]}. $
\ \ \ {\tenbf(} {\tenrm Note:} {\bf d}
$
\Rightarrow \hbox{\tenrm Lk}\!\!\!\!\!\!\! {\lower0.5pt\hbox{${
_{_{\!\hbox{\eightrm
cost}\!\!_{_{\Delta}}\!\!\!{\lower2.4pt\hbox{${^{_{{{\hbox{\eighti
{\char"1C}}}}}}}$}} }}\! }$}} \!{{\delta\!_{_{ }}}} =
\hbox{\tenrm Lk}\!_{_{\Delta}}\!\!{{{\delta}} }$
\ \underbar{iff} \
$\delta {\hbox{\eightsy {\char"5B}}}
\tau\!\notin\!\Delta\!\!$
\ \ or equivalently, \underbar{iff} \
$\tau\lower1pt\hbox{$\ _{\!}\!^{_{_{\setminus}}}\
_{\!}\!\delta$}\notin\ \!\overline{\hbox{\tenrm
{st}}}\!_{_{\Delta}}\!\delta$.
{\bf)}

\smallskip
\noindent %
{\bf e.}\ {If} $\delta\ne\emptyset$ {then}; \
$ \cases
\hbox{\tenbf 1.}\ %
[\big({\hbox{\tenrm cost}\!_{_{\Delta\!\!}}{{\delta}}}\big)
\!^{^{{\hbox{$_{\prime}$}\!}}}
=
{\hbox{\tenrm
cost}\!\!_{_{\Delta\!^{^{_{\prime\!\!}}}}}{{\delta}}}]
\Longleftrightarrow [n\!_{_{{\delta}}}\!\!
=
\!n\!_{_{\Delta}}\ \!\!]
\cr
\hbox{\tenbf 2.}\
[{\hbox{\tenrm cost}\!_{_{\Delta\!\!}}{{\delta}}} \ =\
{\hbox{\tenrm
cost}\!\!_{_{\Delta\!^{^{_{\prime\!\!}}}}}{{\delta}}}]
\Longleftrightarrow [n\!_{_{{\delta}}}\!\! =\!n\!_{_{\Delta}}\
\!\!-1].
\cr
\endcases
$

\smallskip
\noindent
{\bf f.}\ $\ \!\overline{\hbox{\tenrm {st}}}\!_{_{\Delta}}\!\delta
\cap {\hbox{\tenrm cost}\!_{_{\Delta\!\!}}{{\delta}}}
=
\hbox{\tenrm cost}\!\!\!\!\lower3.5pt\hbox{${\phantom{\ }
\over{^{\hbox{\eightrm st}}}}$} {\lower2.4pt\hbox{${
_{_{^{{\!\!\!_{_{\Delta}}\!\!\!{{{{\hbox{\eighti {\char"0E}}}
}}} }}}}\! }$}} \!\!\delta = {\dot{\delta}}\ast {\hbox{\tenrm
Lk}\!_{_{\Delta\!\!}}{{\delta}}}. $
\endproclaim

\demo{Proof}
If $ \Gamma\!\subset\Delta $ then $ \hbox{\tenrm
cost}_{_{\Gamma}}\!{{\gamma}} \!=\! \hbox{\tenrm
cost}\!_{_{\Delta}}\!{{\gamma}} \cap {\Gamma}\ \forall\ \!\gamma
$
and $ \hbox{\tenrm
cost}\!_{_{\Delta}}\!({{\delta\!_{_{1}}}}\!\cup{{\delta\!_{_{2}}}})
\!=\! \{\tau\!\in\!\Delta \ \!\!{|{_{\!}} _{\!\!}|}\ \!
\neg[{{\delta\!_{_{1}}}}\!\cup{{\delta\!_{_{2}}}}\!\!\subset\!\!\tau]\}
\!=\! \{\tau\!\in\!{{\Delta}}\ \!\!\vert \
\!\neg[\delta\!_{_{1}}\!\!\subset\!\!\tau]\lor
\!\neg[\delta\!_{_{2}}\!\!\subset\!\!\tau] \} \!\!= $ $
\hbox{\tenrm cost}\!_{_{\Delta}}\!{{\delta\!_{_{1}}}} \!\cup
\hbox{\tenrm cost}\!_{_{\Delta}}\!{{\delta\!_{_{2}}}}, $ giving
{\bf a} and {\bf b}.
A ``brute force"-check gives {\bf c}, the $``\tau\!\!=\!\!$
$\{$v$\}"$-case of {\bf d}, while
$ {\tau}\!\in\!{\hbox{\tenrm Lk}\!_{_{\Delta}}\!{{\delta}}}\
\!\!\Longrightarrow\!\!{{\{\hbox{\tenrm v}\}}}\!\in\!{\hbox{\tenrm
Lk}\!_{_{\Delta}}\!{{\delta}}} \ \ \!\forall\ \!\hbox{\tenrm
\{v\}}\!\in\!{\tau}\in\Delta \!\Longleftrightarrow
{{\delta}}\!\in\!{\hbox{\tenrm Lk}\!_{_{\Delta}}\!\!{\hbox{\tenrm
v}}}\ \! (\subset\hbox{\tenrm cost}\!_{_{\Delta}}\!\!{\hbox{\tenrm
v}})\ \ \!\! \forall\ \!\hbox{\tenrm v}\!\in\!{\tau}\in\Delta $
gives {\bf d} from {\bf a}, {\bf c} and Proposition\ 2$\ \!\!$a
{\bf i} above.

With $ n\!_{_{\lower1pt\hbox{${^{_{{\varphi}}}}$}} }
\!\!:=\!\dim\hbox{\tenrm
cost}_{_{\!\Delta}}\!{\raise1pt\hbox{$\varphi$}} $ and
$n\!_{_{\Delta}}\!\!:=\!\dim\!\Delta<\infty,$ we get;
$\tau\!\subset\!\delta\Rightarrow n\!_{_{\Delta}}\!\!-\!1\ \!\le
n_{_{\!\raise1pt\hbox{${_{{{\tau}}}}$}} } \!\le
n\!_{_{\raise1pt\hbox{${_{{{\delta}}}}$}} } \!\le n\!_{_{\Delta}}
$
$\forall\ \!\emptyset\!_{_{^o}}\!\!\ne\!\tau\!,\delta$.
Now, for {\bf e};\ ${\hbox{\tenrm
cost}\!_{_{\Delta}}\!\!{{\delta}}}\!$ = [\underbar{iff}
$n\!_{_{{\delta}}}\! \!=\!n\!_{_{\Delta}}\!\!-1\ \!\!]$ =
$\underline{({\hbox{\tenrm cost}\!_{_{\Delta}}\!\!{{\delta}}})
{\hbox{\eightsy{\char"5C}}}
{\Delta\!^{^{\hbox{$_{\prime}$}\!}}}} = {\hbox{\tenrm
cost}\!\!_{_{\Delta\!^{^{_{\prime}}}\!}}\!{{\delta}}}$ $\ \!
\rlap{$^{^{_{\hbox{\sevenrm?}}}}$}{\!=}\ ({\hbox{\tenrm
cost}\!_{_{\Delta\!\!}}{{\delta}}})^{^{{\hbox{$_{\prime}$}\!}}}=$
=\nobreak\ [\underbar{iff} $n\!_{_{{\delta}}}\!
\!=\!n\!_{_{\Delta}}] =\ \! \underline{({\hbox{\tenrm
cost}\!_{_{\Delta}}\!\!{{\delta}}})
{\hbox{\eightsy{\char"5C}}}
{\Delta\!^{^{\hbox{$_{\prime}$}\!}}}},$
and
{\bf f};\ $\ \!\overline{\hbox{\tenrm
{st}}}\!_{_{\Delta}}\!\delta\!=\! {
\bar{\delta}}\!\ast\!\hbox{\tenrm {Lk}}\!_{_{\Delta\!\!}}\delta=
{\dot{\delta}}\!\ast\!\hbox{\tenrm {Lk}}\!_{_{\Delta\!\!}}\delta
\ {\hbox{\eightsy{\char"5B}}}\ 
\{\tau\!\in\!\Delta\ \!{|{_{\!}} _{\!\!}|}\
\!\delta\!\subset\!\tau \}. $
\qed
\enddemo

$ \hbox{\tenrm With}\ \delta\!_{_{^{\ \!\!}}}\!, \tau\!\!_{_{^{\
\!\!}}}\!\in\!\Sigma;$
$\delta\!_{_{^{\ \!\!}}}\cup \tau\!\!_{_{^{\
\!\!}}}\notin\!\Sigma\!_{_{^o}}\!
\!\Longleftrightarrow\!
\delta\!_{_{^{\ \!\!}}}\!\notin \overline{\hbox{\tenrm {st}}}\!\!
\lower1.1pt\hbox{${_{_{\Sigma\!_{_{ }}}}}$}\!\tau
\!\Longleftrightarrow\!
{\hbox{\tenrm {st}}}\!\! \lower1.1pt\hbox{${_{_{\Sigma\!_{_{
}}}}}$}\!\delta\ \!\cap\ \! |\overline{\hbox{\tenrm {st}}}\!\!
\lower1.1pt\hbox{${_{_{\Sigma\!_{_{ }}}}}$}\!\tau|
\!=\!\{\alpha\!_{_{^0}}{_{\!}}\}
\!\Longleftrightarrow
{\hbox{\tenrm {st}}}\!\! \lower1.1pt\hbox{${_{_{\Sigma\!_{_{
}}}}}$}\!\tau
\ \!\cap\ \!
{\hbox{\tenrm {st}}}\!\! \lower1.1pt\hbox{${_{_{\Sigma\!_{_{
}}}}}$}\!\delta
=\!
\{\alpha\!_{_{^0}}{_{\!}}\}\ \hbox{\tenrm and}
$

\noindent
$ |\overline{\hbox{\tenrm {st}}}\!\!
\lower1.1pt\hbox{${_{_{\Sigma\!_{_{ }}}}}$}\!\sigma\!_{_{\!^{\
\!\!}}}|
{\raise1.5pt\hbox{\eightmsbm \char"72}} 
{\hbox{\tenrm {st}}}\!\! \lower1.1pt\hbox{${_{_{\Sigma\!_{_{
}}}}}$}\!\sigma\!
=
|{\dot{\sigma}}\!_{_{\!^{\ \!\!}}}\ast
\hbox{\tenrm {Lk}}\!\! \lower1.1pt\hbox{${_{_{\Sigma\!_{_{
}}}}}$}\!\sigma|
=
|\overline{\hbox{\tenrm {st}}}\!\!
\lower1.1pt\hbox{${_{_{\Sigma\!_{_{ }}}}}$}\!\sigma\cap
{\hbox{\tenrm {cost}}}\!\! \lower1.1pt\hbox{${_{_{\Sigma\!_{_{
}}}}}$}\!\sigma|$,
{\rm by\ \cite{21} %
p.\ 372,\ 62.6\
and Proposition\ 3.{\bf f} above.
}

\noindent
{\rm({\bf III})}
$ \ \!\overline{\hbox{\tenrm {st}}}\!\!
\lower1.1pt\hbox{${_{_{\Sigma\!_{_{ }}}}}$}\!\sigma\!_{_{\!^{\
\!\!}}}
=
\{\tau\!_{_{\!^{\ \!\!}}}\in\Sigma\vert\ \sigma\cup\tau\in\Sigma\}
= {\bar{\sigma}}\!_{_{\!^{\ \!\!}}}\ast
\hbox{\tenrm {Lk}}\!\! \lower1.1pt\hbox{${_{_{\Sigma\!_{_{
}}}}}$}\!\sigma
$
{\rm and} identifying $|{\hbox{\tenrm {cost}}}\!\!
\lower1.1pt\hbox{${_{_{\Sigma}}}$}\!\sigma|$ with its homeomorphic
image in $|{\Sigma}|$ through;
$
|{\hbox{\tenrm {cost}}}\!\!
\lower1.1pt\hbox{${_{_{\Sigma}}}$}\!\sigma|
\simeq
|{\Sigma}| {\raise1.5pt\hbox{\eightmsbm \char"72}} 
{\hbox{\tenrm {st}}}\!\! \lower1.1pt\hbox{${_{_{\Sigma\!_{_{
}}}}}$}\!\sigma,
$
we get
$
{\hbox{\tenrm {st}}}\!\! \lower1.1pt\hbox{${_{_{\Sigma\!_{_{
}}}}}$}\!\sigma
\!=\!
|{\Sigma}| {\raise1.5pt\hbox{\eightmsbm \char"72}} 
|{\hbox{\tenrm {cost}}}\!\!
\lower1.1pt\hbox{${_{_{\Sigma}}}$}\!\sigma|.
$

\medskip
\centerline{------$\ast\ast\ast$------\ \
\------$\ast\ast\ast$------\ \ \
------$\ast\ast\ast$------}

\medskip
\noindent
$\!\!$\FFrame{0.0pt}{0.0pt}
{%
\hsize0.99\hsize
\tenrm
Definitions of {\teni the} {\teni category} {\teni of} {\teni
simplicial} {\teni sets}, {\tenrm former} semi-simplicial
complexes, usually uses {\teni the} {\teni category} {\teni of}
{\teni non}-{\teni empty} {\teni ordinals} but to comply with the
introduction of the categories
\hbox{\tensy D}{\lower1.5pt\hbox{\fivei{\char"7D}}} (\hbox{\tensy
K}{\lower1.5pt\hbox{\fivei o}})
in Ch.\ 2, we have to use {\teni the} {\teni category} {\teni of}
{\teni ordinals} i.e. we include the empty ordinal
{\tenrm{\char"01}}({\tensy{\char"3B}})        
in a consistent way.
Denote an ordered simplicial complex {\tenrm{\char"06}} 
when regarded as a {\teni simplicial} {\teni set} by
\rlap{\raise2.0pt\hbox{\tent{\char"5E}}}{\tenrm$_{\!}${\char"06}},
and let
\rlap{\raise2.0pt\hbox{\tent{\char"5E}}}{{\tenrm$_{\!}${\char"06}}$_{_{^{\!1}}}$}%
\rlap{\raise1.0pt\hbox{\tent{\char"5E}}}{\vbox{\moveleft1.35pt\hbox{\tensy{\char"02}}}} $\!$%
\rlap{\raise2.0pt\hbox{\tent{\char"5E}}}{{\tenrm$_{\!}${\char"06}}$_{_{^{\!2}}}\!$} %
be the semi-simplicial product of \
\rlap{\raise2.0pt\hbox{\tent{\char"5E}}}{{\tenrm$_{\!}${\char"06}}$_{_{^{\!1}}}\!$} 
and
\rlap{\raise2.0pt\hbox{\tent{\char"5E}}}{{\tenrm$_{\!}${\char"06}}$_{_{^{\!2}}}\!$}, %
while
{\rlap{\raise1.5pt\hbox{\tent{\char"5E}}}{\vbox{\moveleft-1.2pt\hbox{\tensy{\char"6A}}}} %
$\!$%
\rlap{\raise2.0pt\hbox{\tent{\char"5E}}}{\tenrm$_{\!}${\char"06}}
\rlap{\raise1.5pt\hbox{\tent{\char"5E}}}{\vbox{\moveleft-1.2pt\hbox{\tensy{\char"6A}}}} %
$\!$%
} is the Milnor realization of
\rlap{\raise2.0pt\hbox{\tent{\char"5E}}}{\tenrm$_{\!}${\char"06}}. 
\cite{9} %
p.\ $\!$160 Prop.\ 4.3.15 \raise1pt\hbox{\sixrm+} p.\ 165
Ex.\ 1\ $\!\!$\raise1pt\hbox{\sixrm+}2 gives;
\ (The same is true also for joins.)
\medskip
\centerline{
{\tensy{\char"6A}}%
{{\tenrm$_{\!}${\char"06}}$_{_{^{\!1}}}\!$}%
{\tensy{\char"6A}} ${\!}$%
\rlap{\raise1.0pt\hbox{\tent{\char"16}}}{\vbox{\moveleft1.35pt\hbox{\tensy{\char"02}}}}
{\tensy{\char"6A}}%
{{\tenrm$_{\!}${\char"06}}$_{_{^{\!2}}}\!$}%
{\tensy{\char"6A}}
{\tensy{\char"27}} 
{%
\rlap{\raise1.5pt\hbox{\tent{\char"5E}}}{\vbox{\moveleft-1.2pt\hbox{\tensy{\char"6A}}}} %
$\!$%
\rlap{\raise2.0pt\hbox{\tent{\char"5E}}}{{\tenrm$_{\!}${\char"06}}$_{_{^{\!1}}}$}\!\!
\rlap{\raise1.5pt\hbox{\tent{\char"5E}}}{\vbox{\moveleft-1.2pt\hbox{\tensy{\char"6A}}}} %
$\!$%
\rlap{\raise1.0pt\hbox{\tent{\char"16}}}{\vbox{\moveleft1.35pt\hbox{\tensy{\char"02}}}}
\rlap{\raise1.5pt\hbox{\tent{\char"5E}}}{\vbox{\moveleft-1.2pt\hbox{\tensy{\char"6A}}}} %
$\!$%
\rlap{\raise2.0pt\hbox{\tent{\char"5E}}}{{\tenrm$_{\!}${\char"06}}$_{_{^{\!2}}}$}\!\!
\rlap{\raise1.5pt\hbox{\tent{\char"5E}}}{\vbox{\moveleft-1.2pt\hbox{\tensy{\char"6A}}}} %
$\!$%
{\tensy{\char"27}}          
\rlap{\raise1.5pt\hbox{\tent{\char"5E}}}{\vbox{\moveleft-1.2pt\hbox{\tensy{\char"6A}}}} %
$\!$%
\rlap{\raise2.0pt\hbox{\tent{\char"5E}}}{{\tenrm$_{\!}${\char"06}}$_{_{^{\!1}}}$}\!%
\rlap{\raise1.0pt\hbox{\tent{\char"5E}}}{\vbox{\moveleft1.35pt\hbox{\tensy{\char"02}}}}%
\rlap{\raise2.0pt\hbox{\tent{\char"5E}}}{{\tenrm$_{\!}${\char"06}}$_{_{^{\!2}}}\!$}%
\rlap{\raise1.5pt\hbox{\tent{\char"5E}}}{\vbox{\moveleft-1.2pt\hbox{\tensy{\char"6A}}}} %
$\!$%
} %
{\tensy{\char"27}}     
{\tensy{\char"6A}}%
{{\tenrm$_{\!}${\char"06}}$_{_{^{\!1\!}}}\!$}%
{\tensy{\char"02}}
{{\tenrm$_{\!}${\char"06}}$_{_{^{\!2}}}\!$}%
{\tensy{\char"6A}} %
}
\medskip
{%
The Milnor realization
\rlap{\raise1.5pt\hbox{\tent{\char"5E}}}{\vbox{\moveleft-1.2pt\hbox{\tensy{\char"6A}}}} %
$\!$%
{\tenrm{\char"04}}%
\rlap{\raise1.5pt\hbox{\tent{\char"5E}}}{\vbox{\moveleft-1.2pt\hbox{\tensy{\char"6A}}}} %
$\!$%
of any simplicial
set {\tenrm{\char"04}} 
is {triangulable} by \cite{9} %
p.\ 209 Cor.\ 4.6.12. E.g; the augmental complex
\tenrm{\char"01}}$\!_{\!}${{\raise4.5pt\hbox{\fivei{\char"7D}}}
$\!$({\teni X})
w.r.t. $\!$any topological space {\teni X},
is a simplicial set and, cf.$^{\!}$ \cite{19} %
\nobreak\ p.$_{_{^{\ \!\!}}}$362 Th.\ 4,
the map %
{\teni j}:%
\rlap{\raise1.5pt\hbox{\tent{\char"5E}}}{\vbox{\moveleft-1.2pt\hbox{\tensy{\char"6A}}}} %
$\!\!$%
\tenrm{\char"01}}$\!_{\!}${{\raise4.5pt\hbox{\fivei{\char"7D}}}
$\!$({\teni X})%
\rlap{\raise1.5pt\hbox{\tent{\char"5E}}}{\vbox{\moveleft-1.2pt\hbox{\tensy{\char"6A}}}} %
$\!$%
{\tensy{\char"21}}{\teni X} %
is a weak  homotopy equivalence i.e. induces isomorphisms in
homotopy groups, and {\teni j} is a true homotopy equivalence if
{\teni X} is of homotopy CW-type,
cf. \cite{9} %
pp.\ 76,\ {189ff},\nobreak\
\underbar{221-2}.}%
}

%


\refstyle{C}


\Refs


\ref\key 1 \by A. Bj\"orner \paper Topological Methods \inbook
Handbook of Combinatorics; R. Graham, M. Gr\"otschel and L.
Lov\'asz eds
%
\publ Elsevier Science B.V. \publaddr   \yr 1995 \pages
1819--1872$.
\ \ \triangleright$\ {\sevenrm 16} %
\endref


\ref\key 2 \by R. Brown \book Topology \publ Ellis Horwood
\publaddr \yr 1988.\ \ $\triangleright$\ {\sevenrm
4,7,8,9,14,15,16%
}
\endref


\ref\key 3 \by W.\ Bruns \& J.\ Herzog \book Cohen-Macaulay rings
\publ
Cambr. Univ. Press \publaddr %
\yr 1998.\ \ $\triangleright$\ {\sevenrm 18,21%
}
\endref


\ref\key 4 \by D.E. Cohen \paper Products and Carrier Theory \jour
Proc. London Math. Soc. \vol VII \yr 1957 \pages
219--248$.\hfill\break \triangleright$\ {\sevenrm 13,14%
}
\endref


\ref\key 5 \by A. Dold \book Lectures on Algebraic Topology \publ
Springer-Verlag. \publaddr %
\yr 1972.\ \
$\triangleright$\ {\sevenrm 8%
}
\endref


\ref\key 6 \by J. Dugundji \book Topology \publ Allyn \& Bacon
Inc. \publaddr Boston \yr 1966.\ \ $\triangleright$\ {\sevenrm
3,4,7,15} %
\endref


\ref\key 7 \by S. Eilenberg \& N. Steenrod \book Foundations of
Algebraic Topology \publ Princeton University Press \publaddr \yr
1952.\ \ $\triangleright$\ {\sevenrm 2,7,14,22} %
\endref


\ref\key 8 \by G. Fors \paper Algebraic Topological Results on
Stanley-Reisner Rings \inbook Commutative Algebra, Trieste, Italy
1992
%
\publ World Scientific Publ. Co. Pte. Ltd. \publaddr \yr 1994
\pages
69-88$.%
\ \ \triangleright$\ {\sevenrm 16,17,21,22}
\endref


\ref\key 9 \by R.\ Fritsch \& R.A.\ Piccinini \book Cellular
Structures in Topology \publ Cambr. Univ. Press. \publaddr \yr
1990. \hfill\break$\triangleright$\ {\sevenrm 5,14,31}
\endref


\ref\key 10 \by H,-G. Gr\"abe \paper \"Uber den Stanley-%
Reisner-Ring von Quasimannigfaltigkeiten \jour Math. Nachr. \vol
117 \yr 1984 \pages 161--174. $\triangleright$\ {\sevenrm
11,16,24,25,27}
\endref


\ref\key 11 \by H,-G. Gr\"abe \paper \"Uber den Rand von
Homologiemannigfaltigkeiten \jour Beitr\"age Algebra. Geom. \vol
22 \yr 1986 \pages 29--37. $\triangleright$\ {\sevenrm 25}
\endref


\ref\key 12 \by T. Hibi \paper Union and Gluing of a Family of
Cohen-Macaulay Partially Ordered Sets \jour Nagoya Math. J. \vol
107 \yr 1987 \pages 91--119. $\triangleright$\ {\sevenrm 20}
\endref




\ref\key 13 \by T. Hibi \paper Level Rings and Algebras with
Straightening Laws \jour J. Algebra \vol 117 \yr 1988 \pages
343--362. $\triangleright$\ {\sevenrm 20}
\endref


\ref\key 14 \by T. Hibi \book Algebraic Combinatorics on Convex
Polytopes \publ Carslaw Publications \publaddr \yr 1992.\ \
$\triangleright$\ {\sevenrm 21,29}
\endref


\ref\key 15 \by P.J.\ Hilton \& S.\ Wylie \book Homology Theory
\publ Cambr. Univ. Press. \publaddr \yr 1960.\ \ $\triangleright$\
{\sevenrm 9}
\endref


\ref\key 16 \by J. Lawson \& B. Madison \paper Comparisons of
Notions of Weak Hausdorffness \inbook Topology Vol {\eightbf 24}.
%
Proc. Memphis State Univ. Conf.,
%
Eds: AS.P. Franklin, B.V. (Smith)Thomas.
%
Lecture Notes in Pure and Appl. Math.
%
\publ Marcel Dekker Inc.  \publaddr New York  \yr 1976 \pages
207-215$.%
\ \ \triangleright$\ {\sevenrm 14}
\endref


\ref\key 17 \by C.R.F. Maunder \book Algebraic Topology  \publ Van
Nostrand Reinhold \publaddr \yr 1970.\ \ $\triangleright$\
{\sevenrm 28,29}
\endref


\ref\key 18 \by J.W. Milnor \paper Construction of Universal
Boundles. {\rm II}\jour Ann. of Math. (2) \vol 63 \yr 1956 \pages
430--436$.\hfill\break \triangleright$\ {\sevenrm 2,7,8}
\endref


\ref\key 19 \by J.W. Milnor \paper The Geometric Realization of a
Semi-Simplicial Complex \jour Ann. of Math. (2) \vol 65 \yr 1957
\pages 357--362$.\ \ \triangleright$\ {\sevenrm 31}
\endref


\ref\key 20 \by J.R. Munkres \book Elementary Differential
Topology {\eightrm Ann. of Math. Stud.} \vol 54 \publ Princeton
Univ. Press, NJ \publaddr \yr 1966$.\ \ \triangleright$\ {\sevenrm
14}
\endref


\ref\key 21 \by J.R. Munkres \book Elements of Algebraic Topology
\publ Benjamin/Cummings \publaddr \yr 1984. %
$\triangleright$\ {\sevenrm 7,9,10,11,13,22, 25,28,31}
\endref


\ref\key 22 \by J.R. Munkres \paper Topological Results in
Combinatorics. \jour Michigan Math. J. \vol 31 \yr 1984 \pages
113--128.\ \ $\triangleright$\ {\sevenrm 11,18,22}
\endref


\ref\key 23 \by A.A. Ranicki \paper On the Hauptvermutung
%
\inbook The Hauptvermutung Book, {\eighti K}-Monogr. Math., 1
%
\publ Kluwer Acad. Publ. \publaddr Dordrecht  \yr 1996 \pages
3--31$.%
\ \ \triangleright$\ {\sevenrm 14}
\endref


\ref\key 24 \by D. Rolfsen \book Knots and Links \publ
Dale Rolfsen/Publish or Perish, Inc \publaddr \yr 1990$.%
\ \ \triangleright$\ {\sevenrm 22} %
\endref


\ref\key 25 \by E.H. Spanier \book Algebraic Topology \publ
McGraw-Hill, Inc \publaddr \yr 1966 $.%
\ \ \triangleright$\ {\sevenrm 1,3,5,7,8,10,11,22,23,24,27}
\endref


\ref\key 26 \by R. Stanley \book Combinatorics and Commutative
Algebra{\eightrm ; 2nd. ed} \publ Birkh\"auser \publaddr \yr 1996
$.%
\ \ \triangleright$\ {\sevenrm 18,19,22}
\endref


\ref\key 27 \by B. St\"uckrad \& W. Vogel \book Buchsbaum Rings
and Applications \publ {\sevenbf VEB}/Springer-Verlag \publaddr
\yr 1986$.\hfill\break\triangleright$\ {\sevenrm 13,17,18,21,28}
\endref


\ref\key 28 \by J.E. Walker \paper Canonical Homeomorphisms of
Posets \jour European J. Combin. \vol 9 \yr 1988 \pages 97--107$.
\hfill\break\triangleright$\ {\sevenrm 14}
\endref


\ref\key 29 \by G.W. Whitehead \paper Homotopy Groups of Joins and
Unions.\jour Trans. Amer. Math. Soc. \vol 83 \yr 1956 \pages
55--69$.\ \triangleright$ {\sevenrm 2,7,14,27}
\endref


\endRefs
\enddocument